# Non-affine Landau-Ginzburg models and intersection cohomology

Thomas Reichelt and Christian Sevenheck

June 22, 2016

**Abstract**

We construct Landau-Ginzburg models for numerically effective complete intersections in toric manifolds as partial compactifications of families of Laurent polynomials. We show a mirror statement saying that the quantum $\mathcal{D}$-module of the ambient part of the cohomology of the submanifold is isomorphic to an intersection cohomology $\mathcal{D}$-module defined from this partial compactification and we deduce Hodge properties of these differential systems.

## Contents

2010 *Mathematics Subject Classification.* 14J33, 14M25, 32S40, 32S60, 14D07, 34Mxx, 53D45

During the preparation of this paper, Th.R. was supported by a postdoctoral fellowship of the "Fondation sciences mathématiques Paris" and by the DFG grant He 2287/2-2; Ch.S. was supported by a DFG Heisenberg fellowship (Se 1114/2-1/2). Both authors acknowledge partial support by the ANR grant ANR-08-BLAN-0317-01 (SEDIGA).

# 1 Introduction

The aim of this paper is the construction of a mirror model for complete intersections in smooth toric varieties. We consider the case where these subvarieties have a numerically effective anticanonical bundle. This includes in particular toric Fano manifolds, whose mirror is usually described by oscillating integrals defined by a family of Laurent polynomials and also the most prominent and classical example of mirror symmetry, namely, that of Calabi-Yau hypersurfaces in toric Fano manifolds. Here the mirror is a family of Calabi-Yau manifolds and the mirror correspondence involves the variation of Hodge structures defined by this family. One interesting feature of our results is that these apparently rather different situations occur as special cases of a general mirror construction, called *non-affine Landau-Ginzburg model*.

It is well-known that quantum cohomology theories admit expressions in terms of certain differential systems, called quantum $\mathcal{D}$-modules. This yields a convenient framework in which mirror symmetry is stated as an equivalence of such systems. Moreover, Hodge theoretic aspects of mirror correspondences can be incorporated using the machinery of (mixed) Hodge modules. However, quantum $\mathcal{D}$-modules have usually irregular singularities, except in the Calabi-Yau case. In our mirror construction, this corresponds to the fact that we let the Fourier-Laplace functor act on various regular $\mathcal{D}$-modules obtained from the Landau-Ginzburg model.

The quantum cohomology of a smooth complete intersection (which in our case is given as the zero locus of a generic section of a vector bundle) can be computed using the so-called *Euler-twisted Gromov-Witten invariants*. Basically, these are integrals over moduli spaces of stable maps of pull-backs of cohomology classes on the variety **and** of the Euler class of the vector bundle. It is well known (see [Kon95, Pan98, Giv98b] and also [Iri11] as well as [MM11] for more recent accounts) that the *ambient part* of the quantum cohomology of the subvariety (consisting of those classes which are induced from cohomology classes of the ambient variety), is given as a quotient of the Euler-twisted quantum cohomology.

From the combinatorial toric data of this vector bundle, we construct in a rather straightforward manner an *affine Landau-Ginzburg model*, which is a family of Laurent polynomials. The Euler-twisted quantum $\mathcal{D}$-module (which encodes the above mentioned Euler-twisted Gromov-Witten invariants) can then be shown to be isomorphic a certain proper FL-transformed *Gauß-Manin system*, namely, the Fourier-Laplace transformation of the top cohomology group of the compactly supported direct image complex (in the sense of $\mathcal{D}$-modules) of this affine Landau-Ginzburg model. On the other hand we show that the $Euler^{-1}$-twisted quantum $\mathcal{D}$-module which encodes the so-called local Gromov-Witten invariants is isomorphic to the usual FL-transformed Gauß-Manin system.

The actual non-affine Landau Ginzburg model is constructed by a certain partial compactification of the affine one, which yields a family of projective varieties. Our main result is Theorem 6.13 (which also contains the above mirror statements on twisted resp. local quantum $\mathcal{D}$-modules), it states that the *ambient quantum $\mathcal{D}$-module* is isomorphic to a Fourier-Laplace transform of the direct image of the intersection cohomology $\mathcal{D}$-module of the total space of this family, notice that this total space is usually not smooth.

One of the big advantages of using this singular variety together with the intersection cohomology $\mathcal{D}$-module is the fact that we do not need any kind of resolutions. In particular, we do not need to construct (or suppose the existence of) crepant resolutions like in [Bat94]. Notice also that [Iri11] discusses Landau-Ginzburg models of a more special class of subvarieties in toric orbifolds (the so-called nef partitions). In that paper, a mirror statement is shown in terms of A- resp. B-periods, but this construction needs a hypothesis on the smoothness of a certain complete intersection (given as the intersection of fibres of several Laurent polynomials, see section 5.2 of loc.cit.). Some more remarks on the nef-partition model and how it relates to our construction can be found in subsection 1.5 below.

We will show that the direct image of the intersection cohomology $\mathcal{D}$-module of the total space is itself (modulo some irrelevant free $\mathcal{O}$-modules) an intersection cohomology $\mathcal{D}$-module with respect to a local system measuring the intersection cohomology of the fibers of the projective family. An important point in our paper is that this intersection cohomology $\mathcal{D}$-module admits a hypergeometric description, that is, it can be derived from so-called GKZ-systems (as defined and studied by Gelfand', Kapranov and Zelevinsky). More precisely, it appears as the image of a morphism between two such GKZ-systems (Theorem 2.16). This result is interesting in its own, as in general there are only very few cases where geometrically interesting intersection cohomology $\mathcal{D}$-modules have an explicit description by differential operators.

Notice that the intersection cohomology $\mathcal{D}$-module mentioned above underlies a pure Hodge module.

From this we can deduce a Hodge-type property of the reduced quantum $\mathcal{D}$-module (see Corollary 6.14). As already mentioned above, it cannot underly itself a Hodge module, as in general it acquires irregular singularities (this never happens for $\mathcal{D}$-modules coming from variation of Hodge structures resp. Hodge modules due to Schmid's theorem). Rather, it is part of a non-commutative Hodge (ncHodge) structure due to a key result by Sabbah ([Sab08]).

There is another important aspect in the paper that has not yet been mentioned. The various quantum $\mathcal{D}$-modules are actually not $\mathcal{D}$-modules in the proper sense, rather, they are families of vector bundles on $\mathbb{P}^1$ together with a connection operator *with poles* along zero and infinity. This is reflected in the fact that we are looking at Fourier-Laplace transforms of certain regular $\mathcal{D}$-modules (like Gauß-Manin systems) together with a given filtration. The filtration induces a lattice structure on the FL-tranformed $\mathcal{D}$-module (i.e., it yields a coherent $\mathcal{O}$-submodule generating the FL-transformed $\mathcal{D}$-module). These lattices can be reconstructed by a twisted logarithmic de Rham complex (in the sense of log geometry) of an intermediate compactification of the family of Laurent polynomials. We show in Corollary 3.20 that this twisted logarithmic de Rham complex can also be explicitly described by hypergeometric equations. Notice that for this result to hold true, we have to restrict to an open subspace of the parameter space, where certain singularities at infinity of these Laurent polynomials are allowed, but not all of them. This situation is different to the one in our earlier paper [RS15] where we had to exclude any singularity at infinity.

The remaining part of this introduction is a rather detailed synopsis of the content of the paper. It can be read as a warm-up, where the main playing characters are introduced together with some examples which illustrates the constructions done later.

Our main case of interest is the following: Let $X_\Sigma$ be an $n$-dimensional smooth projective toric variety. Suppose that $\mathcal{L}_1 = \mathcal{O}_{X_\Sigma}(L_1), \ldots, \mathcal{L}_c = \mathcal{O}_{X_\Sigma}(L_c)$ are ample line bundles on $X_\Sigma$ such that $-K_{X_\Sigma} - \sum_{j=1}^c L_j$ is nef (for many intermediate results, we can actually relax both assumptions and suppose only that the individual bundles $\mathcal{L}_1, \ldots, \mathcal{L}_c$ are nef). Put $\mathcal{E} := \oplus_{j=1}^c \mathcal{L}_j$, then $\mathcal{E}$ is a convex vector bundle. We will be interested in several quantum $\mathcal{D}$-modules, which correspond to twisted Gromov-Witten invariants of $(X_\Sigma, \mathcal{E})$ as well as to Gromov-Witten invariants on the ambient cohomology of the complete intersection $Y := s^{-1}(0)$ defined by a generic section $s \in \Gamma(X_\Sigma, \mathcal{E})$. Let us consider the total space $\mathbb{V}(\mathcal{E}^\vee)$, which is a quasi-projective toric variety with defining fan $\Sigma'$. We set $\Sigma'(1) = \{\mathbb{R}_{\geq 0}\underline{b}_1, \ldots, \mathbb{R}_{\geq 0}\underline{b}_t\}$, where the vectors $\underline{b}_i$ are the primitive integral generators of the rays of $\Sigma'$. From this set of data one can construct Lefschetz fibrations, that is, family of hyperplane sections of some projective toric varieties. The actual Landau-Ginzburg models of the above toric variety (resp. of the complete intersection $Y$) will be obtained by restricting the base of such families to a certain sub-parameter space which is an open subset of the Kähler moduli space of $X_\Sigma$. Actually, we do consider two different situations: Either we start with the data of a toric variety and some line bundles satisfying the above positivity conditions, and construct the set $\{\underline{b}_i\}$ as sketched, or we consider only such a set of vectors, in which case we do not have a Kähler moduli space, and the reduction of the parameter space of the Lefschetz fibration is done using an equivariance property of this fibration with respect to a natural torus action. Nevertheless, many of our constructions also make sense in this more general setup, therefore, the material in sections 2 and 3 below only depend on vectors $\{\underline{b}_i\}$ and do not suppose the existence of $X_\Sigma, \mathcal{L}_1, \ldots, \mathcal{L}_c$.

## 1.1 Lefschetz fibrations

We consider the following situation: Let $B$ be $s \times t$-matrix of integer numbers, written as $B = (\underline{b}_1, \ldots, \underline{b}_t)$. The only assumption we make is that $\sum_{i=1}^t \mathbb{Z}\underline{b}_i = \mathbb{Z}^s$. As just explained, the example the reader should have in mind is when these vectors are the primitive integral generators of the rays of a possibly non-compact toric variety, but most of the constructions below do not depend on this assumption. As a concrete and easy though non-trivial example which will be considered throughout this introduction, let $X_\Sigma = \mathbb{P}^5$, $H \subset \mathbb{P}^5$ a hyperplane, and take the bundles $\mathcal{L}_1 = \mathcal{O}_{\mathbb{P}^5}(2H)$ and $\mathcal{L}_2 = \mathcal{O}_{\mathbb{P}^5}(3H)$. They are obviously ample, and we have $\mathcal{O}_{\mathbb{P}^5}(-K_{\mathbb{P}^5} - 2H - 3H) = \mathcal{O}_{\mathbb{P}^5}(6 - 2 - 3) = \mathcal{O}_{\mathbb{P}^5}(1)$, which is also ample. The defining fan $\Sigma'$ of the total space $\mathbb{V}(\mathcal{L}_1^\vee \oplus \mathcal{L}_2^\vee)$ has rays $\underline{b}_1, \ldots, \underline{b}_8$, and the matrix $B = (\underline{b}_1, \ldots, \underline{b}_8)$

is given by
$$B = \begin{pmatrix} 1 & 0 & 0 & 0 & 0 & -1 & 0 & 0 \\ 0 & 1 & 0 & 0 & 0 & -1 & 0 & 0 \\ 0 & 0 & 1 & 0 & 0 & -1 & 0 & 0 \\ 0 & 0 & 0 & 1 & 0 & -1 & 0 & 0 \\ 0 & 0 & 0 & 0 & 1 & -1 & 0 & 0 \\ 0 & 1 & 1 & 0 & 0 & 0 & 1 & 0 \\ 0 & 0 & 0 & 1 & 1 & 1 & 0 & 1 \end{pmatrix}$$

Let us return to the general setup of a matrix $B \in M(s \times t, \mathbb{Z})$ of rank $s$. Put $S := (\mathbb{C}^*)^s$, and consider the following map
$$g : S \longrightarrow \mathbb{P}^t,$$
$$(y_1, \ldots, y_s) \mapsto (1 : \underline{y}^{b_1}, \ldots, \underline{y}^{b_t}),$$

which is an embedding due to the assumption on the rank of $B$. Here we write $\underline{y}^{b_i}$ for the product $\prod_{k=1}^{s} y_k^{b_{ki}}$, $b_{ki}$ being the entries of $B$. The map $g$ is only locally closed, so we denote by $X$ its closure in $\mathbb{P}^t$. We are interested in a family of hyperplane sections of $X$, constructed in the following way: Consider the incidence variety $Z := \left\{\sum_{i=0}^{t} \lambda_i \cdot w_i\right\} \subset \mathbb{P}^t \times \mathbb{C}^{t+1}$, where $\lambda_0, \ldots, \lambda_t$ are coordinates on $\mathbb{C}^{t+1}$ and where $w_0 : \ldots : w_t$ are homogenous coordinates on $\mathbb{P}^t$. The situation is visualized in the following diagram where $p_1$ resp. $p_2$ is the restriction of the projection to the first resp. second factor

$$\begin{array}{c} Z \\ {}^{p_1}\swarrow \quad \searrow^{p_2} \\ S \subset X \subset \mathbb{P}^t \qquad \mathbb{C}^{t+1} \end{array} \qquad (1)$$

Here we have identified $S$ with its image under $g$. The family of hyperplane sections is by definition the morphism $\Phi := (p_2)_{|p_1^{-1}(X)} : p_1^{-1}(X) \to \mathbb{C}^{t+1}$. It is a projective map, and its restriction $\varphi := \Phi_{|p_1^{-1}(S)}$ is nothing but the family of Laurent polynomials

$$\begin{array}{rcl} S \times \mathbb{C}^t & \longrightarrow & \mathbb{C}_{\lambda_0} \times \mathbb{C}^t \\ (y_1, \ldots, y_s, \lambda_1, \ldots, \lambda_t) & \longmapsto & \left(-\sum_{i=1}^{t} \lambda_t \cdot \underline{y}^{b_i}, \lambda_1, \ldots, \lambda_t\right) \end{array} \qquad (2)$$

For the concrete example from above, the first component of $\varphi$ is given by

$$(y_1, \ldots, y_7, \lambda_1, \ldots, \lambda_8) \mapsto -\lambda_1 \cdot y_1 - \lambda_2 \cdot y_2 y_6 - \lambda_3 \cdot y_3 y_6 - \lambda_4 \cdot y_4 y_7 - \lambda_5 \cdot y_5 y_7 - \lambda_6 \frac{y_7}{y_1 \cdot \ldots \cdot y_5} - \lambda_7 y_6 - \lambda_8 y_7$$

The partial compactification $\Phi$ of this family is easy to calculate as the closure $X$ of $im(g)$ is a hypersurface in $\mathbb{P}^8$, namely, it is given by the binomial equation $w_0 w_7^2 w_8^3 - w_1 w_2 w_3 w_4 w_5 w_6 = 0$. Hence $Z$ is the codimension 2 subvariety of $\mathbb{P}^8 \times \mathbb{C}^9$ cut out by the two equations

$$w_0 w_7^2 w_8^3 - w_1 w_2 w_3 w_4 w_5 w_6 = 0 \quad \text{and} \quad \lambda_0 w_0 + \ldots + \lambda_8 w_8 = 0,$$

and $\Phi$ is the projection from this variety to the space $\mathbb{C}^9$ with coordinates $\lambda_0, \ldots, \lambda_8$.

For various reasons, we will also need to work with the family $\Phi^U$, where in the above diagram (1) the incidence variety $Z$ is replaced by its complement $U := (\mathbb{P}^t \times \mathbb{C}^{t+1}) \backslash Z$. Although geometrically the two morphisms $\Phi$ and $\Phi^U$ behave differently (e.g. $\Phi^U$ is no longer proper), they are strongly related on the cohomological level. The transformation corresponding in cohomology to the geometrical operation of taking the inverse image of $X$ under $p_1$ followed by the projection by $p_2$ is the so-called *Radon transformation* for $\mathcal{D}$-modules (see subsection 2.2 for more details).

The morphism $\varphi$ resp. $\Phi$ can be considered as the maximal family of hyperplane sections of $S$ resp. of its compactification $X$. However, in applications like those presented in section 6 of this paper, we need to restrict these families to some subspace of the parameter space $\mathbb{C}^t$ which is called $\mathcal{KM}^\circ$ in the main part of this article (see the discussion before Definition 6.3). We will not give the precise definition of $\mathcal{KM}^\circ$ here, let us only mention that the torus $S$ acts on $(\mathbb{C}^*)^t$ by $(\underline{y}, \underline{\lambda}) \mapsto (\underline{y}^{-b_1}, \ldots, \underline{y}^{-b_t}) \cdot \underline{\lambda}$ (see formula (26) in subsection 2.4 below). Then we consider the orbit space of this action, which is a torus of dimension $t-s$. The parameter subspace $\mathcal{KM}^\circ$ is a certain open subvariety of this orbit space. We will actually chose an embedding $\mathcal{KM}^\circ \hookrightarrow \mathbb{C}^t$, so that we always see $\mathbb{C}_{\lambda_0} \times \mathcal{KM}^\circ$ as a locally closed subspace of $\mathbb{C}^{t+1}$. In the case where our matrix $B$ is defined by a toric variety $X_\Sigma$ together with a set of line bundles, $\mathcal{KM}^\circ$ is not just an open subset of an abstract torus, but of $(\mathbb{C}^*)^{t-s}$, i.e., it comes with a set of coordinates called $q_1, \ldots, q_{t-s}$. Notice however that the choice of these coordinates is not unique, it depends on the choice of a basis of $H^2(X_\Sigma, \mathbb{Z})$ with good properties.

**Definition 1.1** (see Definition 6.3). Let $X_\Sigma$ be smooth, toric and projective. Let $\mathcal{L}_1 = \mathcal{O}_{X_\Sigma}(L_1), \ldots, \mathcal{L}_c = \mathcal{O}_{X_\Sigma}(L_c)$ be ample line bundles on $X_\Sigma$ such that $-K_{X_\Sigma} - \sum_{j=1}^{c} L_j$ is nef. Let $\Sigma'$ be the defining fan of the total space $\mathbb{V}(\mathcal{E}^\vee)$, where $\mathcal{E} := \oplus_{j=1}^{c} \mathcal{L}_j$ is a convex vector bundle on $X_\Sigma$. Let $\Sigma'(1) = \{\mathbb{R}_{\geq 0}\underline{b}_1, \ldots, \mathbb{R}_{\geq 0}\underline{b}_t\}$, where $\underline{b}_i$ are the primitive integral generators of the rays of $\Sigma'$. Let $\mathcal{KM}^\circ$ be the parameter space described above. Then the restrictions

$$\Pi := \Phi_{|\mathcal{Z}_X^\circ} : \mathcal{Z}_X^\circ := Z \cap p_1^{-1}(X) \cap \left(\mathbb{P}^t \times \mathbb{C}_{\lambda_0} \times \mathcal{KM}^\circ\right) \longrightarrow \mathbb{C}_{\lambda_0} \times \mathcal{KM}^\circ$$

resp.

$$\pi := \varphi_{|Z \cap p_1^{-1}(S) \cap \left(\mathbb{P}^t \times \mathbb{C}_{\lambda_0} \times \mathcal{KM}^\circ\right)} : Z \cap p_1^{-1}(S) \cap \left(\mathbb{P}^t \times \mathbb{C}_{\lambda_0} \times \mathcal{KM}^\circ\right) \longrightarrow \mathbb{C}_{\lambda_0} \times \mathcal{KM}^\circ$$

are called the **non-affine** resp. **affine** Landau-Ginzburg model of $(X_\Sigma, \mathcal{L}_1, \ldots, \mathcal{L}_c)$.

Let us notice that in the main body of this text, the affine Landau-Ginzburg model appears in two versions, called $\pi$ and $\overline{\pi}$. Actually, $\overline{\pi}$ is an intermediate partial compactification of $\pi$ (i.e., the fibres of $\overline{\pi}$ contain those of $\pi$ and are contained in those of $\Pi$).

To illustrate this definition, we discuss the parameter subspace $\mathcal{KM}^\circ$ for the above example of complete intersections of degree $(2, 3)$ in $\mathbb{P}^5$. As $B$ is a $7 \times 8$-matrix in this case, we see from what has been said above that $\mathcal{KM}^\circ$ must be an open subset of $\mathbb{C}^*$. We can choose the embedding $\mathbb{C}^* \hookrightarrow \mathbb{C}^7$, $q \mapsto (1, 1, 1, 1, q, 1, 1)$. The condition for a point $q \in \mathbb{C}^*$ to be in $\mathcal{KM}^\circ$ is then simply that the family $\varphi$, when restricted to $Z \cap p_1^{-1}(S) \cap (\mathbb{P}^t \times \mathbb{C}_{\lambda_0} \times \{q\})$ yields a *non-degenerate* Laurent polynomial, i.e., has no singularities at infinity (see Definition 3.8). One can easily show that the condition that $-K_{\mathbb{P}^5} - L_1 - L_2$ is ample (and not only nef) implies that this is the case for all $q \in \mathbb{C}^*$ (one has to argue along the lines of [RS15, Lemma 2.8]). Hence in this example, we have $\mathcal{KM}^\circ = \mathbb{C}^*$, and therefore the affine and the non-affine Landau-Ginzburg model of $(\mathbb{P}^5, \mathcal{O}_{\mathbb{P}^5}(2), \mathcal{O}_{\mathbb{P}^5}(3))$ are given as

$$\begin{aligned}
\pi : (\mathbb{C}^*)^7 \times \mathbb{C}^* &\longrightarrow \mathbb{C}_{\lambda_0} \times \mathbb{C}^* \\
(y_1, \ldots, y_7, q) &\longmapsto \left(-y_1 - y_2 y_6 - y_3 y_6 - y_4 y_7 - y_5 y_7 - q \frac{y_7}{y_1 \cdot \ldots \cdot y_5} - y_6 - y_7, q\right)
\end{aligned}$$

$$\begin{aligned}
\Pi : \mathcal{Z}_X^\circ &\longrightarrow \mathbb{C}_{\lambda_0} \times \mathbb{C}^* \\
(w_0 : \ldots : w_t, \lambda_0, q) &\longmapsto (\lambda_0, q)
\end{aligned}$$

where the quasi-projetive subvariety $\mathcal{Z}_X^\circ$ of $\mathbb{P}^t \times \mathbb{C}_{\lambda_0} \times \mathbb{C}^*$ is given by

$$\mathcal{Z}_X^\circ = \left\{w_0 w_7^2 w_8^3 - w_1 w_2 w_3 w_4 w_5 w_6 = 0, \lambda_0 w_0 + w_1 + \ldots + w_5 + q w_6 + w_7 + w_8 = 0\right\} \subset \mathbb{P}^t \times \mathbb{C}_{\lambda_0} \times \mathbb{C}^*.$$

## 1.2 GKZ-hypergeometric systems and Fourier-Laplace transformation

The main idea of this paper is that mirror correspondences can be expressed using the language of (filtered) $\mathcal{D}$-modules. For the toric varieties (and possibly non-toric subvarieties of them) that we are concerned with here, these $\mathcal{D}$-modules are of special type, namely they are constructed from the GKZ-system. Let us therefore start by recalling their definition (see Definition 2.7 below). We only treat here a special case which leads to a regular holonomic system. Let as before a matrix $B \in M(s \times t, \mathbb{Z})$ be given. Here we do not even need the condition $\text{rank}(B) = s$. Consider the matrix $\widetilde{B} \in M((s+1) \times (t+1), \mathbb{Z})$ defined by

$$\begin{pmatrix} 1 & 1 & \ldots & 1 \\ \hline 0 & & & \\ \vdots & & B & \\ 0 & & & \end{pmatrix}$$

**Definition 1.2.** *Let $\widetilde{B}$ be as above. Moreover, let $\widetilde{\beta} = (\beta_0, \ldots, \beta_s)$ be an element in $\mathbb{C}^{s+1}$. Write $\mathbb{L}$ for the module of relations among the columns of $\widetilde{B}$, i.e. the kernel of the linear mapping $\mathbb{Z}^{t+1} \to \mathbb{Z}^{s+1}$ given by $\widetilde{B}$. Let $D_{\mathbb{C}^{t+1}}$ be the Weyl algebra in $t+1$ variables, i.e., $D_{\mathbb{C}^{t+1}} := \mathbb{C}[\lambda_0, \lambda_1, \ldots, \lambda_t]\langle \partial_{\lambda_0}, \partial_{\lambda_1}, \ldots, \partial_{\lambda_t}\rangle$ Define*

$$M_{\widetilde{B}}^{\widetilde{\beta}} := D_{\mathbb{C}^{t+1}} / \left( (\Box_{\underline{l}})_{\underline{l} \in \mathbb{L}} + (E_k - \beta_k)_{k=0,\ldots s} \right),$$

*where*

$$\Box_{\underline{l}} := \prod_{i: l_i < 0} \partial_{\lambda_i}^{-l_i} - \prod_{i: l_i > 0} \partial_{\lambda_i}^{l_i}, \quad \underline{l} \in \mathbb{L}$$

$$E_k := \sum_{i=0}^{t} \widetilde{b}_{ki} \lambda_i \partial_{\lambda_i}, \quad k \in \{0, \ldots, s\}$$

*where $\widetilde{B} = \left(\widetilde{b}_{ki}\right)$. Then $M_{\widetilde{B}}^{\beta}$ is called a GKZ-system.*

We will quite often work with the corresponding sheaf of $\mathcal{D}_{\mathbb{C}^{t+1}}$-modules, denoted by $\mathcal{M}_{\widetilde{B}}^{\widetilde{\beta}}$. It is well known (see, e.g., [Ado94, Hot98]) that $\mathcal{M}_{\widetilde{B}}^{\widetilde{\beta}}$ is a regular holonomic $\mathcal{D}_{\mathbb{C}^{t+1}}$-module.

Given the matrix $B$ from the example of the last section (i.e for $(2,3)$-complete intersections in $\mathbb{P}^5$), we have $M_{\widetilde{B}}^{\widetilde{\beta}} = D_{\mathbb{C}^{t+1}}/I$, with

$$\begin{aligned} I &= \Big( \partial_{\lambda_0} \partial_{\lambda_7}^2 \partial_{\lambda_8}^3 - \partial_{\lambda_1} \partial_{\lambda_2} \partial_{\lambda_3} \partial_{\lambda_4} \partial_{\lambda_5} \partial_{\lambda_6}, \lambda_0 \partial_{\lambda_0} + \lambda_1 \partial_{\lambda_1} + \ldots + \lambda_8 \partial_{\lambda_8} - \beta_0, \\ & \quad \lambda_1 \partial_{\lambda_1} - \lambda_6 \partial_{\lambda_6} - \beta_1, \lambda_2 \partial_{\lambda_2} - \lambda_6 \partial_{\lambda_6} - \beta_2, \ldots, \lambda_5 \partial_{\lambda_5} - \lambda_6 \partial_{\lambda_6} - \beta_5, \\ & \quad \lambda_2 \partial_{\lambda_2} + \lambda_3 \partial_{\lambda_3} + \lambda_7 \partial_{\lambda_7} - \beta_6, \lambda_4 \partial_{\lambda_4} + \lambda_5 \partial_{\lambda_5} + \lambda_6 \partial_{\lambda_6} + \lambda_8 \partial_{\lambda_8} - \beta_7 \Big), \end{aligned}$$

where $\tilde{\beta} = (\beta_0, \beta_1, \ldots, \beta_7)$.

Let us describe a basic result from [Rei14] that shows how these $\mathcal{D}$-modules enter into the study of Landau-Ginzburg models. It uses the notion of Gauß-Manin systems, which are differential systems associated to any morphism between smooth algebraic (or analytic) varieties. Intuitivly, solutions of such systems are given by period integrals (at least on the smooth locus of the map). The formal definition requires the notion of direct images of $\mathcal{D}$-modules and is recalled in subsection 2.1 below. With these remarks in mind, we can state the result as follows (in the main part of the text it appears in a more precise version as Theorem 2.11). For simplicity, we also impose the additional assumption of normality, which is discussed in detail in section 5. We write $\mathbb{N}\widetilde{B}$ for the semi-group associated to $\widetilde{B}$, that is, $\mathbb{N}\widetilde{B} := \sum_{i=0}^{t} \mathbb{N}\underline{b}_i \subset \mathbb{Z}^s$.

**Theorem 1.3.** *Let the matrices $B$ and $\widetilde{B}$ be as above. Suppose moreover that the associated semi-group ring $\mathbb{C}[\mathbb{N}\widetilde{B}]$ is normal. Consider the family of Laurent polynomials $\varphi : S \times \mathbb{C}^t \to \mathbb{C}^{t+1}$ defined in equation (2). Then we have an exact sequence of regular holonomic $\mathcal{D}_{\mathbb{C}^{t+1}}$-modules*

$$0 \longrightarrow H^{s-1}(S, \mathbb{C}) \otimes \mathcal{O}_{\mathbb{C}^{t+1}} \longrightarrow \mathcal{H}^0 \varphi_+ \mathcal{O}_{S \times \mathbb{C}^{t+1}} \longrightarrow \mathcal{M}_{\widetilde{B}}^{\underline{0}} \longrightarrow H^s(S, \mathbb{C}) \otimes \mathcal{O}_{\mathbb{C}^{t+1}} \longrightarrow 0$$

*Here the left- and the rightmost terms are vector bundles on $\mathbb{C}^{t+1}$ together with the trivial connection operator which annihilates sections in $H^{s-1}(S, \mathbb{C})$ resp. $H^s(S, \mathbb{C})$, and $\mathcal{H}^0 \varphi_+ \mathcal{O}_{S \times \mathbb{C}^{t+1}}$ is the Gauß-Manin system of $\varphi$ alluded to above.*

An important aspect of the construction in [Rei14] that yields this result is that all the above $\mathcal{D}_{\mathbb{C}^{t+1}}$-modules underly *mixed Hodge modules* and that the exact sequence exists in the abelian category $MHM(\mathbb{C}^{t+1})$. Although Hodge theoretic considerations are one of the main motivations of this paper, we will not use this fact directly, and results on Hodge modules will not come into play until Corollary 6.14.

The theorem above shows that there is a tight connection between the Gauß-Manin system of the morphism $\varphi$ and the GKZ-system associated to the matrix $\widetilde{B}$. However, they are not equal, but their difference (i.e., kernel and cokernel of the morphism $\mathcal{H}^0 (\varphi_+ \mathcal{O}_{S \times \mathbb{C}^{t+1}}) \to \mathcal{M}_{\widetilde{B}}^{\underline{0}}$) are relativly simple. The next construction has the effect of erasing this difference and yields an isomorphism of the two $\mathcal{D}$-modules we are interested in. First we need a certain variant of the Fourier-Laplace transformation for holonomic $\mathcal{D}$-modules. Again we present a simplified version, the actual definition can be found in the next subsection as Definition 2.4.

**Definition 1.4.** *Let $Y$ be a smooth affine variety and let $D_{\mathbb{C}_\lambda \times Y}$ the ring of global algebraic differential operators on $\mathbb{C}_\lambda \times Y$. If $M$ is a $D_{\mathbb{C} \times Y}$-module, we denote by $\mathrm{FL}_Y(M)$ the object which is equal to $M$ as a module over $D_Y$ and where the new variable $\tau$ acts as $\partial_\lambda$ from the left and where $\partial_\tau$ acts as left multiplication by $-\lambda$. In this way $\mathrm{FL}_Y(M)$ becomes a left module over $D_{\mathbb{C}_\tau \times Y}$. Then we define*

$$\mathrm{FL}_Y^{loc}(M) := \mathrm{FL}_Y(M)[\tau^{-1}]$$

*to be the* localized Fourier-Laplace transformation *of $M$. Again we will denote by the same symbol the corresponding functor acting on sheaves of left $\mathcal{D}_{\mathbb{C}_\lambda \times Y}$-modules.*

With this definition at hand, we have the following easy consequence of Theorem 1.3.

**Corollary 1.5.** *Let $B$ and $\widetilde{B}$ be as above. Write $\widehat{\mathbb{C}}^{t+1}$ for the affine space $\operatorname{Spec} \mathbb{C}[\tau, \lambda_1, \ldots, \lambda_t]$. Then there is an isomorphism of $\mathcal{D}_{\widehat{\mathbb{C}}^{t+1}}$-modules*

$$\mathrm{FL}_{\mathbb{C}^t}^{loc}\left(\mathcal{H}^0 \varphi_+ \mathcal{O}_{S \times \mathbb{C}^{t+1}}\right) \cong \mathrm{FL}_{\mathbb{C}^t}^{loc}(\mathcal{M}_{\widetilde{B}}^0)$$

In the above example, we have $\mathrm{FL}_{\mathbb{C}^t}^{loc}(\mathcal{M}_{\widetilde{B}}^0) = D_{\widehat{\mathbb{C}}^{t+1}}/\widehat{I}$, where

$$\begin{aligned}\widehat{I} &= \Big(\tau\partial_{\lambda_7}^2\partial_{\lambda_8}^3 - \partial_{\lambda_1}\partial_{\lambda_2}\partial_{\lambda_3}\partial_{\lambda_4}\partial_{\lambda_5}\partial_{\lambda_6}, -\tau\partial_\tau + \lambda_1\partial_{\lambda_1} + \ldots + \lambda_8\partial_{\lambda_8} - 1, \lambda_1\partial_{\lambda_1} - \lambda_6\partial_{\lambda_6}, \\ &\quad \lambda_2\partial_{\lambda_2} - \lambda_6\partial_{\lambda_6}, \ldots, \lambda_5\partial_{\lambda_5} - \lambda_6\partial_{\lambda_6}, \lambda_2\partial_{\lambda_2} + \lambda_3\partial_{\lambda_3} + \lambda_7\partial_{\lambda_7}, \lambda_4\partial_{\lambda_4} + \lambda_5\partial_{\lambda_5} + \lambda_6\partial_{\lambda_6} + \lambda_8\partial_{\lambda_8}\Big),\end{aligned}$$

The partial compactification $\Phi$ of $\varphi$ is a projective morphism, but its source space $p_1^{-1}(X)$ is usually singular. For that reason, we are more interested in the direct image of the corresponding intersection cohomology $\mathcal{D}$-module. More precisely, consider the regular holonomic $\mathcal{D}_{\mathbb{P}^t}$-module $\mathcal{M}^{IC}(X)$ which correspondes to the intersection complex $IC_X$ of the variety $X$ (recall that $X$ was defined as the closure in $\mathbb{P}^t$ of the image of the embedding $g : S \hookrightarrow \mathbb{P}^t$) under the Riemann-Hilbert correspondence. Formally, $\mathcal{M}^{IC}(X)$ can be defined as the image of the natural morphism $g_\dagger \mathcal{O}_S \to g_+ \mathcal{O}_S$, where $g_\dagger$ is the "direct image with proper support"-functor for holonomic $\mathcal{D}$-modules. It is the minimal (also called intermediate) extension of its restriction to the smooth part of $X$, and as such is an irreducible $\mathcal{D}_{\mathbb{P}^t}$-module. More important, it underlies a pure polarizable algebraic Hodge module, i.e., an object of the category $MH^p(\mathbb{P}^t)$ (see [Sai88]). This last property will play a key role in Hodge theoretic application of our mirror statement (see Corollary 6.14).

In general it is quite hard to describe such intersection cohomology $\mathcal{D}$-modules explicitly, however, this is possible in the current situation. We have the following result (which we state directly in a form involving the functor $\mathrm{FL}_{\mathbb{C}^t}^{loc}$ since this is the result that will be used later)

**Theorem 1.6** (see Theorem 3.6 below). *Suppose that $\mathbb{C}[\mathbb{N}\widetilde{B}]$ is normal, then there is some parameter $\widetilde{\gamma} = (\gamma_0, \gamma_1, \ldots, \gamma_s) \in \mathbb{Z}^{s+1}$ such that*

$$\mathrm{FL}_{\mathbb{C}^t}^{loc}\left(\mathcal{H}^0 p_{2+} p_1^+ \mathcal{M}^{IC}(X)\right) \cong im\left(\mathrm{FL}_{\mathbb{C}^t}^{loc}(\mathcal{M}_{\widetilde{B}}^{\widetilde{\gamma}}) \xrightarrow{D} \mathrm{FL}_{\mathbb{C}^t}^{loc}(\mathcal{M}_{\widetilde{B}}^0)\right) \tag{3}$$

*Here $D$ is the morphism induced from right multiplication by $\tau^{-\gamma_0} \cdot \partial_{\lambda_1}^{g_1} \cdot \ldots \cdot \partial_{\lambda_t}^{g_t}$, where $\underline{g} = (g_1, \ldots, g_t)$ is any element in $\mathbb{Z}^t$ such that $B \cdot \underline{g}^{tr} = -(\gamma_1, \ldots, \gamma_s)$.*

The object on the left hand side of the above isomorphism should be seen (up to the action of the functor $\mathrm{FL}_{\mathbb{C}^t}^{loc}$) as a $\mathcal{D}$-module extending a local system the fibres of which are itself intersection cohomology groups, namely those of the fibres of the morphism $\Phi$. We could also replace the object $p_1^+ \mathcal{M}^{IC}(X)$ by $\mathcal{M}^{IC}(\mathcal{Z}_X^\circ)$, by which we mean the regular holonomic $\mathcal{D}_{\mathbb{P}^t \times \mathbb{C}^{t+1}}$-module corresponding to the intersection complex $IC_{\mathcal{Z}_X^\circ}$ via the Riemann-Hilbert correspondence (so that the complex $p_{2+} p_1^+ \mathcal{M}^{IC}(X) \cong p_{2+} \mathcal{M}^{IC}(\mathcal{Z}_X^\circ)$ corresponds to the topological direct image complex $R\Pi_* IC_{\mathcal{Z}_X^\circ}$).

Notice that the functors $p_1^+$ and $\mathcal{H}^0 p_{2+}$ exist in $MH^p$, hence the object occurring in the last theorem is the Fourier-Laplace transform of a $\mathcal{D}$-module underlying a pure polarizable Hodge module (this is basically the proof of Corollary 6.14)

We would like to explain in an informal way the reason for this theorem to hold true. The main point is that GKZ-systems behave quite well with respect to the duality functor for holonomic $\mathcal{D}$-modules. More precisely, we have the following very nice result of Walther (see [Wal07]).

**Theorem 1.7.** *Let $B$ and $\widetilde{B}$ be as above. Suppose again for simplicity that the semi-group ring $\mathbb{C}[\mathbb{N}\widetilde{B}]$ is normal. Then there is a parameter $\widetilde{\gamma} \in \mathbb{Z}^{s+1}$ such that*
$$\mathbb{D}\mathcal{M}^0_{\widetilde{B}} \cong \mathcal{M}^{\widetilde{\gamma}}_{\widetilde{B}}$$

From this statement we see that the above morphism $D$ can actually be seen (up to some shifts and notational conventions) as a morphism $\mathbb{D}\,\mathrm{FL}^{loc}_{\mathbb{C}^t}(\mathcal{M}^0_{\widetilde{B}}) \to \mathrm{FL}^{loc}_{\mathbb{C}^t}(\mathcal{M}^0_{\widetilde{B}})$. As mentioned above, $\mathcal{M}^{IC}(X)$ is the image of $g_\dagger \mathcal{O}_S \to g_+ \mathcal{O}_S$, notice further that these two $\mathcal{D}$-modules are also dual to each other. Applying the Radon transformation functor to them yields precisely the two GKZ-systems on the right hand side of equation (3) (see Theorem 2.11 below for more details), hence it is plausible that the intersection cohomology module (resp. its Fourier-Laplace transform) on the left hand side of equation (3) can be identified with the image of the morphisms $D$ between these two GKZ-systems.

For our purposes, we need actually a stronger duality statement: We consider the object $(\mathcal{M}^0_{\widetilde{B}}, F^{ord}_\bullet)$ consisting of the regular holonomic $\mathcal{D}_{\mathbb{C}^{t+1}}$-module $\mathcal{M}^0_{\widetilde{B}}$ together with the good filtration by coherent $\mathcal{O}_{\mathbb{C}^{t+1}}$-submodules induced from the filtration by the order of differential operators on $\mathcal{D}_{\mathbb{C}^{t+1}}$. This is an object of M. Saito's category $MF(\mathcal{D}_{\mathbb{C}^{t+1}})$ (see [Sai88, section 2.4]), and there is duality functor on this category extending the duality functor for holonomic $\mathcal{D}$-modules. Then we have (see Theorem 5.4) that $\mathbb{D}(\mathcal{M}^0_{\widetilde{B}}, F^{ord}_\bullet) \cong (\mathcal{M}^{\widetilde{\gamma}}_{\widetilde{B}}, F^{ord}_{\bullet+k})$ for some integer $k$.

For our guiding example, a parameter $\widetilde{\gamma}$ such that $\mathbb{D}\mathcal{M}^0_{\widetilde{B}} = \mathcal{M}^{\widetilde{\gamma}}_{\widetilde{B}}$ can be chosen as
$$\widetilde{\gamma} = (-1, 0, 0, 0, 0, 0, -1, -1)$$
and the map $D$ is induced by right multiplication with $\tau \cdot \partial_{\lambda_7} \cdot \partial_{\lambda_8}$.

Similarly to the considerations of Lefschetz families above, we will need to restrict these $\mathcal{D}$-modules to the parameter subspace $\mathcal{KM}^\circ \subset \mathbb{C}^t$. We will not explain here how to do this in detail, since it is a bit technical (see the presentation in subsection 2.4 and section 6 below). Instead, let us consider again the above example and the embedding

$$\begin{array}{ccc} \mathbb{C}_\tau \times \mathcal{KM}^\circ & \hookrightarrow & \mathbb{C}_\tau \times \mathbb{C}^7 \\ (\tau, q) & \longmapsto & (\tau, 1, 1, 1, 1, 1, q, 1, 1) \end{array}$$

Then we consider the inverse images under this map of $\mathcal{M}^0_{\widetilde{B}}$ and $\mathcal{M}^{\widetilde{\gamma}}_{\widetilde{B}}$ as well as the morphism $D$. For simplicity, we will also set $\tau = 1$, more precisely, we will consider the inverse image under the map $q \mapsto (1, q)$. We will also twist the restriction of $\mathcal{M}^{\widetilde{\gamma}}_{\widetilde{B}}$ by some invertible map (see Definition 6.1) Then the (restriction of the) morphism $D$ is given as

$$\begin{array}{ccc} \mathbb{C}[q^\pm]\langle \partial_q \rangle / (P_1) & \longrightarrow & \mathbb{C}[q^\pm]\langle \partial_q \rangle / (P_2) \\ Q & \longmapsto & Q \cdot (q\partial_q)^2 \end{array} \quad (4)$$

where
$$\begin{aligned} P_1 &= q \cdot (3q\partial_q + 1)(3q\partial_q + 2)(3q\partial_q + 3)(2q\partial_q + 1)(2q\partial_q + 2) + (q\partial_q)^6 \\ &= (q\partial_q)^2 \cdot \underbrace{\left(6q \cdot (3q\partial_q + 1)(3q\partial_q + 2)(2q\partial_q + 1) + (q\partial_q)^4\right)}_{Q^{(2,3)}} =: (q\partial_q)^2 \cdot Q^{(2,3)} \end{aligned}$$

$$\begin{aligned} P_2 &= q \cdot (3q\partial_q)(3q\partial_q + 1)(3q\partial_q + 2)(2q\partial_q)(2q\partial_q + 1) + (q\partial_q)^6 \\ &= \underbrace{\left(6q \cdot (3q\partial_q + 1)(3q\partial_q + 2)(2q\partial_q + 1) + (q\partial_q)^4\right)}_{Q^{(2,3)}} \cdot (q\partial_q)^2 =: Q^{(2,3)} \cdot (q\partial_q)^2 \end{aligned}$$

The map $D$ is obviously well defined, its kernel is generated by $Q^{(2,3)}$ and we see that
$$im(D) \cong \frac{\mathbb{C}[q^\pm]\langle \partial_q \rangle / (P_1)}{ker(D)} \cong \mathbb{C}[q^\pm]\langle \partial_q \rangle / (Q^{(2,3)}).$$

$Q^{(2,3)}$ is an inhomogenous hypergeometric operator with a regular singularity at $q = 0$ and irregular singularity at $q = \infty$.

The following statement summarizes the above calculation and can be seen as an illustration of Theorem 1.6

**Proposition 1.8.** *Consider the example from above. Then we have an isomorphism of left $\mathbb{C}[q^{\pm}]\langle\partial_q\rangle$-modules*
$$\left(\mathrm{FL}^{loc}_{\mathbb{C}^t}\left(\mathcal{H}^0\Pi_+\mathcal{M}^{IC}(\mathcal{Z}^\circ_X))\right)\right)_{|\tau=1} \cong \mathbb{C}[q^{\pm}]\langle\partial_q\rangle/(Q^{(2,3)}).$$

Let us finish this discussion with some remarks on the case of Calabi-Yau complete intersections in toric manifolds. Suppose that instead of the above example we had considered a $(2,4)$-complete intersections in $\mathbb{P}^5$, i.e., the matrix

$$B = \begin{pmatrix} 1 & 0 & 0 & 0 & 0 & -1 & 0 & 0 \\ 0 & 1 & 0 & 0 & 0 & -1 & 0 & 0 \\ 0 & 0 & 1 & 0 & 0 & -1 & 0 & 0 \\ 0 & 0 & 0 & 1 & 0 & -1 & 0 & 0 \\ 0 & 0 & 0 & 0 & 1 & -1 & 0 & 0 \\ 1 & 1 & 0 & 0 & 0 & 0 & 1 & 0 \\ 0 & 0 & 1 & 1 & 1 & 1 & 0 & 1 \end{pmatrix}$$

then the same arguments as above would lead to the operator

$$Q^{(2,4)} := 8q \cdot (2q\partial_q + 1)(4q\partial_q + 1)(4q\partial_q + 2)(4q\partial_q + 3) - (q\partial_q)^4$$

which is regular and homogeneous (with singularities at $q = 0, 2^{-10}, \infty$). In that case, the above statement can be sharpened in the following way.

**Proposition 1.9.** *Consider a $(2,4)$-complete intersection in $\mathbb{P}^5$, then we have the isomorphism of left $\mathbb{C}[q^{\pm}]\langle\partial_q\rangle$-modules*
$$\left(\mathcal{H}^0\Pi_+\mathcal{M}^{IC}(\mathcal{Z}^\circ_X)\right)_{|\lambda_0=1} \cong \mathbb{C}[q^{\pm}]\langle\partial_q\rangle/(Q^{(2,4)})$$

The reason for this to be true is that in the Calabi-Yau case, Fourier-Laplace transformation together with restriction to $\tau = 1$ has basically no effect, i.e., can be identified with restriction to $\lambda_0 = 1$. In particular, the object thus obtained still underlies a pure polarizable Hodge module (whereas in general, we obtain a variation on non-commutative Hodge structures, see 6.14 and Conjecture 6.15 below). This is consistent with classical results on mirror symmetry for Calabi-Yau hypersurfaces like the quintic in $\mathbb{P}^4$. Notice however that in these constructions, one uses certain crepant desingularizations in order to work with ordinary cohomology together with its Hodge structures instead of intersection cohomology as in the present paper. This may introduces a new difference when compared to our construction, which will however disappear when using the functor $\mathrm{FL}^{loc}$. This should basically follow from the decomposition theorem (say, for pure Hodge modules, see [Sai88, Corollaire 3]) when applied to the desingularization map. We will make some more remarks on how our construction is related to known mirror models for complete intersections in the later subsections 1.4 and 1.5 below. However, a thorough treatment of this comparison issue is delicate and will be postponed to a subsequent paper.

## 1.3 Mirror Correspondence

We would like to state here in a slightly informal way the main results of this paper. They can be expressed as isomorphism of two $\mathcal{D}_{\mathbb{C}_\tau \times \mathcal{KM}^\circ}$-modules, one obtained as sketched above (i.e., direct image $\mathcal{D}$-modules under the morphisms $\Pi$ resp. $\pi$), the other one derived from Gromov-Witten theory of the variety $X_\Sigma$ resp. from its subvarieties. The actual picture is considerably more complicated, in the sense that we do not just look at $\mathcal{D}$-modules over $\mathbb{C}_\tau \times \mathcal{KM}^\circ$, but at modules over $\mathcal{D}_{\mathbb{P}^1 \times \mathcal{KM}^\circ}$ together with a structure of $\mathcal{R}_{\mathbb{C}_{\tau^{-1}} \times \mathcal{KM}^\circ}$-modules where $\mathcal{R}_{\mathbb{C}_{\tau^{-1}} \times \mathcal{KM}^\circ}$ is the sheaf of Rees rings for the filtration by orders on differential operators on $\mathcal{D}_{\mathcal{KM}^\circ}$. This corresponds to the fact that the Gauß-Manin-systems as well as the direct image modules of intersection cohomology modules occurring do carry Hodge filtrations, i.e., underly objects of the category *MHM* of (algebraic) mixed Hodge modules (see [Sai90]). This very important additional information can be reformulated as the structure of an $\mathcal{R}$-module, and the latter is conserved by the functor $\mathrm{FL}^{loc}$. Hence our actual statements in section 6 are considerably stronger than what is announced here. In particular, the simplified statement below is basically only an identification of local systems and hence does not take into account the fact that these $\mathcal{D}$-modules have irregular singularities in general. Nevertheless, we think that it is still instructive. It should be seen as a snapshot of what the actual result looks like.

We consider the situation described above, that is, we let $X_\Sigma$ be a smooth, projective toric variety, and $\mathcal{L}_1 = \mathcal{O}_{X_\Sigma}(L_1), \ldots, \mathcal{O}_{X_\Sigma}(L_c)$ ample line bundles such that the class $-K_{X_\Sigma} - \sum_{j=1}^{c} L_j$ is nef. Then for a generic section $s \in \Gamma(X_\Sigma, \mathcal{E})$, the zero locus $Y := s^{-1}(0) \subset X_\Sigma$ is a complete intersection with nef anticanonical class. Put $\mathcal{E} = \oplus_{j=1}^{c} \mathcal{L}_j$, then $\mathcal{E}$ is a convex vector bundle on $X_\Sigma$, and we can consider twisted Gromov-Witten invariants which give rise to the (small) twisted quantum-$\mathcal{D}$-modules QDM$(X_\Sigma, \mathcal{E})$, i.e., a vector bundle on $\mathbb{P}^1 \times \mathcal{KM}^\circ$ with fibre $H^*(X_\Sigma, \mathbb{C})$ with connection operator defined by the twisted quantum product. Moreover, we have the endomorphism $c_{top}(\mathcal{E})$ of $H^*(X_\Sigma, \mathbb{C})$ given by cup product with the Euler class of $\mathcal{E}$, and we put $\overline{H^*(X_\Sigma, \mathbb{C})} := H^*(X_\Sigma, \mathbb{C})/ker(c_{top}(\mathcal{E}))$. Then the reduced or ambient quantum $\mathcal{D}$-modules, denoted by $\overline{\text{QDM}}(X_\Sigma, \mathcal{E})$, is a vector bundle on $\mathbb{P}^1 \times \mathcal{KM}^\circ$ with fibres $\overline{H^*(X_\Sigma, \mathbb{C})}$, and the connection is defined via the quantum product on *ambient* cohomology classes, i.e., classes in the image of the morphism $H^*(X_\Sigma, \mathbb{C}) \to H^*(Y, \mathbb{C})$ (notice that this image is isomorphic to the quotient $\overline{H^*(X_\Sigma, \mathbb{C})}$). Finally, we can also consider moduli spaces of stable maps into the total space $\mathbb{V}(\mathcal{E}^\vee)$ of the vector bundle dual to $\mathcal{E}$ (then $\mathcal{E}^\vee$ is concave), and this yields the so-called *local* Gromov-Witten invariants. The corresponding quantum-$\mathcal{D}$-module is denoted by QDM$(\mathcal{E}^\vee)$. We refer to the sections 4 and 6 for precise definitions of the various quantum $\mathcal{D}$-modules.

With these notions at hand, we have the following results.

**Theorem 1.10** (see theorem 6.13 and theorem 6.16). *Given $X_\Sigma$, $\mathcal{L}_1 = \mathcal{O}_{X_\Sigma}(L_1), \ldots, \mathcal{O}_{X_\Sigma}(L_c)$ as above, but suppose moreover that $\mathcal{L}_1, \ldots, \mathcal{L}_c$ are ample (but $-K_{X_\Sigma} - L_1 - \ldots - L_c$ is still only required to be nef). Consider the matrix $B$ constructed from the rays of $\Sigma'$. Let $\pi$ resp. $\Pi$ the affine resp. the non-affine Landau-Ginzburg models of $(X_\Sigma, \mathcal{L}_1, \ldots, \mathcal{L}_c)$. Then we have*

$$\text{FL}_{\mathbb{C}^t}^{loc}\left(\mathcal{H}^0\pi_\dagger \mathcal{O}_{S \times \mathcal{KM}^\circ}\right)_{|\mathbb{C}_\tau^* \times B_\varepsilon^*} \cong (\text{id}_{\mathbb{C}_\tau^*} \times \text{Mir})^* \left(\text{QDM}(X_\Sigma, \mathcal{E})_{|\mathbb{C}_\tau^* \times B_\varepsilon^*}\right),$$

$$\text{FL}_{\mathbb{C}^t}^{loc}\left(\mathcal{H}^0\pi_+ \mathcal{O}_{S \times \mathcal{KM}^\circ}\right)_{|\mathbb{C}_\tau^* \times B_{\varepsilon'}^*} \cong (\text{id}_{\mathbb{C}_z} \times \text{Mir}')^* \left(\text{QDM}(\mathcal{E}^\vee)_{|\mathbb{C}_\tau^* \times B_{\varepsilon'}^*}\right),$$

$$\text{FL}_{\mathbb{C}^t}^{loc}\left(\mathcal{H}^0\Pi_+ \mathcal{M}^{IC}(\mathcal{Z}_X^\circ)\right)_{|\mathbb{C}_\tau^* \times B_\varepsilon^*} \cong (\text{id}_{\mathbb{C}_\tau^*} \times \text{Mir})^* \left(\overline{\text{QDM}}(X_\Sigma, \mathcal{E})\right)_{|\mathbb{C}_\tau^* \times B_\varepsilon^*}.$$

*Here $B_\varepsilon^*, B_{\varepsilon'}^*$ are some (pointed) convergency neighborhoods of the large volume limit point in $\overline{\mathcal{KM}^\circ}$, Mir is the mirror map (see, e.g., [CG07, MM11]) and Mir' is some other coordinate change (which also involves the mirror map Mir).*

From the pureness property of $\mathcal{M}^{IC}(\mathcal{Z}_X^\circ)$ we can deduce the following corollary, which is (part of) the Hodge theoretic aspect of our mirror correspondence. As mentioned earlier, it relies on the notion of *non-commutative Hodge structures* (see [Sab11] for an overview) which is adapted to the occurrence of irregular singularities in the various quantum $\mathcal{D}$-modules.

**Corollary 1.11** (see corollary 6.14). *Under the assumptions of the last theorem, the ambient quantum $\mathcal{D}$-module $\overline{\text{QDM}}(X_\Sigma, \mathcal{E})$ (or at least its restriction to the convergency neighborhood $B_\varepsilon^*$) is part of a variation of non-commutative Hodge structures.*

We conjecture in 6.15 below that $\overline{\text{QDM}}(X_\Sigma, \mathcal{E})$ is itself a non-commutative Hodge structure, however, the proof of this conjecture would need some additional results on the Hodge filtration of $\mathcal{H}^0 p_{2+} \mathcal{M}^{IC}(\mathcal{Z}_X^\circ)$ which are not yet available.

### 1.4 Givental's mirror model

The aim of the next two subsection is to give some ideas on the relation of our construction to other mirror models for Calabi-Yau resp. nef-complete intersections inside a smooth toric variety $X_\Sigma$. The reader should be warned that a complete comparison of the construction presented in this paper to other models is not yet available, and will be subject to some future work. Nevertheless, we hope that the following remarks indicate that our mirror model can be considered as a unification and generalization of other constructions.

First we consider a construction that can be found (although in a very sketchy form) in [Giv98b, page 10-11]. Let as above $X_\Sigma$ be smooth, projective and toric, and let $\mathcal{L}_1 = \mathcal{O}_{X_\Sigma}(L_1), \ldots, \mathcal{L}_c = \mathcal{O}_{X_\Sigma}(L_c)$ be nef line bundles such that $-K_{X_\Sigma} - \sum_{j=1}^{c} L_j$ is nef. Let again $\underline{b}_1, \ldots, \underline{b}_t$ be the primitive integral generators of the rays of $\Sigma'$. Notice that if $\underline{a}_1, \ldots, \underline{a}_m$ are the generators of the rays of $\Sigma$, then $t = m + c$ (namely, we have $m$ generators $b_i$ projecting to the $a'_i s$, and $c$ generators $b_i$ that projects to zero under $\Sigma' \twoheadrightarrow \Sigma$.) Consider the affine space $\mathbb{C}^{m+c}$ with coordinates $w_1, \ldots, w_{m+c}$. Let $\underline{l}_1, \ldots, \underline{l}_r$ be a basis of the module of relations between the vectors $\underline{b}_1, \ldots, \underline{b}_{m+n}$ (so that $r = m - n$, since the $\underline{b}_i$'s ly in $\mathbb{Z}^{n+c}$). Actually, this basis should not be chosen in an arbitrary way, it is the basis dual to the basis $p_1, \ldots, p_r$ one of $H^2(X_\Sigma, \mathbb{Z})$ chosen in section 6 (see the discussion after the exact sequence (59) below). Write $\underline{l}_a = (l_{a1}, \ldots, l_{at})$ and consider the affine variety

$$E := \left\{ (\underline{w}, \underline{q}) \in (\mathbb{C}^*)^m \times \mathbb{C}^c \times \mathcal{KM}^\circ \mid \left( \prod_{i=1}^{m+c} w_i^{l_{ai}} = q_a \right)_{a=1,\ldots,r} \right\}.$$

Actually, Givental has a slightly different definition as he consider equivariant quantum cohomology, but we ignore this aspect here and concentrate on the case of the non-equivariant limit. Notice also that in Givental's paper appears only the restriction $E_q := E_{|(\mathbb{C}^*)^m \times \mathbb{C}^c \times \{q\}}$.

It can be shown that the closure of $E$ inside $\mathbb{C}^{m+c} \times \mathcal{KM}^\circ$ equals $E$ (see the argument in Proposition 4.8 below), and that the projection $E \to \mathcal{KM}^\circ$; $(\underline{w}, \underline{q}) \mapsto \underline{q}$ is precisely the mapping $\alpha \circ \beta \circ \gamma_1 : \mathcal{Z}^\circ_{X^{\text{aff}}} \to \mathcal{KM}^\circ$ as appearing in the diagram (60) in section 6.

In [Giv98b, page 10-11], Givental very briefly mentions the following oscillating integral

$$\int_{\Gamma \subset E} e^{\tau \cdot \left( \sum_{i=1}^m w_i - \sum_{j=1}^c w_{m+j} \right)} \cdot \frac{d\log(w_1) \wedge \ldots \wedge d\log(w_m) \wedge dw_{m+1} \wedge \ldots \wedge dw_{m+c}}{d\log(q_1) \wedge \ldots \wedge d\log(q_{m-n})} \tag{5}$$

where $\Gamma$ is some *real* non-compact $n + c$-dimensional cycle inside $E$, i.e., a Lefschetz thimble. Notice however that $E$ is singular in general, so that in any case one would need to specify further how to define this cycle. It is claimed in loc.cit (and easily verified) for any relation $\underline{l} = (l_1, \ldots, l_t)$ with $\sum_{i=1}^{m+c} l_i \underline{b}_i = 0$ satisfying $l_1 \geq 0, \ldots, l_m \geq 0$, this integral is annihilated by the differential operators

$$\Delta_{\underline{l}} := \prod_{i=1}^m \prod_{\nu=0}^{l_i - 1} \left( \sum_{a=1}^r l_{ai} q_a \tau \partial_{q_a} - \tau \nu \right) - \prod_{a=1}^r q_a^{\langle p_a, \underline{l} \rangle} \cdot \prod_{i=m+1}^{m+c} \prod_{\nu=1}^{l_i} \left( \sum_{a=1}^r l_{ai} \tau q_a \partial_{q_a} + \tau \nu \right)$$

In order to connect this statement to our construction, one needs to discuss the relation between oscillating integrals and Gauß-Manin systems in some detail. This is a rather classical subject, although there does not seem to exist a general reference covering the present situation. One can find in [Pha83, Pha85] a definition of oscillating integrals for certain polynomial mappings, and in [Sab08, section 1.b] a discussion of the topological Fourier-Laplace transformation which yields a cohomological description of Lefschetz thimbles.

Assuming that this relation between oscillating integrals and Gauß-Manin systems is properly established, one may conjecture that the integral (5) yields a solution of the module

$$\text{FL}^{loc}_{\mathbb{C}^t} \left( \mathcal{H}^0 \pi_\dagger \mathcal{O}_{S \times \mathcal{KM}^\circ} \right)$$

that appeared in Theorem 1.10. However, even if this were proved, it is still unclear whether this integral satisfy a stronger differential equation (this has been noticed by Givental himself in [Giv98b, page 10]), namely, one would like to show that it is even a solution of the system $\text{FL}^{loc}_{\mathbb{C}^t} \left( \mathcal{H}^0 \Pi_+ \mathcal{M}^{IC}(\mathcal{Z}^\circ_X) \right)$. However, we do not have any further evidence at this point for this conjecture.

## 1.5 Nef partition Landau-Ginzburg models

There is a special case of the construction described in the last subsection, for which a more complete description of a mirror model is available in the literature. Namely, assume that the nef line bundles $\mathcal{L}_1, \ldots, \mathcal{L}_c$ are obtained as (line bundles associated to) a sum of some of the torus invariant divisors of $X_\Sigma$. Then already in [Giv98b] is a sketch of how the above mentioned general construction of oscillating integrals can be made more precise. We will use the more recent paper [Iri11] as a references, which also incorporates ideas from [BB96b]. Let us give a very brief reminder of the part of [Iri11] relevant in the present situation.

**Definition 1.12** ([Giv98b, Iri11]). *Let $X_\Sigma$ be a $n$-dimensional, smooth toric variety given by a fan $\Sigma$ with torus-invariant divisors $D_1, \ldots, D_m$. A nef-partition is a partition $\{1, \ldots, m\} = I_0 \sqcup I_1 \ldots \sqcup I_c$ such that $\mathcal{L}_j = \mathcal{O}(\sum_{k \in I_j} D_k)$ is nef for $j = 0, \ldots, c$.*

Notice that Iritani's paper covers a larger domain of applications than this definition since he considers (nef partitions of) toric orbifolds. However, in the main body of our paper we are only concerned with complete intersections in manifolds, so we restrict to this situation here.

To a nef-partition one associates as above a vector bundle $\mathcal{E} = \bigoplus_{j=1}^c \mathcal{L}_j$ (notice that the sum here is running only from 1 to $c$ and does not include the bundle $\mathcal{L}_0$). Choosing a generic section $s \in \Gamma(X_\Sigma, \mathcal{E})$ gives a smooth nef complete intersection $Y \subset X_\Sigma$. Hence we see that the data of a toric manifold with a nef partition are a particular case of the setup considered in the main part of our paper. We will show in remark 1.17 below an example which falls in the scope of this paper but which does not come from a nef partition. In that sense our construction is a true generalization of the Givental-Iritani model. Notice also that in the case $I_0 = \emptyset$ the complete intersection $Y$ is a Calabi-Yau manifold, this is exactly the situation considered by Batyrev and Borisov in [BB96b].

In the approach of Iritani the mirror model of $Y$ is given by a function on a family of complete intersections of Laurent polynomials inside an $n$-dimensional torus $\check{T} = (\mathbb{C}^*)^n$ with coordiantes $t_1, \ldots, t_n$. In order to construct this, one associates to each $I_j$ the following Laurent polynomial

$$W_{\underline{\alpha}}^{(j)} = \sum_{i \in I_j} \alpha_i \underline{t}^{b_i} \qquad \text{for} \quad j = 0, \ldots, c$$

where the $\underline{b}_i$ are again the primitive integral generators of the one-dimensional cones of $\Sigma$. The family of complete intersections $\check{Y}_{\underline{\alpha}}$ over $\underline{\alpha} \in (\mathbb{C}^*)^m$ is then given by

$$\check{Y}_{\underline{\alpha}} = \{\underline{t} \in \check{T} \mid W_{\underline{\alpha}}^{(1)} = \ldots = W_{\underline{\alpha}}^{(c)} = 1\}.$$

**Assumption 1.13** ([Iri11, page 2936]). *In the above situation, we suppose that the affine variety $\check{Y}_{\underline{\alpha}}$ is a smooth complete intersection in $\check{T}$ for generic $\underline{\alpha} \in (\mathbb{C}^*)^m$.*

To the best of our knowledge, there is up to now no result available which would show how restrictive this assumption is. There are some speculations in [Iri11, remark 5.6] that the smoothness of $\check{Y}_{\underline{\alpha}}$ should be related to the smoothness of $Y \subset X_\Sigma$, at least in the case $I_0 = \emptyset$.

**Definition 1.14.** *Let $X_\Sigma$ and a nef partition $\{1, \ldots, m\} = I_0 \sqcup I_1 \ldots \sqcup I_c$ be given. Suppose that assumption 1.13 holds true. If $I_0 \neq \emptyset$, we call the restriction*

$$W_{\underline{\alpha}}^{(0)} : \check{Y}_{\underline{\alpha}} \longrightarrow \mathbb{C}$$

*the* nef-partition Landau-Ginzburg model *of $(X_\Sigma, (I_j)_{j=1,\ldots,m})$. If $I_0 = \emptyset$ we consider the (affine) Calabi-Yau complete intersection $\check{Y}_{\underline{\alpha}} \subset \check{T}$ (smooth by the above assumption) itself as the* nef-partition mirror model.

In the case $I_0 = \emptyset$, Baytrev and Borisov considered a compactification $\mathcal{Y}_{\underline{\alpha}}$ of $\check{Y}_{\underline{\alpha}}$ inside the projetice toric variety $\mathbb{P}_\nabla$, given by the polytope $\nabla = \nabla_1 + \ldots + \nabla_c$ with $\nabla_j = \text{Conv}(\{b_i\}_{i \in I_j})$. Since the Calabi-Yau varieties $\mathcal{Y}_{\underline{\alpha}}$ are usual singular and the ambient variety $\mathbb{P}_\nabla$ does not always admit a crepant resolution Batyrev and Borisov introduced so-called string-theoretic Hodge numbers $h_{st}^{p,q}(\mathcal{Y}_{\underline{\alpha}})$ and could show in [BB96a] that

$$h_{st}^{p,q}(\mathcal{Y}_{\underline{\alpha}}) = h_{st}^{\overline{n}-p,q}(Y) \qquad \text{for} \quad 0 \leq p, q \leq \overline{n}$$

where $\overline{n} = n - c$ is the dimension of $Y$.

In the case where $I_0$ is non-empty Iritani defines (under assumption 1.13) an oscillating integral

$$\int_{\Gamma_\mathbb{R}(\underline{\alpha})} e^{-W_{\underline{\alpha}}^{(0)}(t) \cdot \tau} \Omega_{\underline{\alpha}}$$

where $\Gamma_{\mathbb{R}(\underline{\alpha})} = \check{Y}_{\underline{\alpha}} \cap \check{T}_\mathbb{R}$ is a non-compact cycle in $\check{Y}_{\underline{\alpha}}$ ($\check{T}_\mathbb{R} := (\mathbb{R}_{>0})^n \subset \check{T}$ being the real torus), i.e., a real Lefschetz thimble, and

$$\Omega_{\underline{\alpha}} = \frac{\frac{dt_1}{t_1} \wedge \ldots \wedge \frac{dt_n}{t_n}}{dW_{\underline{\alpha}}^{(1)} \wedge \ldots \wedge dW_{\underline{\alpha}}^{(c)}}$$

is a holomorphic volume form on $\check{Y}_{\underline{\alpha}}$. Notice that the definition of the volume form uses the smoothness assumption 1.13 in an essential way. Notice also that in loc.cit., the variable $z := \frac{1}{\tau}$ is used.

In order to get a mirror theorem, Iritani defines a so-called $A$-period of the complete intersection $Y \subset X_\Sigma$

$$\Pi(1, \mathcal{O}_Y) = (J_Y(q, -z), z^{n - \frac{deg}{2}} z^\rho \hat{\Gamma}_Y)$$

here $J_Y$ is the $J$-function, which is a particular solution of the quantum $\mathcal{D}$-module of $Y$, and $\hat{\Gamma}_Y$ is the Gamma class of $Y$. We refer the reader to [Iri11] for details. The aforementioned mirror theorem of Iritani is the equality (see [Iri11, Theorem 5.7])

$$\Pi_Y(1, \mathcal{O}_Y) = \frac{1}{F(\underline{\alpha})} \int_{\Gamma_\mathbb{R}(\underline{\alpha})} e^{-W_{\underline{\alpha}}^{(0)}(t) \cdot \tau} \Omega_{\underline{\alpha}} \qquad \underline{\alpha} \in \mathcal{KM}^\circ \tag{6}$$

where $F(\underline{\alpha})$ is a certain coordinate change.

From the remarks in the last subsection on the relation between oscillating integrals and Gauß-Manin systems it seems plausible that the oscillating integral on the right hand side of equation (6) gives a solution of the following (Fourier-Laplace transformed) Gauß-Manin system

$$\mathrm{FL}^{loc} \mathcal{H}^0(W_{\underline{\alpha}}^{(0)})_+ \mathcal{O}_{\check{Y}_{\underline{\alpha}}} \in \mathrm{Mod}_{hol}(\mathcal{D}_{\mathbb{C}_\tau})$$

As before, we need to consider families of such differential systems by letting the parameter $\alpha$ vary within $\mathcal{KM}^\circ$, hence, we rather look at the $\mathcal{D}_{\mathbb{C}_\tau \times \mathcal{KM}^\circ}$-module

$$\mathrm{FL}^{loc}_{\mathcal{KM}^\circ} \mathcal{H}^0(W_{\underline{\alpha}}^{(0)}, \mathrm{pr}_{\underline{\alpha}})_+ \mathcal{O}_{Y_{\underline{\alpha}} \times \mathcal{KM}^\circ}$$

Finally, we have to eliminate the asymmetry between $W_{\underline{\alpha}}^{(0)}$ and $W_{\underline{\alpha}}^{(j)}$ for $j = 1, \ldots, c$. This can be easily done in the following way: As $Y_{\underline{\alpha}} = \bigcap_{j=1}^c \left(W_{\underline{\alpha}}^{(j)}\right)^{-1}(1)$, instead of considering the direct image $\mathcal{H}^0(W_{\underline{\alpha}}^{(0)})_+ \mathcal{O}_{Y_{\underline{\alpha}}}$, we can consider the direct image $\mathcal{H}^0\left(W_{\underline{\alpha}}^{(0)}, W_{\underline{\alpha}}^{(1)}, \ldots, W_{\underline{\alpha}}^{(c)}\right)_+ \mathcal{O}_{\check{T}}$ and restrict it to the subspace where the last $c$ coordinates are set to 1 (this follows from the base change theorem for holonomic $\mathcal{D}$-modules, see theorem 2.1 below). Finally, in order to be able to use the Radon transformation functor alluded to above, we need to allow $\underline{\alpha}$ to vary within $\mathbb{C}^m$, and not just in the subspace $\mathcal{KM}^\circ$. This motivates the following definition.

**Definition 1.15.** *Let us be given a nef partition $\{1, \ldots, m\} = I_0 \sqcup I_1 \ldots \sqcup I_c$ on $X_\Sigma$. Then we call the morphism*

$$\Theta : \check{T} \times \mathbb{C}^m \longrightarrow \mathbb{C}^{t+1} = \mathbb{C} \times \mathbb{C}^m \times \mathbb{C}^c$$
$$(\underline{t}, \underline{\alpha}) \mapsto (\lambda_0, \ldots, \lambda_{m+c}) = (-W_{\underline{\alpha}}^{(0)}, \underline{\alpha}, -W_{\underline{\alpha}}^{(1)}, \ldots, -W_{\underline{\alpha}}^{(c)})$$

*the nef partition Landau-Ginzburg model of $(X_\Sigma, (I_j)_{j=0,\ldots,c})$*

We proceed by comparing the nef partition Landau-Ginzburg model to our construction, as outlined before in this introduction. First notice that if we consider the variety $X_\Sigma$, then the collection of line bundles $\mathcal{L}_1 = \mathcal{O}(\sum_{k \in I_1} D_k), \ldots, \mathcal{L}_c = \mathcal{O}(\sum_{k \in I_c} D_k)$ satisfies (almost) the assumptions of our construction: Namely we have that $-K_{X_\Sigma} - \sum_{j=1}^c \left(\sum_{k \in I_j} D_k\right)$ is nef, since it is simply equal to $\sum_{k \in I_0} D_k$ and then its nefness follows from the nef partition assumption. For our main result (Theorem 6.13, see also Theorem 1.10), we need the stronger assumption that $\mathcal{L}_1, \ldots, \mathcal{L}_c$ are ample but this is unnecessary for many intermediate results.

In any case, given a nef partition, we have nef line bundles $\mathcal{L}_1, \ldots, \mathcal{L}_c$ on $X_\Sigma$ such that the bundle $\mathcal{O}_{X_\Sigma}(-K_{X_\Sigma}) \otimes \mathcal{L}_1^{-1} \otimes \ldots \otimes \mathcal{L}_c^{-1}$ is nef, and we can consider the matrix $B$ as constructed in the beginning of this introduction from the primitive integral generators of the fan $\Sigma'$ of the total space $\mathbb{V}(\mathcal{E}^\vee)$ with $\mathcal{E} = \oplus_{j=1}^c \mathcal{L}_j$. If we denote by $B_j$ the matrix with columns $(\underline{b}_i)_{i \in I_j}$ for $j = 0, \ldots, c$, then $B$ is an $(n+c) \times (m+c)$ integer matrix given by

$$B := \left(\begin{array}{ccccc|c} B_0 & B_1 & B_2 & \cdots & B_c & 0 \\ \hline 0 \cdots 0 & 1 \cdots 1 & 0 \cdots 0 & \cdots & 0 \cdots 0 & 1 \\ & 0 \cdots 0 & 1 \cdots 1 & \cdots & 0 \cdots 0 & \\ \vdots & \vdots & & & \vdots & \ddots \\ 0 \cdots 0 & 0 \cdots 0 & & & 1 \cdots 1 & 1 \end{array}\right)$$

As before we get an associated family of Laurent polynomials

$$\varphi : S \times \mathbb{C}^t \longrightarrow \mathbb{C} \times \mathbb{C}^t$$

$$(y_1, \ldots, y_s, \lambda_1, \ldots, \lambda_t) \mapsto -W_{\underline{\lambda}}^{(0)} - \sum_{i=1}^{c}(W_{\underline{\lambda}}^{(i)} + \lambda_{m+i}) \cdot y_{n+i},$$

remember that $S = (\mathbb{C}^*)^s$ where $s := n + c$ and $t := m + c$.

With these definition, we can state the following conjectural relationship of the nef-partion Landau-Ginzburg to our non-affine Landau-Ginzburg model.

**Conjecture 1.16.** *Consider a smooth projective toric variety $X_\Sigma$ with torus invariant divisors $D_1, \ldots, D_m$ and a nef partition $\{1, \ldots, m\} = I_0 \sqcup I_1 \ldots \sqcup I_c$. Then*

1. *There exists a morphism*

$$\mathrm{FL}^{loc}_{\mathbb{C}^t}(\mathcal{H}^0 \Theta_+ \mathcal{O}_{\check{T} \times \mathbb{C}^m}) \longrightarrow \mathrm{FL}^{loc}_{\mathbb{C}^t}(\mathcal{H}^0 \varphi_+ \mathcal{O}_{S \times \mathbb{C}^t}) \simeq \mathrm{FL}^{loc}_{\mathbb{C}^t}(\mathcal{M}^0_{\widetilde{B}}).$$

   *of holonomic $\mathcal{D}_{\widehat{\mathbb{C}}^{t+1}}$-modules.*

2. *This morphism induces an epimorphism $\mathcal{D}_{\widehat{\mathbb{C}}^{t+1}}$-modules*

$$\mathrm{FL}^{loc}_{\mathbb{C}^t}(W_{min}\mathcal{H}^0 \Theta_+ \mathcal{O}_{\check{T} \times \mathbb{C}^m}) \twoheadrightarrow \mathrm{FL}^{loc}_{\mathbb{C}^t}(W_{min}\mathcal{H}^0 \varphi_+ \mathcal{O}_{S \times \mathbb{C}^t}) \simeq im\left(\mathrm{FL}^{loc}_{\mathbb{C}^t}(\mathcal{M}^{\widetilde{\gamma}}_{\widetilde{B}}) \xrightarrow{D} \mathrm{FL}^{loc}_{\mathbb{C}^t}(\mathcal{M}^0_{\widetilde{B}})\right)$$

   *where $W_{min}$ is the minimal step of the weight filtration on the Gauß-Manin system $\mathcal{H}^0 \Theta_+ \mathcal{O}_{\check{T} \times \mathbb{C}^m}$ resp. $\mathcal{H}^0 \varphi_+ \mathcal{O}_{S \times \mathbb{C}^t}$ (which underlies a mixed Hodge modules, i.e., an element of the abelian category $MHM(\mathcal{D}_{\mathbb{C}^{t+1}})$).*

We are actually able to show the first part of this conjecture, but the proof is far to technical to be reproduced here. The second part is still open. Some evidence for this part of the conjecture comes from the fact that the $\mathcal{D}_{\mathbb{C}^{t+1}}$-module $W_{min}\mathcal{H}^0 \varphi_+ \mathcal{O}_{S \times \mathbb{C}^t}$ (which underlies a pure Hodge module) is, as we will see later, irreducible, and hence the above morphism must be surjective if it is not the zero map. Using the mirror symmetry statement of [Iri11] together with our main result, one may even speculate further that this map must be an isomorphism. We postpone a thorough discussion of these matters to a subsequent paper.

**Remark 1.17.** *We give an example of a smooth toric variety $X_\Sigma$ with two nef line bundles $\mathcal{L}_1 = \mathcal{O}_{X_\Sigma}(L_1)$, $\mathcal{L}_2 = \mathcal{O}_{X_\Sigma}(L_2)$ and $-K_{X_\Sigma} - L_1 - L_2 = 0$, which is not representable as a nef partition. Consider the two-dimensional toric variety given by the fan*

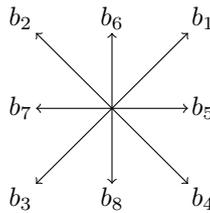

*The primitive generators of the rays give rise to a matrix*

$$A = \begin{pmatrix} 1 & -1 & -1 & 1 & 1 & 0 & -1 & 0 \\ 1 & 1 & -1 & -1 & 0 & 1 & 0 & -1 \end{pmatrix}$$

*The rows of the following matrix*

$$L = \begin{pmatrix} 1 & 0 & 0 & 0 & 0 & 0 & 1 & 1 \\ 0 & 1 & 0 & 0 & 0 & 0 & -1 & 1 \\ 0 & 0 & 1 & 0 & 0 & 0 & -1 & -1 \\ 0 & 0 & 0 & 1 & 0 & 0 & 1 & -1 \\ 0 & 0 & 0 & 0 & 1 & 0 & 1 & 0 \\ 0 & 0 & 0 & 0 & 0 & 1 & 0 & 1 \end{pmatrix}$$

*provide a basis for the module of relations among the columns of A. The well-known sequence*

$$0 \longrightarrow M \xrightarrow{A^t} \bigoplus_{i=1}^{8} \mathbb{Z} D_i \xrightarrow{L^t} H^2(X, \mathbb{Z}) \longrightarrow 0$$

*for torus-invariant Weil divisors endows the free $\mathbb{Z}$-module $H^2(X,\mathbb{Z})$ (which has rank 6) with a basis. The coordinates of the image $[D_i]$ of $D_i$ with respect to this basis are given by the i-th column vector of the matrix L. The closure of the Kähler cone is generated by the vectors*

$$\begin{pmatrix}1\\0\\0\\1\\1\\0\end{pmatrix}, \begin{pmatrix}1\\1\\-1\\1\\1\\1\end{pmatrix}, \begin{pmatrix}2\\0\\0\\2\\2\\1\end{pmatrix}, \begin{pmatrix}1\\1\\0\\1\\1\\1\end{pmatrix}, \begin{pmatrix}2\\0\\0\\0\\1\\1\end{pmatrix}, \begin{pmatrix}2\\1\\0\\0\\1\\1\end{pmatrix}, \begin{pmatrix}2\\0\\0\\1\\1\\1\end{pmatrix}, \begin{pmatrix}2\\2\\0\\0\\1\\2\end{pmatrix}, \begin{pmatrix}2\\1\\-1\\1\\1\\1\end{pmatrix}, \begin{pmatrix}3\\1\\-1\\1\\1\\2\end{pmatrix}, \begin{pmatrix}1\\1\\0\\0\\0\\1\end{pmatrix}, \begin{pmatrix}3\\1\\-1\\1\\2\\1\end{pmatrix}$$

*One easily sees that the vector*

$$[L_1] := \begin{pmatrix}1\\1\\-1\\1\\1\\1\end{pmatrix} \quad \text{and} \quad [L_2] := \begin{pmatrix}2\\0\\0\\0\\1\\1\end{pmatrix},$$

*which lie in the closure of the Kähler cone, are not of the type $[D_{i_1}] + \ldots + [D_{i_s}]$ for $\{i_1, \ldots, i_s\} \subset \{1, \ldots, 8\}$. But*

$$[L_1] + [L_2] = [D_1] + \ldots + [D_8] = [-K_{X_\Sigma}]$$

*Therefore the line bundles $\mathcal{L}_1 := \mathcal{O}(2D_2 + D_6 + D_7)$ and $\mathcal{L}_2 := \mathcal{O}(2D_1 + D_5 + D_6)$ as well as $-K_{X_\Sigma} - L_1 - L_2$ are nef.*

*Notice that although most of the constructions of our article apply to this example, it does not satisfy the assumptions of our main theorem 6.13 simply because the bundles $\mathcal{L}_1$ and $\mathcal{L}_2$ are nef but not ample. We could also give an example which consists of ample line bundles on a toric variety that do not come from a nef partition, but which would even be more complicated. Actually, we need ampleness only to apply results from quantum cohomology (like those presented in [MM11]), whereas for the constructions of the present paper, the nef assumption is sufficient.*

### 1.6 Acknowledgements

We thank Claude Sabbah for continuing support and interest in our work, Hiroshi Iritani, Etienne Mann and Thierry Mignon for useful discussions and for sending us an early version of their paper [IMM16]. We are very grateful to Takuro Mochizuki for his careful reading of the paper and for several very useful remarks. Special thanks go to Jörg Schürmann for pointing us to the reference [KW06].
We thank the anonymous referees for pointing out a number of inaccuracies and for suggesting several improvements of the paper.

## 2 Intersection Cohomology of Lefschetz fibrations

In this section we use the comparison result between Gauß-Manin systems of Laurent polynomials and GKZ-systems from [Rei14] to describe the direct image of the intersection complex of a natural compactification of a generic family of Laurent polynomials. The input data is an integer matrix $B$ of maximal rank and the GKZ-system in question will be defined by a certain homogenized matrix $\widetilde{B}$. The main tool is the Radon transformation resp. the Fourier-Laplace transformation for monodromic $\mathcal{D}$-modules ([Bry86]).

We start by a short remainder on some basic notions from the theory of algebraic $\mathcal{D}$-modules. Then we discuss Gauß-Manin systems, GKZ-systems and intersection cohomology $\mathcal{D}$-modules associated to the above mentioned families. Finally, we show using some facts about *quasi-equivariant* $\mathcal{D}$-modules that most of the objects considered here behave well with respect to a natural torus action on the parameter space of the families of Laurent polynomials resp. of their compactification.

## 2.1 Preliminaries

We review very briefly some basic results from the theory of algebraic $\mathcal{D}$-modules, which will be needed later. Let $\mathcal{X}$ be a smooth algebraic variety (we only consider algebraic varieties defined over $\mathbb{C}$ in the paper) of dimension $n$ and $\mathcal{D}_\mathcal{X}$ be the sheaf of algebraic differential operators on $\mathcal{X}$. We denote by $M(\mathcal{D}_\mathcal{X})$ the abelian category of algebraic $\mathcal{D}_\mathcal{X}$-modules on $\mathcal{X}$ and the abelian subcategory of (regular) holonomic $\mathcal{D}_\mathcal{X}$-modules by $M_h(\mathcal{D}_\mathcal{X})$ (resp. $(M_{rh}(\mathcal{D}_\mathcal{X}))$. The full triangulated subcategory in $D^b(\mathcal{D}_\mathcal{X})$ consisting of objects with (regular) holonomic cohomology is denoted by $D^b_h(\mathcal{D}_\mathcal{X})$ (resp. $D^b_{rh}(\mathcal{D}_\mathcal{X})$).

Let $f : \mathcal{X} \to \mathcal{Y}$ be a map between smooth algebraic varieties. Let $\mathcal{M} \in D^b(\mathcal{D}_\mathcal{X})$ and $\mathcal{N} \in D^b(\mathcal{D}_\mathcal{Y})$ be given, then we denote by $f_+ \mathcal{M} := Rf_*(\mathcal{D}_{\mathcal{Y} \leftarrow \mathcal{X}} \overset{L}{\otimes} \mathcal{M})$ resp. $f^+ \mathcal{N} := \mathcal{D}_{\mathcal{X} \to \mathcal{Y}} \overset{L}{\otimes} f^{-1} \mathcal{N}$ the direct resp. inverse image for $\mathcal{D}$-modules. Notice that the functors $f_+, f^+$ preserve (regular) holonomicity (see e.g., [HTT08, Theorem 3.2.3]). We denote by $\mathbb{D} : D^b_h(\mathcal{D}_\mathcal{X}) \to (D^b_h(\mathcal{D}_\mathcal{X}))^{op}$ the holonomic duality functor. Recall that for a single holonomic $\mathcal{D}_\mathcal{X}$-module $M$, the holonomic dual is also a single holonomic $\mathcal{D}_\mathcal{X}$-module ([HTT08, Proposition 3.2.1]) and that holonomic duality preserves regular holonomicity ([HTT08, Theorem 6.1.10]). For a morphism $f : \mathcal{X} \to \mathcal{Y}$ between smooth algebraic varieties we additionally define the functors $f_\dagger := \mathbb{D} \circ f_+ \circ \mathbb{D}$ and $f^\dagger := \mathbb{D} \circ f^+ \circ \mathbb{D}$.

In [HTT08], the definition of the inverse image functors from above follows a different convention, which is better adapted to the Riemann-Hilbert correspondence. Our functor $f^+$ corresponds to $f^\dagger[\dim(\mathcal{Y}) - \dim(\mathcal{X})]$ from [HTT08, page 31], whereas our functor $f^\dagger$ corresponds to $f^\star[\dim(\mathcal{X}) - \dim(\mathcal{Y})]$ from loc.cit, Definition 3.2.13.

Let $i : \mathcal{Z} \to \mathcal{X}$ be a closed embedding of a smooth subvariety of codimension $d$ and $j : \mathcal{U} \to \mathcal{X}$ be the open embedding of its complement. This gives rise to the following triangles for $\mathcal{M} \in D^b_{rh}(\mathcal{D}_\mathcal{X})$

$$i_+ i^+ \mathcal{M}[-d] \longrightarrow \mathcal{M} \longrightarrow j_+ j^+ \mathcal{M} \overset{+1}{\longrightarrow}, \tag{7}$$

$$j_\dagger j^\dagger \mathcal{M} \longrightarrow \mathcal{M} \longrightarrow i_\dagger i^\dagger \mathcal{M}[d] \overset{+1}{\longrightarrow}. \tag{8}$$

The first triangle is [HTT08, Proposition 1.7.1] and the second triangle follows by dualization. We will often use the following base change theorem.

**Theorem 2.1** ([HTT08, Theorem 1.7.3]). *Consider the following cartesian diagram of algebraic varieties*

$$\begin{array}{ccc} \mathcal{Z} & \overset{f'}{\longrightarrow} & \mathcal{W} \\ {\scriptstyle g'} \downarrow & & \downarrow {\scriptstyle g} \\ \mathcal{Y} & \overset{f}{\longrightarrow} & \mathcal{X} \end{array}$$

*then we have the canonical isomorphism $f^+ g_+[d] \simeq g'_+ f'^+[d']$, where $d := \dim \mathcal{Y} - \dim \mathcal{X}$ and $d' := \dim \mathcal{Z} - \dim \mathcal{W}$.*

**Remark 2.2.** *Notice that by symmetry we have also the canonical isomorphism $g^+ f_+[\tilde{d}] \simeq f'_+ g'^+[\tilde{d}']$ with $\tilde{d} := \dim \mathcal{W} - \dim \mathcal{X}$ and $\tilde{d}' := \dim \mathcal{Z} - \dim \mathcal{Y}$. In the former case we say we are doing a base change with respect to $f$, in the latter case with respect to $g$.*

**Remark 2.3.** *Using the duality functor we get isomorphisms:*

$$f^\dagger g_\dagger[-d] \simeq g'_\dagger f'^\dagger[-d'] \qquad and \qquad g^\dagger f_\dagger[-\tilde{d}] \simeq f'_\dagger g'^\dagger[-\tilde{d}'].$$

In the sequel, we will consider Fourier-Laplace transformations of various $\mathcal{D}$-modules. We give a short reminder on the definition and basic properties of the Fourier-Laplace transformation. Let $\mathcal{X}$ be a smooth algebraic variety, $U$ be a finite-dimensional complex vector space and $U'$ its dual vector space. Denote by $\mathcal{E}'$ the trivial vector bundle $\tau : U' \times \mathcal{X} \to \mathcal{X}$ and by $\mathcal{E}$ its dual. Write $\mathrm{can} : U \times U' \to \mathbb{C}$ for the canonical morphism defined by $\mathrm{can}(a, \varphi) := \varphi(a)$. This extends to a function $\mathrm{can} : \mathcal{E} \times \mathcal{E}' \to \mathbb{C}$.

**Definition 2.4.** *Define $\mathcal{L} := \mathcal{O}_{\mathcal{E}' \times_\mathcal{X} \mathcal{E}} e^{-can}$, this is by definition the free rank one module with differential given by the product rule. Denote by $p_1 : \mathcal{E}' \times_\mathcal{X} \mathcal{E} \to \mathcal{E}'$, $p_2 : \mathcal{E}' \times_\mathcal{X} \mathcal{E} \to \mathcal{E}$ the canonical projections. The Fourier-Laplace transformation is then defined by*

$$FL_\mathcal{X}(\mathcal{M}) := p_{2+}(p_1^+ \mathcal{M} \overset{L}{\otimes} \mathcal{L}) \quad \mathcal{M} \in D_h^b(\mathcal{D}_{\mathcal{E}'}).$$

If the base $\mathcal{X}$ is a point we will simply write FL. In general, the Fourier-Laplace transformation does not preserve regular holonomicity. However, it does preserve regular holonomicity for the derived category of complexes of $\mathcal{D}$-modules the cohomology of which are so-called *monodromic* $\mathcal{D}$-modules. We will give a short reminder on this notion. Let $\chi : \mathbb{C}^* \times \mathcal{E}' \to \mathcal{E}'$ be the natural $\mathbb{C}^*$ action on the fiber $U'$ and let $\theta$ be a coordinate on $\mathbb{C}^*$. We denote the push-forward $\chi_*(\theta \partial_\theta)$ as the Euler vector field $\mathfrak{E}$.

**Definition 2.5.** *[Bry86] A regular holonomic $\mathcal{D}_{\mathcal{E}'}$-module $\mathcal{M}$ is called monodromic, if the Euler field $\mathfrak{E}$ acts locally finite on $\tau_*(\mathcal{M})$, i.e. for a local section $v$ of $\tau_*(\mathcal{M})$ the set $\mathfrak{E}^n(v)$, $(n \in \mathbb{N})$, generates a finite-dimensional vector space. We denote by $D_{mon}^b(\mathcal{D}_{\mathcal{E}'})$ the derived category of bounded complexes of $\mathcal{D}_{\mathcal{E}'}$-modules with regular holonomic and monodromic cohomology.*

**Theorem 2.6.** *[Bry86]*

1. *$FL_\mathcal{X}$ preserves complexes with monodromic cohomology.*

2. *In $D_{mon}^b(\mathcal{D}_{\mathcal{E}'})$ we have*

$$FL_\mathcal{X} \circ FL_\mathcal{X} \simeq Id \quad \text{and} \quad \mathbb{D} \circ FL_\mathcal{X} \simeq FL_\mathcal{X} \circ \mathbb{D}.$$

3. *$FL_\mathcal{X}$ is t-exact with respect to the natural t-structure on $D_{mon}^b(\mathcal{D}_{\mathcal{E}'})$ resp. $D_{mon}^b(\mathcal{D}_\mathcal{E})$.*

*Proof.* The above statements are stated in [Bry86] for constructible monodromic complexes. One has to use the Riemann-Hilbert correspondence [Bry86, Proposition 7.12, Theorem 7.24] to translate the statements. So the first statement is Corollaire 6.12, the second statement is Proposition 6.13 and the third is Corollaire 7.23 in [Bry86]. $\square$

We will make occasionally use of the so-called $\mathcal{R}$-modules. More precisely, let $M$ be a smooth algebraic variety and consider the product of $M$ with the affine line $\mathbb{C}_z$ where $z$ is a fixed coordinate. Then by definition $\mathcal{R}_{\mathbb{C}_z \times M}$ is the $\mathcal{O}_{\mathbb{C}_z \times M}$-subalgebra of $\mathcal{D}_{\mathbb{C}_z \times M}$ locally generated by $z^2 \partial_z$ and by $z \partial_{x_1}, \ldots, z \partial_{x_n}$ where $(x_1, \ldots, x_n)$ are local coordinates on $M$. Notice that $j_M^* \mathcal{R}_{\mathbb{C}_z \times M} \cong \mathcal{D}_{\mathbb{C}_z^* \times M}$, where $j_M : \mathbb{C}_z^* \times M \hookrightarrow \mathbb{C}_z \times M$ is the canonical open embedding.

We will also consider the $\mathcal{O}_{\mathbb{C}_z \times M}$-subalgebra $\mathcal{R}'_{\mathbb{C}_z \times M}$ of $\mathcal{R}_{\mathbb{C}_z \times M}$ which is locally generated by $z \partial_{x_1}, \ldots, z \partial_{x_n}$ only. Sometimes we omit the subscript which denotes the underlying space, so we write $\mathcal{R}$ resp. $\mathcal{R}'$ instead of $\mathcal{R}_{\mathbb{C}_z \times M}$ resp. $\mathcal{R}'_{\mathbb{C}_z \times M}$. The inclusion $\mathcal{R}' \hookrightarrow \mathcal{R}$ induces a functor from the category of $\mathcal{R}$-modules to the category of $\mathcal{R}'$-modules, which we denote by $\text{For}_{z^2 \partial_z}$ ("forgetting the $z^2 \partial_z$-structure").

## 2.2 Gauß-Manin systems, hypergeometric $\mathcal{D}$-modules and the Radon transformation

In this subsection we adapt some results from [Rei14] to our situation. More precisely, for a given generic family of Laurent polynomials, we describe the canonical morphism between its Gauß-Manin-systems with compact support and its usual Gauß-Manin-systems. This mapping can be expressed as a morphism between the corresponding GKZ-systems. We will use this result in the next subsection to describe certain intersection cohomology modules.

We start by fixing our initial data and by introducing the GKZ-hypergeometric $\mathcal{D}$-modules. Let $B$ be a $s \times t$-integer matrix such that the columns of $B$, which we denote by $(\underline{b}_1, \ldots, \underline{b}_t)$, generate $\mathbb{Z}^s$. Consider the torus $S = (\mathbb{C}^*)^s$ and the $t+1$-dimensional vector space $V$ (with coordinates $\lambda_0, \lambda_1, \ldots, \lambda_t$) as well as its dual $V'$ (with coordinates $\mu_0, \mu_1, \ldots, \mu_t$). Define the map

$$\begin{aligned} g : S &\longrightarrow \mathbb{P}(V'), \\ (y_1, \ldots, y_s) &\mapsto (1 : \underline{y}^{\underline{b}_1}, \ldots, \underline{y}^{\underline{b}_t}), \end{aligned} \tag{9}$$

where $\underline{y}^{\underline{b}_i} := \prod_{k=1}^{s} y_k^{b_{ki}}$ for $i \in \{1,\ldots,t\}$. The condition on the columns of the matrix $B$ ensures that this is an embedding. If we denote the closure of the image of $g$ in $\mathbb{P}(V')$ by $X$, then $X$ is a (possibly non-normal) toric variety in the sense of [GKZ08, Chapter 5]. So we have the following sequence of maps

$$S \xrightarrow{j} X \xrightarrow{i} \mathbb{P}(V'), \qquad (10)$$

where $j$ is an open embedding and $i$ a closed embedding.

We will denote the homogeneous coordinates on $\mathbb{P}(V')$ by $(\mu_0 : \ldots : \mu_t)$. Let $Q$ be the convex hull of the elements $\{\underline{b}_0 = 0, \underline{b}_1, \ldots, \underline{b}_t\}$ in $\mathbb{R}^s$. Then by [GKZ08, Chapter 5, Prop 1.9] the projective variety $X$ has a natural stratification by torus orbits $X^0(\Gamma)$, which are in one-to-one correspondence with faces $\Gamma$ of the polytope $Q$. The orbit $X^0(\Gamma)$ is isomorphic to $(\mathbb{C}^*)^{\dim(\Gamma)}$ and is specified inside $X$ by the conditions

$$\mu_i = 0 \quad \text{for all} \quad \underline{b}_i \notin \Gamma, \quad \mu_i \neq 0 \quad \text{for all} \quad \underline{b}_i \in \Gamma. \qquad (11)$$

In particular the torus $S \subset X$ is given by the face $\Gamma = Q$, i.e. by the equations $\mu_i \neq 0$ for all $i \in \{0,\ldots,t\}$.

To this setup we associate the following $\mathcal{D}$-modules. Write $W = \mathbb{C}^t$ with coordinates $\lambda_1,\ldots,\lambda_t$ so that $V = \mathbb{C}_{\lambda_0} \times W$.

**Definition 2.7** ([GKZ90], [Ado94]). *Consider a lattice $\mathbb{Z}^s$ and vectors $\underline{b}_1,\ldots,\underline{b}_t \in \mathbb{Z}^s$. Moreover, let $\beta = (\beta_1,\ldots,\beta_s)$ be an element in $\mathbb{C}^s$. Write $\mathbb{L}$ for the module of relations among the columns of $B$. Any element $\underline{l} \in \mathbb{L}$ will be written as a vector $\underline{l} = (l_1,\ldots,l_t)$ in $\mathbb{Z}^t$. Define*

$$\mathcal{M}_B^\beta := \mathcal{D}_W / \left( (\Box_{\underline{l}})_{\underline{l} \in \mathbb{L}} + (E_k - \beta_k)_{k=1,\ldots s} \right),$$

*where*

$$\Box_{\underline{l}} := \prod_{i:l_i<0} \partial_{\lambda_i}^{-l_i} - \prod_{i:l_i>0} \partial_{\lambda_i}^{l_i}, \quad \underline{l} \in \mathbb{L}$$

$$E_k := \sum_{i=1}^{s} b_{ki} \lambda_i \partial_{\lambda_i}, \quad k \in \{1,\ldots,s\}$$

*where $b_{ki}$ is the $k$-th component of $\underline{b}_i$. The $\mathcal{D}_W$-module $\mathcal{M}_B^\beta$ is called a GKZ-system.*

As GKZ-systems are defined on the affine space $W$, we will often work with the $D_W$-modules of global sections $M_B^\beta := \Gamma(W, \mathcal{M}_B^\beta)$ rather than with the sheaves themselves, where $D_W = \mathbb{C}[\lambda_1,\ldots,\lambda_t]\langle \partial_{\lambda_1},\ldots,\partial_{\lambda_t}\rangle$ is the Weyl algebra.

We will also consider a homogenization of the systems above. Let $\widetilde{B}$ be the $(s+1) \times (t+1)$ integer matrix with columns $\widetilde{\underline{b}}_0 := (1,\underline{0}), \widetilde{\underline{b}}_1 := (1,\underline{b}_1),\ldots,\widetilde{\underline{b}}_t := (1,\underline{b}_t)$.

**Definition 2.8.** *Consider the hypergeometric system $\mathcal{M}_{\widetilde{B}}^{\widetilde{\beta}}$ on $V = \mathbb{C}^{t+1}$ associated to the vectors $\widetilde{\underline{b}}_0, \widetilde{\underline{b}}_1, \ldots, \widetilde{\underline{b}}_t \in \mathbb{Z}^{s+1}$ and $\widetilde{\beta} \in \mathbb{C}^{s+1}$. More explicitly, $\mathcal{M}_{\widetilde{B}}^{\widetilde{\beta}} := \mathcal{D}_V / \mathcal{I}$, where $\mathcal{I}$ is the sheaf of left ideals in $\mathcal{D}_V$ defined by*

$$\mathcal{I} := \mathcal{D}_V(\Box_{\underline{l}})_{\underline{l} \in \mathbb{L}} + \mathcal{D}_V(E_k - \beta_k)_{k=0,\ldots,s},$$

*where*

$$\Box_{\underline{l}} := \partial_{\lambda_0}^{\bar{l}} \cdot \prod_{i:l_i<0} \partial_{\lambda_i}^{-l_i} - \prod_{i:l_i>0} \partial_{\lambda_i}^{l_i} \quad \text{if} \quad \bar{l} \geq 0,$$

$$\Box_{\underline{l}} := \prod_{i:l_i<0} \partial_{\lambda_i}^{-l_i} - \partial_{\lambda_0}^{\bar{l}} \cdot \prod_{i:l_i>0} \partial_{\lambda_i}^{l_i} \quad \text{if} \quad \bar{l} < 0,$$

$$E_k := \sum_{i=1}^{t} b_{ki} \lambda_i \partial_{\lambda_i},$$

$$E_0 := \sum_{i=0}^{t} \lambda_i \partial_{\lambda_i}.$$

The generic rank of the GKZ-systems $\mathcal{M}_B^\beta$ resp. $\mathcal{M}_{\widetilde{B}}^{\widetilde{\beta}}$ may be difficult to predict depending the parameter (see, e.g., [MMW05]), but if we suppose that the matrix $B$ resp. $\widetilde{B}$ satisfies the normality assumption (see Proposition 5.1 and its proof below), then it is known that the rank of both modules equals $s! \cdot \text{vol}(\text{Conv}(\underline{b}_1, \ldots, \underline{b}_t))$ (where vol denotes the normalized volume, i.e. such that the volume of the hypercube is $[0,1]^d$ is one).

Let $h$ be the map given by
$$h : T \longrightarrow V', \tag{12}$$
$$(y_0, \ldots, y_s) \mapsto (\underline{y}^{\widetilde{b}_0}, \ldots, \underline{y}^{\widetilde{b}_t}) = (y_0, y_0\underline{y}^{b_1}, \ldots, y_0\underline{y}^{b_t}),$$

where $T = \mathbb{C}^* \times S = (\mathbb{C}^*)^{s+1}$. Notice that the restriction of $h$ to $\{1\} \times S$ is exactly the map $g$ from formula (9), when seen as a map to the affine chart $\{\mu_0 = 1\} \subset \mathbb{P}(V')$. We will later also need the closure of the image of $h$ in $V'$, which we denote by $Y$. Hence $Y$ is the affine cone over $X$.

As a piece of notation, for any matrix $C = (\underline{c}_1, \ldots, \underline{c}_k)$, we write $\mathbb{N}C$ for the semi-group generated by the columns $\underline{c}_1, \ldots, \underline{c}_k$, that is $\mathbb{N}C := \sum_{i_1}^{k} \mathbb{N}\underline{c}_i$, where we adopt the convention that the set $\mathbb{N}$ of natural numbers contains the element 0. Then we can consider the semi-group ring $\mathbb{C}[\mathbb{N}\widetilde{B}]$, which is naturally $\mathbb{Z}$-graded due to the first line of the matrix $\widetilde{B}$. Hence we can consider the ordinary spectrum of this ring as well as its projective spectrum, and it is clear that we have $Y = \text{Spec}\,\mathbb{C}[\mathbb{N}\widetilde{B}]$ and $X = \text{Proj}\,\mathbb{C}[\mathbb{N}\widetilde{B}]$.

We will now consider natural $D_V$-linear maps between GKZ-systems, which will induce a shift of the parameter. Let $\widetilde{B}$ be as above and consider the map of monoids
$$\rho : \mathbb{N}^{t+1} \longrightarrow \mathbb{N}\widetilde{B} \tag{13}$$
$$e_i \mapsto \widetilde{b}_i$$

where the $e_i$ are the standard generators of $\mathbb{N}^{t+1}$. Let $c \in \mathbb{N}^{t+1}$ be given and put $\widetilde{\gamma} := \rho(c)$. Notice that for every $\widetilde{\beta} \in \mathbb{C}^{s+1}$ the morphism
$$M_{\widetilde{B}}^{\widetilde{\beta}} \longrightarrow M_{\widetilde{B}}^{\widetilde{\beta}+\widetilde{\gamma}}$$
$$P \mapsto P \cdot \partial^c$$

is well-defined. Now let $c_1, c_2 \in \rho^{-1}(\widetilde{\gamma})$. Because $c_1$ and $c_2$ map to the same image, their difference $c_1 - c_2$ is a relation $\underline{l}$ among the columns of the matrix $\widetilde{B}$, thus $\partial^{c_1} - \partial^{c_2} \in (\Box_{\underline{l}})$. This shows that $P \cdot \partial^{c_1} = P \cdot \partial^{c_2}$ in $M_{\widetilde{B}}^{\widetilde{\beta}+\widetilde{\gamma}}$. Thus, we are lead to the following definition.

**Definition 2.9.** *Let $\widetilde{B}$ and $\widetilde{\beta}$ be as above. For every $\widetilde{\gamma} \in \mathbb{N}\widetilde{B}$ define the morphism*
$$\mathcal{M}_{\widetilde{B}}^{\widetilde{\beta}} \xrightarrow{\cdot \partial^{\widetilde{\gamma}}} \mathcal{M}_{\widetilde{B}}^{\widetilde{\beta}+\widetilde{\gamma}}$$

*given by right multiplication with $\partial^c$ for any $c \in \rho^{-1}(\widetilde{\gamma})$.*

In the next lemma, we establish a relation between a direct image under this morphism $h$ and the GKZ-systems just introduced.

**Lemma 2.10.** *There exists a $\delta_{\widetilde{B}} \in \mathbb{N}\widetilde{B}$ such that we have an isomorphism*
$$a : \text{FL}(h_+\mathcal{O}_T) \xrightarrow{\simeq} \mathcal{M}_{\widetilde{B}}^{\widetilde{\beta}} \tag{14}$$

*for every $\widetilde{\beta} \in \delta_{\widetilde{B}} + (\mathbb{R}_{\geq 0}\widetilde{B} \cap \mathbb{Z}^{s+1})$. Furthermore, we have a dual isomorphism*
$$a^\vee : \text{FL}(h_\dagger\mathcal{O}_T) \xleftarrow{\simeq} \mathcal{M}_{\widetilde{B}}^{-\widetilde{\beta}'} \tag{15}$$

*for every $\widetilde{\beta}' \in (\mathbb{R}_{\geq 0}\widetilde{B})^\circ \cap \mathbb{Z}^{s+1}$. For every $\widetilde{\beta}, \widetilde{\beta}'$ as above, the diagram below commutes up to a non-zero constant*

$$\begin{array}{ccc}
\mathcal{M}_{\widetilde{B}}^{-\widetilde{\beta}'} & \xrightarrow{\cdot \partial^{\widetilde{\beta}+\widetilde{\beta}'}} & \mathcal{M}_{\widetilde{B}}^{\widetilde{\beta}} \\
a^\vee \downarrow \simeq & & \simeq \uparrow a \\
\text{FL}(h_\dagger\mathcal{O}_T) & \longrightarrow & \text{FL}(h_+\mathcal{O}_T),
\end{array}$$

*where the lower horizontal morphism is induced by the natural morphism $h_\dagger\mathcal{O}_T \to h_+\mathcal{O}_T$.*

*Proof.* By [SW09, Corollary 3.7] we have the isomorphism $\mathrm{FL}(h_+(\mathcal{O}_T \cdot \underline{y}^{\widetilde{\beta}})) \simeq \mathcal{M}_{\widetilde{B}}^{\widetilde{\beta}}$ for every $\widetilde{\beta} \notin sRes(\widetilde{B})$ where $sRes(\widetilde{B})$ is the set of so-called strongly resonant parameters ([SW09, Definition 3.4]). Here $\mathcal{O}_T \cdot \underline{y}^{\widetilde{\beta}}$ is again the free rank one module with differential given by the product rule. Using [Rei14, Lemma 1.16], which says that there exists an $\delta_{\widetilde{B}} \in \mathbb{N}\widetilde{B}$ such that $\delta_{\widetilde{B}} + (\mathbb{R}_{\geq 0}\widetilde{B} \cap \mathbb{Z}^{s+1}) \cap sRes(\widetilde{B}) = \emptyset$ and the fact that $\mathcal{O}_T \simeq \mathcal{O}_T \cdot \underline{y}^{\widetilde{\gamma}}$ for every $\widetilde{\gamma} \in \mathbb{Z}^{s+1}$, the first statement follows. The second statement follows from taking the holonomic dual of (14), namely, we put

$$a^\vee := \mathbb{D}a : \mathbb{D}\mathcal{M}_{\widetilde{B}}^{\widetilde{\beta}} \xrightarrow{\simeq} \mathbb{D}\,\mathrm{FL}(h_+\mathcal{O}_T) \simeq \mathrm{FL}(\mathbb{D}h_+\mathcal{O}_T) \simeq \mathrm{FL}(h_\dagger \mathcal{O}_T)$$

and then we conclude by applying [Rei14, Proposition 1.23].

The last statement follows from the fact that the only non-zero morphism between $\mathcal{M}_{\widetilde{B}}^{-\widetilde{\beta}'}$ and $\mathcal{M}_{\widetilde{B}}^{\widetilde{\beta}}$ is right multiplication $\partial^{\widetilde{\beta}+\widetilde{\beta}'}$ up to a non-zero constant (cf. [Rei14, Proposition 1.24]). $\square$

We will denote by $Z \subset \mathbb{P}(V') \times V$ the universal hyperplane given by $Z := \{\sum_{i=0}^t \lambda_i \mu_i = 0\}$ and by $U := (\mathbb{P}(V') \times V) \setminus Z$ its complement. Consider the following diagram

$$\begin{array}{c}
U \\
\pi_1^U \swarrow \quad \downarrow j_U \quad \searrow \pi_2^U \\
\mathbb{P}(V') \xleftarrow{\pi_1} \mathbb{P}(V') \times V \xrightarrow{\pi_2} V \\
\pi_1^Z \nwarrow \quad \uparrow i_Z \quad \nearrow \pi_2^Z \\
Z
\end{array}$$

We will use in the sequel several variants of the so-called Radon transformation. These are functors from $D_{rh}^b(\mathcal{D}_{\mathbb{P}(V')})$ to $D_{rh}^b(\mathcal{D}_V)$ given by

$$\mathscr{R}(M) := \pi_{2+}^Z (\pi_1^Z)^+ M \simeq \pi_{2+}\, i_{Z+}\, i_Z^+ \pi_1^+ M\,,$$
$$\mathscr{R}^\circ(M) := \pi_{2+}^U (\pi_1^U)^+ M \simeq \pi_{2+}\, j_{U+}\, j_U^+ \pi_1^+ M\,,$$
$$\mathscr{R}_c^\circ(M) := \pi_{2\dagger}^U (\pi_1^U)^+ M \simeq \pi_{2+}\, j_{U\dagger}\, j_U^+ \pi_1^+ M\,,$$
$$\mathscr{R}_{cst}(M) := \pi_{2+}\, (\pi_1)^+ M\,,$$

The adjunction triangle corresponding to the open embedding $j_U$ and the closed embedding $i_Z$ gives rise to the following triangles of Radon transformations.

$$\mathscr{R}[-1](M) \longrightarrow \mathscr{R}_{cst}(M) \longrightarrow \mathscr{R}^\circ(M) \xrightarrow{+1}, \tag{16}$$

$$\mathscr{R}_c^\circ(M) \longrightarrow \mathscr{R}_{cst}(M) \longrightarrow \mathscr{R}[1](M) \xrightarrow{+1}, \tag{17}$$

where the second triangle is dual to the first.

We can now introduce the generic family of Laurent polynomials mentioned at the beginning of this subsection. It is defined by the columns of the matrix $B$, more precisely, we put

$$\varphi_B : S \times W \longrightarrow V = \mathbb{C}_{\lambda_0} \times W\,, \tag{18}$$

$$(y_1, \ldots, y_s, \lambda_1, \ldots, \lambda_t) \mapsto (-\sum_{i=1}^t \lambda_i \underline{y}^{b_i}, \lambda_1, \ldots, \lambda_t)\,.$$

The following theorem of [Rei14] constructs a morphism between the Gauß-Manin system $\mathcal{H}^0(\varphi_{B,+}\mathcal{O}_{S \times W})$ resp. the its proper version $\mathcal{H}^0(\varphi_{B,\dagger}\mathcal{O}_{S \times W})$ and certain GKZ-hypergeometric systems. For this we apply the triangle (16) to $M = g_\dagger \mathcal{O}_S$ and the triangle (17) to $M = g_+ \mathcal{O}_S$, which gives us the result.

**Theorem 2.11.** *[Rei14, Lemma 1.16, Theorem 2.7] There exists an $\delta_{\widetilde{B}} \in \mathbb{N}\widetilde{B}$ such that for every $\widetilde{\beta} \in \delta_{\widetilde{B}} + \mathbb{R}_{\geq 0}\widetilde{B} \cap \mathbb{Z}^{s+1}$ and every $\widetilde{\beta}' \in (\mathbb{N}\widetilde{B})^\circ = \mathbb{N}\widetilde{B} \cap (\mathbb{R}_{\geq 0}\widetilde{B})^\circ$, the following sequences of $\mathcal{D}_V$-modules are exact and dual to each other:*

$$\begin{array}{ccccccc}
H^{s-1}(S, \mathbb{C}) \otimes \mathcal{O}_V & & \mathcal{H}^0(\varphi_{B,+}\mathcal{O}_{S \times W}) & & \mathcal{M}_{\widetilde{B}}^{\widetilde{\beta}} & & H^s(S, \mathbb{C}) \otimes \mathcal{O}_V \\
\uparrow \simeq & & \uparrow \simeq & & \uparrow \simeq & & \uparrow \simeq \\
0 \longrightarrow \mathcal{H}^{-1}(\mathscr{R}_{cst}(g_+\mathcal{O}_S)) \longrightarrow & \mathcal{H}^0(\mathscr{R}(g_+\mathcal{O}_S)) & \longrightarrow \mathcal{H}^0(\mathscr{R}_c^\circ(g_+\mathcal{O}_S)) & \longrightarrow & \mathcal{H}^0(\mathscr{R}_{cst}(g_+\mathcal{O}_S)) \longrightarrow 0
\end{array}$$

$$\begin{array}{ccccccc}
0 \longleftarrow \mathcal{H}^1(\mathscr{R}_{cst}(g_\dagger\mathcal{O}_S)) \longleftarrow & \mathcal{H}^0(\mathscr{R}(g_\dagger\mathcal{O}_S)) & \longleftarrow \mathcal{H}^0(\mathscr{R}^\circ(g_\dagger\mathcal{O}_S)) & \longleftarrow & \mathcal{H}^0(\mathscr{R}_{cst}(g_\dagger\mathcal{O}_S)) \longleftarrow 0. \\
\downarrow \simeq & & \downarrow \simeq & & \downarrow \simeq & & \downarrow \simeq \\
H_c^{s+1}(S, \mathbb{C}) \otimes \mathcal{O}_V & & \mathcal{H}^0(\varphi_{B,\dagger}\mathcal{O}_{S \times W}) & & \mathcal{M}_{\widetilde{B}}^{-\widetilde{\beta}'} & & H_c^s(S, \mathbb{C}) \otimes \mathcal{O}_V
\end{array}$$

*If moreover $\mathbb{N}\widetilde{B}$ is saturated, then the vector $\delta_{\widetilde{B}}$ can be taken to be $\underline{0} \in \mathbb{N}\widetilde{B}$, in particular, the above statement holds for $\widetilde{\beta} = \underline{0} \in \mathbb{Z}^{s+1}$.*

Thus we get the following exact 4-term sequences which can be connected vertically by the map $\eta : H^0(\mathcal{R}(g_\dagger\mathcal{O}_S)) \to H^0(\mathcal{R}(g_+\mathcal{O}_S))$ induced by the natural morphism $g_\dagger\mathcal{O}_S \to g_+\mathcal{O}_S$. Define $\theta$ to be the composition $\kappa_2 \circ \eta \circ \kappa_1$. The next result gives a concrete description of this morphism:

$$\begin{array}{ccccccccc}
0 & \longrightarrow & H^{s-1}(S, \mathbb{C}) \otimes \mathcal{O}_V & \longrightarrow & \mathcal{H}^0(\mathscr{R}(g_+\mathcal{O}_S)) & \xrightarrow{\kappa_2} & \mathcal{M}_{\widetilde{B}}^{\widetilde{\beta}} & \longrightarrow & H^s(S, \mathbb{C}) \otimes \mathcal{O}_V & \longrightarrow & 0 \\
 & & & & \uparrow \eta & & \uparrow \theta & & \\
0 & \longleftarrow & H_c^{s+1}(S, \mathbb{C}) \otimes \mathcal{O}_V & \longleftarrow & \mathcal{H}^0(\mathscr{R}(g_\dagger\mathcal{O}_S)) & \xleftarrow{\kappa_1} & \mathcal{M}_{\widetilde{B}}^{-\widetilde{\beta}'} & \longleftarrow & H_c^s(S, \mathbb{C}) \otimes \mathcal{O}_V & \longleftarrow & 0.
\end{array}$$

**Lemma 2.12.** *The morphism $\theta$ is induced by right multiplication with $\partial^{\widetilde{\beta}+\widetilde{\beta}'}$ up to a non-zero constant.*

*Proof.* Once we can prove that $\kappa_2 \circ \eta \circ \kappa_1$ is not equal to zero we apply a rigidity result of [Rei14, Proposition 1.24] which says that the only maps between $\mathcal{M}_{\widetilde{B}}^{-\widetilde{\beta}'}$ and $\mathcal{M}_{\widetilde{B}}^{\widetilde{\beta}}$ is right-multiplication with $c \cdot \partial^{\widetilde{\beta}+\widetilde{\beta}'}$ for $c \in \mathbb{C}$. We only have to show that $\kappa_2 \circ \eta \circ \kappa_1$ becomes an isomorphism after microlocalizing with respect to $\partial_0 \cdots \partial_t$. This is sufficient as the microlocalization of the GKZ-systems $\mathcal{M}_{\widetilde{B}}^{\widetilde{\beta}}$ resp. $\mathcal{M}_{\widetilde{B}}^{-\widetilde{\beta}'}$ are not zero for otherwise the sheaves $h_+\mathcal{O}_T$ and $h_\dagger\mathcal{O}_T$ would be supported on the divisor $\{\mu_0 \cdot \mu_1 \cdot \ldots \cdot \mu_t = 0\}$, which is obviously wrong.

It is clear that $\kappa_1$ and $\kappa_2$ become isomorphisms after (micro-)localization with respect to $\partial_0 \cdots \partial_t$ because these maps have $\mathcal{O}_V$-free kernel and cokernel. It remains to prove that $\eta$ is an isomorphism after this micro-localization. To prove this we will use a theorem of [DE03] which compares the Radon transformation with the Fourier-Laplace transformation for $\mathcal{D}$-modules. Consider the following diagram

$$\begin{array}{ccccc}
T & \xrightarrow{h} & V' & \xleftarrow{p} & Bl_0(V') \\
& \searrow^{\tilde{h}} & \uparrow j_0 & \swarrow^q & \\
\pi_T \downarrow & & V' \setminus \{0\} & & \\
& & \downarrow \pi & & \\
S & \xrightarrow{g} & \mathbb{P}(V') & &
\end{array}$$

where $Bl_0(V') \subset \mathbb{P}(V') \times V'$ is the blow-up of $0$ in $V'$ and $q$ is the restriction of the projection to the first component. Notice that the map $h : T \to V'$ from formula (12) factors via $V' \backslash \{0\}$, that is, we have $h = j_0 \circ \widetilde{h}$, where $j_0 : V' \backslash \{0\} \hookrightarrow V'$ is the canonical inclusion.

It follows from [DE03, Proposition 1] that we have the following isomorphism
$$\mathscr{R}(g_+\mathcal{O}_S) \simeq \mathrm{FL}(p_+q^+g_+\mathcal{O}_S) \tag{19}$$

and its holonomic dual
$$\mathscr{R}(g_\dagger\mathcal{O}_S) \simeq \mathrm{FL}(p_+q^+g_\dagger\mathcal{O}_S), \tag{20}$$

where we have used $\mathscr{R} \circ \mathbb{D} \cong \mathbb{D} \circ \mathscr{R}$, $\mathrm{FL} \circ \mathbb{D} \cong \mathbb{D} \circ \mathrm{FL}$, $p_+ \circ \mathbb{D} \cong \mathbb{D} \circ p_+$ ($p$ is proper) and $q^+ \circ \mathbb{D} \cong \mathbb{D} \circ q^+$ ($q$ is smooth). Recall that we want to show that the morphism

$$\mathcal{H}^0(\mathscr{R}(g_\dagger\mathcal{O}_S)) \xrightarrow{\eta} \mathcal{H}^0(\mathscr{R}(g_+\mathcal{O}_S)),$$

becomes an isomorphism after localization with respect to $\partial_{\lambda_0} \cdots \partial_{\lambda_t}$. Using the isomorphisms (19) and (20) and the fact that FL is an exact functor and that it exchanges the action of $\mu_i$ and $\partial_{\lambda_i}$ we see that it is enough to show that
$$\mathcal{H}^0(p_+q^+g_\dagger\mathcal{O}_S) \longrightarrow \mathcal{H}^0(p_+q^+g_+\mathcal{O}_S) \tag{21}$$

becomes an isomorphism after localization with respect to $\mu_0 \cdots \mu_t$. In other words, we have to show that the kernel and the cokernel of the morphism (21) are supported on $\{\mu_0 \cdots \mu_t = 0\} \subset V'$. Obviously, we have $\{0\} \subset \{\mu_0 \cdot \ldots \cdot \mu_t = 0\}$ and hence $V' \backslash \{\mu_0 \cdot \ldots \cdot \mu_t = 0\} \subset V' \backslash \{0\}$. It is thus sufficient to show that kernel and cokernel of the restriction of the morphism (21) to $V' \backslash \{0\}$ are supported on $\{\mu_0 \cdot \ldots \cdot \mu_t = 0\} \backslash \{0\}$. Notice that the restriction of $\mathcal{H}^0(p_+q^+g_\dagger\mathcal{O}_S)$ resp. $\mathcal{H}^0(p_+q^+g_+\mathcal{O}_S)$ to $V' \backslash \{0\}$ is isomorphic to $\mathcal{H}^0(\pi^+g_\dagger\mathcal{O}_S)$ resp. $\mathcal{H}^0(\pi^+g_+\mathcal{O}_S)$. Thus the kernel and the cokernel of (21) are supported on $\{\mu_0 \cdots \mu_t = 0\}$ if and only if kernel and cokernel of

$$\mathcal{H}^0(\pi^+g_\dagger\mathcal{O}_S) \longrightarrow \mathcal{H}^0(\pi^+g_+\mathcal{O}_S)$$

are supported on $\{\mu_0 \cdots \mu_t = 0\} \backslash \{0\}$. The map $\pi$ is smooth and therefore $\pi^+$ is an exact functor. It is therefore enough to show that kernel and cokernel of

$$\mathcal{H}^0(g_\dagger\mathcal{O}_S) \longrightarrow \mathcal{H}^0(g_+\mathcal{O}_S)$$

are supported on $\{\mu_0 \ldots \mu_t = 0\} \subset \mathbb{P}(V')$. But this follows from the description of the map $g$, namely, by the remark right after equation (11) the support of the cone of the morphism $g_\dagger\mathcal{O}_S \to g_+\mathcal{O}_S$ is contained in $\{\mu_0 \ldots \mu_t = 0\}$. $\square$

### 2.3 Intersection cohomology $\mathcal{D}$-modules

As mentioned in the beginning of this section, our aim is to describe a $\mathcal{D}_V$-module derived from the intersection complex of a natural compactification of the family of Laurent polynomials $\varphi_B$ as defined in formula (18). This module will actually appear as the Radon transformation of the ($\mathcal{D}$-module corresponding to the) intersection complex of the variety $X \subset \mathbb{P}(V')$.

We start by fixing some notations concerning these $\mathcal{D}$-modules. Let $\mathcal{P}$ be a smooth variety and $\mathcal{U} \subset \mathcal{P}$ be a smooth subvariety, write $\mathcal{X}$ for the closure of $\mathcal{U}$ inside $\mathcal{P}$, $j_\mathcal{U} : U \hookrightarrow \mathcal{X}$ for the open embedding of $\mathcal{U}$ in $\mathcal{X}$ and $i_\mathcal{X} : \mathcal{X} \to \mathcal{P}$ for the closed embedding of the closure of $\mathcal{X}$ in $\mathcal{P}$. Consider the abelian category $\mathrm{Perv}(\mathcal{P})$ of perverse sheaves on $\mathcal{P}$ (with respect to middle perversity). For a reference about the definition and basic properties of perverse sheaves, see [Dim04]. Recall that the simple objects in $\mathrm{Perv}(\mathcal{P})$ are the objects $(i_\mathcal{X})_!IC(\mathcal{X}, \mathcal{L})$ where $\mathcal{L}$ is an irreducible local system on $\mathcal{U}$ and $IC(\mathcal{X}, \mathcal{L})$ is the intersection complex of $\mathcal{X}$ with coefficient in $\mathcal{L}$, that is the image of the morphism ${}^p\mathcal{H}^0((j_\mathcal{U})_!\mathcal{L}) \to {}^p\mathcal{H}^0((Rj_\mathcal{U})_*\mathcal{L})$ in $Perv(\mathcal{X})$. We will denote the corresponding $\mathcal{D}$-module on $\mathcal{P}$ by $\mathcal{M}^{IC}(\mathcal{X}, \mathcal{L})$. If $\mathcal{L}$ is the constant sheaf $\mathbb{C}_\mathcal{U}$ we will simply write $\mathcal{M}^{IC}(\mathcal{X})$. The $p$-th intersection cohomology group of $\mathcal{X}$ (see [GM83]) is denoted by $IH^p(\mathcal{X})$ and is obtained from the intersection complex by the formula $IH^p(\mathcal{X}) = \mathbb{H}^{p-\dim(\mathcal{X})}(IC(\mathcal{X}, \mathbb{C}_\mathcal{U}))$.

We will apply this formalism to the special situation where $\mathcal{U} = g(S)$ (where $g$ is the embedding defined by formula (9)), $\mathcal{X} = X$ and $\mathcal{P} = \mathbb{P}(V')$. The module $\mathcal{M}^{IC}(X)$ is the image of the morphism $g_\dagger\mathcal{O}_S \to g_+\mathcal{O}_S$. In the next result, we will compute the Radon transformation of this module.

**Proposition 2.13.** *In the above situation, we have the following (non-canonical) isomorphism of $\mathcal{D}_V$-modules*
$$\mathcal{H}^0 \mathscr{R}(\mathcal{M}^{IC}(X)) \simeq \mathcal{M}^{IC}(X^\circ, \mathcal{L}) \oplus (IH^{s-1}(X) \otimes \mathcal{O}_V),$$

*and*
$$\mathcal{H}^i \mathscr{R}(\mathcal{M}^{IC}(X)) \simeq IH^{i+s+1}(X) \otimes \mathcal{O}_V \quad \text{for } i > 0,$$
$$\mathcal{H}^i \mathscr{R}(\mathcal{M}^{IC}(X)) \simeq IH^{i+s-1}(X) \otimes \mathcal{O}_V \quad \text{for } i < 0,$$

*where $X^\circ$ is some subvariety of $V$ and $\mathcal{L}$ some local system on some smooth open subset of $X^\circ$.*

*Proof.* Using the comparison isomorphism between the Radon transformation and the Fourier-Laplace transformation (equation (19)) from above, we have
$$\begin{aligned}\mathcal{H}^i \mathscr{R}(\mathcal{M}^{IC}(X)) &\simeq \mathcal{H}^i \operatorname{FL}(p_+ q^+ \mathcal{M}^{IC}(X)) \\ &\simeq \operatorname{FL} \mathcal{H}^i(p_+ q^+ \mathcal{M}^{IC}(X)) \\ &\simeq \operatorname{FL} \mathcal{H}^i(p_+ \mathcal{M}^{IC}(q^{-1}(X))),\end{aligned}$$

where the second isomorphism follows from the exactness of FL and the last isomorphism follows from the smoothness of $q$. We now apply the decomposition theorem [Sai88, corollaire 3, equation 0.12] which gives
$$\mathcal{H}^i(p_+ \mathcal{M}^{IC}(q^{-1}(X))) \simeq \bigoplus_k \mathcal{M}^{IC}(Y_k^i, \mathcal{L}_k^i) \tag{22}$$

for some subvarieties $Y_k^i \subset V'$ and some local systems $\mathcal{L}_k^i$ on a Zariski open subset of $Y_k^i$. Notice that
$$\begin{aligned}j_0^+ \mathcal{H}^i(p_+ \mathcal{M}^{IC}(q^{-1}(X))) &\simeq j_0^+ \mathcal{H}^i(p_+ q^+ \mathcal{M}^{IC}(X)) \\ &\simeq \mathcal{H}^i(j_0^+ p_+ q^+ \mathcal{M}^{IC}(X)) \\ &\simeq \mathcal{H}^i(\pi^+ \mathcal{M}^{IC}(X)) \\ &\simeq \mathcal{H}^i(\mathcal{M}^{IC}(\pi^{-1}(X))),\end{aligned}$$

which is equal to 0 for $i \neq 0$ and equal to $\mathcal{M}^{IC}(Y \setminus \{0\})$ for $i = 0$ (recall from subsection 2.2, more precisely, from the discussion before Lemma 2.10, that $Y$ is the cone of $X$ in $V'$). Thus the decomposition from (22) becomes
$$\mathcal{H}^0(p_+ \mathcal{M}^{IC}(q^{-1}(X))) \simeq \mathcal{M}^{IC}(Y) \oplus \mathcal{S}_0,$$

resp.
$$\mathcal{H}^i(p_+ \mathcal{M}^{IC}(q^{-1}(X))) \simeq \mathcal{S}_i \quad i \neq 0,$$

where the $\mathcal{S}_i$ are $\mathcal{D}$-modules with support at 0, i.e. $\mathcal{S}_i \simeq i_{0+} S_i$, where the $S_i$ are finite-dimensional vector spaces and $i_0 : \{0\} \to V'$ is the natural embedding. We now use the fact that FL is an equivalence of categories, which means that it transforms simple object to simple objects, so we set
$$\mathcal{M}^{IC}(X^\circ, \mathcal{L}) := \operatorname{FL}(\mathcal{M}^{IC}(Y)). \tag{23}$$

It also transforms $\mathcal{D}$-modules with support at 0 to free $\mathcal{O}$-modules, i.e. $\operatorname{FL}(\mathcal{S}_i) \simeq S_i \otimes \mathcal{O}_V$. In order to show the claim, we have to compute the $S_i$. Recall that we have
$$p_+ q^+ \mathcal{M}^{IC}(X) \simeq \bigoplus_j \mathcal{H}^j(p_+ q^+ \mathcal{M}^{IC}(X))[-j] \simeq \bigoplus_{j \neq 0} \mathcal{S}_j[-j] \oplus \mathcal{S}_0 \oplus \mathcal{M}^{IC}(Y), \tag{24}$$

where the first isomorphism is non-canonical. We compute
$$H^i(a_{V'})_+ p_+(q^+ \mathcal{M}^{IC}(X)) \simeq H^i(a_{\mathbb{P}})_+ q_+(q^+ \mathcal{M}^{IC}(X)) \simeq H^i(a_{\mathbb{P}})_+ \mathcal{M}^{IC}(X)[1] \simeq IH^{i+s+1}(X)$$

(here $a_{V'} : V' \to \{pt\}$ resp. $a_{\mathbb{P}} : \mathbb{P}(V') \to \{pt\}$ are the projections to a point), where the second isomorphism follows from [KS94, Corollary 2.7.7 (iv)] and the Riemann-Hilbert correspondence. For the right hand side of Equation (24) we get

$$H^i(a_{V'})_+ \left(\bigoplus_{j\neq 0} \mathcal{S}_j[-j] \oplus \mathcal{S}_0 \oplus \mathcal{M}^{IC}(Y)\right) \simeq S_i \qquad \text{for } i \geq 0,$$

$$H^i(a_{V'})_+ \left(\bigoplus_{j\neq 0} \mathcal{S}_j[-j] \oplus \mathcal{S}_0 \oplus \mathcal{M}^{IC}(Y)\right) \simeq S_i \oplus IH^{i+s+1}(Y) \simeq S_i \oplus IH_p^{i+s+1}(X) \qquad \text{for } i < 0,$$

where $IH_p^{i+s+1}(X)$ is the primitive part of $IH^{i+s+1}(X)$ and where the last isomorphism follows from [KW06, Chapter 4.10]. Therefore we have

$$S_i \simeq IH^{i+s+1}(X) \quad \text{for } i \geq 0,$$
$$S_i \simeq L(IH^{i+s-1}(X)) \simeq IH^{i+s-1}(X) \quad \text{for } i < 0,$$

where $L : IH^{i+s-1}(X) \to IH^{i+s+1}(X)$ is the Lefschetz operator which is injective for $i \leq 0$. $\square$

In the next proposition we show that at a generic point $\underline{\lambda} \in V$ the Radon transformation $\mathscr{R}(\mathcal{M}^{IC}(X))$ of $\mathcal{M}^{IC}(X)$ measures the intersection cohomology of $X \cap H_{\underline{\lambda}}$, where $H_{\underline{\lambda}}$ is the hyperplane in $\mathbb{P}(V')$ corresponding to $\underline{\lambda}$.

**Proposition 2.14.** *Let $\underline{\lambda}$ be a generic point of $V$ and denote by $i_{\underline{\lambda}} : \{\underline{\lambda}\} \longrightarrow V$ its embedding. We have the following isomorphism*

$$i_{\underline{\lambda}}^+ \mathscr{R}(\mathcal{M}^{IC}(X)) \simeq R\Gamma(X \cap H_{\underline{\lambda}}, IC_{X \cap H_{\underline{\lambda}}}),$$

*in particular*

$$\mathcal{H}^j(i_{\underline{\lambda}}^+ \mathscr{R}(\mathcal{M}^{IC}(X))) \simeq IH^{j+s-1}(X \cap H_{\underline{\lambda}}).$$

*Proof.* Consider the following diagram where all squares are cartesian

$$\begin{array}{ccccc}
X & \xleftarrow{\pi_1^X} & Z_X & \xleftarrow{i_X} & X \cap H_{\underline{\lambda}} \\
{\scriptstyle i}\downarrow & & {\scriptstyle \eta}\downarrow & & {\scriptstyle \eta^H}\downarrow \\
\mathbb{P}(V') & \xleftarrow{\pi_1^Z} & Z & \xleftarrow{i_H} & H_{\underline{\lambda}} \\
& & {\scriptstyle \pi_2^Z}\downarrow & & {\scriptstyle \pi^H}\downarrow \\
& & V & \xleftarrow{i_{\underline{\lambda}}} & \{\underline{\lambda}\}
\end{array}$$

We have

$$\begin{aligned}
DR(i_{\underline{\lambda}}^+ \mathscr{R}(\mathcal{M}^{IC}(X))) &\simeq i_{\underline{\lambda}}^! R\pi_{2*}^Z (\pi_1^Z)^! i_! IC(X)[1] \\
&\simeq i_{\underline{\lambda}}^! R\pi_{2*}^Z R\eta_* \pi_1^{X!} IC(X)[1] \\
&\simeq R\pi_*^H i_H^! R\eta_* \pi_1^{X!} IC(X)[1] \\
&\simeq R\pi_*^H R\eta_*^H i_X^! \pi_1^{X!} IC(X)[1] \\
&\simeq R(\pi^H \circ \eta^H)_* (\pi_1^X \circ i_X)^! IC(X)[1] \\
&\simeq R(\pi^H \circ \eta^H)_* IC(X \cap H_{\underline{\lambda}}) \\
&\simeq R\Gamma(X \cap H_{\underline{\lambda}}, IC(X \cap H_{\underline{\lambda}})),
\end{aligned}$$

where the first isomorphism follows from $DR \circ i_{\underline{\lambda}}^+ = i_{\underline{\lambda}}^! \circ DR[t+1]$ and $DR \circ (\pi_1^Z)^+ \simeq (\pi_1^Z)^! \circ DR[-t]$ (see e.g. [HTT08, Theorem 7.1.1]), the second, third and fourth isomorphism follows from base change (see e.g. [Dim04, Theorem 3.2.13(ii)] and the sixth isomorphism follows from [GM83, Section 5.4.1] (notice that their $IC(X)$ is our $IC(X)[n]$ where $n = \dim_{\mathbb{C}}(X)$) and the fact that for a generic $\underline{\lambda}$ the hyperplane $H_{\underline{\lambda}}$ is transversal to a given Whitney stratification of $X$. The first claim now follows from the fact that the de Rham functor $DR$ is the identity on a point. The second claim follows from $\mathbb{H}^{j-s+1}(X \cap H_{\underline{\lambda}}, IC(X \cap H_{\underline{\lambda}})) \simeq IH^j(X \cap H_{\underline{\lambda}})$. $\square$

**Remark 2.15.** *Combining Proposition 2.13 and Proposition 2.14 we see that we have the following decomposition for generic $\underline{\lambda} \in V$:*

$$IH^{s-1}(X \cap H_{\underline{\lambda}}) \simeq \mathcal{H}^0(i_{\underline{\lambda}}^+ \mathscr{R}(\mathcal{M}^{IC}(X))) \simeq i_{\underline{\lambda}}^+ \mathcal{H}^0(\mathscr{R}(\mathcal{M}^{IC}(X))) \simeq i_{\underline{\lambda}}^+ \mathcal{M}^{IC}(X^\circ, \mathcal{L}) \oplus IH^{s-1}(X).$$

*This is the intersection cohomology analogon of the decomposition of the cohomology of a smooth hyperplane section of a smooth projective variety into its vanishing part and the ambient part.*

We will now show that $\mathcal{M}^{IC}(X^\circ, \mathcal{L})$ can expressed as an image of a morphism between GKZ-systems.

**Theorem 2.16.** *Let $\widetilde{\beta}, \widetilde{\beta}'$ be as in Theorem 2.11, then $\mathcal{M}^{IC}(X^\circ, \mathcal{L}) \simeq im(\mathcal{M}_{\widetilde{B}}^{-\widetilde{\beta}'} \stackrel{\cdot \partial^{\widetilde{\beta}+\widetilde{\beta}'}}{\longrightarrow} \mathcal{M}_{\widetilde{B}}^{\widetilde{\beta}})$.*

*Proof.* First recall that we have shown in the proof of Proposition 2.13. that $\mathcal{M}^{IC}(X^\circ, \mathcal{L}) \simeq \text{FL}(\mathcal{M}^{IC}(Y))$. On the other hand, as $Y$ is the closure in $V'$ of the image of the morphism $h$, the module $\mathcal{M}^{IC}(Y)$ is isomorphic to the image of $h_\dagger \mathcal{O}_T \to h_+ \mathcal{O}_T$. As the Fourier-Laplace transformation is exact we can conclude that $\mathcal{M}^{IC}(X^\circ, \mathcal{L})$ is isomorphic to the image of $\text{FL}(h_\dagger \mathcal{O}_T) \to \text{FL}(h_+ \mathcal{O}_T)$.

By Lemma 2.10 we know that $\text{FL}(h_+ \mathcal{O}_T)$ is isomorphic to $\mathcal{M}_{\widetilde{B}}^{\widetilde{\beta}}$ for every $\widetilde{\beta} \in \delta_{\widetilde{B}} + (\mathbb{R}_{\geq 0} \widetilde{B} \cap \mathbb{Z}^{s+1})$ and that $\text{FL}(h_\dagger \mathcal{O}_T)$ is isomorphic to $\mathcal{M}_{\widetilde{B}}^{-\widetilde{\beta}'}$ for every $\widetilde{\beta}' \in (\mathbb{R}_{\geq 0} \widetilde{B})^\circ \cap \mathbb{Z}^{s+1}$. It follows now from the last statement of Lemma 2.10, that the induced morphism between $\mathcal{M}_{\widetilde{B}}^{-\widetilde{\beta}'}$ and $\mathcal{M}_{\widetilde{B}}^{\widetilde{\beta}}$ is equal to $\cdot \partial^{\widetilde{\beta}+\widetilde{\beta}'}$ up to some non-zero constant. □

In general it is quite difficult to make any precise statement on the variety $X^\circ$ and the local system $\mathcal{L}$ which define the module $\mathcal{M}^{IC}(X^\circ, \mathcal{L})$. Nevertheless, if we restrict our attention to the situation where the matrix $B$ defining the embedding $g : S \hookrightarrow \mathbb{P}(V')$ is given by the primitive integral generators of a toric manifold which is given by a total bundle $\mathbb{V}(\mathcal{E}^\vee) \twoheadrightarrow X_\Sigma$, where $\mathcal{E}$ is a split convex vector bundle over another toric manifold $X_\Sigma$ such that the zero locus of a generic section is a nef complete intersection (i.e., the situation considered from section 4 on, see also the introduction in section 1), then we expect that $X^\circ = V$. In order to show this, one would need to prove that if we restrict the morphism $\mathcal{M}_{\widetilde{B}}^{-\widetilde{\beta}'} \stackrel{\cdot \partial^{\widetilde{\beta}+\widetilde{\beta}'}}{\longrightarrow} \mathcal{M}_{\widetilde{B}}^{\widetilde{\beta}}$ from the last theorem to a generic point of $V$, then it is not the zero map. It is well known (see, e.g., [Ado94]) that the restriction of $\mathcal{M}_{\widetilde{B}}^{\widetilde{\beta}}$ to a generic point is the quotient of $\mathbb{C}[\mathbb{N}\widetilde{B}]$ by the ideal generated by the Euler vector fields $\sum_{i=0}^{t} \widetilde{b}_{ki} \lambda_i \partial_{\lambda_i}$ $(k = 0, \ldots, s)$, where now $(\lambda_0, \ldots, \lambda_t)$ are the components of the generic point we restrict to. Hence one needs to show that the monomial $\partial^{\widetilde{\beta}+\widetilde{\beta}'}$ does not lie in this ideal. Nevertheless, we do not have, at this moment, any further evidence for this to be true.

Even under the above restrictive assumptions on $B$ and even if we suppose that $X^\circ = V$, it is not easy to predict the rank of $\mathcal{L}$. What we expect is that the generic rank of the module $\widehat{\mathcal{M}}^{IC}(X^\circ, \mathcal{L})$ from theorem 3.6 below can be identified with the dimension of the image of the map

$$H^*(X_\Sigma, \mathbb{C}) \stackrel{\cup c_{top}(\mathcal{E})}{\longrightarrow} H^*(X_\Sigma, \mathbb{C}).$$

However, the module $\mathcal{M}^{IC}(X^\circ, \mathcal{L})$ resp. the local system $\mathcal{L}$ may contain constant subobjects that vanish after localized Fourier-Laplace transformation (see section 3 below for more details). Hence its rank may be different from that of $\widehat{\mathcal{M}}^{IC}(X^\circ, \mathcal{L})$.

For applications like those in the last section, we need a description of $\mathcal{M}^{IC}(X^\circ, \mathcal{L})$ as a quotient of a GKZ-system, rather than submodule of it. For this purpose, denote by $\mathcal{K}$ the kernel of the morphism $\mathcal{M}_{\widetilde{B}}^{-\widetilde{\beta}'} \stackrel{\cdot \partial^{\widetilde{\beta}+\widetilde{\beta}'}}{\longrightarrow} \mathcal{M}_{\widetilde{B}}^{\widetilde{\beta}}$, then $\mathcal{M}^{IC}(X^\circ, \mathcal{L})$ is isomorphic to the quotient $\mathcal{M}_{\widetilde{B}}^{-\widetilde{\beta}'}/\mathcal{K}$ in the abelian category of regular holonomic $\mathcal{D}$-modules. The next result gives a concrete description of $\mathcal{K}$ as a submodule of $\mathcal{M}_{\widetilde{B}}^{-\widetilde{\beta}'}$.

First, we define a sub-$D_V$-module $\Gamma_{\partial,c}(M_{\widetilde{B}}^{-\widetilde{\beta}'})$ of $M_{\widetilde{B}}^{-\widetilde{\beta}'}$, where $c \in \rho^{-1}(\widetilde{\beta} + \widetilde{\beta}')$ (cf. Equation (13)) :

$$\Gamma_{\partial,c}(M_{\widetilde{B}}^{-\widetilde{\beta}'}) := \{m \in M_{\widetilde{B}}^{-\widetilde{\beta}'} \mid \exists n \in \mathbb{N} \text{ with } (\partial^c)^n \cdot m = 0\}$$

Recall that two elements $\partial^{c_1}$ and $\partial^{c_2}$ with $c_1, c_2 \in \rho^{-1}(\widetilde{\beta} + \widetilde{\beta}')$ differ by some element $P \cdot \Box_{\underline{l}}$, where $P \in \mathbb{C}[\partial_0, \ldots, \partial_s]$ and $\underline{l} = c_1 - c_2$. Any element $m \in M_{\widetilde{B}}^{-\widetilde{\beta}'}$ is eliminated by left multiplication with some high enough power of $P \cdot \square_{\underline{l}}$. This shows that $\Gamma_{\partial,c}(M_{\widetilde{B}}^{-\widetilde{\beta}'})$ is actually independent of the chosen $c \in \rho^{-1}(\widetilde{\beta} + \widetilde{\beta}')$. Thus we denote it just by $\Gamma_{\partial}(M_{\widetilde{B}}^{-\widetilde{\beta}'})$ and the corresponding sub-$\mathcal{D}_V$-module of $\mathcal{M}_{\widetilde{B}}^{-\widetilde{\beta}'}$ by $\Gamma_{\partial}(\mathcal{M}_{\widetilde{B}}^{-\widetilde{\beta}'})$.

**Proposition 2.17.** *Let $\widetilde{\beta}, \widetilde{\beta}'$ be as in Theorem 2.11 and let $\mathcal{K}$ be the kernel of $\mathcal{M}_{\widetilde{B}}^{-\widetilde{\beta}'} \overset{\cdot \partial^{\widetilde{\beta}+\widetilde{\beta}'}}{\longrightarrow} \mathcal{M}_{\widetilde{B}}^{\widetilde{\beta}}$. Then*
$$\mathcal{K} \simeq \Gamma_{\partial}(\mathcal{M}_{\widetilde{B}}^{-\widetilde{\beta}'}), \quad \text{in particular} \quad \mathcal{M}^{IC}(X^{\circ}, \mathcal{L}) \simeq \mathcal{M}_{\widetilde{B}}^{-\widetilde{\beta}'}/\Gamma_{\partial}(\mathcal{M}_{\widetilde{B}}^{-\widetilde{\beta}'}).$$

*Proof.* Recall that the morphism $\mathcal{M}_{\widetilde{B}}^{-\widetilde{\beta}'} \overset{\cdot \partial^{\widetilde{\beta}+\widetilde{\beta}'}}{\longrightarrow} \mathcal{M}_{\widetilde{B}}^{\widetilde{\beta}}$ is induced by the morphism $\mathrm{FL}(h_{\dagger}\mathcal{O}_T) \to \mathrm{FL}(h_+\mathcal{O}_T)$, where we used the isomorphisms $\mathcal{M}_{\widetilde{B}}^{-\widetilde{\beta}'} \simeq \mathrm{FL}(h_{\dagger}\mathcal{O}_T)$ and $\mathcal{M}_{\widetilde{B}}^{\widetilde{\beta}} \simeq \mathrm{FL}(h_+\mathcal{O}_T)$. Applying the Fourier-Laplace transformation again and using $\mathrm{FL} \circ \mathrm{FL} = Id$, we see that the morphism $\mathrm{FL}(\mathcal{M}_{\widetilde{B}}^{-\widetilde{\beta}'}) \overset{\cdot w^{\widetilde{\beta}+\widetilde{\beta}'}}{\longrightarrow} \mathrm{FL}(\mathcal{M}_{\widetilde{B}}^{\widetilde{\beta}})$ is induced by the morphism $h_{\dagger}\mathcal{O}_T \longrightarrow h_+\mathcal{O}_T$. We will calculate the kernel of $\mathrm{FL}(\mathcal{M}_{\widetilde{B}}^{-\widetilde{\beta}'}) \overset{\cdot w^{\widetilde{\beta}+\widetilde{\beta}'}}{\longrightarrow} \mathrm{FL}(\mathcal{M}_{\widetilde{B}}^{\widetilde{\beta}})$. First notice that the map $h$ can be factorized as $h = k \circ l$, where $k$ is the canonical inclusion of $(\mathbb{C}^*)^{t+1} \to V'$ and the map $l$ is given by

$$l : T \longrightarrow (\mathbb{C}^*)^{t+1},$$
$$(y_0, \ldots, y_r) \longmapsto (\underline{y}^{\widetilde{b}_0}, \ldots, \underline{y}^{\widetilde{b}_t}) = (y_0, y_0\underline{y}^{b_1}, \ldots, y_0\underline{y}^{b_t}).$$

This shows that $\mathrm{FL}(\mathcal{M}_{\widetilde{B}}^{\widetilde{\beta}}) \simeq k_+ l_+ \mathcal{O}_T$ is localized along $V' \setminus (\mathbb{C}^*)^{t+1}$, i.e. $\mathrm{FL}(\mathcal{M}_{\widetilde{B}}^{\widetilde{\beta}}) \simeq k_+ k^+ \mathrm{FL}(\mathcal{M}_{\widetilde{B}}^{\widetilde{\beta}})$. Let $D_1 = \{w^{\widetilde{\beta}+\widetilde{\beta}'} = 0\}$, set $U_1 := V' \setminus D_1$ and denote by $j_1 : U_1 \to V'$ the canonical inclusion. Because $(\mathbb{C}^*)^{t+1} \subset U_1$, the $\mathcal{D}$-module $\mathrm{FL}(\mathcal{M}_{\widetilde{B}}^{\widetilde{\beta}})$ is also localized along $D_1$, i.e, $\mathrm{FL}(\mathcal{M}_{\widetilde{B}}^{\widetilde{\beta}}) \simeq j_{1+}j_1^+ \mathrm{FL}(\mathcal{M}_{\widetilde{B}}^{\widetilde{\beta}})$. Notice that the induced morphism $j_1^+ \mathrm{FL}(\mathcal{M}_{\widetilde{B}}^{-\widetilde{\beta}'}) \to j_1^+ \mathrm{FL}(\mathcal{M}_{\widetilde{B}}^{\widetilde{\beta}})$ is an isomorphism, because $w^{\widetilde{\beta}+\widetilde{\beta}'}$ is invertible on $U_1$. Therefore we can conclude that $j_{1+}j_1^+ \mathrm{FL}(\mathcal{M}_{\widetilde{B}}^{-\widetilde{\beta}'}) \to j_{1+}j_1^+ \mathrm{FL}(\mathcal{M}_{\widetilde{B}}^{\widetilde{\beta}}) \simeq \mathrm{FL}(\mathcal{M}_{\widetilde{B}}^{\widetilde{\beta}})$ is an isomorphism. It is therefore enough to calculate the kernel of $\mathrm{FL}(\mathcal{M}_{\widetilde{B}}^{-\widetilde{\beta}'}) \to j_{1+}j_1^+ \mathrm{FL}(\mathcal{M}_{\widetilde{B}}^{-\widetilde{\beta}'})$. On the level of global sections this is $H^0_{D_1}(\mathrm{FL}(M_{\widetilde{B}}^{-\widetilde{\beta}'}))$ (cf. [HTT08, Proposition 1.7.1]) which is given by

$$H^0_{D_1}(\mathrm{FL}(M_{\widetilde{B}}^{-\widetilde{\beta}'})) = \{m \in \mathrm{FL}(M_{\widetilde{B}}^{-\widetilde{\beta}'}) \mid \exists n \in \mathbb{N} \text{ with } (w^{\widetilde{\beta}+\widetilde{\beta}'})^n \cdot m = 0\}.$$

Applying the Fourier-Laplace transformation to this kernel shows the claim. □

### 2.4 The equivariant setting

In this section we show that the $\mathcal{D}$-modules discussed above are quasi-equivariant with respect to a natural torus action. We review the definition of an quasi-equivariant $\mathcal{D}$-modules from [Kas08, Chapter 3] and prove some simple statements for these.

Let $\mathcal{X}$ be a smooth, complex, quasi-projective variety and $G$ be a complex affine algebraic group, which acts on $\mathcal{X}$. Denote by $\nu : G \times \mathcal{X} \to \mathcal{X}$ the action of $G$ on $\mathcal{X}$ and by $p_2 : G \times \mathcal{X} \to \mathcal{X}$ the second projection. A $\mathcal{D}_{\mathcal{X}}$-module $\mathcal{M}$ is called quasi-$G$-equivariant if it satisfies $\nu^+\mathcal{M} \simeq p_2^+\mathcal{M}$ as $\mathcal{O}_G \boxtimes \mathcal{D}_{\mathcal{X}}$-modules together with an associative law (cf. [Kas08, Definition 3.1.3]). We denote the abelian category of quasi-$G$-equivariant $\mathcal{D}_{\mathcal{X}}$-modules by $M(\mathcal{D}_{\mathcal{X}}, G)$ and the subcategories of coherent, holonomic and regular holonomic quasi-$G$-equivariant $\mathcal{D}_{\mathcal{X}}$-modules by $M_{coh}(\mathcal{D}_{\mathcal{X}}, G)$ resp. $M_h(\mathcal{D}_{\mathcal{X}}, G)$ resp. $M_{rh}(\mathcal{D}_{\mathcal{X}}, G)$. The corresponding bounded derived categories are denoted by $D^b_*(\mathcal{D}_{\mathcal{X}}, G)$ for $* = \emptyset, coh, h, rh$.

A $\mathcal{O}_{\mathcal{X}}$-module $\mathcal{F}$ is called $G$-equivariant if $\nu^*\mathcal{F} \simeq pr^*\mathcal{F}$ as $\mathcal{O}_{G \times \mathcal{X}}$-modules and if it satisfies an associative law (expressed as the commutativity of a certain diagram, see [Kas08, Definition 3.1.2]). We denote by $Mod(\mathcal{O}_{\mathcal{X}}, G)$ the category of $G$-equivariant $\mathcal{O}_{\mathcal{X}}$-modules and by $Mod_{coh}(\mathcal{O}_{\mathcal{X}}, G)$ the subcategory of coherent $G$-equivariant $\mathcal{O}_{\mathcal{X}}$-modules.

Let $f : \mathcal{X} \to \mathcal{Y}$ be a $G$-equivariant map. Then the direct image resp. the inverse image functors preserve quasi-$G$-equivariance (cf. [Kas08, Equation (3.4.1), Equation (3.5.2)]).

We will now show that the duality functor preserves quasi-$G$-equivariance.

**Proposition 2.18.** *Let $M \in D^b_{coh}(\mathcal{D}_\mathcal{X}, G)$ then $\mathbb{D}(M) \in D^b_{coh}(\mathcal{D}_\mathcal{X}, G)^{opp}$.*

*Proof.* By a dévissage we may assume that $M$ is a single degree complex, i.e. $M \in Mod_{coh}(\mathcal{D}_\mathcal{X}, G)$. By [Kas08, Lemma 3.3.2] for every $N \in Mod_{coh}(\mathcal{O}_\mathcal{X}, G)$ there exists a $G$-equivariant locally-free $\mathcal{O}_\mathcal{X}$ module $L$ of finite rank and a surjective $G$-equivariant morphism $L \twoheadrightarrow N$. Notice that there exists a $G$-equivariant coherent $\mathcal{O}_\mathcal{X}$-submodule $K$ of $M$ with $\mathcal{D}_\mathcal{X} \otimes K = M$. This enables us to construct a locally-free, $G$-equivariant resolution

$$\cdots \to L_2 \to L_1 \to L_0 \to K \to 0$$

of $K$ in $Mod_{coh}(\mathcal{O}_\mathcal{X}, G)$, which gives rise to a resolution of $M$

$$\cdots \to \mathcal{D}_\mathcal{X} \otimes L_2 \to \mathcal{D}_\mathcal{X} \otimes L_1 \to \mathcal{D}_\mathcal{X} \otimes L_0 \to M \to 0$$

in $Mod_{coh}(\mathcal{D}_\mathcal{X}, G)$ by the exactness of $\mathcal{D}_\mathcal{X} \otimes_{\mathcal{O}_\mathcal{X}}$. We have

$$\begin{aligned}\mathbb{D}M &= R\mathcal{H}om_{\mathcal{D}_\mathcal{X}}(M, \mathcal{D}_\mathcal{X}) \otimes \Omega_\mathcal{X}^{\otimes -1}[dim\mathcal{X}] \\ &\simeq \mathcal{H}om_{\mathcal{D}_\mathcal{X}}(\mathcal{D}_\mathcal{X} \otimes L_\bullet, \mathcal{D}_\mathcal{X}) \otimes \Omega_\mathcal{X}^{\otimes -1}[dim\mathcal{X}] \\ &\simeq (\mathcal{H}om_{\mathcal{O}_\mathcal{X}}(L_\bullet, \mathcal{O}_\mathcal{X}) \otimes \mathcal{D}_\mathcal{X}) \otimes \Omega_\mathcal{X}^{\otimes -1}[dim\mathcal{X}] \\ &\simeq \mathcal{D}_\mathcal{X} \otimes \mathcal{H}om_{\mathcal{O}_\mathcal{X}}(L_\bullet, \mathcal{O}_\mathcal{X})[dim\mathcal{X}].\end{aligned}$$

But $\mathcal{H}om_{\mathcal{O}_\mathcal{X}}(L_\bullet, \mathcal{O}_\mathcal{X})$ is again a complex in $Mod_{coh}(\mathcal{O}_\mathcal{X}, G)$, which can be easily seen by the local-freeness of the $L_i$. Thus we can conclude that $\mathbb{D}M \in D^b_{coh}(\mathcal{D}_\mathcal{X}, G)^{opp}$. $\square$

**Corollary 2.19.** *Let $f : \mathcal{X} \to \mathcal{Y}$ be a $G$-equivariant map. Then the proper direct image and the exceptional inverse image functor preserves quasi-$G$-equivariance.*

*Proof.* This follows from $f_\dagger = \mathbb{D} \circ f_+ \circ \mathbb{D}$ and $f^\dagger = \mathbb{D} \circ f^+ \circ \mathbb{D}$. $\square$

In the next proposition we will show that the characteristic variety of a quasi-$G$-equivariant $\mathcal{D}$-module is $G$-invariant. For that purpose, we will consider the action induced by $\nu$ on the cotangent bundle $T^*\mathcal{X}$. More precisely, consider the differential $d\nu$ of the action map, which is a map of vector bundles $d\nu : \nu^* T^*\mathcal{X} \to T^*(G \times \mathcal{X}) = T^*G \boxtimes T^*\mathcal{X}$ over $G \times \mathcal{X}$, or, equivalently, a map $d\nu : (G \times \mathcal{X}) \times_\mathcal{X} T^*\mathcal{X} \to T^*G \times T^*\mathcal{X}$ of smooth complex varieties. Notice that

$$\begin{aligned} t: G \times T^*\mathcal{X} &\longrightarrow (G \times \mathcal{X}) \times_\mathcal{X} T^*\mathcal{X} = \{((g,x),v) \mid \pi(v) = \nu(g,x) \in \mathcal{X}\}, \\ (g,v) &\longmapsto (g, \nu(g^{-1}, \pi(v)), v)\end{aligned}$$

is an isomorphism, with inverse map sending $((g,x),v)$ to $(g,v)$. Now consider the composition $\xi : \widetilde{p}_2 \circ d\nu \circ t : G \times T^*\mathcal{X} \to T^*\mathcal{X}$, where $\widetilde{p}_2 : T^*G \times T^*\mathcal{X} \to T^*\mathcal{X}$ is the second projection. One easily checks that we have $\xi(g_1 \cdot g_2, x) = \xi(g_1, \xi(g_2, x))$, i.e., that we obtain an action of $G$ on $T^*\mathcal{X}$. Notice that for any $g \in G$, the map $\xi(g, -) : T^*\mathcal{X} \to T^*\mathcal{X}$ is nothing but the differential $d\nu_g$ of the map $\nu_g : \mathcal{X} \to \mathcal{X}$ where $\nu_g(x) := \nu(g, x)$. Notice that for $M \in D^b(\mathcal{D}_\mathcal{X}, G)$ one has $\nu_g^+ M \simeq M$ by the quasi-$G$-equivariance of $M$.

**Proposition 2.20.** *Let $M \in D^b_{coh}(\mathcal{D}_\mathcal{X}, G)$, then the characteristic variety $\mathrm{char}(M)$ of $M$ is invariant under the $G$-action on $T^*\mathcal{X}$ given by $\xi$. Moreover, if $G$ is irreducible then the irreducible components of $\mathrm{char}(M)$ are also $G$-invariant.*

*Proof.* For both statements it is sufficient to show invariance under the morphism $\nu_g$ for any $g \in G$. We are going to use the following fact (cf. [HTT08, Lemma 2.4.6(iii)]). Let $f : \mathcal{X} \to \mathcal{Y}$ be a morphism between smooth algebraic varieties. One has the natural morphisms

$$T^*\mathcal{X} \xleftarrow{\rho_f} \mathcal{X} \times_\mathcal{Y} T^*\mathcal{Y} \xrightarrow{\omega_f} T^*\mathcal{Y}.$$

Let $M \in Mod_{coh}(\mathcal{D}_\mathcal{Y})$. If $f$ is non-characteristic then $\operatorname{char}(f^+M) \subset \rho_f \omega_f^{-1}(\operatorname{char}(M))$.

We want to apply this to the case $f = \nu_g$. Notice that in this case the maps $\rho_{\nu_g}$ and $\omega_{\nu_g}$ are isomorphisms and $\rho_{\nu_g} \circ \omega_{\nu_g}^{-1} = d\nu_g$. Thus we have

$$\operatorname{char}(M) = \operatorname{char}(\nu_g^+ M) \subset d\nu_g(\operatorname{char}(M)).$$

Repeating the argument with $\nu_{g^{-1}}$ gives $\operatorname{char}(M) \subset d\nu_{g^{-1}}(\operatorname{char}(M))$. Now applying $d\nu_g$ to both sides of the latter inclusion shows the first claim.

Now assume that $G$ is irreducible and let $C_i$ be an irreducible component of $Ch(M)$. Notice that $G \times C_i$ is irreducible. Consider the scheme-theoretic image $I$ of $G \times C_i$ under the induced action map $\xi : G \times \operatorname{char}(M) \to \operatorname{char}(M)$. Then $\overline{\xi} : G \times C_i \to I$ is a dominant morphism. We want to show that $I$ is irreducible. Let $U \subset I$ be an affine open set. The restriction $\overline{\xi}^{-1}(U) \to U$ is still dominant and induces an injective ring homomorphism $\mathcal{O}_I(U) \to \mathcal{O}_{G \times C_i}(\overline{\xi}^{-1}(U))$. As $G \times C_i$ is irreducible and reduced the ring $\mathcal{O}_{G \times C_i}(\overline{\xi}^{-1}(U))$ is a domain. Thus $\mathcal{O}_I(U)$ is also a domain and because $U$ was chosen arbitrary we conclude that $I$ is irreducible. Notice that we have $C_i \subset I \subset \operatorname{char}(M)$ and therefore $C_i = I$, which shows the claim. $\square$

The proposition above enables us to prove that a section of a quotient map of a free action is non-characteristic with respect to quasi-$G$-equivariant $\mathcal{D}$-modules.

**Lemma 2.21.** *Let $G \times \mathcal{X} \to \mathcal{X}$ be a free action and $\pi_G : \mathcal{X} \to \mathcal{X}/G$ a geometric quotient. Let $i_G : \mathcal{X}/G \to \mathcal{X}$ be a section of $\pi_G$, then $i_G$ is non-characteristic with respect to every $M \in D^b_{rh}(G, \mathcal{D}_\mathcal{X})$.*

*Proof.* We consider $\mathcal{X}/G$ as smooth subvariety of $\mathcal{X}$. Notice that $\mathcal{X}/G$ is transversal to the orbits of the $G$-action on $\mathcal{X}$ given by $\nu$. Let $\operatorname{char}(M) = \bigcup_{i \in I} C_i$ be the decomposition into irreducible components and put $\mathcal{X}_i := \pi(X_i)$ so that $C_i = T^*_{\mathcal{X}_i}\mathcal{X}$. From Proposition 2.20 we know that $C_i$ is invariant under the action given by $\xi$, and hence a union of orbits of this $G$-action. On the other hand, the image under the projection $\pi : T^*\mathcal{X} \to \mathcal{X}$ of such an orbit is necessarily an orbit of the original action given by $\nu$. Hence $\mathcal{X}_i$ is a union $\bigcup_j \mathcal{X}_i^{(j)}$ of $G$-orbits, more precisely, these orbits form a Whitney stratification of $\mathcal{X}_i$ (see, [Dim92, Proposition 1.14]). Whitney's condition A then implies that $T^*_{\mathcal{X}_i}\mathcal{X} \subset \bigcup_j T^*_{\mathcal{X}_i^{(j)}}\mathcal{X}_i$. Transversality of $\mathcal{X}/G$ and the orbits $\mathcal{X}_i^{(j)}$ means that $T^*_{\mathcal{X}/G}\mathcal{X} \cap T^*_{\mathcal{X}_i^{(j)}}\mathcal{X}_i \subset T^*_\mathcal{X}\mathcal{X}$, from which we deduce that $T^*_{\mathcal{X}/G}\mathcal{X} \cap T^*_{\mathcal{X}_i}\mathcal{X}_i \subset T^*_\mathcal{X}\mathcal{X}$ and hence $T^*_{\mathcal{X}/G}\mathcal{X} \cap \operatorname{char}(M) \subset T^*_\mathcal{X}\mathcal{X}$. Thus $i_G$ non-characteristic with respect to $M$ as required. $\square$

Let $V^* = \mathbb{C} \times (\mathbb{C}^*)^t$ and let $j_{V^*} : V^* \to V$ be the canonical embedding. Consider the following diagram

$$\begin{array}{ccccc}
S & \xleftarrow{\pi_1^S} & \Gamma & \xleftarrow{j_{\Gamma^*}} & \Gamma^* \\
{\scriptstyle j}\downarrow & & {\scriptstyle \theta}\downarrow & & {\scriptstyle \zeta}\downarrow \\
X & \longleftarrow & Z_X & \xleftarrow{j_{Z_X^*}} & Z_X^* \\
{\scriptstyle i}\downarrow & & {\scriptstyle \eta}\downarrow & & {\scriptstyle \varepsilon}\downarrow \\
\mathbb{P}(V') & \xleftarrow{\pi_1^Z} & Z & \xleftarrow{j_{Z^*}} & Z^* \\
& & {\scriptstyle \pi_Z^2}\downarrow & & {\scriptstyle \delta}\downarrow \\
& & V & \xleftarrow{j_{V^*}} & V^*
\end{array} \qquad (25)$$

where the varieties $Z^*, Z_X^*, \Gamma^*$ together with the maps $j_{Z^*}, j_{Z_X^*}, j_{\Gamma^*}$ and $\delta, \varepsilon, \zeta$ are induced by the base change $j_{V^*}$. Thus all squares in the diagram above are cartesian.

We now specify to the case $G = (\mathbb{C}^*)^s$. We let $G$ act on $S$ and $V$ by

$$G \times S \longrightarrow S, \qquad (26)$$
$$(g_1, \ldots, g_s, y_1, \ldots, y_s) \mapsto (g_1 y_1, \ldots g_s y_s),$$
$$G \times V \longrightarrow V,$$
$$(g_1, \ldots, g_s, \lambda_0, \ldots, \lambda_t) \mapsto (\lambda_0, \underline{g}^{-\underline{b}_1} \lambda_1, \ldots, \underline{g}^{-\underline{b}_t} \lambda_t).$$

We also define the following $G$-action on $\mathbb{P}(V')$:

$$G \times \mathbb{P}(V') \longrightarrow \mathbb{P}(V'), \qquad (27)$$
$$(g_1, \ldots, g_s, (\mu_0 : \ldots : \mu_t)) \mapsto (\mu_0 : \underline{g}^{\underline{b}_1} \mu_1 : \ldots : \underline{g}^{\underline{b}_t} \mu_t).$$

This makes map $g = i \circ j : S \to \mathbb{P}(V')$ $G$-equivariant. There is a natural action of $G$ on $\mathbb{P}(V') \times V$ resp. $S \times V$ which leaves the subvarieties $Z = \{\sum_{i=0}^{t} \lambda_i \mu_i = 0\}$ resp. $\Gamma = \{\lambda_0 + \sum_{i=1}^{t} \lambda_i y^{\underline{b}_i} = 0\}$ invariant. It is now easy to see, using the induced actions on $\Gamma$ resp. $Z$, that the maps $\pi_1^Z, \pi_2^Z, \pi_1^S$ as well as $\eta$ and $\theta$ are $G$-equivariant.

Notice that $G$ leaves $V^*$ invariant and acts freely on it, but this shows that $G$ acts also freely on $Z^*$, $Z_X^*$ and $\Gamma^*$. Therefore also the maps $\delta, \varepsilon, \zeta$ are $G$-equivariant. Notice that the action of $G$ on $\mathbb{P}(V')$ as defined in formula (27) is not free, there are orbits of dimension strictly smaller dimension than $s = \dim(G)$. Because we have $\mathbb{Z} B = \mathbb{Z}^s$, there exist matrices $N_1 \in Gl(s \times s, \mathbb{Z})$ and $N_2 \in Gl(t \times t, \mathbb{Z})$ such that

$$B = N_1 \cdot (I_s \mid 0_{s \times r}) \cdot N_2, \qquad (28)$$

where $r := t - s$. Define matrices

$$L := N_2^{-1} \cdot \begin{pmatrix} 0_{s \times r} \\ I_r \end{pmatrix}, \quad M := (0_{r \times s} \mid I_r) \cdot N_2, \quad C := N_2^{-1} \cdot \begin{pmatrix} I_s \\ 0_{r \times s} \end{pmatrix} \cdot N_1^{-1}, \quad D := (C \cdot B)^t,$$

whose entries we denote by $l_{ij}$, $m_{ji}$, $c_{ik}$ and $d_{il}$, respectively. Then $M \cdot L = I_r$, $B \cdot C = I_s$, $B \cdot L = 0$, $M \cdot C = 0$ and

$$C \cdot B + L \cdot M = I_t. \qquad (29)$$

Consider the following map, where $F := (\mathbb{C}^*)^s$ and $\mathcal{KM} := (\mathbb{C}^*)^s$:

$$T_\mathbb{P} : \mathbb{P}(V') \times \mathbb{C} \times F \times \mathcal{KM} \longrightarrow \mathbb{P}(V') \times V^*,$$
$$((\mu_0 : \ldots : \mu_t), \lambda_0, f_1, \ldots, f_s, q_1, \ldots, q_r) \mapsto ((\mu_0 : \underline{f}^{-\underline{b}_1} \mu_1 : \ldots : \underline{f}^{-\underline{b}_t} \mu_t), \lambda_0, \underline{f}^{\underline{b}_1} \cdot \underline{q}^{\underline{m}_1}, \ldots, \underline{f}^{\underline{b}_t} \cdot \underline{q}^{\underline{m}_t})$$

with $\underline{f}^{\underline{b}_i} = \prod_{k=1}^{s} f_k^{b_{ki}}$, $\underline{q}^{\underline{m}_i} = \prod_{j=1}^{r} q_j^{m_{ji}}$ and inverse

$$T_\mathbb{P}^{-1} : \mathbb{P}(V') \times V^* \longrightarrow \mathbb{P}(V') \times \mathbb{C} \times F \times \mathcal{KM},$$
$$((\mu_0 : \ldots : \mu_t), \lambda_0, \ldots, \lambda_t) \mapsto ((\mu_0 : \underline{\lambda}^{\underline{d}_1} \mu_1 : \ldots : \underline{\lambda}^{\underline{d}_t} \mu_t), \lambda_0, \underline{\lambda}^{\underline{c}_1}, \ldots, \underline{\lambda}^{\underline{c}_s}, \underline{\lambda}^{\underline{l}_1}, \ldots, \underline{\lambda}^{\underline{l}_r})$$

with $\underline{\lambda}^{\underline{c}_k} := \prod_{i=1}^{t} \lambda_i^{c_{ik}}$, $\underline{\lambda}^{\underline{l}_j} = \prod_{i=1}^{t} \lambda_i^{l_{ij}}$ and $\underline{\lambda}^{\underline{d}_l} := \prod_{i=1}^{t} \lambda_i^{d_{il}} = \prod_{i=1}^{t} \lambda_i^{\sum_k c_{ik} b_{kl}}$.

Notice that the space $\mathcal{KM}$ will reappear in section 6 (see the explanations after the exact sequence (59)), where it similarly denotes the $r$-dimensional torus $(\mathbb{C}^*)^r$. There is however a difference: in the present section, our input data is the matrix $B$, and the map $T_\mathbb{P}$ and its inverse $T_\mathbb{P}^{-1}$ are defined by the choice of the matrices $N_1$ and $N_2$ which have to satisfy only equation 28. In section 6, we work with a toric variety (and the matrix $B$ is given by the primitiv integral generators of its rays), and here these choices have to satisfy much finer conditions. Nevertheless, we will use the same symbol in order to avoid overloading the notation too much.

Recall the following $G$-action on $\mathbb{P}(V') \times V^*$

$$G \times (\mathbb{P}(V') \times V^*) \longrightarrow \mathbb{P}(V') \times V^*,$$
$$(g_1, \ldots, g_s, (\mu_0 : \ldots : \mu_t), \lambda_0, \ldots, \lambda_t) \mapsto ((\mu_0 : \underline{g}^{\underline{b}_1} \mu_1 : \ldots : \underline{g}^{\underline{b}_t} \mu_t), \lambda_0, \underline{g}^{-\underline{b}_1} \lambda_1, \ldots, \underline{g}^{-\underline{b}_t} \lambda_t).$$

Consider the following $G$-action on $\mathbb{P}(V') \times \mathbb{C} \times F \times \mathcal{KM}$

$$G \times (\mathbb{P}(V') \times \mathbb{C} \times F \times \mathcal{KM}) \longrightarrow \mathbb{P}(V') \times \mathbb{C} \times F \times \mathcal{KM},$$
$$(g_1, \ldots, g_s, (\mu_0 : \ldots : \mu_t), \lambda_0, f_1, \ldots, f_s, q_1, \ldots, q_r) \mapsto ((\mu_0 : \mu_1 : \ldots : \mu_t), \lambda_0, g_1^{-1} f_1, \ldots, g_s^{-1} f_s, q_1, \ldots, q_r).$$

It is easy to see that $T_\mathbb{P}$ resp. $T_\mathbb{P}^{-1}$ is $G$-equivariant with respect to the $G$-actions above.

Consider the map

$$T_S : S \times \mathbb{C} \times F \times \mathcal{KM} \longrightarrow S \times V^*,$$
$$(y_1, \ldots, y_s, \lambda_0, f_1, \ldots, f_s, q_1, \ldots, q_r) \mapsto (f_1^{-1} y_1, \ldots, f_s^{-1} y_s, \lambda_0, \underline{f}^{\underline{b}_1} \cdot \underline{q}^{\underline{m}_1}, \ldots, \underline{f}^{\underline{b}_t} \cdot \underline{q}^{\underline{m}_t})$$

and its inverse

$$T_S^{-1} : S \times V^* \longrightarrow S \times \mathbb{C} \times F \times \mathcal{KM},$$
$$(y_1, \ldots, y_s, \lambda_0, \ldots, \lambda_t) \mapsto (\underline{\lambda}^{\underline{c}_1} y_1, \ldots, \underline{\lambda}^{\underline{c}_s} y_s, \lambda_0, \underline{\lambda}^{\underline{c}_1}, \ldots, \underline{\lambda}^{\underline{c}_s}, \underline{\lambda}^{\underline{l}_1}, \ldots, \underline{\lambda}^{\underline{l}_r}),$$

where one has to use (29).
Recall the $G$-action on $S \times V^*$

$$G \times (S \times V^*) \longrightarrow S \times V^*,$$
$$(g_1, \ldots, g_s, \lambda_0, \ldots, \lambda_t) \mapsto (g_1 y_1, \ldots, g_s y_s, \lambda_0, \underline{g}^{-\underline{b}_1} \lambda_1, \ldots, \underline{g}^{-\underline{b}_t} \lambda_t)$$

and consider the following $G$-action on $S \times \mathbb{C} \times F \times \mathcal{KM}$

$$G \times (S \times \mathbb{C} \times F \times \mathcal{KM}) \longrightarrow S \times \mathbb{C} \times F \times \mathcal{KM},$$
$$(g_1, \ldots, g_s, y_1, \ldots, y_s, \lambda_0, f_1, \ldots, f_s, q_1, \ldots, q_r) \mapsto (y_1, \ldots, y_s, \lambda_0, g_1^{-1} f_1, \ldots, g_s^{-1} f_s, q_1, \ldots, q_r).$$

It is again easy to see that $T_S$ resp. $T_S^{-1}$ is $G$-equivariant with respect to the $G$-actions above.

The subvarieties $Z^*$ resp. $\Gamma^*$ are then given by $\lambda_0 \mu_0 + \sum_{i=1}^t \mu_i \cdot \underline{q}^{\underline{m}_i} = 0$ resp. $\lambda_0 + \sum_{i=1}^t \underline{y}^{\underline{b}_i} \cdot \underline{q}^{\underline{m}_i} = 0$.

Finally consider the maps

$$\begin{aligned}
T : \mathbb{C} \times F \times \mathcal{KM} &\longrightarrow V^*, \\
(\lambda_0, f_1, \ldots, f_s, q_1, \ldots, q_r) &\mapsto (\lambda_0, \underline{f}^{\underline{b}_1} \cdot \underline{q}^{\underline{m}_1}, \ldots \underline{f}^{\underline{b}_t} \cdot \underline{q}^{\underline{m}_t}), \\
T^{-1} : V^* &\longrightarrow \mathbb{C} \times F \times \mathcal{KM}, \\
(\lambda_0, \ldots, \lambda_t) &\mapsto (\lambda_0, \underline{\lambda}^{\underline{c}_1}, \ldots, \underline{\lambda}^{\underline{c}_s}, \underline{\lambda}^{\underline{l}_1}, \ldots, \underline{\lambda}^{\underline{l}_r}),
\end{aligned} \quad (30)$$

which are $G$-equivariant with respect to the $G$-action on $V^*$ and the following $G$-action on $\mathbb{C} \times F \times \mathcal{KM}$

$$G \times (\mathbb{C} \times F \times \mathcal{KM}) \longrightarrow \mathbb{C} \times F \times \mathcal{KM},$$
$$(g_1, \ldots, g_s, \lambda_0, f_1, \ldots, f_s, q_1, \ldots, q_r) \mapsto (\lambda_0, g_1^{-1} f_1, \ldots, g_s^{-1} f_s, q_1, \ldots, q_r).$$

The $G$-equivariant isomorphisms above show that the geometric quotients of $V^*$, $Z^*$ and $\Gamma^*$ by $G$ exist and are given by $\mathbb{C} \times \mathcal{KM}$,

$$\mathcal{Z} := \{\lambda_0 \mu_0 + \sum_{i=1}^t \underline{q}^{\underline{m}_i} \mu_i = 0\} \subset \mathbb{P}(V') \times \mathbb{C} \times \mathcal{KM}$$

and

$$\mathcal{G} := \{\lambda_0 + \sum_{i=1}^t \underline{q}^{\underline{m}_i} y_{\underline{b}_i} = 0\} \subset S \times \mathbb{C} \times \mathcal{KM},$$

respectively. We denote the corresponding quotient maps by $\pi_G^{V^*}, \pi_G^{Z^*}$ and $\pi_G^{\Gamma^*}$.

Notice that we have a natural section $i_G^{V^*}$ to $\pi_G^{V^*}$, which is induced by the inclusion

$$\mathbb{C} \times \mathcal{KM} \longrightarrow \mathbb{C} \times F \times \mathcal{KM},$$
$$(\lambda_0, q_1, \ldots, q_r) \mapsto (\lambda_0, 1, \ldots, 1, q_1, \ldots, q_r)$$

and the isomorphism above. This gives also rise to sections $i_G^{Z^*}$ and $i_G^{\Gamma^*}$ of $\pi_G^{Z^*}$ resp. $\pi_G^{\Gamma^*}$. Consider the following diagram

$$\begin{array}{c}
S \xleftarrow{\pi_1^S} \Gamma \xleftarrow{j_{\Gamma^*}} \Gamma^* \xleftarrow[\pi_G^{\Gamma^*}]{i_G^{\Gamma^*}} \mathcal{G} \\
\downarrow j \quad \downarrow \theta \quad \downarrow \zeta \quad \downarrow \gamma \\
X \xleftarrow{} Z_X \xleftarrow{j_{Z_X^*}} Z_X^* \xleftarrow{i_G^{Z_X^*}} \mathcal{Z}_X \\
\downarrow i \quad \downarrow \eta \quad \downarrow \varepsilon \quad \downarrow \beta \\
\mathbb{P}(V') \xleftarrow{\pi_1^Z} Z \xleftarrow{j_{Z^*}} Z^* \xleftarrow[\pi_G^{Z^*}]{i_G^{Z^*}} \mathcal{Z} \\
\downarrow \pi_Z^2 \quad \downarrow \delta \quad \downarrow \alpha \\
V \xleftarrow{j_{V^*}} V^* \xleftarrow[\pi_G^{V^*}]{i_G^{V^*}} \mathbb{C} \times \mathcal{KM}
\end{array} \quad (31)$$

Notice also that all squares are cartesian.

**Proposition 2.22.** *Let $i_G^{Z^*} : \mathcal{Z} \to Z^*$ resp. $i_G^{V^*} : \mathbb{C} \times \mathcal{KM} \to V^*$ be the sections constructed above.*

1. *The $\mathcal{D}_{Z^*}$-modules*
$$(\varepsilon \circ \zeta)_\dagger \mathcal{O}_{\Gamma^*}, \quad (\varepsilon \circ \zeta)_+ \mathcal{O}_{\Gamma^*} \quad \text{and} \quad \mathcal{M}^{IC}(Z_X^*)$$
   *are quasi-$G$-equivariant and non-characteristic with respect to $i_G^{Z^*}$.*

2. *The $\mathcal{D}_{V^*}$-modules*
$$\mathcal{H}^0(\varphi_{B,\dagger} \mathcal{O}_{S \times W^*}) \quad \text{and} \quad \mathcal{H}^0(\varphi_{B,+} \mathcal{O}_{S \times W^*})$$
   *are quasi-$G$-equivariant and non-characteristic with respect to $i_G^{V^*}$.*

3. *We have*
$$(i_G^{Z^*})^+ \mathcal{M}^{IC}(Z_X^*) \simeq \mathcal{M}^{IC}(\mathcal{Z}_X).$$

*In particular we have*
$$\alpha_+ \mathcal{M}^{IC}(\mathcal{Z}_X) \simeq i_{\mathcal{KM}}^+ \mathscr{R}\left(\mathcal{M}^{IC}(X)\right), \quad (32)$$

*where $i_{\mathcal{KM}} := j_{V^*} \circ i_G^{V^*}$ is non-characteristic with respect to $\mathscr{R}(\mathcal{M}^{IC}(X))$.*

*Proof.* 1. First notice that because the map $(i \circ j) : S \to \mathbb{P}(V')$ is affine and this property is preserved by base change, the map $(\varepsilon \circ \zeta)$ is also affine. Thus the direct image as well as the proper direct image of $\mathcal{O}_{\Gamma^*}$ is a single $\mathcal{D}_{Z^*}$-module. The closure of $\Gamma^*$ in $Z^*$ is $Z_X^*$, therefore we have

$$\mathcal{M}^{IC}(Z_X^*) = im((\varepsilon \circ \zeta)_\dagger \mathcal{O}_{\Gamma^*} \to (\varepsilon \circ \zeta)_+ \mathcal{O}_{\Gamma^*}) \in Mod_{rh}(\mathcal{D}_{Z^*}). \quad (33)$$

To show the first claim, it is enough by Lemma 2.21 to show that the corresponding $\mathcal{D}$-modules are quasi-$G$-equivariant. First recall that $\Gamma^* \subset S \times V^*$ and denote by $\iota : \Gamma^* \to S$ the restriction of the projection to the first factor. Notice that $\iota$ is $G$-equivariant and $\mathcal{O}_{\Gamma^*} \simeq \iota^+ \mathcal{O}_S$. Therefore $\mathcal{O}_{\Gamma^*}$ is a quasi-$G$-equivariant $\mathcal{D}$-module. Because $\varepsilon, \zeta$ is $G$-equivariant we see that $(\varepsilon \circ \zeta)_\dagger \mathcal{O}_{\Gamma^*}$ and $(\varepsilon \circ \zeta)_+ \mathcal{O}_{\Gamma^*}$ are quasi-$G$-equivariant. Furthermore, because of Equation (33) and the fact that $Mod(G, \mathcal{D}_{Z^*})$ is an abelian category the $\mathcal{D}$-module $\mathcal{M}^{IC}(Z_X^*)$ is quasi-$G$-equivariant.

2. For the second point, consider the action of $G$ on $W^* = (\mathbb{C}^*)^t$ which is given by
$$G \times W^* \longrightarrow W^*,$$
$$(g_1, \ldots, g_s, \lambda_1, \ldots, \lambda_t) \mapsto (\underline{g}^{-\underline{b}_1}\lambda_1, \ldots, \underline{g}^{-\underline{b}_t}\lambda_t).$$

This action together with the action (26) induces a $G$-action on $S \times W^*$. It is easy to see that $\varphi_{B|S \times W^*}$ is $G$-equivariant. Thus the $\mathcal{D}_{V^*}$-modules $\mathcal{H}^0(\varphi_{B,\dagger}\mathcal{O}_{S \times W^*})$ and $\mathcal{H}^0(\varphi_{B,+}\mathcal{O}_{S \times W^*})$ are quasi-$G$-equivariant. The fact that $i_G^{V^*}$ is non-characteristic with respect to these $\mathcal{D}_{V^*}$-modules follows now again from Lemma 2.21.

3. To show the third claim, consider the following isomorphisms
$$\begin{aligned}\mathcal{M}^{IC}(\mathcal{Z}_X) &\simeq im\left((\beta \circ \gamma)_\dagger \mathcal{O}_\mathcal{G} \to (\beta \circ \gamma)_+ \mathcal{O}_\mathcal{G}\right) \\ &\simeq im\left((\beta \circ \gamma)_\dagger (i_G^{\Gamma^*})^\dagger \mathcal{O}_{\Gamma^*} \to (\beta \circ \gamma)_+ (i_G^{\Gamma^*})^+ \mathcal{O}_{\Gamma^*}\right) \\ &\simeq im\left((i_G^{Z^*})^\dagger (\varepsilon \circ \zeta)_\dagger \mathcal{O}_{\Gamma^*} \to (i_G^{Z^*})^+ (\varepsilon \circ \zeta)_+ \mathcal{O}_{\Gamma^*}\right) \\ &\simeq (i_G^{Z^*})^+ im\left((\varepsilon \circ \zeta)_\dagger \mathcal{O}_{\Gamma^*} \to (\varepsilon \circ \zeta)_+ \mathcal{O}_{\Gamma^*}\right) \\ &\simeq (i_G^{Z^*})^+ \mathcal{M}^{IC}(Z_X^*),\end{aligned}$$

where the second isomorphism follows from $(i_G^{\Gamma^*})^+ \mathcal{O}_{\Gamma^*} \simeq \mathcal{O}_\mathcal{G}$, the fact that $\mathcal{O}_{\Gamma^*}$ is non-characteristic for $i_G^{\Gamma^*}$ and [HTT08, Theorem 2.7.1(ii)]. The third isomorphism follows by base change and the fourth isomorphism follows from the fact that $i_G^{Z^*}$ is non-characteristic with respect to $(\varepsilon \circ \zeta)_\dagger \mathcal{O}_{\Gamma^*}$ and $(\varepsilon \circ \zeta)_+ \mathcal{O}_{\Gamma^*}$.

For the last claim consider the following diagram
$$\begin{array}{ccccc} \mathcal{Z} & \xleftarrow{j_{Z^*}} & Z^* & \xleftarrow{i_G^{Z^*}} & \mathcal{Z} \\ {\scriptstyle \pi_Z^2}\downarrow & & \downarrow{\scriptstyle \delta} & & \downarrow{\scriptstyle \alpha} \\ V & \xleftarrow{j_{V^*}} & V^* & \xleftarrow{i_G^{V^*}} & \mathbb{C} \times \mathcal{KM} \end{array}$$

We have the following isomorphisms
$$\begin{aligned}\alpha_+ \mathcal{M}^{IC}(\mathcal{Z}_X) &\simeq \alpha_+ (i_G^{Z^*})^+ \mathcal{M}^{IC}(Z_X^*) \\ &\simeq \alpha_+ (i_G^{Z^*})^+ j_{Z^*}^+ \mathcal{M}^{IC}(Z_X) \\ &\simeq (i_G^{V^*})^+ j_{V^*}^+ \pi_{2+}^Z \mathcal{M}^{IC}(Z_X) \\ &\simeq i_{\mathcal{KM}}^+ \pi_{2+}^Z \mathcal{M}^{IC}(Z_X) \\ &\simeq i_{\mathcal{KM}}^+ \pi_{2+}^Z (\pi_1^Z)^+ \mathcal{M}^{IC}(X) \\ &\simeq i_{\mathcal{KM}}^+ \mathscr{R}(\mathcal{M}^{IC}(X)).\end{aligned}$$

The non-characteristic property of $i_{\mathcal{KM}} = j_{V^*} \circ i_G^{V^*}$ follows from Lemma 2.21 and the fact that $j_{V^*}^+ \mathscr{R}(\mathcal{M}^{IC}(X))$ is quasi-$G$-equivariant. □

# 3 Fourier-Laplace transformation and lattices

In this section we apply the Fourier-Laplace transformation functor $FL_W$ to the various $\mathcal{D}$-modules considered in section 2. For the families of Laurent polynomials resp. compactifications thereof that appear in mirror symmetry, we obtain $\mathcal{D}$-modules that can eventually be matched with the differential systems defined by quantum cohomology. They have in general irregular singularities, and this is reflected in the fact that although the modules considered in section 2 were monodromic on $V$, they do not have necessarily that property with respect to the vector bundle $V = \mathbb{C}_{\lambda_0} \times W \to W$. Hence the functor $\mathrm{FL}_W$ will in general not preserve regularity.

In the second part of this section, we study a lattice in the Fourier-Laplace transformation of the Gauß-Manin system of the family of Laurent polynomials $\varphi_B$. It is given by a so-called twisted de Rham complex, however, in order to obtain a good hypergeometric description of it, we have to introduce a certain intermediate compactification of $\varphi_B$ and replace this de Rham complex by a logarithmic version. Moreover, the parameters of the family $\varphi_B$ have to be restricted to a Zariski open set excluding certain (but not all) singularities at infinity. Then we can show the necessary finiteness and freeness of the lattice. It will later correspond to the twisted quantum $\mathcal{D}$-module (see section 4), seen as a family of algebraic vector bundles over $\mathbb{C}_z$ (not only over $\mathbb{C}_z^*$) with connection operator which is meromorphic along $\{z = 0\}$.

## 3.1 Localized Fourier-Laplace Transform

We discuss here a partial localized Fourier-Laplace transform of the Gauß-Manin systems of $\varphi_B$ and of the $\mathcal{D}$-module $\mathcal{M}^{IC}(X^\circ, \mathcal{L})$.

Consider the product decomposition $V = \mathbb{C}_{\lambda_0} \times W$, where $W$ is the hyperplane given by $\lambda_0 = 0$. We interpret $V$ as a rank one bundle with base $W$ and consider the Fourier-Laplace transformation with respect to the base $W$ as in Definition 2.4, where we denote the coordinate on the dual fiber by $\tau$. Set $z = 1/\tau$ and denote by $j_\tau : \mathbb{C}_\tau^* \times W \hookrightarrow \mathbb{C}_\tau \times W$ and $j_z : \mathbb{C}_\tau^* \times W \hookrightarrow \hat{V} := \mathbb{C}_z \times W = \mathbb{P}_\tau^1 \setminus \{\tau = 0\} \times W$ the canonical embeddings. Let $\mathcal{N}$ be a $\mathcal{D}_V$-module, the partial, localized Fourier-Laplace transformation is defined by

$$\mathrm{FL}_W^{loc}(\mathcal{N}) := j_{z+} j_\tau^+ \, \mathrm{FL}_W(\mathcal{N}) \, .$$

The localized Fourier-Laplace transformations of the Gauß-Manin systems are denoted by

$$\mathcal{G}^+ := \mathrm{FL}_W^{loc}(\mathcal{H}^0(\varphi_{B,+} \mathcal{O}_{S \times W})) \, , \tag{34}$$

$$\mathcal{G}^\dagger := \mathrm{FL}_W^{loc}(\mathcal{H}^0(\varphi_{B,\dagger} \mathcal{O}_{S \times W})) \, . \tag{35}$$

We also consider the partial, localized Fourier-Laplace transform of the $\mathcal{D}$-modules $\mathcal{M}_{\widetilde{B}}^{\widetilde{\beta}}$. The following notation will be useful.

**Definition 3.1.** *Let $\widehat{M}_B^{(\beta_0, \beta)}$ be the $D_{\widehat{V}}$-module $D_{\widehat{V}}[z^{-1}]/I$, where $I$ is the left ideal generated by the operators $\widehat{\Box}_{\underline{l}}$, $\widehat{E}_k - \beta_k z$ and $\widehat{E} - \beta_0 z$, which are defined by*

$$\widehat{\Box}_{\underline{l}} := \prod_{i : l_i < 0} (z \cdot \partial_{\lambda_i})^{-l_i} - \prod_{i : l_i > 0} (z \cdot \partial_{\lambda_i})^{l_i}, \quad \underline{l} \in \mathbb{L}_B$$

$$\widehat{E} := z^2 \partial_z + \sum_{i=1}^t z \lambda_i \partial_{\lambda_i} \, ,$$

$$\widehat{E}_k := \sum_{i=1}^t b_{ki} z \lambda_i \partial_{\lambda_i} \, , \quad k = 1, \ldots, s.$$

*We denote the corresponding $\mathcal{D}_{\widehat{V}}$-module by $\widehat{\mathcal{M}}_B^{(\beta_0, \beta)}$.*

**Lemma 3.2.** *We have the following isomorphism*

$$\mathrm{FL}_W^{loc}(\mathcal{M}_{\widetilde{B}}^{\widetilde{\beta}}) \simeq \widehat{\mathcal{M}}_B^{(\beta_0 + 1, \beta)}$$

*for every $\widetilde{\beta} = (\beta_0, \beta) \in \mathbb{Z}^{s+1}$.*

*Proof.* This is an easy calculation, using the substitution

$$\lambda_0 \to -\partial_\tau = z^2 \partial_z \quad \text{and} \quad \partial_{\lambda_0} \to \tau = 1/z$$

and the fact that $\widehat{\mathcal{M}}_B^{(\beta_0, \beta)}$ is localized along $z = 0$. □

Notice that in the lemma above we used the subscript $\widetilde{B}$ for the GKZ-system on the left hand side and the subscript $B$ for its localized Fourier-Laplace transform on the right hand side. This notation takes into account the fact that the properties of the system $\mathcal{M}_{\widetilde{B}}^{\widetilde{\beta}}$ are governed by the geometry of the semi-group $\mathbb{N}\widetilde{B}$, whereas the properties of its localized Fourier-Laplace transform $\widehat{\mathcal{M}}^{(\beta_0+1,\beta)}$ depend on the geometry of $\mathbb{N}B$. This explains the different sets of allowed parameters in Proposition 3.3 resp. Theorem 3.6 in contrast to Theorem 2.11 resp. Theorem 2.16 and Proposition 2.17.

Notice that under the normality assumption on the semi-group $\mathbb{N}\widetilde{B}$, the rank of $\widehat{\mathcal{M}}_B^{(\beta_0,\beta)}$ is also equal to $d! \cdot \mathrm{vol}(\mathrm{Conv}(\underline{b}_1,\ldots,\underline{b}_t))$ (this can be shown by an argument similar to [RS15, Proposition 2.7]).

The following proposition gives an isomorphism between the localized partial Fourier-Laplace transform of the Gauß-Manin systems $\mathcal{G}^+$ and $\mathcal{G}^\dagger$ and the hypergeometric systems $\widehat{\mathcal{M}}_B^{(\beta_0,\beta)}$ introduced above.

**Proposition 3.3.** *There exists a $\delta_B \in \mathbb{N}B$ such that we have an isomorphism*

$$\mathcal{G}^+ \simeq \widehat{\mathcal{M}}_B^{(\beta_0,\beta)}$$

*for every $\beta_0 \in \mathbb{Z}$ and $\beta \in \delta_B + (\mathbb{R}_{\geq 0}B \cap \mathbb{Z}^s)$. If $\mathbb{N}B$ is saturated, then $\delta_B$ can be taken to be $\underline{0} \in \mathbb{N}B$ (in particular, the statement holds for $(\beta_0,\beta) = (\beta_0,\underline{0}) \in \mathbb{Z}^{1+s}$).*

*Furthermore, we have an isomorphism*

$$\mathcal{G}^\dagger \simeq \widehat{\mathcal{M}}_B^{(\beta_0',-\beta')}$$

*for every $\beta_0' \in \mathbb{Z}$ and $\beta' \in (\mathbb{R}_{\geq 0}B)^\circ \cap \mathbb{Z}^s$.*

*Proof.* We construct the isomorphisms by applying the Fourier-Laplace transform $\mathrm{FL}_W$ to the exact sequences in Theorem 2.11. First notice that the first and last term in the exact sequences are free $\mathcal{O}_V$-modules, thus their Fourier-Laplace transform has support on $\tau = 0$, i.e. their localized Fourier-Laplace transform is 0. Thus there is some $\delta_{\widetilde{B}} \in \mathbb{N}\widetilde{B}$ such that we have the following isomorphisms

$$\mathcal{G}^+ = \mathrm{FL}_W^{loc}(\mathcal{H}^0(\varphi_{B,+}\mathcal{O}_{S \times W})) \simeq \mathrm{FL}_W^{loc}(\mathcal{M}_{\widetilde{B}}^{\widetilde{\beta}})$$

and

$$\mathcal{G}^\dagger = \mathrm{FL}_W^{loc}(\mathcal{H}^0(\varphi_{B,\dagger}\mathcal{O}_{S \times W})) \simeq \mathrm{FL}_W^{loc}(\mathcal{M}_{\widetilde{B}}^{-\widetilde{\beta}'})$$

for any $\widetilde{\beta} \in \delta_{\widetilde{B}} + (\mathbb{R}_{\geq 0}\widetilde{B} \cap \mathbb{Z}^{s+1})$ and any $\widetilde{\beta}' \in (\mathbb{R}_{\geq 0}\widetilde{B})^\circ \cap \mathbb{Z}^{s+1}$. Write $\delta_{\widetilde{B}} = (\delta_0, \delta_B)$ with $\delta_B \in \mathbb{Z}^s$. Now given any $(\beta_0, \beta) \in \mathbb{Z} \times (\delta_B + (\mathbb{R}_{\geq 0}B \cap \mathbb{Z}^s))$ resp. $(\beta_0', \beta') \in \mathbb{Z} \times ((\mathbb{R}_{\geq 0}B)^\circ \cap \mathbb{Z}^s)$ we can find a $\gamma_0, \gamma_0' \in \mathbb{Z}$ such that $(\gamma_0, \beta) \in \delta_{\widetilde{B}} + (\mathbb{R}_{\geq 0}\widetilde{B} \cap \mathbb{Z}^{s+1})$ resp. $(\gamma_0', \beta') \in (\mathbb{R}_{\geq 0}\widetilde{B})^\circ \cap \mathbb{Z}^{s+1}$. It remains to show that there are isomorphism

$$\widehat{\mathcal{M}}_B^{(\beta_0,\beta)} \simeq \widehat{\mathcal{M}}_B^{(\gamma_0,\beta)} \tag{36}$$

for $(\beta_0, \beta) \in \mathbb{Z} \times (\delta_B + (\mathbb{R}_{\geq 0}B \cap \mathbb{Z}^s))$ and $(\gamma_0, \beta) \in \delta_{\widetilde{B}} + (\mathbb{R}_{\geq 0}\widetilde{B} \cap \mathbb{Z}^{s+1}))$ resp.

$$\widehat{\mathcal{M}}_B^{(\beta_0',-\beta')} \simeq \widehat{\mathcal{M}}_B^{(-\gamma_0',-\beta)} \tag{37}$$

for $(\beta_0', \beta') \in \mathbb{Z} \times ((\mathbb{R}_{\geq 0}B)^\circ \cap \mathbb{Z}^s)$ and $(-\gamma_0', -\beta') \in ((\mathbb{R}_{\geq 0}\widetilde{B})^\circ \cap \mathbb{Z}^{s+1})$. Notice that $\widehat{\mathcal{M}}_B^{(\beta_0,\beta)}$ is localized along $z = 0$ for all $(\beta_0, \beta) \in \mathbb{Z}^{s+1}$ by Lemma (3.2). Therefore the morphism given by right multiplication with $z$

$$\widehat{\mathcal{M}}_B^{(\beta_0,\beta)} \xrightarrow{\cdot z} \widehat{\mathcal{M}}_B^{(\beta_0-1,\beta)} \tag{38}$$

is an isomorphism, which shows (36) and (37).

Concerning the last statement, suppose that $\mathbb{N}B$ is saturated. Let $\beta \in \mathbb{N}B = (\mathbb{R}_{\geq 0}B \cap \mathbb{Z}^s)$ and let $\beta_0 \in \mathbb{Z}$ be arbitrary. By [Rei14, Lemma 1.17] we have $\beta \notin sRes(B)$, where $sRes(B) \subset \mathbb{C}^s$ is the set of strongly resonant values (cf. [SW09, Definition 3.4]). Using [Rei14, Lemma 1.19] there exists a $\gamma_0 \in \mathbb{Z}$ such that $(\gamma_0, \beta) \notin sRes(\widetilde{B})$. Now we argue as above, i.e. by [Rei14, Theorem 2.7] we have $\mathcal{G}^+ = \mathrm{FL}_W^{loc}(\mathcal{H}^0(\varphi_{B,+}\mathcal{O}_{S \times W})) \simeq \mathrm{FL}_W^{loc}(\mathcal{M}_{\widetilde{B}}^{(\gamma_0,\beta)})$ which in turn is isomorphic to $\widehat{\mathcal{M}}_B^{(\beta_0,\beta)}$. □

If the semigroup $\mathbb{N}B$ is saturated, we will compute the isomorphism above explicitly for $(\beta_0, \beta) = (0, \underline{0})$. For this we will need a direct description of the localized, partial Fourier-Laplace transformed Gauß-Manin system $\mathcal{G}^+$.

**Lemma 3.4.** *Write $\varphi_B = (F, pr)$, where $F : S \times W \to \mathbb{C}$, $(\underline{y}, \underline{\lambda}) \mapsto -\sum_{i=1}^{t} \lambda_i \underline{y}^{b_i}$ and $pr : S \times W \to W$ is the projection. Recall from formula (34) that we denote by $\mathcal{G}^+$ the localized Fourier-Laplace transformation of the Gauß-Manin system of the morphism $\varphi_B$. Write $G^+ := H^0(\widehat{V}, \mathcal{G}^+)$ for its module of global sections. Then there is an isomorphism of $D_{\widehat{V}}$-modules*

$$G^+ \cong H^0\left(\Omega^{\bullet+s}_{S\times W/W}[z^{\pm}], d - z^{-1} \cdot d_y F \wedge\right),$$

*where $d$ is the differential in the relative de Rham complex $\Omega^{\bullet}_{S\times W/W}$. The structure of a $\mathcal{D}_{\widehat{V}}$-module on the right hand side is defined as follows*

$$\partial_z(\omega \cdot z^i) = i \cdot \omega \cdot z^{i-1} + F \cdot \omega \cdot z^{i-2},$$
$$\partial_{\lambda_i}(\omega \cdot z^i) := \partial_{\lambda_i}(\omega) \cdot z^i - \partial_{\lambda_i} F \cdot \omega \cdot z^{i-1} = \partial_{\lambda_i}(\omega) \cdot z^i + \underline{y}^{b_i} \cdot \omega \cdot z^{i-1},$$

*where $\omega \in \Omega^s_{S \times W/W}$.*

*Proof.* The expression for the module $G^+$ as well as the formulas for the $\mathcal{D}_{\widehat{V}}$-structure are an immediate consequence of the definition of the direct image functor. See, e.g. [Rei14, equations 2.0.18, 2.0.19], from which the desired formulas can be easily obtained. $\square$

Using the description of $G^+$ via relative differential forms, we find a distinguished element, which is (the class of) the volume form on $S$, that is

$$\omega_0 := \frac{dy_1}{y_1} \wedge \ldots \wedge \frac{dy_s}{y_s}.$$

In the next lemma we compute the image of $\omega_0$ under the isomorphisms in Proposition 3.3 under the assumption of normality of $\mathbb{N}B$.

**Lemma 3.5.** *Let $\mathbb{N}B$ be a saturated semigroup, then the isomorphism from Proposition 3.3*

$$\Phi : \mathcal{G}^+ \xrightarrow{\simeq} \widehat{\mathcal{M}}^{(0,\underline{0})}_B$$

*maps $\omega_0$ to $1$.*

*Proof.* Recall from the proof of Proposition 3.3, that there exists a $\gamma_0 \in \mathbb{Z}$ such that $(\gamma_0, \underline{0}) \notin sRes(\widetilde{B})$ (notice that here we only assume that $\mathbb{N}B$ is saturated which does not imply that $\mathbb{N}\widetilde{B}$ is saturated). Denote by

$$\psi_{(\gamma_0, \underline{0})} : \Gamma(V, \mathcal{H}^0(\varphi_{B,+} \mathcal{O}_{S \times W})) \to M^{(\gamma_0, \underline{0})}_{\widetilde{B}}$$

the morphism from Theorem 2.11. We first compute the image of $\omega_0$ under the morphism $\psi_{(\gamma_0, \underline{0})}$ using the description of $\mathcal{H}^0(\varphi_{B,+} \mathcal{O}_{S \times W})$ by relative differential forms (see e.g. [Rei14, Equation 2.0.17]). We will use the following two facts of loc. cit. Proposition 2.8 whose proofs extend directly to our slightly more general situation (there it was assumed that $\mathbb{N}\widetilde{B}$ is saturated). Namely first, that there exists a non-zero morphism $M^{(-1,\underline{0})}_{\widetilde{B}} \to \Gamma(V, \mathcal{H}^0(\varphi_{B,+} \mathcal{O}_{S \times W}))$ which sends $1$ to $\omega_0$ and second that $\psi_{(\gamma_0, \underline{0})}(\omega_0) \neq 0$. Concatenating this morphism with $\psi_{(\gamma_0, \underline{0})}$ gives a non-zero morphism $M^{(-1,\underline{0})}_{\widetilde{B}} \to M^{(\gamma_0, \underline{0})}_{\widetilde{B}}$, where $1 \in M^{(-1,\underline{0})}_{\widetilde{B}}$ is sent to the image of $\omega_0$ under $\psi_{(\gamma_0, \underline{0})}$. By [Rei14, Proposition 1.24] this morphism is uniquely given by right multiplication with $\partial^{\gamma_0+1}_{\lambda_0}$ (up to a non-zero constant). Applying now the partial localized Fourier-Laplace transform to the morphism $\psi_{(\gamma_0, \underline{0})}$, we see that $\psi_{(\gamma_0, \underline{0})}(\omega_0) = z^{-\gamma_0-1}$. Using the isomorphism $\widehat{\mathcal{M}}^{(\gamma'_0, \underline{0})}_B \xrightarrow{\cdot z} \widehat{\mathcal{M}}^{(\gamma'_0-1, \underline{0})}_B$, which holds for any $\gamma'_0 \in \mathbb{Z}$, shows the claim. $\square$

By Proposition 2.16, we can now give a concrete description of the partial, localized Fourier-Laplace transform $\widehat{\mathcal{M}^{IC}}(X^\circ, \mathcal{L}) := \mathrm{FL}^{loc}_W(\mathcal{M}^{IC}(X^\circ, \mathcal{L}))$ of the intersection cohomology $\mathcal{D}$-module $\mathcal{M}^{IC}(X^\circ, \mathcal{L})$.

**Theorem 3.6.** *Let $\beta \in \delta_B + (\mathbb{R}_{\geq 0} B \cap \mathbb{Z}^s)$, $\beta' \in (\mathbb{R}_{\geq 0} B)^\circ \cap \mathbb{Z}^s$ and $\beta_0, \beta_0' \in \mathbb{Z}$, then we have the following isomorphisms*

$$\widehat{\mathcal{M}^{IC}}(X^\circ, \mathcal{L}) \simeq im\left( \widehat{\mathcal{M}}_B^{(\beta_0', -\beta')} \xrightarrow{\cdot z^{\beta_0' - \beta_0} \partial^{\beta + \beta'}} \widehat{\mathcal{M}}_B^{(\beta_0, \beta)} \right),$$

*resp.*

$$\widehat{\mathcal{M}^{IC}}(X^\circ, \mathcal{L}) \simeq \widehat{\mathcal{M}}_B^{(\beta_0', -\beta')} / \widehat{\Gamma}_\partial\left( \widehat{\mathcal{M}}_B^{(\beta_0', -\beta')} \right),$$

*where $\widehat{\Gamma}_\partial\left( \widehat{\mathcal{M}}_B^{(\beta_0', -\beta')} \right)$ is the sub-$\mathcal{D}$-module corresponding to the sub-$D$-module*

$$\widehat{\Gamma}_\partial\left( \widehat{M}_B^{(\beta_0', -\beta')} \right) := \{ m \in \widehat{M}_B^{(\beta_0', -\beta')} \mid \exists n \in \mathbb{N} \text{ with } \left( \partial^{\beta + \beta'} \right)^n \cdot m = 0 \}.$$

*Furthermore, if $\mathbb{N}B$ is saturated, then $\delta_B$ can be taken to be $\underline{0} \in \mathbb{N}B$ (so that, similarly to Proposition 3.3, the statement holds true for $(\beta_0, \beta) = (\beta_0, \underline{0}) \in \mathbb{Z}^{1+s}$).*

*Proof.* Using the isomorphism

$$\widehat{\mathcal{M}}_B^{(\beta_0, \beta)} \xrightarrow{\cdot z} \widehat{\mathcal{M}}_B^{(\beta_0 - 1, \beta)}, \tag{39}$$

which holds for every $(\beta_0, \beta) \in \mathbb{Z}^{s+1}$, we can assume that $(\beta_0 + 1, \beta) \in \delta_{\widetilde{B}} + (\mathbb{R}_{\geq 0}\widetilde{B} \cap \mathbb{Z}^{s+1})$ resp. $(\beta_0' + 1, \beta') \in (\mathbb{R}_{\geq 0}\widetilde{B})^\circ \cap \mathbb{Z}^{s+1}$. Then the first isomorphism follows by applying the functor $\mathrm{FL}_W^{loc}$ to the isomorphism in Theorem 2.16 and Lemma 3.2.

For the second isomorphism we can assume again that $(\beta_0' + 1, \beta') \in (\mathbb{R}_{\geq 0}\widetilde{B})^\circ \cap \mathbb{Z}^{s+1}$. Now the desired statement is obtained by applying $\mathrm{FL}_W^{loc}$ to the second isomorphism in Proposition 2.17 and the fact that $\widehat{\Gamma}_\partial(\widehat{\mathcal{M}}_B^{(\beta_0', -\beta')})$ is stable under left multiplication with $z$.

Now assume that $\mathbb{N}B$ is saturated and let $\beta \in \mathbb{N}B$. Arguing as in the last part of the proof of Proposition 3.3 we can find a $\gamma_0 \in \mathbb{Z}$ such that $(\gamma_0, \beta) \notin sRes(\widetilde{\beta})$. By [SW09, Corollary 3.7] we have an isomorphism $\mathrm{FL}(h_+ \mathcal{O}_T) \simeq \mathcal{M}_{\widetilde{B}}^{(\gamma_0, \beta)}$. Now the proof of Theorem 2.16 shows that

$$\mathcal{M}^{IC}(X^\circ, \mathcal{L}) \simeq im(\mathcal{M}_{\widetilde{B}}^{-\widetilde{\beta}'} \xrightarrow{\cdot \partial^{(\gamma_0, \beta) + \widetilde{\beta}'}} \mathcal{M}_{\widetilde{B}}^{(\gamma_0, \beta)}).$$

Now applying the functor $F_W^{loc}$ and using the isomorphism (39) shows the claim in the saturated case. $\square$

## 3.2 Tameness and Lattices

In this section we define a natural lattice in the Fourier-Laplace transformed Gauß-Manin system $\mathcal{G}^+$ outside some bad locus where the Laurent polynomial acquires singularities at infinity. For this we need to study the characteristic variety of the Gauß-Manin system of $\varphi_B$ and the corresponding GKZ system $\mathcal{M}_{\widetilde{B}}^{\widetilde{\beta}}$. Throughout this section we assume that $\mathbb{N}B$ is a saturated semigroup. Recall the embedding of the torus $S$ in the projective space from formula (10)

$$S \xrightarrow{j} X \xrightarrow{i} \mathbb{P}(V').$$

The projective variety $X$ serves as a convenient ambient space to compactify fibers of the family of Laurent polynomials $\varphi_B$. However, we will also need an intermediate partial compactification of $S$, which is still an affine variety.

**Definition 3.7.** *The restriction of $X$ to the affine chart of $\mathbb{P}(V')$ given by $\mu_0 = 1$ is called $X^{aff}$, in other words, $X^{aff}$ is the closure of the map*

$$g_B : S \longrightarrow \mathbb{C}^t,$$
$$(y_1, \ldots, y_s) \mapsto (\underline{y}^{b_1}, \ldots, \underline{y}^{b_t}),$$

*and thererfore isomorphic to $\mathrm{Spec}\,(\mathbb{C}[\mathbb{N}B])$.*

Consider the following diagram, which is a refinement of a part of diagram (25):

$$\begin{array}{ccccccccc}
\Gamma & \xrightarrow{\theta_2} & Z_{X^{aff}} & \xrightarrow{\theta_1} & Z_X & \xrightarrow{\eta} & Z & \xrightarrow{\pi_2^Z} & V \\
\downarrow & & \downarrow & & \downarrow & & \downarrow \pi_1^Z & & \\
S & \xrightarrow{j_2} & X^{aff} & \xrightarrow{j_1} & X & \xrightarrow{i} & \mathbb{P}(V') & &
\end{array} \qquad (40)$$

where $j_1$ and $j_2$ are the canonical inclusions and the three squares are cartesian. Recall that $Z \subset \mathbb{P}(V') \times V$ was given by the incidence relation $\sum_{i=0}^{t} \lambda_i \mu_i = 0$ and the composed map $g = i \circ j = i \circ j_1 \circ j_2$ was defined by formula (9). Thus $\Gamma$ resp. $Z_{X^{aff}}$ is the subvariety of $S \times V = S \times \mathbb{C}_{\lambda_0} \times W$ resp. $X^{aff} \times V$ given by the equation $\lambda_0 + \sum_{i=1}^{r} \lambda_i \underline{y}^{\underline{b}_i} = 0$. It follows from the definition that $\Gamma$ is the graph of $\varphi_B$. Therefore the maps

$$\pi_{Z_X} := \pi_2^Z \circ \eta : Z_X \longrightarrow V$$

resp.

$$\pi_{Z_{X^{aff}}} := \pi_2^Z \circ \eta \circ \theta_1 : Z_{X^{aff}} \longrightarrow V$$

provide natural (partial) compactifications of the family of Laurent polynomials $\varphi_B$. Putting $H_{\widetilde{\underline{\lambda}}} := \{\sum_{i=0}^{t} \lambda_i \mu_i = 0\} \subset \mathbb{P}(V')$ for any $\widetilde{\underline{\lambda}} \in V$, we see that the fiber $\pi_{Z_X}^{-1}(\widetilde{\underline{\lambda}})$ resp. $\pi_{Z_{X^{aff}}}^{-1}(\widetilde{\underline{\lambda}})$ is given by $X \cap H_{\widetilde{\underline{\lambda}}}$ resp. $\{\lambda_0 + \sum_{i=1}^{t} \lambda_i \underline{y}^{\underline{b}_i} = 0\} \subset X^{aff}$.

Recall that the toric variety $X$ has a natural stratification by torus orbits $X^0(\Gamma)$, which are in one-to-one correspondence with the faces $\Gamma$ of the polytope $Q$, which is the convex hull of the elements $\{\underline{b}_0 := \underline{0}, \underline{b}_1, \ldots, \underline{b}_t\}$. Notice that the stratification $\mathcal{S} := \{X^0(\Gamma)\}$ is a Whitney stratification of $X$ (see e.g. [Dim92, Proposition 1.14].

By [GKZ08, Chapter 5, Prop 1.9] the orbit $X^0(\Gamma) \simeq (\mathbb{C}^*)^{\dim(\Gamma)}$ is the image of the map

$$g_\Gamma : S \longrightarrow \mathbb{P}(V'),$$
$$(y_1, \ldots, y_s) \mapsto (\varepsilon_0 1 : \varepsilon_1 \underline{y}^{\underline{b}_1} : \ldots : \varepsilon_t \underline{y}^{\underline{b}_t}),$$

where $\varepsilon_i = 0$ if $\underline{b}_i \notin \Gamma$ and $\varepsilon_i = 1$ if $\underline{b}_i \in \Gamma$. It is easy to see that

$$X^{aff} = \bigcup_{\Gamma \mid 0 \in \Gamma} X^0(\Gamma)$$

and this induces a Whitney stratification of $X^{aff}$.

The preimage of $X^0(\Gamma) \cap H_{\widetilde{\underline{\lambda}}}$ under $g_\Gamma$ is given by

$$\{(y_1, \ldots, y_s) \in S \mid \sum_{\underline{b}_i \in \Gamma} \lambda_i \underline{y}^{\underline{b}_i} = 0\}.$$

It follows from [GKZ08, Chapter 5.D] that the morphism $g_\Gamma : S \longrightarrow X^0(\Gamma) \simeq (\mathbb{C}^*)^{\dim(\Gamma)}$ is a trivial fibration with fiber being isomorphic to $(\mathbb{C}^*)^{d-\dim(\Gamma)}$.

Denote by $S_\Gamma^{crit, \widetilde{\underline{\lambda}}}$ the set

$$\left\{ (y_1, \ldots, y_s) \in S \mid \sum_{\underline{b}_i \in \Gamma} \widetilde{\lambda}_i \underline{y}^{\underline{b}_i} = 0; \; y_k \partial_{y_k} (\sum_{\underline{b}_i \in \Gamma} \widetilde{\lambda}_i \underline{y}^{\underline{b}_i}) = 0 \quad \text{for all} \quad k \in \{1, \ldots, s\} \right\}. \qquad (41)$$

Then its image under $g_\Gamma$ is exactly the singular set $sing(X^0(\Gamma) \cap H_{\widetilde{\underline{\lambda}}})$ of $X^0(\Gamma) \cap H_{\widetilde{\underline{\lambda}}}$. This motivates the following definition.

**Definition 3.8.** *Let* $\widetilde{\underline{\lambda}} \in V$

1. The fiber $\pi_{Z_X}^{-1}(\widetilde{\underline{\lambda}})$ has stratified singularities in $X^0(\Gamma)$ if $X^0(\Gamma) \cap H_{\widetilde{\underline{\lambda}}}$ is singular, i.e. $S_\Gamma^{crit,\widetilde{\underline{\lambda}}} \neq 0$.

2. The set
$$\Delta_B := \{\widetilde{\underline{\lambda}} \in V \mid S_Q^{crit,\widetilde{\underline{\lambda}}} \neq \emptyset\}$$
$$= \{\widetilde{\underline{\lambda}} \in V \mid \varphi_B^{-1}(\widetilde{\underline{\lambda}}) \text{ is singular}\}$$
is called the discriminant of $\varphi_B$.

3. The fiber $\varphi_B^{-1}(\widetilde{\underline{\lambda}})$ has **singularities at infinity** if there exists a proper face $\Gamma$ of the Newton polyhedron $Q$ so that $S_\Gamma^{crit,\widetilde{\underline{\lambda}}} \neq \emptyset$. The set
$$\Delta_B^\infty := \{\widetilde{\underline{\lambda}} \in V \mid \exists\, \Gamma \neq Q \text{ so that } S_\Gamma^{crit,\widetilde{\underline{\lambda}}} \neq \emptyset\}$$
is called the non-tame locus of $\varphi_B$.

4. The fiber $\varphi_B^{-1}(\widetilde{\underline{\lambda}})$ has **bad singularities at infinity** if there exists a proper face $\Gamma$ of the Newton polyhedron $Q$ not containing the origin such that $S_\Gamma^{crit,\widetilde{\underline{\lambda}}} \neq \emptyset$. The set
$$\Delta_B^{bad} := \{\widetilde{\underline{\lambda}} \in V \mid \exists\, \Gamma \neq Q, 0 \notin \Gamma \text{ so that } S_\Gamma^{crit,\widetilde{\underline{\lambda}}} \neq \emptyset\} \subset \Delta_B^\infty$$
is called the bad locus of $\varphi_B$.

**Remark 3.9.** Notice that $\Delta_B^{bad}$ is independent of $\lambda_0$. We denote its projection to $W$ by $W^{bad}$. Let $W^* = W \setminus \{\lambda_1 \ldots \lambda_t = 0\}$ and define
$$W^\circ := W^* \setminus W^{bad},$$
which we call the set of good parameters for $\varphi_B$.

Recall that $X^{aff}$ is isomorphic to $\mathrm{Spec}\,(R_B)$ with $R_B := \mathbb{C}[\mathbb{N}B]$. Let $\underline{\lambda} \in W$ and set $f_{\underline{\lambda}}(\bullet) := \varphi_B(\bullet, \underline{\lambda})$. Notice that the Laurent polynomials $f_{\underline{\lambda}}$ and $y_k \partial f_{\underline{\lambda}}/\partial y_k$ for $k = 1, \ldots, s$, which were defined on $S$ before are actually elements of $R_B$ and can thus naturally be considered as functions on $X^{aff}$.

**Lemma 3.10.** Let $\underline{\lambda} \in W^\circ$ be a good parameter, then
$$\dim_\mathbb{C} \left( R_B/(y_k \partial f_{\underline{\lambda}}/\partial y_k)_{k=1,\ldots,s} \right) = \mathrm{vol}(Q),$$
where the volume of a hypercube $[0,1]^s \subset \mathbb{R}^s$ is normalized to $s!$. Moreover, we have
$$supp(R_B/(y_k \partial f_{\underline{\lambda}}/\partial y_k)_{k=1,\ldots,s}) = \bigcup_{\lambda_0 \in \mathbb{C}} sings_\mathcal{S}(\pi_{Z_X}^{-1}(\lambda_0, \underline{\lambda})),$$
where we see $\pi_{Z_X}^{-1}(\lambda_0, \underline{\lambda})$ as a subset of $X \subset \mathbb{P}(V')$ and where $sings_\mathcal{S}(\pi_{Z_X}^{-1}(\lambda_0, \underline{\lambda}))$ denotes the stratified singular locus with respect to the stratification $\mathcal{S}$ of $X$ by torus orbits defined above.

*Proof.* For the first claim consider the following increasing filtration on $R_B$. Let as above $Q$ be the convex hull of $\underline{b}_1, \ldots, \underline{b}_t$ and $0$ in $\mathbb{R}^s$. Let $u \in \mathbb{N}B$ then the weight of $\underline{y}^u$ is defined by $\inf\{\lambda \in \mathbb{R}_{\geq 0} \mid u \in \lambda \cdot Q\}$. It is easy to see that there is an integer $e$ so that all weights lie in $e^{-1}\mathbb{N}$. Denote by $R_B^{\frac{k}{e}}$ the elements in $R_B$ with weight $\leq k/e$. Let $gr R_B$ be the graduated ring with respect to this filtration. By [Ado94, Equation 5.12] we have
$$\dim_\mathbb{C} gr(R_B)/(\overline{y_k \partial f_{\underline{\lambda}}/\partial y_k})_{k=1,\ldots,s} = \mathrm{vol}(Q),$$
where $\overline{y_k \partial f_{\underline{\lambda}}/\partial y_k}$ is the image of $y_k \partial f_{\underline{\lambda}}/\partial y_k$ in $gr(R_B)$. It remains to show that
$$\dim_\mathbb{C} gr(R_B)/(\overline{y_k \partial f_{\underline{\lambda}}/\partial y_k})_{k=1,\ldots,s} = \dim_\mathbb{C} R_B/(y_k \partial f_{\underline{\lambda}}/\partial y_k)_{k=1,\ldots,s}\,.$$

The proof of this equality is an easy adaptation of the proof of [Ado94, Theorem 5.4].

For the proof of the second statement we notice first that

$$sing_{\mathcal{S}}(\pi_{Z_X}^{-1}(\lambda_0, \underline{\lambda})) = \bigcup_{\Gamma | 0 \in \Gamma} sing(X^0(\Gamma) \cap H_{(\lambda_0, \underline{\lambda})})$$

because the fiber over $(\lambda_0, \underline{\lambda})$ has no bad singularities at infinity.

Define the following $r$ hyperplanes $H_{\underline{\lambda}}^k$ for $k \in \{1, \ldots, s\}$ and $\underline{\lambda} \in W^\circ$:

$$H_{\underline{\lambda}}^k := \{(\mu_0 : \ldots : \mu_t) \in \mathbb{P}(V') \mid \sum_{i=1}^{t} b_{ki} \lambda_i \mu_i = 0\}.$$

We have $sing(X^0(\Gamma) \cap H_{(\lambda_0, \underline{\lambda})}) = X^0(\Gamma) \cap H_{(\lambda_0, \underline{\lambda})} \cap (\bigcap_{k=1}^{s} H_{\underline{\lambda}}^k)$ by equation (41) and therefore

$$sing_{\mathcal{S}}(\pi_{\mathcal{Z}_X}^{-1}(\lambda_0, \underline{\lambda})) = X^{aff} \cap H_{(\lambda_0, \underline{\lambda})} \cap (\bigcap_{k=1}^{s} H_{\underline{\lambda}}^k).$$

Notice that

$$\bigcup_{\lambda_0 \in \mathbb{C}} (X^{aff} \cap H_{(\lambda_0, \underline{\lambda})} \cap (\bigcap_{k=1}^{s} H_{\underline{\lambda}}^k)) = \bigcup_{\lambda_0 \in \mathbb{C}} supp(R_B/R_B(f_{\underline{\lambda}} - \lambda_0) + R_B(\partial f_{\underline{\lambda}}/\partial y_k)_{k=1,\ldots,s})$$
$$= supp(R_B/R_B(\partial f_{\underline{\lambda}}/\partial y_k)_{k=1,\ldots,s}),$$

which shows the claim. □

Let $\widetilde{B}$ be the $(s+1) \times (t+1)$-matrix as introduced before Definition 2.8. Let $\widetilde{Q}$ be the convex hull of $\underline{\widetilde{b}}_0, \ldots, \underline{\widetilde{b}}_t$ in $\mathbb{R}^{s+1}$. Notice that $\widetilde{Q} \subset \{1\} \times \mathbb{R}^s$ and therefore no face $\widetilde{\Gamma}$ of $\widetilde{Q}$ contains the origin. Adolphson characterized the characteristic variety $char(\mathcal{M}_{\widetilde{B}}^{\widetilde{\beta}})$ of the GKZ system $\mathcal{M}_{\widetilde{B}}^{\widetilde{\beta}}$ as follows. Let $T^*V \simeq V \times V'$ be the holomorphic cotangent bundle with coordinates $(\lambda_0, \ldots, \lambda_t, \mu_0, \ldots, \mu_t)$. Define the following Laurent polynomials on $(\mathbb{C}^*)^{s+1}$

$$\widetilde{f}_{\underline{\widetilde{\lambda}}}(\underline{y}) := \widetilde{f}_{\underline{\widetilde{\lambda}}, \widetilde{Q}}(\underline{y}) := \sum_{i=0}^{t} \lambda_i \underline{y}^{\widetilde{b}_i},$$
$$\widetilde{f}_{\underline{\widetilde{\lambda}}, \widetilde{\Gamma}}(\underline{y}) := \sum_{\underline{\widetilde{b}}_i \in \widetilde{\Gamma}} \lambda_i \underline{y}^{\widetilde{b}_i},$$

where we define $\underline{y}^{\widetilde{b}_i} := \prod_{k=0}^{r} y_k^{\widetilde{b}_{ki}}$.

**Lemma 3.11** ([Ado94] Lemma 3.2, Lemma 3.3)**.**

1. For each $(\underline{\widetilde{\lambda}}^{(0)}, \underline{\widetilde{\mu}}^{(0)}) \in char(\mathcal{M}_{\widetilde{A}}^{\widetilde{\beta}})$ there exists a (possibly empty) face $\widetilde{\Gamma}$ such that $\widetilde{\mu}_j^{(0)} \neq 0$ if and only if $\underline{\widetilde{b}}_j \in \widetilde{\Gamma}$.

2. If $\underline{\widetilde{\lambda}}^{(0)}$ is a singular point of $\mathcal{M}_{\widetilde{B}}^{\widetilde{\beta}}$ and $\widetilde{\Gamma}$ the corresponding (non-empty) face, then the Laurent polynomials $\partial \widetilde{f}_{\underline{\widetilde{\lambda}}^{(0)}, \widetilde{\Gamma}}/\partial y_0, \ldots, \partial \widetilde{f}_{\underline{\widetilde{\lambda}}^{(0)}, \widetilde{\Gamma}}/\partial y_s$ have a common zero in $(\mathbb{C}^*)^{s+1}$.

We can use this result in the next lemma to compute the singular locus of the $\mathcal{D}$-modules we are interested in.

**Lemma 3.12.** The singular locus of $\mathcal{M}_{\widetilde{B}}^{\widetilde{\beta}}$ as well as the singular locus of the modules $\mathcal{H}^0(\varphi_{B+}\mathcal{O}_{S \times W})$ resp. $\mathcal{H}^0(\varphi_{B\dagger}\mathcal{O}_{S \times W})$ is given by

$$\Delta_S := \Delta_B \cup \Delta_B^\infty.$$

*Proof.* Notice that the polytope $\widetilde{Q} \subset \{1\} \times \mathbb{R}^s$ is just the shifted polytope $Q \subset \mathbb{R}^s$ defined above. One easily sees that the Laurent polynomials $\partial \widetilde{f}_{\underline{\widetilde{\lambda}}^{(0)}, \widetilde{Q}}/\partial y_0, \ldots, \partial \widetilde{f}_{\underline{\widetilde{\lambda}}^{(0)}, \widetilde{Q}}/\partial y_s$ have a common zero in $(\mathbb{C}^*)^{s+1}$ if and only if $\varphi_B^{-1}(\underline{\widetilde{\lambda}}^{(0)})$ is singular, i.e. the set of $\underline{\widetilde{\lambda}}^{(0)}$'s which satisfy this condition is exactly the discriminant $\Delta_B$ of $\varphi_B$. If there exists a proper face $\widetilde{\Gamma}$ of $\widetilde{Q}$ such that the Laurent polynomials $\partial \widetilde{f}_{\underline{\widetilde{\lambda}}^{(0)}, \widetilde{\Gamma}}/\partial y_0, \ldots, \partial \widetilde{f}_{\underline{\widetilde{\lambda}}^{(0)}, \widetilde{\Gamma}}/\partial y_s$ have a common zero in $(\mathbb{C}^*)^{s+1}$, then then fiber $\varphi_B^{-1}(\underline{\widetilde{\lambda}}^{(0)})$ has a singularity at infinity, i.e. its compactification has a singularity in $X^0(\Gamma)$, where $\Gamma$ is the corresponding face of $Q$. $\square$

**Lemma 3.13.** *The restriction of the discriminant $\Delta_S$ to $\mathbb{C} \times W^\circ \subset V$ is finite over $W^\circ \subset W$.*

*Proof.* We will first show quasi-finiteness of the map $p : \Delta_{S|\mathbb{C} \times W^\circ} \to W^\circ$. First notice that we have $\Delta_{S|\mathbb{C} \times W^\circ} = (\Delta_S \setminus \Delta_B^{bad})_{|\mathbb{C} \times W^\circ}$. Fix some $\underline{\lambda} \in W^\circ$. We have to show that $\Delta_{S|\mathbb{C} \times \{\underline{\lambda}\}}$ is a finite set. By the definition of $\Delta_S$ it is enough to show that $sing_{\mathcal{S}}(\pi_{Z_X}^{-1}(\lambda_0, \underline{\lambda}))$ is a finite set, but this is Lemma 3.10.

To prove finiteness of the map $p : \Delta_{S|\mathbb{C} \times W^\circ} \to W^\circ$ it remains to show that it is proper. Let $K$ be any compact subset of $W^\circ$. Suppose that $p^{-1}(K)$ is not compact, then it must be unbounded in $V \simeq \mathbb{C}^{t+1}$ for the standard metric. Hence there is a sequence $(\lambda_0^{(i)}, \underline{\lambda}^{(i)}) \in p^{-1}(K)$ with $\lim_{i \to \infty} |\lambda_0^{(i)}| = \infty$, as $K$ is closed and bounded in $W^\circ \subset W = \mathbb{C}^t$.

In order to construct a contradiction, we use the partial compactification of the family $\varphi_B$ from above. Recall the spaces $Z := \{\sum_{i=0}^t \lambda_i \cdot \mu_i = 0\} \subset \mathbb{P}(V') \times V$ and $Z_X := (X \times V) \cap Z$. Introduce the spaces $Z_k := \{\sum_{i=1}^t b_{ki} \lambda_i \mu_i = 0\}$ for $k \in \{1, \ldots, t\}$. Then $Z_X \cap (\bigcap_{k=1}^d Z_k)$ is the stratified critical locus $crit_{\mathcal{S}}(\pi_{Z_X})$ of the family $\pi_{Z_X}$, where we denote by abuse of notation by $\mathcal{S}$ also the stratification on $Z_X$ induced from the torus stratification on $X$ used above.

Because the projection from the stratified critical locus $crit_{\mathcal{S}}(\pi_{Z_X})$ of $\pi_{Z_X}$ to $\Delta_S$ is onto, there is a sequence $((\mu_0^{(i)} : \underline{\mu}^{(i)}), (\lambda_0^{(i)}, \underline{\lambda}^{(i)})) \in X^{aff} \times p^{-1}(K)$ projecting under $\pi_{Z_X|X^{aff} \times p^{-1}(K)}$ to $(\lambda_0^{(i)}, \underline{\lambda}^{(i)})$ (Notice that we consider here $X^{aff}$ as a subset of $\mathbb{P}(V')$ under the embedding $i \circ j_1$). Consider the first component of the sequence $((\mu_0^{(i)} : \underline{\mu}^{(i)}), (\lambda_0^{(i)}, \underline{\lambda}^{(i)}))$, then this is a sequence $(\mu_0^{(i)} : \underline{\mu}^{(i)})$ in $X$ which converges (after possibly passing to a subsequence) to a limit $(0 : \mu_1^{\lim} : \ldots : \mu_t^{\lim})$ (this is forced by the incidence relation $\sum_{i=0}^t \lambda_i \mu_i$). In other words this limit lies in $X \setminus X^{aff}$ by the definition of $X^{aff}$ (see Definition 3.7). But because $(X \times V) \cap Z \cap \bigcap_{k=1}^d Z_k = Z_X \cap \bigcap_{k=1}^d Z_k$ is closed, the point $((0 : \mu_1^{\lim} : \ldots : \mu_t^{\lim}), (\lambda_0^{\lim}, \underline{\lambda}^{\lim}))$ lies in $((X \setminus X^{aff}) \times p^{-1}(K)) \cap Z \cap \bigcap_{k=1}^d Z_k$, i.e. $\pi_{Z_X}^{-1}(\lim_{i \to \infty}(\lambda_0^{(i)}, \underline{\lambda}^{(i)}))$ has a bad singularity at infinity, which is a contradiction by the definition of $W^\circ$. $\square$

We can now prove the following regularity property of $\widehat{\mathcal{M}}_B^{(\beta_0, \beta)}$, which is essentially the same proof as in [RS15, Lemma 4.4].

**Lemma 3.14.** *Consider $\widehat{\mathcal{M}}_B^{(\beta_0, \beta)}$ as a $\mathcal{D}_{\mathbb{P}^1 \times \overline{W}}$-module, where $\overline{W}$ is a smooth projective compactification of $W$. Then $\widehat{\mathcal{M}}_B^{(\beta_0, \beta)}$ is regular outside $(\{z = 0\} \times W) \cup (\mathbb{P}^1_z \times (\overline{W} \setminus W^\circ))$ and smooth on $\mathbb{C}_z^* \times W^\circ$.*

*Proof.* It suffices to show that any $\underline{\lambda} = (\lambda_1, \ldots, \lambda_t) \in W^\circ$ has a small analytic neighborhood $W_{\underline{\lambda}}^\circ \subset W^{\circ an}$ such that the partial analytization $\mathcal{O}_{W_{\underline{\lambda}}^\circ}^{an}[\tau, \tau^{-1}] \otimes_{\mathcal{O}_{\mathbb{C}_\tau^* \times W^\circ}} \widehat{\mathcal{M}}_B^{(\beta_0, \beta)}$ is regular on $\mathbb{C}_\tau \times W_{\underline{\lambda}}^\circ$ (but not at $\tau = \infty$). This is precisely the statement of [DS03, Theorem 1.11 (1)], taking into account the regularity of $\mathcal{M}_{\widetilde{B}}^{\widetilde{\beta}}$ (c.f. [Hot98, section 6]), the fact that the singular locus of $\mathcal{M}_{\widetilde{B}}^{\widetilde{\beta}}$ coincides with $\Delta_S$ (see Lemma 3.12) as well as the last lemma (notice that the non-characteristic assumption in [DS03, Theorem 1.11 (1)] is satisfied, see, e.g., [Pha79, page 281]). $\square$

The next step is to study several natural lattices in $\widehat{\mathcal{M}}_B^{(\beta_0, \beta)}$. They are defined in terms of $\mathcal{R}$-modules, see the end of subsection 2.1.

**Definition 3.15.** *1. Consider the left ideal*

$$\mathcal{I} := \mathcal{D}_{\mathbb{C}_z \times W^*}(\widehat{\square}_{\underline{l}})_{\underline{l} \in \mathbb{L}} + \mathcal{D}_{\mathbb{C}_z \times W^*}(\widehat{E}_k - z \cdot \beta_k)_{k=1,\ldots,r} + \mathcal{D}_{\mathbb{C}_z \times W^*}(\widehat{E} - z \cdot \beta_0)$$

*in $\mathcal{D}_{\mathbb{C}_z \times W^*}$ and write ${}^*\widehat{\mathcal{M}}_B^{(\beta_0,\beta)}$ for the cyclic $\mathcal{D}$-module $\mathcal{D}_{\mathbb{C}_z \times W^*}/\mathcal{I}$. Here the operators $\widehat{\square}_{\underline{l}}$, $\widehat{E}_k$ and $\widehat{E}$ are those from Definition 3.1.*

2. Consider the left ideal $_0\mathcal{I} := \mathcal{R}_{\mathbb{C}_z \times W^*}(\widehat{\square}_{\underline{l}})_{\underline{l} \in \mathbb{L}} + \mathcal{R}_{\mathbb{C}_z \times W^*}(\widehat{E}_k - z \cdot \beta_k)_{k=1,\ldots,r} + \mathcal{R}_{\mathbb{C}_z \times W^*}(\widehat{E} - z \cdot \beta_0)$ in $\mathcal{R}_{\mathbb{C}_z \times W^*}$ and write $_0^*\widehat{\mathcal{M}}_B^{(\beta_0,\beta)}$ for the cyclic $\mathcal{R}$-module $\mathcal{R}_{\mathbb{C}_z \times W^*}/\mathcal{I}$.

3. Consider the open inclusions $W^\circ \subset W^* \subset W$ and define $^\circ\mathcal{R} := \mathcal{R}_{|\mathbb{C}_z \times W^\circ}$ with ring of global sections $^\circ R$. Define the $\mathcal{D}_{\mathbb{C}_z \times W^\circ}$-module

$$^\circ\widehat{\mathcal{M}}_B^{(\beta_0,\beta)} := \left(\widehat{\mathcal{M}}_B^{(\beta_0,\beta)}\right)_{|\mathbb{C}_z \times W^\circ}$$

and the $\mathcal{R}_{\mathbb{C}_z \times W^\circ}$-module

$$_0^\circ\widehat{\mathcal{M}}_B^{(\beta_0,\beta)} := \left(_0^*\widehat{\mathcal{M}}_B^{(\beta_0,\beta)}\right)_{|\mathbb{C}_z \times W^\circ}.$$

**Remark 3.16.**

1. We have $\mathcal{D}_{\mathbb{C}_z \times W^*} \otimes_{\mathcal{R}_{\mathbb{C}_z \times W^*}} {}_0^*\widehat{\mathcal{M}}_B^{(\beta_0,\beta)} = \widehat{\mathcal{M}}_B^{(\beta_0,\beta)}{}_{|\mathbb{C}_z \times W^*}$.

2. The restriction of $_0^*\widehat{\mathcal{M}}_B^{(\beta_0,\beta)}$ to $\mathbb{C}_z^* \times W^*$ equals the restriction of $\widehat{\mathcal{M}}_B^{(\beta_0,\beta)}$ to $\mathbb{C}_z^* \times W^*$.

3. For $_{z^2\partial_z}(_0^*\widehat{\mathcal{M}}_B^{(\beta_0,\beta)}) = \mathcal{R}'/_0\mathcal{I}'$, where $_0\mathcal{I}'$ is given by

$$_0\mathcal{I}' := \mathcal{R}'(\widehat{\square}_{\underline{l}})_{\underline{l} \in \mathbb{L}} + \mathcal{R}'(\widehat{E}_k - z \cdot \beta_k)_{k=1,\ldots,r}.$$

**Lemma 3.17.** *The quotient $_0^*\widehat{\mathcal{M}}_B^{(\beta_0,\beta)}/z \cdot {}_0^*\widehat{\mathcal{M}}_B^{(\beta_0,\beta)}$ is the sheaf of commutative $\mathcal{O}_{W^*}$-algebras associated to*

$$\frac{\mathbb{C}[\lambda_1^\pm,\ldots,\lambda_t^\pm,\kappa_1,\ldots,\kappa_t]}{(\prod_{l_i<0}\kappa_i^{-l_i} - \prod_{l_i>0}\kappa_i^{l_i})_{\underline{l}\in\mathbb{L}} + (\sum_{i=1}^t b_{ki}\lambda_i\kappa_i)_{k=1,\ldots,s}} \simeq \frac{\mathbb{C}[\mathbb{N}B][\lambda_1^\pm,\ldots,\lambda_t^\pm]}{y_k\partial f_{\underline{\lambda}}/\partial y_k}, \quad (42)$$

*where $y_k\partial f_{\underline{\lambda}}/\partial y_k = \sum_{i=1}^t b_{ki}\lambda_i \underline{y}^{b_i}$.*

*Proof.* Let $\kappa_i$ be the class of $z\partial_{\lambda_i}$. Because the commutator $[\kappa_i, \lambda_i]$ is zero we see that $_0^*\widehat{\mathcal{M}}_B^{(\beta_0,\beta)}/z \cdot {}_0^*\widehat{\mathcal{M}}_B^{(\beta_0,\beta)}$ is a commutative algebra and isomorphic to the module on the left hand side of equation (42). To show the isomorphism (42), consider the $\mathbb{C}[\lambda_1^\pm,\ldots,\lambda_t^\pm]$-linear morphism

$$\psi : \mathbb{C}[\lambda_1^\pm,\ldots,\lambda_t^\pm,\kappa_1,\ldots,\kappa_t] \longrightarrow \mathbb{C}[\mathbb{N}B][\lambda_1^\pm,\ldots,\lambda_t^\pm],$$
$$\kappa_i \mapsto \underline{y}^{b_i}$$

which is surjective by the definition of $\mathbb{C}[\mathbb{N}B]$. The kernel of this map is equal to $(\prod_{l_i<0}\kappa_i^{-l_i} - \prod_{l_i>0}\kappa_i^{l_i})_{\underline{l}\in\mathbb{L}}$ by [MS05, Theorem 7.3]. Finally notice that $\psi(\sum_{i=1}^t b_{ki}\lambda_i\kappa_i) = y_k\partial f_{\underline{\lambda}}/\partial y_k$, which shows the claim. $\square$

We need the following result saying that the GKZ-system $\mathcal{M}_B^\beta$ is isomorphic to the restriction of the Fourier-Laplace transformed GKZ system $\widehat{\mathcal{M}}_B^{(\beta_0,\beta)}$.

**Lemma 3.18.** *Let $i_1 : \{1\} \times W \longrightarrow \hat{V} = \mathbb{C}_z \times W$ be the canonical inclusion. Then*

$$\mathcal{H}^0\left(i_1^+\widehat{\mathcal{M}}_B^{(\beta_0,\beta)}\right) \simeq \mathcal{M}_B^\beta.$$

*Proof.* During the proof we will work with modules of global sections rather with the $\mathcal{D}$-modules itself. Recall that the left ideal defining the quotient $\widehat{M}_B^{(\beta_0,\beta)}$ is generated by the operators $\widehat{\square}_{\underline{l}}, \widehat{E}_k - \beta_k z$ and $\widehat{E} - \beta_0 z$, where

$$\widehat{\square}_{\underline{l}} := \prod_{i: l_i < 0}(z \cdot \partial_{\lambda_i})^{-l_i} - \prod_{i: l_i > 0}(z \cdot \partial_{\lambda_i})^{l_i},$$

$$\widehat{E} := z^2\partial_z + \sum_{i=1}^t z\lambda_i\partial_{\lambda_i},$$

$$\widehat{E}_k := \sum_{i=1}^t b_{ki}z\lambda_i\partial_{\lambda_i}.$$

The presence of $z^{-2}(\widehat{E}_0 - \beta_0 z)$ in this ideal show that have the an isomorphism of $\mathbb{C}[z^\pm, \lambda_1, \ldots, \lambda_n]\langle \partial_{\lambda_1}, \ldots \partial_{\lambda_n}\rangle$-modules

$$\widehat{M} \simeq \mathbb{C}[z^\pm, \lambda_1, \ldots, \lambda_n]\langle \partial_{\lambda_1}, \ldots \partial_{\lambda_n}\rangle / \mathbb{C}[z^\pm, \lambda_1, \ldots, \lambda_n]\langle \partial_{\lambda_1}, \ldots \partial_{\lambda_n}\rangle \widehat{I} \tag{43}$$

where the left $\mathbb{C}[z^\pm, \lambda_1, \ldots, \lambda_n]\langle \partial_{\lambda_1}, \ldots \partial_{\lambda_n}\rangle$-ideal $\widehat{I}$ is generated by $\widehat{\square}_{\underline{l} \in \mathbb{L}}$ and $\widehat{E}_k - \beta_k$ for $k \in \{1, \ldots, d\}$. The $D_W$-module corresponding to $\mathcal{H}^0\left(i_1^+ \widehat{\mathcal{M}}\right)$ is given by $\widehat{M}/(z-1)\widehat{M}$. Using the isomorphism (43) one easily sees that

$$\widehat{M}/(z-1)\widehat{M} \simeq M_B^\beta,$$

which shows the claim. $\square$

**Proposition 3.19.** *The $\mathcal{O}_{\mathbb{C}_z \times W^\circ}$-module ${}_0^\circ\widehat{\mathcal{M}}_B^{(\beta_0, \beta)}$ is locally-free of rank* vol$(Q)$.

*Proof.* Notice that it is sufficient to show that ${}_0^\circ\widehat{\mathcal{M}}_B^{(\beta_0, \beta)}$ is $\mathcal{O}_{\mathbb{C} \times W^\circ}$-coherent. Namely, ${}_0^\circ\widehat{\mathcal{M}}_B^{(\beta_0, \beta)}/z \cdot {}_0^\circ\widehat{\mathcal{M}}_B^{(\beta_0, \beta)}$ is $\mathcal{O}_{W^\circ}$-locally free of rank vol$(Q)$ by Lemma 3.10. Moreover, the restriction of ${}_0^\circ\widehat{\mathcal{M}}_B^{(\beta_0, \beta)}$ to $\mathbb{C}_z^* \times W^\circ$ is a locally-free $\mathcal{O}_{\mathbb{C}_z^* \times W^\circ}$-module by Lemma 3.14. Its restriction to $\{1\} \times W^\circ$ is isomorphic to the restriction of $\mathcal{M}_B^\beta$ to $W^\circ$ by Lemma 3.18 which is locally free of rank vol$(Q)$. Now we use the fact that a coherent $\mathcal{O}$-module which has everywhere the same rank is locally-free.

It is actually sufficient to show the coherence of $\mathcal{N} := \text{For}_{z^2 \partial_z}({}_0^\circ\widehat{\mathcal{M}}_B^{(\beta_0, \beta)})$, as this is the same as ${}_0^\circ\widehat{\mathcal{M}}_B^{(\beta_0, \beta)}$ when considered as an $\mathcal{O}_{\mathbb{C}_z \times W^\circ}$-module. Let us denote by $F_\bullet$ the natural filtration on $\mathcal{R}'_{\mathbb{C}_z \times W^\circ}$ defined by

$$F_k \mathcal{R}'_{\mathbb{C}_z \times W^\circ} := \left\{ P \in \mathcal{R}'_{\mathbb{C}_z \times W^\circ} \mid P = \sum_{|\alpha| \leq k} g_\alpha(z, \underline{\lambda})(z\partial_{\lambda_1})^{\alpha_1} \cdot \ldots \cdot (z\partial_{\lambda_t})^{\alpha_t} \right\}.$$

This filtration induces a filtration $F_\bullet$ on $\mathcal{N}$ which satisfies $F_k \mathcal{R}'_{\mathbb{C}_z \times W^\circ} \cdot F_l \mathcal{N} = F_{k+l} \mathcal{N}$. Obviously, for any $k$, $F_k \mathcal{N}$ is $\mathcal{O}_{\mathbb{C}_z \times W^\circ}$-coherent, so that it suffices to show that the filtration $F_\bullet$ becomes eventually stationary. Let $P = \sum_{|\alpha| \leq k} g_\alpha(z, \underline{\lambda})(z\partial_{\lambda_1})^{\alpha_1} \cdot \ldots \cdot (z\partial_{\lambda_t})^{\alpha_t}$ then its symbol is defined as

$$\sigma_k(P) := \sum_{|\alpha| = k} g_\alpha(z, \underline{\lambda})(\kappa_1)^{\alpha_1} \cdot \ldots \cdot (\kappa_t)^{\alpha_t} \in \mathcal{O}_{\mathbb{C}_z \times W^\circ}[\kappa_1, \ldots, \kappa_t],$$

which is a function on $\mathbb{C}_z \times T^* W^\circ$ with fiber variables $\kappa_1, \ldots, \kappa_t$. Let $\mathcal{I}$ be the radical ideal of the ideal generated by the symbols of $\widehat{\square}_{\underline{l} \in \mathbb{L}}$ and $\widehat{E}_k - z \cdot \beta_k$ for $k = 1, \ldots, t$. Then the vanishing locus of $\mathcal{I}$ is the $\mathcal{R}'_{\mathbb{C}_z \times W^\circ}$-characteristic variety of $\mathcal{N}$. Notice that $\mathcal{N}$ is $\mathcal{O}_{\mathbb{C}_z \times W^\circ}$-coherent if and only if its $\mathcal{R}'_{\mathbb{C}_z \times W^\circ}$-characteristic variety is a subset of $\mathbb{C}_z \times T^*_{W^\circ} W^\circ$. The proof of this fact is completely parallel to the $\mathcal{D}$-module case (see e.g. [Pha79, Proposition 10.3]).

To compute the $\mathcal{R}'_{\mathbb{C}_z \times W^\circ}$-characteristic variety, notice that the symbols of $\widehat{\square}_{\underline{l} \in \mathbb{L}}$ and $\widehat{E}_k - z \cdot \beta_k$ are independent of $z$. Thus it is enough to compute its restriction to $\{1\} \times W^\circ$. Now notice that the generators of the ideal corresponding to the GKZ-system $\mathcal{M}_B^\beta$ have exactly the same symbols as the operators above. Thus it is enough to show that the restriction of the GKZ-system $\mathcal{M}_B^\beta$ to $W^\circ$ is $\mathcal{O}_{W^\circ}$-coherent. But this follows from [Ado94, Lemma 3.2 and 3.3] and the definition of $W^\circ$ (see Definition 3.8 and Lemma 3.12). $\square$

**Corollary 3.20.** *The natural map ${}_0^\circ\widehat{\mathcal{M}}_B^{(\beta_0, \beta)} \to {}^\circ\widehat{\mathcal{M}}_B^{(\beta_0, \beta)}$ which is induced by the inclusion $\mathcal{R}_{\mathbb{C}_z \times W^*} \to \mathcal{D}_{\mathbb{C}_z \times W^*}$ is injective.*

*Proof.* Recall that $\mathcal{D}_{\mathbb{C}_z \times W^*} \otimes_{\mathcal{R}} {}_0\widehat{\mathcal{M}}_B^{(\beta_0, \beta)} \simeq \widehat{\mathcal{M}}_B^{(\beta_0, \beta)}{}_{|\mathbb{C}_z \times W^*}$ and $D_{\mathbb{C}_z \times W^*} \simeq R[z^\pm]$. Thus the kernel of ${}_0\widehat{\mathcal{M}}_B^{(\beta_0, \beta)} \to \widehat{\mathcal{M}}_B^{(\beta_0, \beta)}{}_{|\mathbb{C}_z \times W^*}$ has $z$-torsion. On the open set $\mathbb{C}_z \times W^\circ \subset \mathbb{C}_z \times W^*$ the module ${}_0^\circ\widehat{\mathcal{M}}_B^{(\beta_0, \beta)} = {}_0\widehat{\mathcal{M}}_B^{(\beta_0, \beta)}{}_{|\mathbb{C}_z \times W^\circ}$ is $\mathcal{O}_{\mathbb{C}_z \times W^\circ}$-locally free. In particular it has no $z$-torsion, but this shows the claim. $\square$

In order to do this, consider once again the affine toric variety $X^{\text{aff}} = \text{Spec}\,(\mathbb{C}[\mathbb{N}B])$ from Definition 3.7, which contains the torus $g_B(S) \cong S$ as an open subset. Denote by $D$ the complement of $S$ in $X^{\text{aff}}$. We will consider $X^{\text{aff}}$ as a log scheme in the sense of logarithmic geometry (see, e.g., [Gro11]). More precisely, we endow $X^{\text{aff}}$ with divisorial log structure induced by $D$ and $W^*$ with the trivial log structure. We consider the relative log de Rham complex $\Omega^\bullet_{X^{\text{aff}} \times W^*/W^*}(\log D)$ ([Gro11, section 3.3]). We have isomorphisms $\Omega^k_{X^{\text{aff}} \times W^*/W^*}(\log D) \cong \mathcal{O}_{X^{\text{aff}} \times W^*} \otimes_{\mathbb{Z}} \bigwedge^k \mathbb{Z}^r$.

**Proposition 3.21.** *Let $\mathbb{N}B$ be a saturated semigroup. There exists the following $R_{\mathbb{C}_z \times W^\circ}$-linear isomorphism*
$$H^0\left(\Omega^{\bullet+s}_{X^{aff} \times W^\circ/W^\circ}(\log D)[z], zd - d_yF\wedge\right) \cong {}_0\widehat{M}^{(0,\underline{0})}_B,$$
*which maps $\omega_0$ to 1.*

*Proof.* We first define the $R_{\mathbb{C}_z \times W}$-linear morphism
$$\psi : {}_0\widehat{M}^{(0,\underline{0})}_B \longrightarrow H^0\left(\Omega^{\bullet+s}_{X^{aff} \times W^*/W^*}(\log D)[z], zd - d_yF\wedge\right),$$
$$1 \mapsto \omega_0,$$

which is well-defined by 3.4. Let
$$\omega = \sum_{\alpha,\gamma,\delta} c_{\alpha\gamma\delta} \lambda_1^{\gamma_1} \ldots \lambda_t^{\gamma_t} z^\delta \underline{y}^{\alpha_1 \cdot \underline{b}_1} \ldots \underline{y}^{\alpha_t \cdot \underline{b}_t} \omega_0$$

be a general element in $\Omega^s_{X^{aff} \times W^*/W^*}(\log D)[z]$ with $\alpha \in \mathbb{N}^t$, $\gamma \in \mathbb{Z}^t$ and $\delta \in \mathbb{N}$. Then
$$\sum_{\alpha,\gamma,\delta} c_{\alpha\gamma\delta} \lambda_1^{\gamma_1} \ldots \lambda_t^{\gamma_t} z^\delta (z\partial_{\lambda_1})^{\alpha_1} \ldots (z\partial_{\lambda_t})^{\alpha_t}$$

is a preimage, which shows that the map $\psi$ is surjective. Notice that the restricted map
$${}^\circ\psi : {}_0\widehat{\circ M}^{(0,\underline{0})}_B \longrightarrow H^0(\Omega^{\bullet+s}_{X^{aff} \times W^\circ/W^\circ}(\log D)[z], zd - d_yF\wedge)$$

is also surjective. Consider the following commutative diagram

$$\begin{array}{ccc}
\widehat{\circ M}^{(0,\underline{0})}_B & \xrightarrow{\simeq} & H^0(\Omega^{\bullet+s}_{S \times W^\circ/W^\circ}[z^\pm], zd - d_yF\wedge) \\
\uparrow & & \uparrow \\
{}_0\widehat{\circ M}^{(0,\underline{0})}_B & \xrightarrow{{}^\circ\psi} & H^0(\Omega^{\bullet+s}_{X^{aff} \times W^\circ/W^\circ}(\log D)[z], zd - d_yF\wedge)
\end{array}$$

where the upper horizontal map is an isomorphism by Proposition 3.3 and Lemma 3.4, the left vertical map is injective by Corollary 3.20 and the right vertical map is induced by the morphism
$$\Omega^s_{X^{aff} \times W^\circ/W^\circ}(\log D)[z] \longrightarrow \Omega^s_{X^{aff} \times W^\circ/W^\circ}(*D)[z^\pm] = \Omega^s_{S \times W^\circ/W^\circ}[z^\pm].$$

But this shows that ${}^\circ\psi$ is also injective, which shows the claim. Notice that as a by-product, we also obtain that the morphism
$$H^0(\Omega^{\bullet+s}_{X^{aff} \times W^\circ/W^\circ}(\log D)[z], zd - d_yF\wedge) \longrightarrow H^0(\Omega^{\bullet+s}_{S \times W^\circ/W^\circ}[z^\pm], zd - d_yF\wedge)$$

is injective. $\square$

# 4 Quantum cohomology of toric complete intersections

We recall in this section some rather well known notations and results concerning twisted Gromov-Witten invariants on the one hand, and basic constructions from toric geometry for smooth complete intersections in toric varieties on the other hand. Any of the statements of this section can be found in either the original articles like [Kon95], [Giv98b, Giv98a], [CG07] (for twisted Gromov-Witten invariants), the references [Ful93], [CLS11] and [CK99], (for facts on toric geometry of complete intersections) but also in the more recent paper [MM11], from which we borrow some of the notation. By collecting the material we need later here we hope to make this paper more self-contained.

## 4.1 Twisted and reduced quantum $\mathcal{D}$-modules

A smooth complete intersection inside a smooth projective variety can be described as the zero locus of a generic section of a split vector bundle on that variety. Associated to such a bundle are the **twisted Gromov-Witten invariants**, which we describe first. They give rise to the twisted quantum product, and to the twisted quantum-$\mathcal{D}$-module. From this one can derive (basically by dividing by the kernel of the multiplication by the first Chern classes of the factors of the vector bundle) the **reduced quantum $\mathcal{D}$-module**, which corresponds to the ambient part of the quantum cohomology of the subvariety. We also discuss this reduced module here, and we define pairings (coming from the Poincaré pairing on the ambient variety) on both the twisted and the reduced quantum $\mathcal{D}$-module.

Let $\mathcal{X}$ be a smooth projective $n$-dimensional variety. Let $\mathcal{L}_1, \ldots, \mathcal{L}_c$ be line bundles on $\mathcal{X}$ which are globally generated and define $\mathcal{E} := \bigoplus_{i=1}^c \mathcal{L}_i$. We are going to recall the construction of the so-called twisted quantum $\mathcal{D}$-module $\mathrm{QDM}(\mathcal{X}, \mathcal{E})$ and the reduced quantum $\mathcal{D}$-module $\overline{\mathrm{QDM}}(\mathcal{X}, \mathcal{E})$. Our notation follows the exposition in [MM11, Chapter 2.5].

For $l \in \mathbb{N}$ and $d \in H_2(\mathcal{X}, \mathbb{Z})$ we denote by $\overline{\mathcal{M}}_{0,l,d}(\mathcal{X})$ the moduli space of stable maps of degree $d$ from curves of genus 0 with $l$ marked points to $\mathcal{X}$. Denote by $e_i : \overline{\mathcal{M}}_{0,l,d}(\mathcal{X}) \longrightarrow \mathcal{X}$ the evaluation at the $i$ marked point for $i \in \{1, \ldots, l\}$ and denote by $\pi : \overline{\mathcal{M}}_{0,l+1,d}(\mathcal{X}) \longrightarrow \overline{\mathcal{M}}_{0,l,d}(\mathcal{X})$ the map which forgets the last marked point. Let $\mathcal{E}_{0,l,d}$ be the locally free sheaf $R^0 \pi_* e_{l+1}^* \mathcal{E}$ and for any $j \in \{1, \ldots, l\}$, let $\mathcal{E}_{0,l,d}(j)$ be the kernel of the surjective morphism $\mathcal{E}_{0,l,d} \longrightarrow e_j^* \mathcal{E}$ which evaluates a section at the $j$-marked point. For $i \in \{1, \ldots, l\}$ denote by $\mathcal{N}_i$ the line bundle on $\overline{\mathcal{M}}_{0,l,d}(\mathcal{X})$ whose fiber at a point $(C, x_1, \ldots, x_l, f : C \to \mathcal{X})$ is the cotangent space $T^*_{x_i} C$. Put $\psi_i := c_1(\mathcal{N}_i) \in H^2(\overline{\mathcal{M}}_{0,l,d}(\mathcal{X}))$.

**Definition 4.1.** *Let $l \in \mathbb{N}$, $(m_1, \ldots, m_l \in \mathbb{N}^l)$, $\gamma_1, \ldots, \gamma_l \in H^{2*}(\mathcal{X})$ and $d \in H_2(\mathcal{X}, \mathbb{Z})$. The $j$-th twisted Gromov-Witten invariant with descendants is denoted by*

$$\langle \tau_{m_1}(\gamma_1), \ldots, \widetilde{\tau_{m_j}(\gamma_j)}, \ldots, \tau_{m_l}(\gamma_l) \rangle_{0,l,d} := \int_{[\overline{\mathcal{M}}_{0,l,d}(\mathcal{X})]^{vir}} c_{top}(\mathcal{E}_{0,l,d}(j)) \prod_{i=1}^l \psi_i^{m_i} e_i^* \gamma_i \, ,$$

*where $[\overline{\mathcal{M}}_{0,l,d}(\mathcal{X})]^{vir}$ is the virtual fundamental class of $\overline{\mathcal{M}}_{0,l,d}(\mathcal{X})$. An invariant $\langle \ldots, \gamma_k, \ldots \rangle_{0,l,d}$ has to be understood as $\langle \ldots, \tau_0(\gamma_k), \ldots \rangle_{0,l,d}$. Below we will actually use only such non-descendant (i.e., with all $m_k = 0$) invariants.*

Let $(T_0, T_1, \ldots, T_h)$ be a homogeneous basis of $H^{2*}(\mathcal{X})$ such that $T_0 = 1$ and $T_1, \ldots, T_r$ is a basis of $H^2(\mathcal{X}, \mathbb{Z})$ modulo torsion which lies in the Kähler cone of $\mathcal{X}$. Let $T$ be the torus $H^2(\mathcal{X}, \mathbb{C})/2\pi i H^2(\mathcal{X}, \mathbb{Z})$. Then the basis $T_1, \ldots, T_r$ of $H^2(\mathcal{X}, \mathbb{Z})$ gives rise to coordinates $q = (q_1, \ldots, q_r)$ on $T$.

**Definition 4.2.** *Let $\gamma_1, \ldots, \gamma_2 \in H^{2*}(\mathcal{X}, \mathbb{C})$ and $q \in T$. The twisted small quantum product is defined by*

$$\gamma_1 \bullet^{tw}_q \gamma_2 := \sum_{a=1}^h \sum_{d \in H_2(\mathcal{X}, \mathbb{Z})} q^d \langle \gamma_1, \gamma_2, \widetilde{T}_a \rangle_{0,3,d} T^a \, .$$

where $q^d$ is shorthand notation for $q_1^{\langle T_1, d \rangle}, \ldots, q_r^{\langle T_r, d \rangle}$. Notice that $\langle \gamma_1, \gamma_2, \widetilde{T}_a \rangle_{0,3,d} \neq 0$ only if $d$ lies in the semigroup of effective classes, i.e. $d \in \mathrm{Eff}_\mathcal{X} \subset H_2(\mathcal{X}, \mathbb{Z})$. Hence, by our assumption on the basis $T_1, \ldots, T_r$, only positive powers of the $q_i$ appear in the formula above. Let $\bar{T} = \mathbb{C}^r$ be a partial compactification of $T$ with respect to the coordinates $q_1, \ldots q_r$. In the following we assume that there exists an open neighborhood $\bar{U}$ of $0 \in \bar{T}$ such that the twisted quantum product is convergent on $\bar{U}$ as a power series in $q_1, \ldots, q_r$. The twisted quantum product is associative, commutative and has $T_0$ as a unit (see, e.g., [MM11, Proposition 2.14]).

Put $U := \bar{U} \cap T$. In analogy to the untwisted case one defines a trivial vector bundle $F$ on $H^0(\mathcal{X}) \times \mathbb{C}_z \times U$ with fiber $H^{2*}(\mathcal{X})$ together with a flat meromorphic connection

$$\nabla_{\partial_{t_0}} := \partial_{t_0} + \frac{1}{z} T_0 \bullet^{tw}_q, \quad \nabla_{q_a \partial_{q_a}} := q_a \partial_{q_a} + \frac{1}{z} T_a \bullet^{tw}_q, \quad \nabla_{z \partial_z} := z \partial_z - \frac{1}{z} E \bullet^{tw}_q + \mu,$$

where $\mu$ is the diagonal morphism defined by $\mu(T_a) := \frac{1}{2}(deg(T_a) - (\dim_\mathbb{C} \mathcal{X} - rk\mathcal{E}))T_a$ and $E := t_0 T_0 + c_1(\mathcal{T}_\mathcal{X}) - c_1(\mathcal{E})$ is the so-called Euler field.

Define a twisted pairing on $H^{2*}(\mathcal{X})$ by:
$$(\gamma_1, \gamma_2)^{tw} := \int_{\mathcal{X}} \gamma_1 \cup \gamma_2 \cup c_{top}(\mathcal{E}) \quad \text{for } \gamma_1, \gamma_2 \in H^{2*}(\mathcal{X}).$$

This pairing is degenerate with kernel equal to $\ker m_{c_{top}}$, where $m_{c_{top}}$ is defined by
$$m_{c_{top}} : H^{2*}(\mathcal{X}) \longrightarrow H^{2*}(\mathcal{X}),$$
$$\alpha \mapsto c_{top}(\mathcal{E}) \cup \alpha$$

and satisfies the Frobenius relation:
$$(\gamma_1 \bullet_q^{tw} \gamma_2, \gamma_3)^{tw} = (\gamma_1, \gamma_2 \bullet_q^{tw} \gamma_3)^{tw} \quad \text{for } \gamma_1, \gamma_2, \gamma_3 \in H^{2*}(\mathcal{X}).$$

Denote by $\mathcal{F}$ the sheaf sections of $F$ and define an involution $\iota$ by
$$\iota : H^0(\mathcal{X}) \times \mathbb{C}_z \times U \longrightarrow H^0(\mathcal{X}) \times \mathbb{C}_z \times U,$$
$$(t_0, z, q) \mapsto (t_0, -z, q).$$

We define a $\nabla$-flat sesquilinear pairing
$$S : \iota^*(\mathcal{F}) \times \mathcal{F} \longrightarrow \mathcal{O},$$
$$(s_1, s_2) \mapsto S(s_1, s_2)(t_0, z, q) = (s_1(t_0, -z, q), s_2(t_0, z, q))^{tw}.$$

We call $\overline{H^{2*}(\mathcal{X})} := H^{2*}(\mathcal{X}) / \ker m_{c_{top}}$ the reduced cohomology ring of $(\mathcal{X}, \mathcal{E})$. For $\gamma \in H^{2*}(\mathcal{X})$ denote by $\overline{\gamma}$ its class in $\overline{H^{2*}(\mathcal{X})}$. The pairing $(\cdot, \cdot)^{tw}$ gives rise to a pairing $(\cdot, \cdot)^{red}$ on $\overline{H^{2*}(\mathcal{X})}$ by
$$(\overline{\gamma}_1, \overline{\gamma}_2)^{red} := (\gamma_1, \gamma_2)^{tw} \quad \text{for } \gamma_1, \gamma_2 \in H^{2*}(\mathcal{X}).$$

Because the kernel of $(\cdot, \cdot)^{tw}$ is $\ker m_{c_{top}}$ this pairing is well-defined and non-degenerate. Denote by $\overline{F}$ the trivial bundle on $H^0(\mathcal{X}) \times \mathbb{C}_z \times U$ with fiber $\overline{H^{2*}(\mathcal{X})}$. The pairing $S$ induces a pairing $\overline{S}$ on $\overline{F}$ by
$$\overline{S}(\overline{s}_1, \overline{s}_2) := S(s_1, s_2),$$

which is non-degenerate.

Notice that $\overline{H^{2*}(\mathcal{X})}$ is naturally graded because $m_{c_{top}}$ is a graded morphism. Let $(\phi_0, \ldots, \phi_{s'})$ be a homogeneous basis of $\overline{H^{2*}(\mathcal{X})}$ and denote by $(\phi^0, \ldots, \phi^{s'})$ its dual basis w.r.t. $(\cdot, \cdot)^{red}$. The reduced Gromov-Witten invariants are defined by
$$\langle \overline{\gamma}_1, \ldots, \overline{\gamma}_l \rangle_{0,l,d}^{red} := \langle \gamma_1, \ldots, \widetilde{c_{top}(\mathcal{E})\gamma_l} \rangle_{0,l,d}$$

and the reduced quantum product is
$$\overline{\gamma}_1 \bullet_q^{red} \overline{\gamma}_2 := \sum_{a=0}^{s'} \sum_{d \in H_2(\mathcal{X}, \mathbb{Z})} q^d \langle \overline{\gamma}_1, \overline{\gamma}_2, \phi_a \rangle_{0,3,d}^{red} \phi^a,$$

where the restriction is compatible with the multiplication, i.e.
$$\overline{\gamma_1 \bullet_q^{tw} \gamma_2} = \overline{\gamma}_1 \bullet_q^{red} \overline{\gamma}_2.$$

The bundle $\overline{F}$ carries the following connection:
$$\overline{\nabla}_{\partial_{t_0}} := \partial_{t_0} + \frac{1}{z}\overline{T}_0 \bullet_q^{red}, \quad \overline{\nabla}_{q_a \partial_{q_a}} + \frac{1}{z}\overline{T}_a \bullet_q^{red}, \quad \overline{\nabla}_{z\partial_z} := z\partial_z - \frac{1}{z}\overline{E} \bullet_q^{red} + \overline{\mu},$$

where $\overline{\mu}$ is the diagonal morphism defined by $\overline{\mu}(\phi_A) := \frac{1}{2}(deg(\phi_a) - (dim_{\mathbb{C}}\mathcal{X} - rk\mathcal{E}))\phi_a$ and $\overline{E} := t_0\overline{T}_0 + \overline{c_1(\mathcal{T}_\mathcal{X})} - \overline{c_1(\mathcal{E})}$. One can show that $\overline{\nabla}$ is flat and $\overline{S}$ is $\overline{\nabla}$-flat.

**Definition 4.3.** *Consider the above situation of a smooth projective variety $\mathcal{X}$ and globally generated line bundles $\mathcal{L}_1, \ldots, \mathcal{L}_c$.*

1. *The triple $(F, \nabla, S)$ is called the twisted quantum $\mathcal{D}$-module* $\mathrm{QDM}(\mathcal{X}, \mathcal{E})$.

2. *The triple $(\overline{F}, \overline{\nabla}, \overline{S})$ is called the reduced quantum $\mathcal{D}$-module* $\overline{\mathrm{QDM}(\mathcal{X}, \mathcal{E})}$.

## 4.2 Toric geometry of complete intersection subvarieties

In this subsection we consider the case where the variety $\mathcal{X}$ from above is toric. It will be denoted by $X_\Sigma$, where $\Sigma$ is the defining fan (see below). We recall some well-known results on the toric description of the total space of the bundle $\mathcal{E}$ resp. its dual, on Picard groups, Kähler cones etc. All this is needed in section 6 below.

Let, as usual, $N$ be a free abelian group of rank $n$ for which we choose once and for all a basis which identifies it with $\mathbb{Z}^n$. Let $\Sigma$ be a complete smooth fan in $N_\mathbb{R} := N \otimes \mathbb{R}$ and $X_\Sigma$ the associated toric variety, which is compact and smooth. We recall the toric description of the Kähler resp. the nef cone of $\Sigma$. Let $\Sigma(1) = \{\mathbb{R}_{\geq 0}\underline{a}_1, \ldots, \mathbb{R}_{\geq 0}\underline{a}_m\}$ be the rays of $\Sigma$, where $\underline{a}_i \in N \cong \mathbb{Z}^n$ are the primitive integral generators of the rays of $\Sigma$. Then we have an exact sequence

$$0 \longrightarrow \mathbb{L}_A \longrightarrow \mathbb{Z}^{\Sigma(1)} = \mathbb{Z}^m \longrightarrow N \longrightarrow 0, \tag{44}$$

where the morphism $\mathbb{Z}^m \twoheadrightarrow N$ is given by the matrix (henceforth called $A$) having the vectors $\underline{a}_1, \ldots, \underline{a}_m$ as columns. $\mathbb{L}_A$ is the module of relations between these vectors. We also consider the dual sequence

$$0 \longrightarrow M \longrightarrow (\mathbb{Z}^{\Sigma(1)})^\vee = \mathbb{Z}^m \longrightarrow \mathbb{L}_A^\vee \longrightarrow 0,$$

where $M := N^\vee$ is the dual lattice. It is well known that as $X_\Sigma$ is smooth and compact, we have

$$H^2(X_\Sigma, \mathbb{Z}) \simeq \mathrm{Pic}(X_\Sigma) \cong \mathbb{L}_A^\vee,$$

moreover, the group $(\mathbb{Z}^{\Sigma(1)})^\vee$ is the free abelian group generated by the torus invariant divisors on $X_\Sigma$. We denote these generators by $D_1, \ldots, D_m$. Its images in $\mathbb{L}_A^\vee$ (called $\overline{D}_i$) are thus the cohomology classes which are Poincaré dual to these divisors, and they generate the Picard group.

Any element in $(\mathbb{Z}^{\Sigma(1)})^\vee \otimes \mathbb{R}$ can be considered as a function on $N_\mathbb{R}$ (actually on the support of $\Sigma$, but this equals $N_\mathbb{R}$ by completeness), which is linear on each cone of $\Sigma$, these are called piecewise linear functions with respect to $\Sigma$. For a given divisor $D \in \mathrm{Div}(X_\Sigma) \cong (\mathbb{Z}^{\Sigma(1)})^\vee$, we denote the piecewise linear function it corresponds to by $\psi_D^\Sigma$. Inside $(\mathbb{Z}^{\Sigma(1)})^\vee \otimes \mathbb{R}$ we have the cone of convex functions, which are those piecewise linear functions $\psi$ having the property that for any cone $\sigma \in \Sigma$ and for any $n \in N_\mathbb{R}$, we have $\psi(n) \leq \psi_\sigma(n)$, where $\psi_\sigma$ is the extension to a linear function on all of $N_\mathbb{R}$ of the restriction $\psi_{|\sigma}$. The interior of the cone of convex functions are those which are strictly convex, that is, those such that the above inequality is strict on $N_\mathbb{R} \backslash \sigma$. Notice that any linear function on $N$ is piecewise linear and this inclusion is precisely given by $M_\mathbb{R} \hookrightarrow (\mathbb{Z}^{\Sigma(1)})^\vee \otimes \mathbb{R}$. We define the nef cone $\mathcal{K}_{X_\Sigma}$ of $X_\Sigma$ to be the image of the cone of convex functions in $(\mathbb{Z}^{\Sigma(1)})^\vee \otimes \mathbb{R}$ under the projection $(\mathbb{Z}^{\Sigma(1)})^\vee \otimes \mathbb{R} \twoheadrightarrow \mathbb{L}_A^\vee \otimes \mathbb{R}$. Its interior is the Kähler cone $\mathcal{K}_{X_\Sigma}^\circ$ of $\Sigma$. We assume that $\mathcal{K}_{X_\Sigma}^\circ$ is non-empty, which amounts to say that $X_\Sigma$ is projective. Let us recall the following description of the cone $\mathcal{K}_{X_\Sigma}$, the proof of this fact can be found, e.g., in [CK99, section 3.4.2].

**Lemma 4.4.** *For any cone $\sigma \in \Sigma$, put*

$$J_\sigma := \{i \in \{1, \ldots, m\} \,|\, \mathbb{R}_{\geq 0}\underline{a}_i \notin \sigma\}$$

*and define*

$$\check{\sigma} := \sum_{i \in J_\sigma} \mathbb{R}_{\geq 0} \overline{D}_i \subset (\mathbb{L}_A^\vee)_\mathbb{R}.$$

*We call $\check{\sigma}$ the anticone associated to $\sigma$. Then we have $\mathcal{K}_{X_\Sigma} = \bigcap_{\sigma \in \Sigma} \check{\sigma} \subset (\mathbb{L}_A^\vee)_\mathbb{R}$.*

We proceed by considering the toric analogue of the situation from subsection 4.1. More precisely, let $\mathcal{L}_1 = \mathcal{O}_{X_\Sigma}(L_1), \ldots, \mathcal{L}_c = \mathcal{O}_{X_\Sigma}(L_c)$ be line bundles on $X_\Sigma$ with $L_1, \ldots, L_c \in Div(X_\Sigma)$. We suppose that the following two properties hold

**Assumption 4.5.**

1. For all $j = 1, \ldots, c$, the line bundle $\mathcal{L}_j$ is nef. Notice that according to [Ful93, Section 3.4], on a toric variety, $\mathcal{L}_j$ is nef iff it is globally generated.

2. Let $-K_{X_\Sigma}$ be the anti-canonical divisor of $X_\Sigma$. Then we assume that $-K_{X_\Sigma} - \sum_{j=1}^{c} L_j$ is nef.

Put $\mathcal{E} := \oplus_{j=1}^{c} \mathcal{L}_j$ and consider the dual bundle $\mathcal{E}^\vee := \mathcal{H}om_{\mathcal{O}_{X_\Sigma}}(\mathcal{E}, \mathcal{O}_{X_\Sigma})$. We have the following fact.

**Definition-Lemma 4.6.** *The total space $\mathbb{V}(\mathcal{E}^\vee) := \mathbf{Spec}_{\mathcal{O}_{X_\Sigma}}\left(\mathrm{Sym}_{\mathcal{O}_{X_\Sigma}}(\mathcal{E})\right)$ of $\mathcal{E}^\vee$, is a smooth toric variety defined by a fan $\Sigma'$ which is described in the following way. First we define the set of rays $\Sigma'(1)$: For this, we choose divisors $D_{m+j} = \sum_{i=1}^{m} d_{ji} D_i$ with $d_{ji} \geq 0$ and $\mathcal{O}(D_{m+j}) = \mathcal{L}_j$. This choice is possible due to Lemma 4.4 as all $\mathcal{L}_j$ are nef. Write $\underline{d}_i := (d_{1i}, \ldots, d_{ci}) \in \mathbb{Z}^c$ and put $\underline{a}_i' := (\underline{a}_i, \underline{d}_i) \in N' := N \times \mathbb{Z}^c \cong \mathbb{Z}^{n+c}$. Moreover, letting $e_{n+1}, \ldots, e_{n+c}$ be the last $c$ standard generators of $\mathbb{Z}^{n+c}$, we put $\underline{a}_{m+j}' := e_{n+j}$. Then we let $\Sigma'(1) := \{\mathbb{R}_{\geq 0} \underline{a}_1', \ldots, \mathbb{R}_{\geq 0} \underline{a}_{m+c}'\}$ and we group, as before, the column vectors $\underline{a}_1', \ldots, \underline{a}_{m+c}'$ in a matrix $A' \in \mathrm{Mat}\,((n+c) \times (m+c), \mathbb{Z})$. This means that*

$$A' = \left(\begin{array}{c|c} A & 0_{n,c} \\ \hline (d_{ji}) & \mathrm{Id}_c \end{array}\right). \tag{45}$$

*The fan $\Sigma'$ is now defined as follows: For any set of vectors $\underline{b}_1, \ldots, \underline{b}_r \in \mathbb{R}^k$ define $\langle \underline{b}_1, \ldots, \underline{b}_r \rangle := \sum_{j=1}^{r} \mathbb{R}_{\geq 0} \underline{b}_j$. Then we put*

$$\Sigma' := \left\{ \langle \underline{a}_{i_1}', \ldots, \underline{a}_{i_k}', \underline{e}_{j_1}, \ldots, \underline{e}_{j_t} \rangle \subset N_\mathbb{R}' \mid \langle \underline{a}_{i_1}, \ldots, \underline{a}_{i_k} \rangle \in \Sigma(k), \{j_1, \ldots, j_t\} \subset \{n+1, \ldots, n+c\}\right\}.$$

*In other words, considering the canonical projection $\pi : N_\mathbb{R}' \to N_\mathbb{R}$ which forgets the last $c$ components, we have that $\sigma' \in \Sigma'$ iff $\pi(\sigma') \in \Sigma$.*

In the following proposition, we list some rather obvious properties of the cohomology (resp. its toric description) of the space $\mathbb{V}(\mathcal{E}^\vee)$.

**Proposition 4.7.** *Let $X_\Sigma$, $\mathcal{L}_1, \ldots, \mathcal{L}_c$ and the sum $\mathcal{E}$ resp. its dual $\mathcal{E}^\vee$ be as above.*

1. *The projection map $p : \mathbb{V}(\mathcal{E}^\vee) \twoheadrightarrow X_\Sigma$ induces an isomorphism $p^* : H^*(X_\Sigma, \mathbb{Z}) \cong H^*(\mathbb{V}(\mathcal{E}^\vee), \mathbb{Z})$.*

2. *Consider the analogue of sequence (44) for the matrix $A'$, that is, the sequence*

$$0 \longrightarrow \mathbb{L}_{A'} \longrightarrow \mathbb{Z}^{\Sigma'(1)} = \mathbb{Z}^{m+c} \longrightarrow N' \longrightarrow 0, \tag{46}$$

   *then we have an isomorphism*

$$\begin{array}{rcl} \mathbb{L}_A & \longrightarrow & \mathbb{L}_{A'} \\ \underline{l} = (l_1, \ldots, l_m) & \longmapsto & \underline{l}' := (l_1, \ldots, l_m, l_{m+1}, \ldots, l_{m+c}), \end{array} \tag{47}$$

   *where $l_{m+j} := -\sum_{i=1}^{m} l_i d_{ji} = -\langle c_1(\mathcal{L}_j), \underline{l} \rangle$ for all $j = 1, \ldots, c$, and where $\langle -, - \rangle$ is the non-degenerate intersection product between $\mathbb{L} \cong H_2(X_\Sigma, \mathbb{Z})$ and $\mathrm{Pic}(X_\Sigma)$. Notice that in the definition of this isomorphism we consider $\mathbb{L}_A$ resp. $\mathbb{L}_{A'}$ as embedded into $\mathbb{Z}^m$ resp. $\mathbb{Z}^{m+c}$.*

3. *The scalar extension $H^2(X_\Sigma, \mathbb{R}) \xrightarrow{\cong} H^2(\mathbb{V}(\mathcal{E}^\vee), \mathbb{R})$ of the isomorphism $p^*$ from above identifies the Kähler cones (resp. the nef cones) $\mathcal{K}_{X_\Sigma}^\circ$ and $\mathcal{K}_{\mathbb{V}(\mathcal{E}^\vee)}^\circ$ (resp. $\mathcal{K}_{X_\Sigma}$ and $\mathcal{K}_{\mathbb{V}(\mathcal{E}^\vee)}$).*

4. *The manifold $\mathbb{V}(\mathcal{E}^\vee)$ is nef. Moreover, if $s \in \Gamma(X_\Sigma, \mathcal{E})$ is generic, and $Y := s^{-1}(0)$ is the zero locus of this section, then also $Y$ is smooth and also nef.*

*Proof.* The first point follows from the fact that $\mathbb{V}(\mathcal{E}^\vee)$ and $X_\Sigma$ are homotopy equivalent. The second point follows from a direct calculation. For the third point notice that the isomorphism $p^*$ restricted to $H^2(X_\Sigma)$ is given by

$$p^* : H^2(X_\Sigma) \simeq \bigoplus_{i=1}^{m} \mathbb{Z} D_i / (\sum_{i=1}^{m} a_{ki} D_i)_{k=1,\ldots,n} \longrightarrow \bigoplus_{i=1}^{m+c} \mathbb{Z} D_i' / (\sum_{i=1}^{m+c} a_{ki}' D_i')_{k=1,\ldots,n+c} \simeq H^2(\mathbb{V}(\mathcal{E}^\vee)),$$

$$\sum_{i=1}^{m} d_i D_i \mapsto \sum_{i=1}^{m} d_i D_i'.$$

We first prove $p^*(\mathcal{K}_{X_\Sigma}) \subset \mathcal{K}_{\mathbb{V}(\mathcal{E}^\vee)}$. Let $D = \sum_{i=1}^{m} d_i D_i$ be a divisor in $X_\Sigma$ with $\overline{D} \in \mathcal{K}_{X_\Sigma}$. Then $\psi_D^\Sigma$ is given on a maximal cone $\sigma \in \Sigma(n)$ by $u_\sigma^\Sigma \in M \simeq \mathbb{Z}^n$ which is defined by $\langle u_\sigma^\Sigma, \underline{a}_i \rangle = -d_i$ for $\underline{a}_i \in \sigma$. The PL-function $\psi_D^\Sigma$ is convex if and only if for all $\sigma \in \Sigma(n)$ the following inequalities hold $\langle u_\sigma^\Sigma, \underline{a}_i \rangle \geq -d_i$ for all $i \in \{1, \ldots, m\}$. Now consider the corresponding PL-function $\psi_{p^*(D)}^{\Sigma'}$ for $p^*(D)$. Let $\sigma' \in \Sigma'(n+c)$ be a maximal cone in $\Sigma'$ with $\sigma' = \langle \underline{a}'_{i_1}, \ldots, \underline{a}'_{i_n}, \underline{e}_{n+1}, \ldots, \underline{e}_{n+c} \rangle$, where $\{i_1, \ldots, i_n\} \subset \{1, \ldots, m\}$. Then $u_{\sigma'}^{\Sigma'} \in M' \simeq \mathbb{Z}^{n+c}$ is defined by

$$\langle u_{\sigma'}^{\Sigma'}, \underline{a}'_i \rangle = -d_i \quad \text{for } i \in \{i_1, \ldots, i_n\}$$

and

$$\langle u_{\sigma'}^{\Sigma'}, \underline{e}_i \rangle = 0 \quad \text{for } i \in \{n+1, \ldots, n+c\}. \tag{48}$$

But because of equation (48) we have

$$\langle u_{\sigma'}^{\Sigma'}, \underline{a}'_i \rangle = \langle u_\sigma^\Sigma, \underline{a}_i \rangle \geq -d_i \quad \text{for } i \in \{1, \ldots, m\},$$

which shows that $\psi_{p^*(D)}^{\Sigma'}$ is convex, i.e. $p^*(\overline{D}) \in \mathcal{K}_{\mathbb{V}(\mathcal{E}^\vee)}$. Now assume $\overline{D}' \in \mathcal{K}_{\mathbb{V}(\mathcal{E}^\vee)}$. Because $p^*$ is an isomorphism, we can assume that $D'$ has a presentation $\sum_{i=1}^{m+c} d'_i D'_i$ in which $d'_{m+j} = 0$ for $j \in \{1, \ldots, c\}$, i.e. $\overline{D}' = p^*(\overline{D})$ with $D = \sum_{i=1}^m d'_i D_i$. Let $\sigma \in \Sigma(n)$ and $\sigma' \in \Sigma(n+c)$ be maximal cones with $\pi(\sigma') = \sigma$. Because of the presentation of $D'$ we have $\langle u_{\sigma'}^{\Sigma'}, \underline{e}_i \rangle = 0$ for $i \in \{n+1, \ldots, n+c\}$. Therefore we have

$$\langle u_\sigma^\Sigma, \underline{a}_i \rangle = \langle u_{\sigma'}^{\Sigma'}, \underline{a}'_i \rangle \geq -d_i,$$

which shows that $\psi_D^\Sigma$ is convex, i.e. $\overline{D} \in \mathcal{K}_{X_\Sigma}$. The statement for the open parts follows from the fact that $p^*$ is a homeomorphism.

For the fourth point recall that $\mathbb{V}(\mathcal{E}^\vee)$ is nef, i.e. has a nef anticanonical divisor, if the class of the divisor

$$-K_{\mathbb{V}(\mathcal{E}^\vee)} = \sum_{i=1}^m D'_i + \sum_{j=1}^c D'_{m+j}$$

lies in $\mathcal{K}_{\mathbb{V}(\mathcal{E}^\vee)}$. Because of 3. it is enough to show that $(p^*)^{-1}(-K_{\mathbb{V}(\mathcal{E}^\vee)})$ lies in $\mathcal{K}_{X_\Sigma}$. But we have

$$(p^*)^{-1}(-\overline{K}_{\mathbb{V}(\mathcal{E}^\vee)}) = \sum_{i=1}^m \overline{D}_i - \sum_{j=1}^c \sum_{i=1}^m d_{ji} \overline{D}_i = -\overline{K}_{X_\Sigma} - \sum_{j=1}^c c_1(\mathcal{L}_j)$$

and the term on the right hand side lies in $\mathcal{K}_{X_\Sigma}$ by Assumption 4.5 2. Let $s \in \Gamma(X_\Sigma, \mathcal{E})$ be a generic section, then one can show that $Y = s^{-1}(0)$ is smooth by repeatedly applying Bertini's theorem. The nefness of $Y$ is obtained by repeatedly applying the adjunction formula and Assumption 4.5 2. . $\square$

We finish this section by the following remark, which will not be explicitly used in the sequel, but which helps to understand the geometry of the torus embedding considered in the beginning of section 2. More precisely, let $S := \text{Spec } \mathbb{C}[\mathbb{Z}^{n+c}]$ and denote again by $g : S \to \mathbb{P}^{m+c}$ the map defined by $(y_1, \ldots, y_{m+c}) \longmapsto (1 : \underline{y}^{\underline{a}'_1} : \ldots : \underline{y}^{\underline{a}'_{m+c}})$. In section 2 we considered the factorization $g : S \xrightarrow{j} X \xrightarrow{i} \mathbb{P}^{m+c}$ (with $X := \overline{Im(g)}$) where $j$ is an open embedding and $i$ is a closed embedding. However, we will also need to consider some other factorization, namely, we write $g = g^{(2)} \circ g^{(1)}$, where $g^{(1)} : S \longrightarrow \mathbb{C}^m \times (\mathbb{C}^*)^c$ sends $\underline{y}$ to $(\underline{y}^{\underline{a}'_i})_{i=1, \ldots, m+c}$ and $g^{(2)}$ is the composition of the two open embeddings $\mathbb{C}^m \times (\mathbb{C}^*)^c \hookrightarrow \mathbb{C}^{m+c}$ and $\mathbb{C}^{m+c} \hookrightarrow \mathbb{P}^{m+c}$. Now we have the following fact.

**Proposition 4.8.** *The morphism $g^{(1)}$ is a closed embedding. Hence, we have*

$$X \backslash Im(g) \subset \{\mu_0 \cdot \mu_{m+1} \cdot \ldots \cdot \mu_{m+c} = 0\},$$

*where we use $(\mu_0 : \mu_1 : \ldots : \mu_{m+c})$ as homogeneous coordinates on $\mathbb{P}^{m+c}$ and $\mu_1, \ldots, \mu_m$ as coordinates on $\mathbb{C}^{m+c}$ (resp. on $(\mathbb{C}^*)^{m+c}$, $\mathbb{C}^m \times (\mathbb{C}^*)^c$ etc).*

*Proof.* It suffices obviously to show the first statement. We will use a method similar to the proof of [RS15, Proposition 2.1]. First notice that the embedding $\alpha : S \hookrightarrow (\mathbb{C}^*)^{m+c}$ sending $\underline{y}$ to $(\underline{y}^{\underline{a}'_i})_{i=1, \ldots, m+c}$ is obviously closed, so that it suffices to show that $\overline{\text{im}(g^{(1)})} \cap (\mathbb{C}^m \backslash (\mathbb{C}^*)^m) \times (\mathbb{C}^*)^c = \emptyset$. Recall that $\overline{\text{im}(g^{(1)})}$ is the closed subvariety of $\mathbb{C}^m \times (\mathbb{C}^*)^c$ defined by the binomial equations

$$\prod_{i: l'_i > 0} \mu_i^{l'_i} - \prod_{i: l_i < 0} \mu_i^{-l'_i}$$

for any $l' \in \mathbb{L}_{A'}$ (these equations form the toric ideal of $A'$). It was shown in loc.cit. that due to the compactness of $X_\Sigma$, there is some $\underline{l}$ lying in $\mathbb{L}_A \cap \mathbb{Z}_{\geq 0}^m$. Hence, the image $\underline{l}'$ of $\underline{l}$ under the isomorphism (47) lies in $\mathbb{Z}_{>0}^m \times \mathbb{Z}_{<0}^c$, as the coefficients $d_{ji}$ appearing in formula (47) are non-negative (see Definition 4.6) and moreover, for fixed $j$, not all $d_{ji}$ can be zero. It follows that the toric ideal of $A'$ contains an equation

$$\prod_{i=1}^{m} \mu_i^{l'_i} - \prod_{i=m+1}^{m+c} \mu_i^{-l'_i}, \tag{49}$$

where none of the exponents is zero. Now suppose that there is a point $x = (x_1, \ldots, x_m, x_{m+1}, \ldots, x_{m+c}) \in \overline{\text{im}(g^{(1)})} \cap (\mathbb{C}^m \backslash (\mathbb{C}^*)^m) \times (\mathbb{C}^*)^c$, that is, we have $x_i = 0$ for some $i \in \{1, \ldots, m\}$, then as equation (49) vanishes on $x$, we must have some $j \in \{1, \ldots, c\}$ with $x_{m+j} = 0$, which contradicts the assumption that $x \in (\mathbb{C}^m \backslash (\mathbb{C}^*)^m) \times (\mathbb{C}^*)^c$. Hence the intersection $\overline{\text{im}(g^{(1)})} \cap (\mathbb{C}^m \backslash (\mathbb{C}^*)^m) \times (\mathbb{C}^*)^c$ is indeed empty from which it follows that $g^{(1)} : S \hookrightarrow \mathbb{C}^m \times (\mathbb{C}^*)^c$ is a closed embedding. $\square$

**Remark:** The GKZ-systems (see Definition 2.8) associated to the matrix $A'$ is not necessary regular, as the vectors $\underline{a}'_1, \ldots, \underline{a}'_{m+c}$ do not necessarily lie on an affine hyperplane in $\mathbb{Z}^{m+c}$ (see [Hot98] for this regularity criterion). The situation is similar to that considered in our earlier paper [RS15], and for the same reasons as in loc.cit., we will work with the extended matrix $A'' \in \text{Mat}((1 + n + c) \times (1 + m + c), \mathbb{Z})$ with columns $\underline{a}''_0, \underline{a}''_1, \ldots, \underline{a}''_{m+c}$, where $\underline{a}''_i := (1, \underline{a}'_i)$ and $\underline{a}''_0 := (1, \underline{0}, \underline{0})$. In particular we have $\underline{a}''_{m+j} = (1, \underline{e}_{n+j}) \in \mathbb{Z}^{n+c+1}$ for $j = 1, \ldots, c$ where $e_{n+j}$ is the $n+j$-th standard vector in $\mathbb{C}^{n+c}$. We write $\mathbb{L}_{A''}$ for the module of relations between the columns of $A''$, obviously we have an isomorphism $\mathbb{L}_{A'} \to \mathbb{L}_{A''}$ sending $\underline{l} = (l_1, \ldots, l_{m+c})$ to $(-\sum_{i=1}^{m+c} l_i, \underline{l})$. As a matter of notation, we will often write the parameter of the GKZ-systems defined by the matrix $A''$, which are vectors in $\mathbb{C}^{1+m+c}$ by definition, as $(\alpha, \underline{\beta}, \underline{\gamma}) \in \mathbb{C}^{1+m+c}$, where $\alpha \in \mathbb{C}$, $\underline{\beta} \in \mathbb{C}^m$ and $\underline{\gamma} \in \mathbb{C}^c$.

## 5 Euler-Koszul homology and duality of GKZ-systems

In this section, we show a duality result for the GKZ-systems associated to the toric situation just described. We will explain how to calculate the holonomic dual of the system $\mathcal{M}_{A''}^\beta$ for some specific $\beta$, this is used to get a more precise description of the various $\mathcal{D}$-module considered in sections 2 and 3. The methods used here somehow similar [RS15, section 2.3], but we have to take into account the non-compactness of the toric varieties involved.

**Proposition 5.1.** *Let $X_\Sigma$ be smooth, toric and projective and suppose that $\mathcal{L}_1 = \mathcal{O}_{X_\Sigma}(L_1), \ldots, \mathcal{L}_c = \mathcal{O}_{X_\Sigma}(L_c)$ are nef line bundles on $X_\Sigma$. However, we do not make any assumption on the positivity of $-K_{X_\Sigma} - \sum_{j=1}^c L_j$. Let $A'$ be the matrix from in Definition 4.6 (i.e. with columns the primitive integral generator of the fan of $\mathbb{V}(\mathcal{E}^\vee))$. Then the semi-group ring $\mathbb{C}[\mathbb{N}A']$ is normal and Cohen-Macaulay. The map*

$$\Psi : \mathbb{N}A' \longrightarrow (\mathbb{N}A')^\circ,$$

$$\underline{m} \longmapsto \underline{m} + \underline{a}'_{m+1} + \ldots + \underline{a}'_{m+c}$$

*is a bijection. Hence, $\mathbb{C}[\mathbb{N}A']$ is a Gorenstein ring where the generator of the canonical module $\omega_{\mathbb{C}[\mathbb{N}A']}$ is given by the monomial $\underline{y}^{\underline{a}'_{m+1} + \ldots + \underline{a}'_{m+c}}$.*

We can deduce the following immediate corollary.

**Corollary 5.2.** *In the situation of the last proposition, suppose moreover that $-K_{X_\Sigma} - \sum_j^c L_j$ is nef. Let $A''$ be the extension considered at the end of section 4. Then also the semi-group $\mathbb{N}A''$ is normal and we have*

$$(\mathbb{N}A'')^\circ = \underline{a}''_0 + \underline{a}''_{m+1} + \ldots + \underline{a}''_{m+c} + \mathbb{N}A''.$$

*Hence $\mathbb{C}[\mathbb{N}A'']$ is a normal, Cohen-Macaulay and Gorenstein ring, with*

$$\omega_{\mathbb{C}[\mathbb{N}A'']} \cong \mathbb{C}[\mathbb{N}A''] \cdot \underline{y}^{\underline{a}_0'' + \underline{a}_{m+1}'' + \ldots + \underline{a}_{m+c}''}.$$

*Proof.* This follows directly by applying proposition 5.1 to the toric variety $X_\Sigma$ and the collection of nef line bundles $\mathcal{L}_1, \ldots, \mathcal{L}_c, \mathcal{L}_{c+1} := \mathcal{O}_{X_\Sigma}(-K_{X_\Sigma} - \sum_{j=1}^{x} L_j)$. □

The following lemma is a rather obvious consequence of the nefness condition of the bundles $\mathcal{L}_1, \ldots, \mathcal{L}_c$.

**Lemma 5.3.** *Let as before $X_\Sigma$ be toric and let $\mathcal{L}_1, \ldots, \mathcal{L}_c$ be nef line bundles. Consider the fan $\Sigma'$ of the space $\mathbb{V}(\mathcal{E}^\vee)$, where $\mathcal{E} = \oplus_{j=1}^c \mathcal{L}_j$. Then the support $\mathrm{supp}(\Sigma')$ is convex. As a consequence, we have the following equality*

$$\mathrm{supp}(\Sigma') = \mathbb{R}_{\geq 0} A' \tag{50}$$

*where $\mathbb{R}_{\geq 0} A' := \sum_{i=1}^{m+c} \mathbb{R}_{\geq 0} \underline{a}_i'$.*

*Proof.* This is obvious from the construction of $\Sigma'$ as presented in definition 4.6. Namely, for any $j \in \{1, \ldots, c\}$, the functions $\psi_{D_{m+j}}^\Sigma = \sum_{i=1}^m d_{ji} \psi_{D_i}^\Sigma$ are convex due to the nefness of $\mathcal{L}_j$ (remember that $\mathcal{O}(D_{m+j}) = \mathcal{L}_j$), and one can describe the set $\mathrm{supp}(\Sigma')$ as

$$\mathrm{supp}(\Sigma') = \Big\{(x_1, \ldots, x_n, x_{n+1}, \ldots, x_{n+c}) \mid (x_1, \ldots, x_n) \in \mathrm{supp}(\Sigma) = \mathbb{R}^n,$$
$$x_{n+j} \geq -\psi_{D_{m+j}}^\Sigma(x_1, \ldots, x_n) \quad \forall j = 1, \ldots, c\Big\}.$$

Then the convexity of the set $\mathrm{supp}(\Sigma')$ is precisely the convexity condition on the functions $\psi_{D_{m+j}}^\Sigma$. For the second statement, notice that the inclusion $\mathrm{supp}(\Sigma') \subset \mathbb{R}_{\geq 0} A'$ is trivial (and does not depend on the convexity of $\mathrm{supp}(\Sigma')$). On the other hand, if $\mathrm{supp}(\Sigma')$ is convex, then we have the inclusion

$$\mathrm{supp}(\Sigma') \supset \mathrm{Conv}(\underline{a}_1', \ldots, \underline{a}_{m+c}'), \tag{51}$$

where $\mathrm{Conv}(\underline{a}_1', \ldots, \underline{a}_{m+c}')$ denotes the convex hull of the vectors $\underline{a}_1', \ldots, \underline{a}_{m+c}'$, since the left hand side must contain the convex hull of any of its subsets. On the other hand, we obviously have that

$$\mathbb{R}_{\geq 0} A = \left\{\lambda \cdot \underline{x} \mid \underline{x} \in \mathrm{Conv}(\underline{a}_1', \ldots, \underline{a}_{m+c}'), \lambda \in \mathbb{R}_{\geq 0}\right\},$$

so that the desired inclusion $\mathrm{supp}(\Sigma') \supset \mathbb{R}_{\geq 0} A'$ follows from equation (51) and the fact that the set $\mathrm{Supp}(\Sigma')$ is conical, i.e., for all $\underline{x} \in \mathrm{Supp}(\Sigma')$ and all $\lambda \in \mathbb{R}_{\geq 0}$ we have that $\lambda \cdot \underline{x} \in \mathrm{Supp}(\Sigma')$. □

*Proof of the proposition.* We first show the normality of $\mathbb{N}A'$: Given any vector $\underline{x}' \in \mathbb{R}_{\geq 0} A' \cap N'$, then by equation (50) there is some maximal cone $\langle \underline{a}_{i_1}, \ldots, \underline{a}_{i_n} \rangle \in \Sigma$ such that $\underline{x}' \in \langle \underline{a}_{i_1}', \ldots, \underline{a}_{i_n}', \underline{a}_{i_{n+1}}', \ldots, \underline{a}_{i_{n+1}}' \rangle \in \Sigma'$ (recall that $\underline{a}_{i_{n+j}}' = \underline{a}_{m+j}' = \underline{e}_{n+j}$). Hence we have an equation

$$\underline{x}' = \sum_{k=1}^{n+c} \lambda_k \underline{a}_{i_k}' \tag{52}$$

with $\lambda_k \in \mathbb{R}_{\geq 0}$. We know that $(\underline{a}_{i_1}', \ldots, \underline{a}_{i_{n+c}}') = (\underline{a}_{i_1}', \ldots, \underline{a}_{i_n}', \underline{e}_{m+1}, \ldots, \underline{e}_{m+c})$ is a $\mathbb{Z}$-basis of $N'$ as $\langle \underline{a}_{i_1}', \ldots, \underline{a}_{i_{n+c}}' \rangle$ is a smooth $n+c$-dimensional cone in $\Sigma'$. Hence $\lambda_k \in \mathbb{N}$ for $k = 1, \ldots, n+c$, and $\underline{x}' \in \mathbb{N}A'$, which is the defining property of normality of $\mathbb{N}A'$. It follows that $\mathbb{C}[\mathbb{N}A']$ is Cohen-Macaulay by Hochster's theorem ([Hoc72, Theorem 1]).

It remains to show the second statement concerning the characterization of the interior points of $\mathbb{N}A'$. We will actually show the following

**Claim:** Let $\underline{x}' \in \mathbb{N}A'$. Consider the representation (52) of $\underline{x}'$ as an element of $\sum_{j=1}^{n+c} \mathbb{R}_{\geq 0} \underline{a}_{i_j}'$, that is, an equation $x' = \sum_{i=1}^{m+c} \lambda_i \underline{a}_i' \in \mathbb{N}A'$, where $\lambda_k = 0$ if $k \in \{1, \ldots, m\} \setminus \{i_1, \ldots, i_n\}$. Then $x'$ lies in $(\mathbb{N}A')^\circ$ iff $\lambda_i > 0$ for $i \in \{m+1, \ldots, m+c\} = \{i_{n+1}, \ldots, i_{n+c}\}$.

Notice that a representation as in the claim is unique, if there are two maximal cones of $\Sigma(n)$ such that $x'$ is contained in both of the cones generated by the corresponding column vectors of $A'$, then it lies on a common boundary, and the two expressions (52) are equal.

The claim implies that the map $\Psi$ from the proposition is well-defined and surjective, and it is obviously injective. In order to show the claim, notice that

$$(\mathbb{N}A')^\circ = (\mathbb{R}_{\geq 0}A'\backslash\partial(\mathbb{R}_{\geq 0}A')) \cap N' = (\mathbb{R}_{\geq 0}A' \cap N') \setminus (\partial(\mathbb{R}_{\geq 0}A') \cap N') = \mathbb{N}A'\setminus(\partial(\mathbb{R}_{\geq 0}A') \cap N'),$$

so that we have to show that the points in $\partial(\mathbb{R}_{\geq 0}A') \cap N'$ are precisely those from $\mathbb{N}A'$ where in the above representation (52) there is at least one index $i \in \{m+1,\ldots,m+c\}$ with $\lambda_i = 0$. From Formula (50) we deduce that

$$\partial(\mathbb{R}_{\geq 0}A') \subset \bigcup_{\langle \underline{a}_{i_1},\ldots,\underline{a}_{i_n}\rangle \in \Sigma_A(n)} \partial\langle \underline{a}'_{i_1},\ldots,\underline{a}'_{i_n}, \underline{a}'_{m+1},\ldots,\underline{a}'_{m+c}\rangle.$$

More precisely, for each $\langle \underline{a}_{i_1},\ldots,\underline{a}_{i_n}\rangle \in \Sigma(n)$ the cone $\langle \underline{a}'_{i_1},\ldots,\underline{a}'_{i_n}, \underline{a}'_{m+1},\ldots,\underline{a}'_{m+c}\rangle$ has two types of facets: those that are facets of $\partial(\mathbb{R}_{\geq 0}A')$ (call them "outer boundary") and those which are not ("inner boundary"). The union (over all $n$-dimensional cones of $\Sigma$) of the outer boundaries is the set $\partial(\mathbb{R}_{\geq 0}A')$ we are interested in.

The fan $\Sigma'$ is smooth, in particular simplicial, this implies that for any cone $\langle \underline{a}'_{i_1},\ldots,\underline{a}'_{i_n}, \underline{a}'_{m+1},\ldots,\underline{a}'_{m+c}\rangle \in \Sigma'$ we have the following description of its boundary.

$$\partial\langle \underline{a}'_{i_1},\ldots,\underline{a}'_{i_{n+c}}\rangle = \partial\langle \underline{a}'_{i_1},\ldots,\underline{a}'_{i_n}, \underline{e}_{n+1},\ldots,\underline{e}_{n+c}\rangle$$

$$\stackrel{!}{=} \bigcup_{k=1}^{n} \langle \underline{a}'_{i_1},\ldots,\widehat{\underline{a}}'_{i_k},\ldots,\underline{a}'_{i_n}, \underline{e}_{n+1},\ldots,\underline{e}_{n+c}\rangle \cup \bigcup_{l=1}^{c} \langle \underline{a}'_{i_1},\ldots,\underline{a}'_{i_n}, \underline{e}_{m+1},\ldots,\widehat{\underline{e}}_{m+l},\ldots,\underline{e}_{m+c}\rangle.$$

The facet $\langle \underline{a}'_{i_1},\ldots,\widehat{\underline{a}}'_{i_k},\ldots,\underline{a}'_{i_n},\underline{e}_{n+1},\ldots,\underline{e}_{n+c}\rangle$ is an inner boundary, i.e., it is not contained in $\partial(\mathbb{R}_{\geq 0}A')$. This is a consequence of the completeness of $\Sigma$, namely, there is some other cone $\langle \underline{a}_{j_1},\ldots,\underline{a}_{j_n}\rangle \in \Sigma$ having $\langle \underline{a}_{i_1},\ldots,\widehat{\underline{a}}_{i_k},\ldots,\underline{a}_{i_n}\rangle$ as a facet, and then similarly the cone $\langle \underline{a}'_{i_1},\ldots,\widehat{\underline{a}}'_{i_k},\ldots,\underline{a}'_{i_n},\underline{e}_{n+1},\ldots,\underline{e}_{n+c}\rangle$ is a facet of both $\langle \underline{a}'_{i_1},\ldots,\underline{a}'_{i_n},\underline{e}_{n+1},\ldots,\underline{e}_{n+c}\rangle$ and $\langle \underline{a}'_{j_1},\ldots,\underline{a}'_{j_n},\underline{e}_{n+1},\ldots,\underline{e}_{n+c}\rangle$, hence it is not contained in $\partial(\mathbb{R}_{\geq 0}A')$. However, the facet $\langle \underline{a}'_{i_1},\ldots,\underline{a}'_{i_n},\underline{e}_{n+1},\ldots,\widehat{\underline{e}}_{n+l},\ldots,\underline{e}_{n+c}\rangle$ (for $l = 1,\ldots,c$) is an outer boundary, i.e., a facets of $\mathbb{R}_{\geq 0}A'$. We conclude that

$$\partial(\mathbb{R}_{\geq 0}A') = \bigcup_{\langle \underline{a}_{i_1},\ldots,\underline{a}_{i_n}\rangle \in \Sigma(n)} \left[\bigcup_{l=1}^{c} \langle \underline{a}'_{i_1},\ldots,\underline{a}'_{i_n}, \underline{e}_{n+1},\ldots,\widehat{\underline{e}}_{n+l},\ldots,\underline{e}_{n+c}\rangle\right].$$

We see that for any point $\partial(\mathbb{R}_{\geq 0}A') \cap N'$, there must be some $l \in \{1,\ldots,c\}$ such that in the representation (52) the coefficient $\lambda_{m+l}$ is zero. This shows the claim, and proves that the map $\Psi$ is an isomorphism. Finally, it follows from standard arguments about semigroup rings (see, e.g. [BH93, corollary 6.3.8]) that $\mathbb{C}[\mathbb{N}A']$ is Gorenstein, and that the generator of the canonical module $\omega_{\mathbb{C}[\mathbb{N}A']}$ is as claimed. □

As a consequence, we obtain the following duality result for those GKZ systems that we will be interested in the sequel.

**Theorem 5.4.** *et $A''$ be as above, that is, suppose that its columns $(\underline{a}''_0, \underline{a}''_1, \ldots, \underline{a}''_{m+c})$ are of the form $\underline{a}''_i = (1, \underline{a}'_i)$ where $\underline{a}''_0 = (1, \underline{0})$ and where $\underline{a}'_i$ ($i = 1,\ldots,m+c$) are the integral primitive generator of the fan of $\mathbb{V}(\mathcal{E}^\vee)$. For $\beta \in \mathbb{Z}^{1+m+c}$, consider the GKZ-system $\mathcal{M}^\beta_{A''}$ as in Definition 2.7.*

1. *There is an isomorphism*

$$\mathbb{D}(\mathcal{M}^{(0,\underline{0},\underline{0})}_{A''}) \cong \mathcal{M}^{-(c+1,\underline{0},\underline{1})}_{A''} = \mathcal{M}^{-\underline{a}''_0-\sum_{l=1}^{c}\underline{a}''_{m+l}}_{A''}.$$

2. *Consider the natural good filtration $F_\bullet \mathcal{M}^\beta_{A''}$ induced by the order filtration on $\mathcal{D}$. Let $\mathbb{D}(\mathcal{M}^\beta_{A''}, F_\bullet)$ be the dual filtered module in the sense of [Sai88, section 2.4], i.e., $\mathbb{D}(\mathcal{M}^\beta_{A''}, F_\bullet) = (\mathbb{D}\mathcal{M}^\beta_{A''}, F^\mathbb{D}_\bullet)$ where $F^\mathbb{D}_\bullet(\mathbb{D}\mathcal{M}^\beta_{A''})$ is the filtration dual to $F_\bullet \mathcal{M}^\beta_{A''}$. Then we have*

$$\mathbb{D}\left(\mathcal{M}^{-\underline{a}''_0-\sum_{l=1}^{c}\underline{a}''_{m+l}}_{A''}, F_\bullet\right) \cong (\mathcal{M}^{(0,\underline{0},\underline{0})}_{A''}, F_{\bullet+n-(m+c+1)}).$$

*Proof.* 1. The proof is parallel to [Wal07, Proposition 4.1] or [RS15, Theorem 2.15 and Proposition 2.18], so that we only sketch it here, referring to loc.cit. for details. First one has to define the so-called Euler-Koszul complex resp. co-complex (see [MMW05]). Its global sections complex $K_\bullet(T, E - \beta)$ is a complex of free $D_V \otimes_R T$-modules where $R = \mathbb{C}[\partial_0, \partial_1, \ldots, \partial_{m+c}]$ and where $T$ is a so-called *toric R-module*. A particular case is $T = \mathbb{C}[\mathbb{N}A'']$. Notice that the terms of $K_\bullet(T, E - \beta)$ are not free over $D_V$. However, for $T = \mathbb{C}[\mathbb{N}A'']$, this complex is a resolution by left $D_V$-modules of the modules $M_{A''}^\beta$. The differentials of $K_\bullet(T, E - \beta)$ are defined by the operators $E$ and $Z_k$ entering in the definition of $M_{A''}^\beta$. From a resolution of the toric ring $\mathbb{C}[\mathbb{N}A'']$ by free $\mathbb{C}[\partial_0, \partial_1, \ldots, \partial_{m+c}]$-modules one can also construct a resolution of $M_{A''}^{(0,\underline{0},\underline{0})}$ by free $D_V$-modules. Applying $Hom_{D_V}(-, D_V)$ yields basically the same complex, but where the parameters in the differentials are changed, and where the toric module is now the canonical module of the ring $\mathbb{C}[\mathbb{N}A'']$. Now from the Gorenstein property of $\mathbb{C}[\mathbb{N}A'']$ with the precise description of the interior ideal from Corollary 5.2 we obtain the desired result by taking the cohomology of the two complexes, that is, we can show the identification of the holonomic dual of $\mathcal{M}_{A''}^{(0,\underline{0},\underline{0})}$ with $\mathcal{M}_{A''}^{-(c+1,\underline{0},\underline{1})}$.

2. The proof is literally the same as in [RS15, Proposition 2.19, 2.] with the indices shifted appropriately.

□

As a consequence, we can make more specific statements on the parameter vectors of the various GKZ-systems occurring in the results of the previous sections.

**Corollary 5.5.** *Consider the situation in section 2 where the matrix $B$ is $A'$, i.e., given by the primitive integral generators of the fan of $\mathbb{V}(\mathcal{E}^\vee)$, in particular, both $\mathbb{N}B = \mathbb{N}A'$ and $\mathbb{N}\widetilde{B} = \mathbb{N}A''$ are normal semigroups. Then*

1. *The statements of Theorem 2.11, Theorem 2.16 and of Proposition 2.17 hold true for the parameter values $\widetilde{\beta} = (0, \underline{0}, \underline{0}), \widetilde{\beta}' = (c+1, \underline{0}, \underline{1}) \in \mathbb{Z}^{1+n+c}$.*

2. *The statements of Proposition 3.3 and of Theorem 3.6 hold true for the parameter values $\beta = (\underline{0}, \underline{0}), \beta' = (\underline{0}, \underline{1}) \in \mathbb{Z}^{n+c}$ and for any $\beta_0, \beta'_0 \in \mathbb{Z}$.*

For later use, we introduce the following piece of notation.

**Definition 5.6.** *In the situation of Theorem 5.4, we call the map*
$$\phi : \mathcal{M}_{A''}^{-(c+1,\underline{0},\underline{1})} \longrightarrow \mathcal{M}_{A''}^{(0,\underline{0},\underline{0})}$$
*induced by right multiplication by $\partial_0 \cdot \partial_{m+1} \cdot \ldots \cdot \partial_{m+c}$ the duality morphism. For any $\beta_0 \in \mathbb{Z}$, we obtain an induced morphism*
$$\widehat{\phi} : \widehat{\mathcal{M}}_{A'}^{(\beta_0,\underline{0},-\underline{1})} \longrightarrow \widehat{\mathcal{M}}_{A'}^{(\beta_0+c,\underline{0},\underline{0})}$$
*given by right multiplication with $\partial_{m+1} \cdot \ldots \cdot \partial_{m+c}$ (see 3.1 for the definition of the modules $\widehat{\mathcal{M}}^\beta$). The case $\beta_0 = -2c$ will be particularly important, and we will also call the map*
$$\widehat{\phi} : \widehat{\mathcal{M}}_{A'}^{-(2c,\underline{0},\underline{1})} \longrightarrow \widehat{\mathcal{M}}_{A'}^{(-c,\underline{0},\underline{0})}$$
*the duality morphism.*

**Remark:** In our previous paper [RS15], we obtained from a similar construction a non-degenerate pairing on the Fourier-Laplace transformed GKZ-system (see [RS15, corollary 2.20], where this system was called $\widehat{\mathcal{M}}_{\widetilde{A}}$). It was given by an isomorphism of $\widehat{\mathcal{M}}_{\widetilde{A}}$ to its holonomic dual (which is isomorphic to its meromorphic dual, see also the proof of Lemma 6.6 below). The fact that in the current situation, we only have a morphism $\widehat{\phi} : \widehat{\mathcal{M}}_{A'}^{-(2c,\underline{0},\underline{1})} \longrightarrow \widehat{\mathcal{M}}_{A'}^{(-c,\underline{0},\underline{0})}$ which is not an isomorphism unless $c = 0$ (in which case we are exactly in the situation of [RS15], see the remark at the end of section 6 of this paper) corresponds to the fact that the pairing $S$ on the twisted quantum $\mathcal{D}$-module as introduced in Definition 4.3 is degenerate. As we have seen in the definition of the reduced quantum-$\mathcal{D}$-module, it becomes non-degenerate when we divide out the kernel of the cup product with the first Chern classes of the line bundles $\mathcal{L}_j$. We will show below in corollary 6.14 that the reduced quantum $\mathcal{D}$-module is part of a *non-commutative* Hodge structure, which implies in particular that it carries a non-degenerate pairing like the one from [RS15].

# 6 Mirror correspondences

In this section we combine the results obtained so far with the GKZ-type description of the ambient resp. reduced quantum $\mathcal{D}$-modules from [MM11] for the toric case. We obtain a mirror statement which identifies them with $\mathcal{D}$-modules constructed from our Landau-Ginzburg models. The results from section 2 will be applied for the case where the matrix $B$ (used for the construction of GKZ-systems and of families of Laurent polynomials) is given by $A'$ (see Definition 4.6) the columns of which are the primitive integral generators of the fan of the total bundle $\mathbb{V}(\mathcal{E}^\vee)$. Recall also (remark at the end of section 4) that we denote by $A''$ the matrix constructed from $A'$ by adding 1 as an extra component to all columns and by adding $(1,\underline{0})$ as extra column. Hence, if $B$ is equal to $A'$, then the matrix $\widetilde{B}$ used in section 2 is exactly the matrix $A''$. Recall also that the parameter of the GKZ-systems of the matrix $A''$ is written as $(\alpha, \underline{\gamma}, \underline{\delta}) \in \mathbb{C}^{1+m+c}$ with $\alpha \in \mathbb{C}$, $\underline{\gamma} \in \mathbb{C}^m$ and $\underline{\delta} \in \mathbb{C}^c$.

The starting point for our discussion here is the duality morphism from the last section. We need to consider a slight variation of it, which is defined only outside the boundary $\lambda_i = 0$ and only outside the bad parameter locus as defined in subsection 3.2. Recall that $V = \mathbb{C}_{\lambda_0} \times W$, and that this bad parameter locus of the family $\varphi_{A'}$ was called $W^{bad} \subset W$. The complement of this locus outside the boundary $\lambda_i = 0$ was called $W^\circ$, that is, $W^\circ := W^* \backslash W^{bad}$.

**Definition-Lemma 6.1.** *For any $\beta = (\beta_0, \beta_1, \ldots, \beta_m, \beta_{m+1}, \ldots, \beta_{n+c}) \in \mathbb{Z}^{1+n+c}$, consider the restricted, Fourier-Laplace transformed GKZ-system ${}^*\widehat{\mathcal{M}}_{A'}^{(\beta_0, \beta)}$ we have*

$$
{}^*\widehat{\mathcal{M}}_{A'}^{(\beta_0, \beta)} = \frac{\mathcal{D}_{\mathbb{C}_z \times W^*}[z^{-1}]}{\mathcal{D}_{\mathbb{C}_z \times W^*}[z^{-1}](\widehat{\square}_{\underline{l}}^*)_{\underline{l} \in \mathbb{L}_{A'}} + \mathcal{D}_{\mathbb{C}_z \times W^*}[z^{-1}](\widehat{E}_k - z\beta_k)_{k=0,\ldots,n+c}},
$$

*where*

$$
\widehat{\square}_{\underline{l}}^* := \prod_{i \in \{1,\ldots,m+c\}: l_i > 0} \lambda_i^{l_i}(z \cdot \partial_i)^{l_i} - \prod_{i=1}^{m+c} \lambda_i^{l_i} \cdot \prod_{i \in \{1,\ldots,m+c\}: l_i < 0} \lambda_i^{-l_i}(z \cdot \partial_i)^{-l_i}, \qquad \underline{l} \in \mathbb{L}_{A'}
$$

$$
\widehat{E}_0 := z^2 \partial_z + \sum_{i=1}^{m+c} \lambda_i \cdot z \partial_i,
$$

$$
\widehat{E}_k := \sum_{i=1}^{m+c} a'_{ki} \lambda_i \cdot z \partial_i \qquad k = 1, \ldots, n+c
$$

*and moreover, ${}_0^*\widehat{\mathcal{M}}_{A'}^\beta$ is the $\mathcal{R}_{\mathbb{C}_z \times W^*}$-subalgebra generated by $[1]$, and we have*

$$
{}_0^*\widehat{\mathcal{M}}_{A'}^\beta = \frac{\mathcal{R}_{\mathbb{C}_z \times W^*}}{\mathcal{R}_{\mathbb{C}_z \times W^*}(\widehat{\square}_{\underline{l}}^*) + \mathcal{R}_{\mathbb{C}_z \times W^*}(\widehat{E}_k - z\beta_k)_{k=0,\ldots,n+c}}.
$$

*Moreover, we define the modules ${}^*\widehat{\mathcal{N}}_{A'}^{(\beta_0, \beta)}$ as the cyclic quotients of $\mathcal{D}_{\mathbb{C}_z \times W^*}[z^{-1}]$ by the left ideal generated by $\widetilde{\square}_{\underline{l}}$ for $\underline{l} \in \mathbb{L}_{A'}$ and $\widehat{E}_k - z\beta_k$ for $k = 0, \ldots, n+c$, where*

$$
\widetilde{\square}_{\underline{l}} := \prod_{i \in \{1,\ldots,m\}: l_i > 0} \lambda_i^{l_i}(z \cdot \partial_i)^{l_i} \prod_{i \in \{m+1,\ldots,m+c\}: l_i > 0} \prod_{\nu=1}^{l_i} (\lambda_i(z \cdot \partial_i) - z \cdot \nu)
$$

$$
- \prod_{i=1}^{m+c} \lambda_i^{l_i} \cdot \prod_{i \in \{1,\ldots,m\}: l_i < 0} \lambda_i^{-l_i}(z \cdot \partial_i)^{-l_i} \prod_{i \in \{m+1,\ldots,m+c\}: l_i < 0} \prod_{\nu=1}^{-l_i} (\lambda_i(z \cdot \partial_i) - z \cdot \nu).
$$

*Consider the morphism*

$$
\Psi : {}^*\widehat{\mathcal{N}}_{A'}^{(0,\underline{0},\underline{0})} \longrightarrow {}^*\widehat{\mathcal{M}}_{A'}^{-(2c,\underline{0},\underline{1})} \tag{53}
$$

*given by right multiplication with $z^c \cdot \prod_{i=m+1}^{m+c} \lambda_i$. As it is obviously invertible, the two modules ${}^*\widehat{\mathcal{N}}_{A'}^{(0,\underline{0},\underline{0})}$ and ${}^*\widehat{\mathcal{M}}_{A'}^{-(2c,\underline{0},\underline{1})}$ are isomorphic. We define $\widetilde{\phi}$ to be the composition $\widetilde{\phi} := \widehat{\phi} \circ \Psi$, where $\widehat{\phi}$ is the duality morphism introduced in Definition 5.6. In concrete terms, we have:*

$$
\widetilde{\phi} : {}^*\widehat{\mathcal{N}}_{A'}^{(0,\underline{0},\underline{0})} \longrightarrow {}^*\widehat{\mathcal{M}}_{A'}^{(-c,\underline{0},\underline{0})},
$$

$$
m \longmapsto \widehat{\phi}(m \cdot z^c \cdot \lambda_{m+1} \cdot \ldots \cdot \lambda_{m+c}) = m \cdot (z\lambda_{m+1}\partial_{m+1}) \cdot \ldots \cdot (z\lambda_{m+c}\partial_{m+c}).
$$

*In view of corollary 5.5, 2. (see also Theorem 3.6) we obtain*

$$\mathrm{im}(\widetilde{\phi}) \cong \mathrm{im}(\widehat{\phi}) \cong (\mathrm{id}_{\mathbb{C}_z} \times j)^+ \widehat{\mathcal{M}^{IC}}(X^\circ, \mathcal{L}). \tag{54}$$

*For any $\beta \in \mathbb{Z}^{1+n+c}$, consider the $\mathcal{R}_{\mathbb{C}_z \times W^*}$-subalgebra of*

$$\mathcal{D}_{\mathbb{C}_z \times W^*}[z^{-1}] \Big/ \mathcal{D}_{\mathbb{C}_z \times W^*}[z^{-1}]\left((\widetilde{\Box}_{\underline{l}})_{\underline{l} \in \mathbb{L}_{A'}}\right) + \mathcal{D}_{\mathbb{C}_z \times W^*}[z^{-1}]\left(\widehat{E}_k - z\beta_k\right)_{k=0,\ldots,n+c}$$

*generated by the element [1] and denote its restriction to $\mathbb{C}_z \times W^\circ$ by ${}_0\widehat{\mathcal{N}}_{A'}^{(\beta_0, \beta)}$. Similarly to corollary 3.20, we have*

$${}_0\widehat{\mathcal{N}}_{A'}^{(\beta_0, \beta)} = \left[\frac{\mathbb{C}[z, \lambda_1^\pm, \ldots, \lambda_{m+c}^\pm]\langle z^2\partial_z, z\partial_{\lambda_1}, \ldots z\partial_{\lambda_{m+c}}\rangle}{\left((\widetilde{\Box}_{\underline{l}})_{\underline{l} \in \mathbb{L}_{A'}}, (\widehat{E}_k - z \cdot \beta_k)_{k=0,\ldots,n+c}\right)}\right]_{|\mathbb{C}_z \times W^\circ}.$$

In the next lemma we want to describe the restriction of the $\mathcal{D}$-module $\widehat{\mathcal{M}^{IC}}(X^\circ, \mathcal{L})$ to $\mathbb{C}_z \times W^*$.

**Lemma 6.2.** *Consider the morphism $\widehat{\phi} : \widehat{\mathcal{M}}_{A'}^{-(2c,0,1)} \longrightarrow \widehat{\mathcal{M}}_{A'}^{-(c,0,0)}$ from Definition 5.6 and the isomorphisms $\widehat{\mathcal{M}^{IC}}(X^\circ, \mathcal{L}) \simeq \mathrm{im}(\widehat{\phi}) \simeq \widehat{\mathcal{M}}_{A'}^{-(2c,0,1)}/\ker(\widehat{\phi})$ from Corollary 5.5 (see also Theorem 3.6). We have the following isomorphism*

$$(id_{\mathbb{C}_z} \times j)^+ \widehat{\mathcal{M}^{IC}}(X^\circ, \mathcal{L}) \simeq {}^*\widehat{\mathcal{M}}_{A'}^{-(2c,0,1)}/\widehat{\mathcal{K}}_{\mathcal{M}} \simeq {}^*\widehat{\mathcal{N}}_{A'}^{(0,0,0)}/\widehat{\mathcal{K}}_{\mathcal{N}},$$

*where $\widehat{\mathcal{K}}_{\mathcal{M}}$ resp. $\widehat{\mathcal{K}}_{\mathcal{N}}$ are the sub-$\mathcal{D}$-modules associated to the sub-D-modules*

$$\{m \in {}^*\widehat{M}^{-(2c,0,1)} \mid \exists p \in \mathbb{Z}, k \in \mathbb{N} \text{ such that } (\lambda\partial + p)\ldots(\lambda\partial + p + k)m = 0\}$$

*resp.*

$$\{n \in {}^*\widehat{N}^{(0,0,0)} \mid \exists p \in \mathbb{Z}, k \in \mathbb{N} \text{ such that } (\lambda\partial + p)\ldots(\lambda\partial + p + k)n = 0\}$$

*with $(\lambda\partial + i) := \prod_{j=m+1}^{m+c}(\lambda_j\partial_j + i)$ for $i \in \mathbb{Z}$.*

*Proof.* We will first compute the restriction of $\mathcal{M}^{IC}(X^\circ, \mathcal{L})$ to $V^* = \mathbb{C}_{\lambda_0} \times W^*$. Recall the morphism $\phi : M_{A''}^{-(c+1,0,1)} \longrightarrow M_{A''}^{(0,0,0)}$ from Definition 5.6. We know from Theorem 2.16 and from Propostion 2.17 that $M^{IC}(X^\circ, \mathcal{L}) \simeq M_{A''}^{-(c+1,0,1)}/\ker(\phi)$ where $\ker(\phi)$ is given by

$$\{m \in M_{A''}^{-(c+1,0,1)} \mid \exists\, n \in \mathbb{N} \text{ such that } (\partial_0 \cdot \partial_{m+1} \cdots \partial_{m+c})^n m = 0\}.$$

Notice that $\mathbb{C}[\lambda^\pm] \otimes_{\mathbb{C}[\lambda]} M^{IC}(X^\circ, \mathcal{L}) \simeq {}^*M_{A''}^{-(c+1,0,1)}/(\mathbb{C}[\lambda^\pm] \otimes_{\mathbb{C}[\lambda]} \ker(\phi))$, where ${}^*M_{A''}^{-(c+1,0,1)}$ is the module of global sections of ${}^*\mathcal{M}_{A''}^{-(c+1,0,1)}$. The notation $\mathbb{C}[\lambda]$ is shorthand for $\mathbb{C}[\lambda_0, \ldots, \lambda_m, \lambda_{m+1}, \ldots, \lambda_{m+c}]$, and the notation $\mathbb{C}[\lambda^\pm]$ is shorthand for $\mathbb{C}[\lambda_0, \ldots, \lambda_m, \lambda_{m+1}^\pm, \ldots, \lambda_{m+c}^\pm]$ (and not, as it is usual, shorthand for $\mathbb{C}[\lambda_0^\pm, \ldots, \lambda_m^\pm, \lambda_{m+1}^\pm, \ldots, \lambda_{m+c}^\pm]$).

We want to characterize $\mathbb{C}[\lambda^\pm] \otimes_{\mathbb{C}[\lambda]} \ker(\phi)$ inside ${}^*M_{A''}^{-(c+1,0,1)} = \mathbb{C}[\lambda^\pm] \otimes_{\mathbb{C}[\lambda]} M_{A''}^{-(c+1,0,1)}$. For this we define the following submodule in ${}^*M_{A''}^{-(c+1,0,1)}$:

$$K := \{m \in {}^*M_{A''}^{-(c+1,0,1)} \mid \exists p \in \mathbb{Z}, k \in \mathbb{N} \text{ such that } \partial_0^{k+1}(\lambda\partial + p)\ldots(\lambda\partial + p + k)m = 0\}.$$

Consider the following element of $\mathbb{C}[\lambda^\pm] \otimes_{\mathbb{C}[\lambda]} \ker(\phi)$:

$$\lambda_{m+1}^{-p_1}\ldots\lambda_{m+c}^{-p_c} \otimes m \quad \text{with} \quad p_1, \ldots, p_c \in \mathbb{N}, \tag{55}$$

i.e. there exists an $n \in \mathbb{N}$ such that $(\partial_0 \cdot \partial_{m+1}\ldots\partial_{m+c})^{n+1}m = 0$. Therefore we have

$$\begin{aligned}
0 &= \lambda_{m+1}^{-p_1}\ldots\lambda_{m+c}^{-p_c} \otimes (\partial_0 \cdot \partial_{m+1}\ldots\partial_{m+c})^{n+1}m \\
&= \lambda_{m+1}^{-p_1}\ldots\lambda_{m+c}^{-p_c} \otimes (\lambda_{m+1}\ldots\lambda_{m+c})^{n+1}(\partial_0 \cdot \partial_{m+1}\ldots\partial_{m+c})^{n+1}m \\
&= \partial_0^{n+1} \cdot (\lambda_{m+1}^{-p_1}\ldots\lambda_{m+c}^{-p_c} \otimes (\lambda\partial)\ldots(\lambda\partial - n)m) \\
&= \partial_0^{n+1}(\lambda\partial + p_{max})\ldots(\lambda\partial + p_{min} - n) \cdot (\lambda_{m+1}^{-p_1}\ldots\lambda_{m+c}^{-p_c} \otimes m) \\
&= \partial_0^{k+1}(\lambda\partial + p)\ldots(\lambda\partial + p + k) \cdot (\lambda_{m+1}^{-p_1}\ldots\lambda_{m+c}^{-p_c} \otimes m),
\end{aligned}$$

where $p_{max} := \max\{p_i\}, p_{min} := \min\{p_i\}$, $p := p_{min} - n$ and $k := p_{max} - p_{min} + n$. Because $\mathbb{C}[\lambda^\pm] \otimes_{\mathbb{C}[\lambda]} ker(\phi)$ is generated by elements of the form (55), we see that $\mathbb{C}[\lambda^\pm] \otimes_{\mathbb{C}[\lambda]} ker(\phi) \subset K$. Therefore we have a surjective morphism

$$\mathbb{C}[\lambda^\pm] \otimes_{\mathbb{C}[\lambda]} M^{IC}(X^\circ, \mathcal{L}) \simeq {}^*M_{A''}^{-(c+1,\underline{0},\underline{1})}/(\mathbb{C}[\lambda^\pm] \otimes_{\mathbb{C}[\lambda]} ker(\phi)) \twoheadrightarrow {}^*M_{A''}^{-(c+1,\underline{0},\underline{1})}/K.$$

Because $\mathbb{C}[\lambda^\pm] \otimes_{\mathbb{C}[\lambda]} M^{IC}(X^\circ, \mathcal{L})$ corresponds to the restriction of the simple $\mathcal{D}$-module $\mathcal{M}^{IC}(X^\circ, \mathcal{L})$ to the open subset $V^*$, it is itself simple. Thus, ${}^*M_{A''}^{-(c+1,\underline{0},\underline{1})}/K$ is either equal to 0 or is isomorphic to $\mathbb{C}[\lambda^\pm] \otimes_{\mathbb{C}[\lambda]} M^{IC}(X^\circ, \mathcal{L})$.

In order to prove the lemma, we are going to show that $K \subsetneq {}^*M_{A''}^{-(c+1,\underline{0},\underline{1})}$. Denote by $F_\bullet {}^*M_{A''}^{-(c+1,\underline{0},\underline{1})}$ the good filtration on ${}^*M_{A''}^{-(c+1,\underline{0},\underline{1})}$ which is induced by the order filtration on $D_{V^*}$. Notice that we have

$$K \subsetneq {}^*M_{A''}^{-(c+1,\underline{0},\underline{1})} \quad \Longleftrightarrow \quad gr^F K \subsetneq gr^F {}^*M_{A''}^{-(c+1,\underline{0},\underline{1})} \tag{56}$$

In order to show that $gr^F K \subsetneq gr^F {}^*M_{A''}^{-(c+1,\underline{0},\underline{1})}$, we first remark that

$$gr^F K \subset \{\overline{m} \in gr^F {}^*M_{A''}^{-(c+1,\underline{0},\underline{1})} \mid \exists k \in \mathbb{N} \text{ such that } \mu_0^{k+1} \underline{\lambda}^{k+1} \underline{\mu}^{k+1} \overline{m} = 0\},$$

where $\underline{\lambda} = (\lambda_{m+1} \cdots \lambda_{m+c})$, $\underline{\mu} = (\mu_{m+1} \cdots \mu_{m+c})$ and $\mu_i$ is the symbol $\sigma(\partial_{\lambda_i})$.

Thus, in order to show the right hand side of (56), it is enough to show that $\text{char}({}^*M_{A''}^{-(c+1,\underline{0},\underline{1})}) = supp(gr^F {}^*M_{A''}^{-(c+1,\underline{0},\underline{1})}) \subset T^*(V^*)$ is not contained in $\{\mu_0 \cdot \underline{\mu} \cdot \underline{\lambda} = 0\}$.

Therefore it is enough to find a vector $(\underline{\mu}', \underline{\lambda}') \in \text{char}({}^*M_{A''}^{-(c+1,\underline{0},\underline{1})}) \subset T^*(V^*)$ with $\mu_0' \cdot \underline{\mu}' \cdot \underline{\lambda}' \neq 0$, resp. a vector $(\underline{\mu}', \underline{\lambda}') \in \text{char}(M_{A''}^{-(c+1,\underline{0},\underline{1})}) \subset T^*(V)$ with $\mu_0' \cdot \underline{\mu}' \neq 0$ and $\lambda_i' \neq 0$ for $i = 1, \ldots, m+c$. Notice that we have

$$\text{char}(M_{A''}^{-(c+1,\underline{0},\underline{1})}) = \text{char}(M_{A''}^{(0,\underline{0},\underline{0})}) = \text{char}(\text{FL}(M_{A''}^{(0,\underline{0},\underline{0})})) = \text{char}(h_+\mathcal{O}_T),$$

where the first equality follows from [GKZ89, Theorem 4], the second equality follows e.g. from [Bry86, Corollaire 7.25] and the third equality follows from [SW09, Corollary 3.7]. Recall that the coordinates on $V'$ are denoted by $\mu_i$ for $i = 0, \ldots, m+c$ and the symbols of $\partial_{\mu_i}$ are denoted by $\lambda_i$. We now compute the fiber of $\text{char}(h_+\mathcal{O}_T) \to V'$ over the point $\underline{\mu} = (1, \ldots, 1)$. Recall that the map

$$h : T \longrightarrow V',$$
$$(y_0, \ldots, y_{n+c}) \mapsto (y^{\underline{a}_0''}, \ldots, y^{\underline{a}_{m+c}''})$$

can be factored into a closed embedding $h' : T \to (\mathbb{C}^*)^{m+c+1}$ and an open embedding $(\mathbb{C}^*)^{m+c+1} \to V'$. Therefore the fiber of the characteristic variety over $(1, \ldots, 1)$ is just the fiber of the conormal bundle of $h'(T)$ in $(\mathbb{C}^*)^{m+c+1}$. The tangent space of $h'(T)$ at $(1, \ldots, 1)$ is generated by

$$\sum_{i=0}^{m+c} a_{ki}'' \partial_{\mu_i} \quad \text{for} \quad k = 0, \ldots, n+c.$$

Therefore $(\underline{1}, \underline{\lambda}')$ lies in $\text{char}(h_+\mathcal{O}_T)$ if and only if $\sum_{i=0}^{m+c} a_{ki}'' \lambda_i' = 0$ for all $k = 0, \ldots, n+c$. So it remains to prove that there exists such a $\underline{\lambda}'$ with $\lambda_i' \neq 0$ for $i = 1, \ldots, m+c$. First notice that it is enough to construct a $(\lambda_1^\circ, \ldots, \lambda_{m+c}^\circ)$ with

$$\sum_{i=1}^{m+c} a_{ki}' \lambda_i^\circ = 0 \tag{57}$$

for all $k = 1, \ldots, n+c$ and $\lambda_i^\circ \neq 0$ for all $i = 1, \ldots, m+c$. Recall the structure of the matrix $A'$:

$$A' = \left( \begin{array}{c|c} A & 0_{n,c} \\ \hline (d_{ji}) & \text{Id}_c \end{array} \right), \tag{58}$$

where $d_{ji} \geq 0$ and the columns $\underline{a}_i$ of the matrix $A$ are the primitive integral generators of the rays of the fan $\Sigma$ corresponding to a complete, smooth toric variety $X_\Sigma$. This ensures the existence of $(\lambda_1^\circ, \ldots, \lambda_m^\circ) \in \mathbb{Z}_{>0}^m$ with $\sum_{i=1}^m \lambda_i^\circ \underline{a}_i = 0$. Setting $\lambda_{m+j}^\circ := -\sum_{i=1}^m d_{ji}\lambda_i^\circ$, we have constructed an element $(\lambda_1^\circ, \ldots, \lambda_{m+c}^\circ)$ with $\lambda_j^\circ \neq 0$ and satisfying $\sum_{i=1}^{m+c} a'_{ki}\lambda_i^\circ = 0$. Summarizing, this shows that $K \subsetneq {}^*M_{A''}^{-(c+1,\underline{0},1)}$, i.e.

$$\mathbb{C}[\lambda^\pm] \otimes_{\mathbb{C}[\lambda]} M^{IC}(X^\circ, \mathcal{L}) \simeq {}^*M_{A''}^{-(c+1,\underline{0},1)}/K \,.$$

Applying the localized Fourier-Laplace transformation to this isomorphism, we obtain the first isomorphism in the statement of the lemma. The second isomorphism follows from the $\mathcal{D}$-linearity of the isomorphism ${}^*\widehat{\mathcal{M}}_{A'}^{-(2c,\underline{0},1)} \simeq {}^*\widehat{\mathcal{N}}_{A'}^{(0,\underline{0},0)}$. $\square$

As in [RS15, section 3], we proceed by studying the restriction of the modules ${}^*\mathcal{M}_{A''}^\beta$, ${}^*\widehat{\mathcal{M}}_{A''}^\beta$ and ${}^*\widehat{\mathcal{N}}_{A''}^\beta$ to the Kähler moduli space of $\mathbb{V}(\mathcal{E}^\vee)$ as described in the second part of section 4 (see Lemma 4.4 and Proposition 4.7). The following construction has some overlap with the considerations in subsection 2.4 on which we comment later.

We apply $Hom_\mathbb{Z}(-, \mathbb{C}^*)$ to the exact sequence (46) to obtain the following exact sequence

$$1 \longrightarrow (\mathbb{C}^*)^{n+c} \longrightarrow (\mathbb{C}^*)^{m+c} \longrightarrow \mathbb{L}_{A'}^\vee \otimes \mathbb{C}^* \longrightarrow 1. \qquad (59)$$

We will identify the middle torus with $\operatorname{Spec} \mathbb{C}[\lambda_1^\pm, \ldots, \lambda_{m+c}^\pm]$, this space was called $W^*$ in section 2. Choose a basis $(p_1, \ldots, p_r)$ of $\mathbb{L}_{A'}^\vee$ with the following properties

1. $p_a \in \mathcal{K}_{\mathbb{V}(\mathcal{E}^\vee)} = \mathcal{K}_{X_\Sigma}$ for all $a = 1, \ldots, r$,

2. $\sum_{i=1}^{m+c} \overline{D}_i \in \sum_{a=1}^r \mathbb{R}_{\geq 0} p_a$.

Using the basis $(p_a)_{a=1,\ldots,r}$, we identify $\mathbb{L}_{A'}^\vee \otimes \mathbb{C}^*$ with $(\mathbb{C}^*)^r$ and obtain coordinates $q_1, \ldots, q_r$ on this space. We will write $\mathcal{KM}$ for this space and call it complexified Kähler moduli space. Notice that the choice of coordinates is considered as part of the data of $\mathcal{KM}$, that is, we really have $\mathcal{KM} = (\mathbb{C}^*)^r$ and not only $\mathcal{KM} = \mathbb{L}_{A'}^\vee \otimes \mathbb{C}^*$. Notice that this space already occurred in section 2.4 in a slightly more general context (which is consistent with the situation considered here, see the explanations after formula (30)).

Consider the embedding $\mathbb{L}_{A'} \hookrightarrow \mathbb{Z}^{m+c}$, which is given by a matrix $L \in \operatorname{Mat}((m+c) \times r, \mathbb{Z})$ with respect to the basis $p_a^\vee$ of $\mathbb{L}_{A'}$ and the natural basis of $\mathbb{Z}^{m+c}$. Choose a section $\mathbb{Z}^{m+c} \to \mathbb{L}_{A'}$ of this inclusion, which is given by a matrix $M \in \operatorname{Mat}(r \times (m+c), \mathbb{Z})$. This defines a section on the dual lattices, i.e. a section $\mathbb{L}_{A'}^\vee \to \mathbb{Z}^{m+c}$ of the projection $\mathbb{Z}^{m+c} \to \mathbb{L}_{A'}^\vee$ and a closed embedding $\varrho' : \mathcal{KM} = (\mathbb{C}^*)^r \hookrightarrow W^*$. We will need to consider a slight twist of this morphism. Let $\iota : W^* \to W^*$ be the involution given by $\iota(\lambda_i) := (-1)^{\varepsilon(i)}\lambda_i$ with $\varepsilon(i) = 0$ for $i = 1, \ldots, m$ and $\varepsilon(i) = 1$ for $i = m+1, \ldots, m+c$.

We will further restrict our objects of study to that part of the complexified Kähler moduli space which maps to the set of good parameters in $W = \mathbb{C}^{m+c}$ as discussed in subsection 3.2. Hence we put $\mathcal{KM}^\circ := (\iota \circ \varrho')^{-1}(W^\circ) \subset \mathcal{KM}$, and write

$$\varrho := \iota \circ \varrho' : \mathcal{KM}^\circ \hookrightarrow W^* \,.$$

We can now define the main object of study of this paper. We are going to use the constructions of the subsections 2.4 and 3.2, in particular, the diagrams (25), (31) and (40). We consider the composed morphism $\alpha \circ \beta : \mathcal{Z}_X \to \mathbb{C}_{\lambda_0} \times \mathcal{KM}$ as defined by diagram (31). Let $\mathcal{Z}_X^\circ := (\alpha \circ \beta)^{-1}(\mathbb{C}_{\lambda_0} \times \mathcal{KM}^\circ) \subset \mathcal{Z}_X$ be the subspace which is parameterized by the good parameter locus $\mathcal{KM}^\circ$ inside $\mathcal{KM}$.

For future reference, let us collect the relevant morphisms once again in a diagram, in which the spaces $\mathcal{Z}^\circ$, $\mathcal{Z}_{X^{\text{aff}}}^\circ$ and $\mathcal{Z}_X^\circ$ are defined by the requirement that all squares are cartesian. For simplicity of the notation, we denote by $\alpha, \beta, \gamma_1$ and $\gamma_2$ also the corresponding restrictions above $\mathbb{C}_{\lambda_0} \times \mathcal{KM}^\circ$.

$$\begin{array}{ccccccccc}
S & \xleftarrow{\pi_1^S} & \Gamma \cong S \times W & \longleftarrow & \Gamma^* \cong S \times W^* & \longleftarrow & S \times \mathcal{KM} & \longleftarrow & S \times \mathcal{KM}^\circ \\
{\scriptstyle j_2}\downarrow & & {\scriptstyle \theta_2}\downarrow & & {\scriptstyle \zeta_2}\downarrow & & {\scriptstyle \gamma_2}\downarrow & & {\scriptstyle \gamma_2}\downarrow \\
X^{\mathrm{aff}} & \longleftarrow & Z_{X^{\mathrm{aff}}} \cong X^{\mathrm{aff}} \times W & \longleftarrow & Z^*_{X^{\mathrm{aff}}} \cong X^{\mathrm{aff}} \times W^* & \longleftarrow & \mathcal{Z}_{X^{\mathrm{aff}}} \cong X^{\mathrm{aff}} \times \mathcal{KM} & \longleftarrow & \mathcal{Z}^\circ_{X^{\mathrm{aff}}} \cong X^{\mathrm{aff}} \times \mathcal{KM}^\circ \\
{\scriptstyle j_1}\downarrow & & {\scriptstyle \theta_1}\downarrow & & {\scriptstyle \zeta_1}\downarrow & & {\scriptstyle \gamma_1}\downarrow & & {\scriptstyle \gamma_1}\downarrow \\
X & \longleftarrow & Z_X & \longleftarrow & Z_X^* & \longleftarrow & \mathcal{Z}_X & \longleftarrow & \mathcal{Z}_X^\circ \\
{\scriptstyle i}\downarrow & & {\scriptstyle \eta}\downarrow & & {\scriptstyle \varepsilon}\downarrow & & {\scriptstyle \beta}\downarrow & & {\scriptstyle \beta}\downarrow \\
\mathbb{P}(V') & \xleftarrow{\pi_1^Z} & Z & \longleftarrow & Z^* & \longleftarrow & \mathcal{Z} & \longleftarrow & \mathcal{Z}^\circ \\
& & {\scriptstyle \pi_2^Z}\downarrow & & {\scriptstyle \delta}\downarrow & & {\scriptstyle \alpha}\downarrow & & {\scriptstyle \alpha}\downarrow \\
& & V & \longleftarrow & V^* & \longleftarrow & \mathbb{C}_{\lambda_0} \times \mathcal{KM} & \longleftarrow & \mathbb{C}_{\lambda_0} \times \mathcal{KM}^\circ
\end{array}$$
$$\mathrm{id}_{\mathbb{C}_{\lambda_0}} \times \varrho$$

(60)

**Definition 6.3.** *The **non-affine Landau-Ginzburg model** associated to $(X_\Sigma, \mathcal{L}_1, \ldots, \mathcal{L}_c)$ is the morphism*

$$\Pi : \mathcal{Z}_X^\circ \longrightarrow \mathbb{C}_{\lambda_0} \times \mathcal{KM}^\circ,$$

*which is by definition the restriction of the universal family of hyperplane sections of $X$, i.e, of the morphism $\pi_2^Z \circ \eta : Z_X \to V$ to the parameter space $\mathcal{KM}^\circ$. We recall once again that $X$ is defined as the closure of the embedding $g : S \to \mathbb{P}(V')$ sending $(y_1, \ldots, y_{n+c})$ to $(1 : \underline{y}^{\underline{a}'_1} : \ldots : \underline{y}^{\underline{a}'_{m+c}})$ where $\underline{a}'_i$ are the columns of the matrix $A'$ from Definition 4.6.*

*We also consider the restrictions $\overline{\pi} = \alpha \circ \beta \circ \gamma_1 : \mathcal{Z}^\circ_{X^{\mathrm{aff}}} \cong X^{\mathrm{aff}} \times \mathcal{KM}^\circ \to \mathbb{C}_{\lambda_0} \times \mathcal{KM}^\circ$ resp. $\pi = \alpha \circ \beta \circ \gamma_1 \circ \gamma_2 : S \times \mathcal{KM}^\circ :\to \mathbb{C}_{\lambda_0} \times \mathcal{KM}^\circ$. These are nothing but the family of Laurent polynomials*

$$(\underline{y}, \underline{q}) \longmapsto \left( -\sum_{i=1}^m \underline{q}^{\underline{m}_i} \cdot \underline{y}^{\underline{a}'_i} + \sum_{i=m+1}^{m+c} \underline{q}^{\underline{m}_i} \cdot \underline{y}^{\underline{a}'_i}, \underline{q} \right),$$

*where the monomial $\underline{y}^{\underline{a}'_i}$ is seen as an element of $\mathcal{O}_{X^{\mathrm{aff}}}$ in the first case and as an element of $\mathcal{O}_S$ in the second case. Here $\underline{m}_i$ is the $i$'th column of the matrix $M \in \mathrm{Mat}(r \times (m+c), \mathbb{Z})$ from above. Notice that the first component of $\pi$ has been split in two sums with opposite signs of each summand due to the action of the involution $\iota$ entering in the definition of the morphism $\varrho : \mathcal{KM}^\circ \hookrightarrow W^*$. Both morphisms $\pi$ and $\overline{\pi}$ are called the **affine Landau-Ginzburg model** of $(X_\Sigma, \mathcal{L}_1, \ldots, \mathcal{L}_c)$.*

As we will see later, the affine Landau-Ginzburg model is related to the twisted quantum $\mathcal{D}$-module $\mathrm{QDM}(X_\Sigma, \mathcal{E})$ whereas the reduced quantum $\mathcal{D}$-module $\overline{\mathrm{QDM}}(X_\Sigma, \mathcal{E})$ can be obtained from the non-affine Landau-Ginzburg model $\Pi : \mathcal{Z}_X^\circ \to \mathbb{C}_{\lambda_0} \times \mathcal{KM}^\circ$. The next results are parallel to [RS15, corollary 3.3. and corollary 3.4]. They show that the calculation of the Gauß-Manin system resp. the intersection cohomology $\mathcal{D}$-module from section 2 can be used to describe the corresponding objects for the morphism $\Pi$.

We consider, as in subsection 3.1, the localized partial Fourier-Laplace transformation, this time with base $\mathcal{KM}^\circ$, that is, let $j_\tau : \mathbb{C}_\tau^* \times \mathcal{KM}^\circ \hookrightarrow \mathbb{C}_\tau \times \mathcal{KM}^\circ$, $j_z : \mathbb{C}_\tau^* \times \mathcal{KM}^\circ \hookrightarrow \mathbb{C}_z \times \mathcal{KM}^\circ$ then we put $\mathrm{FL}^{loc}_{\mathcal{KM}^\circ} := j_{z,+} j_\tau^+ \mathrm{FL}_{\mathcal{KM}^\circ}$.

**Lemma 6.4.** *We have*

$$\mathrm{FL}^{loc}_{\mathcal{KM}^\circ} \left( \mathcal{H}^0 \pi_+ \mathcal{O}_{S \times \mathcal{KM}^\circ} \right) \cong (\mathrm{id}_{\mathbb{C}_z} \times \varrho)^{+\,*} \widehat{\mathcal{M}}^{(-c, \underline{0}, \underline{0})}_{A'}.$$

*Similarly, the isomorphism*

$$\mathrm{FL}^{loc}_{\mathcal{KM}^\circ} \left( \mathcal{H}^0 \pi_{\dagger} \mathcal{O}_{S \times \mathcal{KM}^\circ} \right) \cong (\mathrm{id}_{\mathbb{C}_z} \times \varrho)^{+\,*} \widehat{\mathcal{N}}^{(0, \underline{0}, \underline{0})}_{A'}$$

*holds.*

Notice that the embedding $(\mathrm{id}_{\mathbb{C}_z} \times \varrho)$ is obviously non-characteristic for both of the modules ${}^*\widehat{\mathcal{M}}_{A''}^{(-c,\underline{0},\underline{0})}$ and ${}^*\widehat{\mathcal{N}}_{A''}^{(0,\underline{0},\underline{0})}$ as their singular locus is contained in

$$(\{0,\infty\} \times \mathcal{KM}^\circ) \cup \left(\mathbb{P}^1_z \times (W^* \backslash \mathcal{KM}^\circ)\right).$$

Hence, the complexes $(\mathrm{id}_{\mathbb{C}_z} \times \varrho)^+ \, {}^*\widehat{\mathcal{M}}_{A'}^{(-c,\underline{0},\underline{0})}$ and $(\mathrm{id}_{\mathbb{C}_z} \times \varrho)^+ \, {}^*\widehat{\mathcal{N}}_{A'}^{(0,\underline{0},\underline{0})}$ have cohomology only in degree zero.

*Proof.* The proof of the first isomorphism is the same as [RS15, corollary 3.3]: Consider the cartesian diagram (which is part of the diagram (60))

$$\begin{array}{ccc} S \times \mathcal{KM}^\circ & \longrightarrow & \Gamma^* \cong S \times W^* \\ \pi \downarrow & & \downarrow \varphi \\ \mathbb{C}_{\lambda_0} \times \mathcal{KM}^\circ & \xrightarrow{\mathrm{id}_{\mathbb{C}_{\lambda_0}} \times \varrho} & V^* = \mathbb{C}_{\lambda_0} \times W^* \end{array} \qquad (61)$$

then the base change property (Theorem 2.1) and the commutation of $\mathrm{FL}^{loc}$ with inverse images shows that

$$\mathrm{FL}^{loc}_{\mathcal{KM}^\circ}(\mathcal{H}^0 \pi_+ \mathcal{O}_{S \times \mathcal{KM}^\circ}) \cong (\mathrm{id}_{\mathbb{C}_z} \times \varrho)^+ \mathcal{G}^+_{|V^*},$$

where $\mathcal{G}^+$ is the $\mathcal{D}_{\mathbb{C}_{\lambda_0} \times W}$-module introduced in subsection 3.1, and then one concludes using Proposition 3.3.

Concerning the second isomorphism, we use base change (with respect to the morphism $\mathrm{id}_{\mathbb{C}_{\lambda_0}} \times \varrho$ in diagram (61)) for proper direct images and exceptional inverse images. However, the latter ones equal ordinary inverse images if the horizontal morphisms in the above diagram are non-characteristic for the modules in question. This is the case by Proposition 2.22, 2., so that we obtain

$$\mathrm{FL}^{loc}_{\mathcal{KM}^\circ}\left(\mathcal{H}^0 \pi_\dagger \mathcal{O}_{S \times \mathcal{KM}^\circ}\right) \cong (\mathrm{id}_{\mathbb{C}_z} \times \varrho)^+ \mathrm{FL}^{loc}_W \left(\mathcal{H}^0 \varphi_{B,\dagger} \mathcal{O}_{S \times W}\right)_{|V^*} = (\mathrm{id}_{\mathbb{C}_z} \times \varrho)^+ \mathcal{G}^\dagger_{|V^*}.$$

The second part of corollary 5.5 (and the second part of Proposition 3.3) tells us that $\mathcal{G}^\dagger \cong \widehat{\mathcal{M}}_{A'}^{-(c,\underline{0},\underline{1})}$. However, the isomorphism $\Psi : {}^*\widehat{\mathcal{N}}_{A'}^{(0,\underline{0},\underline{0})} \longrightarrow {}^*\widehat{\mathcal{M}}_{A'}^{-(2c,\underline{0},\underline{1})}$ given by right multiplication with $z^c \cdot \lambda_{m+1} \cdots \lambda_{m+c}$ (see equation (53)) shows that

$$(\mathrm{id}_{\mathbb{C}_z} \times \varrho)^+ \, {}^*\widehat{\mathcal{M}}_{A'}^{-(2c,\underline{0},\underline{1})} \cong (\mathrm{id}_{\mathbb{C}_z} \times \varrho)^+ \, {}^*\widehat{\mathcal{N}}_{A'}^{(0,\underline{0},\underline{0})}$$

so that finally we arrive at the desired equality

$$\mathrm{FL}^{loc}_{\mathcal{KM}^\circ}\left(\mathcal{H}^0 \pi_\dagger \mathcal{O}_{S \times \mathcal{KM}^\circ}\right) \cong (\mathrm{id}_{\mathbb{C}_z} \times \varrho)^+ \, {}^*\widehat{\mathcal{N}}_{A'}^{(0,\underline{0},\underline{0})}.$$

□

Next we show the analog of Proposition 3.21 for the morphism $\overline{\pi}$.

**Lemma 6.5.** *Let $\widetilde{F} : X^{aff} \times \mathcal{KM}^\circ \to \mathbb{C}_{\lambda_0}$ be the first component of the morphism $\overline{\pi}$, then we have the following isomorphism of $\mathcal{R}_{\mathbb{C}_z \times \mathcal{KM}^\circ}$-modules*

$$z^{-c} H^{n+c}(\Omega^\bullet_{X^{aff} \times \mathcal{KM}^\circ / \mathcal{KM}^\circ}(\log D)[z], zd - d\widetilde{F}) \cong (\mathrm{id}_{\mathbb{C}_z} \times \varrho)^* \left({}_0\widehat{\mathcal{M}}_{A'}^{(-c,\underline{0},\underline{0})}\right). \qquad (62)$$

*Proof.* In order to show the statement, notice that by definition $H^{n+c}(\Omega^\bullet_{X^{aff} \times W^*/W^*}(\log D)[z], zd - dF)$ is the cokernel of

$$\Omega^{n+c-1}_{X^{aff} \times W^*/W^*}(\log D)[z] \xrightarrow{zd-dF} \Omega^{n+c}_{X^{aff} \times W^*/W^*}(\log D)[z],$$

that is, the cokernel of an $\mathcal{O}_{\mathbb{C}_z \times W^*}$-linear morphism between free (though not coherent) $\mathcal{O}_{\mathbb{C}_z \times W^*}$-modules. Hence tensoring with $\mathcal{O}_{\mathbb{C}_z \times \mathcal{KM}^\circ}$ yields the exact sequence

$$\Omega^{n+c-1}_{X^{aff} \times \mathcal{KM}^\circ / \mathcal{KM}^\circ}(\log D)[z] \xrightarrow{zd-d\widetilde{F}} \Omega^{n+c}_{X^{aff} \times \mathcal{KM}^\circ / \mathcal{KM}^\circ}(\log D)[z] \longrightarrow$$

$$\mathcal{O}_{\mathbb{C}_z \times \mathcal{KM}^\circ} \otimes_{\mathcal{O}_{\mathbb{C}_z \times W^*}} H^{n+c}(\Omega^\bullet_{X^{aff} \times W^*/W^*}(\log D)[z], zd - dF) \longrightarrow 0$$

from which we conclude that

$$H^{n+c}(\Omega^\bullet_{X^{aff}\times\mathcal{KM}^\circ/\mathcal{KM}^\circ}(\log D)[z], zd-d\widetilde{F}) =$$

$$\mathcal{O}_{\mathbb{C}_z\times\mathcal{KM}^\circ}\otimes_{\mathcal{O}_{\mathbb{C}_z\times W^*}} H^{n+c}(\Omega^\bullet_{X^{aff}\times W^*/W^*}(\log D)[z], zd-dF).$$

Notice that the restriction functor $(\mathcal{O}_{\mathbb{C}_z\times\mathcal{KM}^\circ}\otimes_{\mathcal{O}_{\mathbb{C}_z\times W^*}} -)$ is defined via the embedding $\varrho:\mathcal{KM}^\circ\hookrightarrow W^*$, and hence involves the involution $\iota$. Therefore the function $\widetilde{F}$ appears on the left hand side of the last formula, whereas on the right hand side we have to put $F$.

We know by Proposition 3.21 that

$$z^{-c}H^{n+c}(\Omega^\bullet_{X^{aff}\times W^\circ/W^\circ}(\log D)[z], zd-dF) \cong z^{-c}{}_0\widehat{\mathcal{M}}^{(0,\underline{0},\underline{0})}_{A'}.$$

On the other hand, we know from equation (38) that right multiplication by $z^c$ induces an isomorphism

$$\left(\widehat{\mathcal{M}}^{(0,\underline{0},\underline{0})}_{A'}\right)_{|\mathbb{C}_z\times W^\circ} \xrightarrow{\cdot z^c} \left(\widehat{\mathcal{M}}^{(-c,\underline{0},\underline{0})}_{A'}\right)_{|\mathbb{C}_z\times W^\circ}$$

which maps $z^{-c}{}_0\widehat{\mathcal{M}}^{(0,\underline{0},\underline{0})}_{A'} \subset \left(\widehat{\mathcal{M}}^{(0,\underline{0},\underline{0})}_{A'}\right)_{|\mathbb{C}_z\times W^\circ}$ isomorphically to ${}_0\widehat{\mathcal{M}}^{(-c,\underline{0},\underline{0})}_{A'} \subset \left(\widehat{\mathcal{M}}^{(0,\underline{0},\underline{0})}_{A'}\right)_{|\mathbb{C}_z\times W^\circ}$. The desired statement, i.e., Formula (62) follows as the restriction map $\varrho:\mathcal{KM}^\circ\hookrightarrow W^*$ factors by definition over $W^\circ$. □

Similarly to the last statement, we now give a geometric interpretation of the (restriction to $\mathbb{C}_z\times\mathcal{KM}^\circ$ of the) modules ${}_0^*\widehat{\mathcal{N}}^{(0,\underline{0},\underline{0})}_{A'}$ using the twisted relative logarithmic de Rham complex on $X^{aff}\times\mathcal{KM}^\circ$. We need some preliminary notations. Denote by $(-)'$ the duality functor in the category of locally free $\mathcal{O}_{\mathbb{C}_z\times\mathcal{KM}^\circ}$-modules with meromorphic connection with poles along $\{0\}\times\mathcal{KM}^\circ$, that is, if $(\mathcal{F},\nabla)$ is an object of this category, we put $(\mathcal{F},\nabla)' := (\mathcal{H}om_{\mathcal{O}_{\mathbb{C}_z\times\mathcal{KM}^\circ}}(\mathcal{F},\mathcal{O}_{\mathbb{C}_z\times\mathcal{KM}^\circ}),\nabla')$, where $\nabla'$ is the dual connection. Notice that the $\mathcal{R}_{\mathbb{C}_z\times\mathcal{KM}^\circ}$-modules from isomorphism (62) are actually objects of this category. Notice also that the duality functor in the category of $\mathcal{R}_{\mathbb{C}_z\times\mathcal{KM}^\circ}$-modules (i.e., the functor $\mathcal{E}xt^{r+1}_{\mathcal{R}_{\mathbb{C}_z\times\mathcal{KM}^\circ}}(-,\mathcal{R}_{\mathbb{C}_z\times\mathcal{KM}^\circ})$) restricts to $(-)'$ on the subcategory described above (this follows from [DS03, Lemma A.12]).

As a piece of notation, for any complex manifold $M$ we denote by $\sigma$ the involution of $\mathbb{C}_z\times M$ defined by $(z,\underline{x})\mapsto(-z,\underline{x})$.

**Lemma 6.6.** *There is an isomorphism of $\mathcal{R}_{\mathbb{C}_z\times\mathcal{KM}^\circ}$-modules*

$$\sigma^*z^n\left(H^{n+c}(\Omega^\bullet_{X^{aff}\times\mathcal{KM}^\circ/\mathcal{KM}^\circ}(\log D)[z], zd-d\widetilde{F})\right)' \xrightarrow{\cong} (\mathrm{id}_{\mathbb{C}_z}\times\varrho)^*\left({}_0^*\widehat{\mathcal{N}}^{(0,\underline{0},\underline{0})}_{A'}\right).$$

*Proof.* Consider the filtration on $\mathcal{D}_{\mathbb{C}_z\times W}$ resp. on $\mathcal{D}_{\mathbb{C}_z\times W^*}$ which extends the order filtration on $\mathcal{D}_W$ (resp. on $\mathcal{D}_{W^*}$) and for which $z$ has degree $-1$ and $\partial_z$ has degree 2. Denote by $G_\bullet$ the induced filtrations on the modules ${}^*\widehat{\mathcal{N}}^{(0,\underline{0},\underline{0})}_{A'}$ and $\widehat{\mathcal{M}}^{(-c,\underline{0},\underline{0})}_{A'}$ resp. on ${}^*\widehat{\mathcal{M}}^{(-c,\underline{0},\underline{0})}_{A'}$, in particular, we have $G_0\left({}^*\widehat{\mathcal{N}}^{(0,\underline{0},\underline{0})}_{A'}\right) = {}_0^*\widehat{\mathcal{N}}^{(0,\underline{0},\underline{0})}_{A'}$ and $G_0\left(\widehat{\mathcal{M}}^{(-c,\underline{0},\underline{0})}_{A'}\right) = {}_0\widehat{\mathcal{M}}^{(-c,\underline{0},\underline{0})}_{A'}$ resp. $G_0\left({}^*\widehat{\mathcal{M}}^{(-c,\underline{0},\underline{0})}_{A'}\right) = {}_0^*\widehat{\mathcal{M}}^{(-c,\underline{0},\underline{0})}_{A'}$.

Similar to the proof of [RS15, Proposition 2.18, 3.], we consider the saturation of the filtration $F_\bullet$ on $\mathcal{M}^\beta_{A''}$ by $\partial^{-1}_{\lambda_0}$. More precisely, we first notice that Lemma 3.2 can be reformulated by saying that for any $\beta' = (\beta'_0,\beta'_1,\ldots,\beta'_{n+c})\in\mathbb{Z}^{1+n+c}$, we have

$$\widehat{\mathcal{M}}^\beta_{A'} = \mathrm{FL}_W\left(\mathcal{M}^{\beta'}_{A''}[\partial^{-1}_{\lambda_0}]\right),$$

where $\beta_0 = \beta'_0 + 1$ and $\beta_i = \beta'_i$ for $i=1,\ldots,n+c$ and where we write $\mathcal{M}^{\beta'}_{A''}[\partial^{-1}_{\lambda_0}] := \mathcal{D}_V[\partial^{-1}_{\lambda_0}]\otimes_{\mathcal{D}_V}\mathcal{M}^{\beta'}_{A''}$. Now we consider the natural localization morphism $\widehat{\mathrm{loc}}:\mathcal{M}^{\beta'}_{A''}\to\mathcal{M}^{\beta'}_{A''}[\partial^{-1}_{\lambda_0}]$ and we put

$$F_k\mathcal{M}^{\beta'}_{A''}[\partial^{-1}_{\lambda_0}] := \sum_{j\geq 0}\partial^{-j}_{\lambda_0}\widehat{\mathrm{loc}}\left(F_{k+j}\mathcal{M}^{\beta'}_{A''}\right).$$

As we have

$$F_k\mathcal{M}^{\beta'}_{A''}[\partial^{-1}_{\lambda_0}] = im\left(\partial^k_{\lambda_0}\mathbb{C}[\lambda_0,\lambda_1,\ldots,\lambda_{m+c}]\langle\partial^{-1}_{\lambda_0},\partial^{-1}_{\lambda_0}\partial_{\lambda_1},\ldots,\partial^{-1}_{\lambda_0}\partial_{\lambda_{m+c}}\rangle\right) \text{ in } \mathcal{M}^{\beta'}_{A''}[\partial^{-1}_{\lambda_0}],$$

the filtration induced by $F_k\mathcal{M}_{A''}^{\beta'}[\partial_{\lambda_0}^{-1}]$ on $\widehat{\mathcal{M}}_{A'}^{\beta}$ is precisely $G_k\widehat{\mathcal{M}}_{A'}^{\beta}$. We conclude from [Sai89, formula 2.7.5] and from the fact that Fourier-Laplace transformation commutes with the duality functor up to the action of $\sigma$ that

$$(G_l\widehat{\mathcal{M}}_{A'}^{-(c,\underline{0},\underline{1})})' = \mathcal{H}om_{\mathcal{O}_{\mathbb{C}_z \times W}}\left(G_l\widehat{\mathcal{M}}_{A'}^{-(c,\underline{0},\underline{1})}, \mathcal{O}_{\mathbb{C}_z \times W}\right) \stackrel{!}{=} \sigma^* G_{l+(m+c+2)}^{\mathbb{D}} \widehat{\mathcal{M}}_{A'}^{(1,\underline{0},\underline{0})},$$

where $G_\bullet^{\mathbb{D}} \widehat{\mathcal{M}}_{A'}^{(1,\underline{0},\underline{0})}$ is the filtration induced by the saturation of the filtration on $\mathcal{M}_{A''}^{(0,\underline{0},\underline{0})}$ dual to the order filtration $F_\bullet$ on $\mathcal{M}_{A''}^{-(c+1,\underline{0},\underline{1})}$. By Theorem 5.4, 2. and by restriction to $\mathbb{C}_z \times W^*$ we obtain

$$G_\bullet^{\mathbb{D}} {}^*\widehat{\mathcal{M}}_{A'}^{(1,\underline{0},\underline{0})} = G_{\bullet+n-(m+c+1)} {}^*\widehat{\mathcal{M}}_{A'}^{(1,\underline{0},\underline{0})}.$$

Hence

$$(G_l {}^*\widehat{\mathcal{M}}_{A'}^{-(c,\underline{0},\underline{1})})' = \sigma^* G_{l+n+1} {}^*\widehat{\mathcal{M}}_{A'}^{(1,\underline{0},\underline{0})}.$$

Now we use the fact that for any $k \in \mathbb{Z}$, the isomorphism (see equation (38))

$$\cdot z^k : \widehat{\mathcal{M}}_{A'}^{(\beta_0,\beta)} \stackrel{\cong}{\longrightarrow} \widehat{\mathcal{M}}_{A'}^{(\beta_0-k,\beta)}$$

sends $G_k\widehat{\mathcal{M}}_{A'}^{(\beta_0,\beta)} = z^{-k}{}_0\widehat{\mathcal{M}}_{A'}^{(\beta_0,\beta)}$ to $G_0\widehat{\mathcal{M}}_{A'}^{(\beta_0-k,\beta)} = {}_0\widehat{\mathcal{M}}_{A'}^{(\beta_0-k,\beta)}$. Therefore (setting $l = 0$) we have

$$(G_0\widehat{\mathcal{M}}_{A'}^{-(c,\underline{0},\underline{1})})' \cong \sigma^* G_{n+1} {}^*\widehat{\mathcal{M}}_{A'}^{(1,\underline{0},\underline{0})} = \sigma^* G_n {}^*\widehat{\mathcal{M}}_{A'}^{(0,\underline{0},\underline{0})}.$$

which implies

$$G_0 {}^*\widehat{\mathcal{M}}_{A'}^{-(c,\underline{0},\underline{1})} \cong \left(\sigma^* G_n {}^*\widehat{\mathcal{M}}_{A'}^{(0,\underline{0},\underline{0})}\right)'$$

The isomorphism $\Psi$ from Formula (53) satisfies

$$\Psi : {}_0^*\widehat{\mathcal{N}}_{A'}^{(0,\underline{0},\underline{0})} \stackrel{\cong}{\longrightarrow} z^c \cdot {}_0^*\widehat{\mathcal{M}}_{A'}^{-(2c,\underline{0},\underline{1})} \cong {}_0^*\widehat{\mathcal{M}}_{A'}^{-(c,\underline{0},\underline{1})}$$

In conclusion, we obtain

$${}_0^*\widehat{\mathcal{N}}_{A'}^{(0,\underline{0},\underline{0})} \cong \left(\sigma^* z^{-n} \cdot {}_0^*\widehat{\mathcal{M}}_{A'}^{(0,\underline{0},\underline{0})}\right)' \cong \sigma^* z^n \cdot \left({}_0^*\widehat{\mathcal{M}}_{A'}^{(0,\underline{0},\underline{0})}\right)',$$

and then the statement follows from Proposition 3.21 as the inverse image under $\mathrm{id}_{\mathbb{C}_z} \times \varrho^*$ commutes with the functor $(-)'$. $\square$

Now we can construct a $\mathcal{D}_{\mathbb{C}_z \times \mathcal{KM}^\circ}$-module from the non-affine Landau-Ginzburg model $\Pi : \mathcal{Z}_X^\circ \longrightarrow \mathbb{C}_{\lambda_0} \times \mathcal{KM}^\circ$ that will ultimately give us the reduced quantum $\mathcal{D}$-module. It will consist in a minimal extension of the local system of intersection cohomologies of the fibres of $\Pi$. Recall that $\mathcal{M}^{IC}(\mathcal{Z}_X^\circ)$ is the intersection cohomology $\mathcal{D}$-module of $\mathcal{Z}_X^\circ$, that is, the unique regular singular $\mathcal{D}_{\mathcal{Z}^\circ}$-module supported on $\mathcal{Z}_X^\circ$ which corresponds to the intermediate extension of the constant sheaf on the smooth part of $\mathcal{Z}_X^\circ$.

**Proposition 6.7.** *1. Consider the local system $\mathcal{L}$ from Proposition 2.13. Then*

$$\mathcal{H}^0\alpha_+\mathcal{M}^{IC}(\mathcal{Z}_X^\circ) \cong (\mathrm{id}_{\mathbb{C}_{\lambda_0}} \times \varrho)^+ \left(\mathcal{M}^{IC}(X^\circ, \mathcal{L}) \oplus \left(IH^{n+c-1}(X) \otimes \mathcal{O}_V\right)\right)_{|V^*}.$$

*Using the Riemann-Hilbert correspondence, the above isomorphism can be expressed in terms of the morphism $\Pi$ as*

$${}^p\mathcal{H}^0 R\Pi_* IC(\mathcal{Z}_X^\circ) \cong (\mathrm{id}_{\mathbb{C}_{\lambda_0}} \times \varrho)^{-1} \left((j_{X^\circ})_! IC(X^\circ, \mathcal{L}) \oplus \underline{IH}^{n+c-1}(X)\right)_{|V^*},$$

*where ${}^p\mathcal{H}$ denotes the perverse cohomology functor, where $j_{X^\circ} : X_0 \hookrightarrow V$ is the canonical closed embedding and where $\underline{IH}^{n+c-1}(X)$ is the constant sheaf on $V$ with fibre $IH^{n+c-1}(X)$.*

*2. We have isomorphisms of $\mathcal{D}_{\mathbb{C}_z \times \mathcal{KM}^\circ}$-modules*

$$\mathrm{FL}_{\mathcal{KM}^\circ}^{loc}\left(\mathcal{H}^0\alpha_+\mathcal{M}^{IC}(\mathcal{Z}_X^\circ)\right) \cong (\mathrm{id}_{\mathbb{C}_z} \times \varrho)^+ \widehat{\mathcal{M}^{IC}}(X^\circ, \mathcal{L})_{|V^*} \cong (\mathrm{id}_{\mathbb{C}_z} \times \varrho)^+ \mathrm{im}(\widetilde{\phi}),$$

*where $\widetilde{\phi} : {}^*\widehat{\mathcal{N}}_{A'}^{(0,\underline{0},\underline{0})} \longrightarrow {}^*\widehat{\mathcal{M}}_{A'}^{(-c,\underline{0},\underline{0})}$ is the morphism introduced in Definition 6.1.*

*Proof.* 1. As the inclusion $\mathcal{Z}^\circ \hookrightarrow \mathcal{Z}$ is open and hence non-characteristic for any $\mathcal{D}_{\mathcal{Z}}$-module, the assertion to be shown follows from Proposition 2.22 (more precisely, from Formula (32)) and Proposition 2.13.

2. The first isomorphism is a direct consequence of the last point, using again the commutation of $\mathrm{FL}^{loc}$ with the inverse image and the fact that $\mathcal{O}_V$-free modules are killed by $\mathrm{FL}^{loc}_W$. The second isomorphism follows from equation (54).

$\square$

For future use, we give names to the $\mathcal{D}$-modules on the Kähler moduli space considered above. We also define natural lattices inside them.

**Definition 6.8.** *Define the following $\mathcal{D}_{\mathbb{C}_z \times \mathcal{KM}^\circ}$-modules:*
$$\mathcal{QM}_{A'} := (\mathrm{id}_{\mathbb{C}_z} \times \varrho)^+ \left({}^*\widehat{\mathcal{N}}_{A'}^{(0,0,0)}\right) \qquad \text{and} \qquad \mathcal{QM}_{A'}^{IC} := (\mathrm{id}_{\mathbb{C}_z} \times \varrho)^+ \left(\mathrm{im}(\widetilde{\phi})\right).$$

*Define moreover*
$$_0\mathcal{QM}_{A'} := (\mathrm{id}_{\mathbb{C}_z} \times \varrho)^* \left({}_0^*\widehat{\mathcal{N}}_{A'}^{(0,0,0)}\right) \qquad \text{and} \qquad _0\mathcal{QM}_{A'}^{IC} := (\mathrm{id}_{\mathbb{C}_z} \times \varrho)^* \left(\widetilde{\phi}\left({}_0^*\widehat{\mathcal{N}}_{A'}^{(0,0,0)}\right)\right),$$

*where here the functor $(\mathrm{id}_{\mathbb{C}_z} \times \varrho)^*$ is the inverse image in the category of holomorphic vector bundles on $\mathbb{C}_z \times \mathcal{KM}^\circ$ with meromorphic connection (meromorphic along $\{0\} \times \mathcal{KM}^\circ$).*

We proceed by comparing the objects $\mathcal{QM}_{A'}$ and $\mathcal{QM}_{A'}^{IC}$ just introduced to the twisted and the reduced quantum $\mathcal{D}$-module from section 4. For the readers convenience, let us recall one of the main results from [MM11] which concerns the toric description of the twisted resp. reduced quantum $\mathcal{D}$-modules.

**Theorem 6.9** ([MM11, Theorem 5.10]). *Let $X_\Sigma$ be as before, and suppose that $\mathcal{L}_1 = \mathcal{O}_{X_\Sigma}(L_1), \ldots, \mathcal{L}_c = \mathcal{O}_{X_\Sigma}(L_c)$ are ample line bundles on $X_\Sigma$ such that $-K_{X_\Sigma} - \sum_{j=1}^c L_j$ is nef. Put again $\mathcal{E} := \oplus_{j=1}^c \mathcal{L}_j$. For any $\mathcal{L} \in \mathrm{Pic}(X_\Sigma)$ with $c_1(\mathcal{L}) = \sum_{a=1}^r d_a p_a \in \mathbb{L}_A^\vee$, we put $\widehat{\mathcal{L}} = \sum_{a=1}^r z d_a q_a \partial_{q_a} \in R_{\mathbb{C}_z \times \mathcal{KM}}$. Define the left ideal $J$ of $R_{\mathbb{C}_z \times \mathcal{KM}}$ by*
$$J := R_{\mathbb{C}_z \times \mathcal{KM}} (Q_{\underline{l}})_{\underline{l} \in \mathbb{L}_{A'}} + R_{\mathbb{C}_z \times \mathcal{KM}} \cdot \widehat{E},$$

*where*
$$Q_{\underline{l}} := \prod_{i \in \{1,\ldots,m\}: l_i > 0} \prod_{\nu=0}^{l_i-1} \left(\widehat{\mathcal{D}}_i - \nu z\right) \prod_{j \in \{1,\ldots,c\}: l_{m+j} > 0} \prod_{\nu=1}^{l_{m+c}} \left(\widehat{\mathcal{L}}_j + \nu z\right)$$

$$- \underline{q}^{\underline{l}} \cdot \prod_{i \in \{1,\ldots,m\}: l_i < 0} \prod_{\nu=0}^{-l_i-1} \left(\widehat{\mathcal{D}}_i - \nu z\right) \prod_{j \in \{1,\ldots,c\}: l_{m+j} < 0} \prod_{\nu=1}^{-l_{m+c}} \left(\widehat{\mathcal{L}}_j + \nu z\right),$$

$$\widehat{E} := z^2 \partial_z - \widehat{K}_{\mathbb{V}(\mathcal{E}^\vee)}.$$

*Here we write $\mathcal{D}_i \in \mathrm{Pic}(X_\Sigma)$ for a line bundle associated to the torus invariant divisor $D_i$, where $i = 1, \ldots, m$. Notice that the ideal $J$ was called $\mathbb{G}$ in [MM11, Definition 4.3].*

*Moreover, let Quot be the left ideal in $R_{\mathbb{C}_z \times \mathcal{KM}}$ generated by the following set*
$$G := \{P \in R_{\mathbb{C}_z \times \mathcal{KM}} \,|\, \widehat{c}_{\mathrm{top}} \cdot P \in J\},$$

*where $\widehat{c}_{\mathrm{top}} := \prod_{j=1}^c \widehat{\mathcal{L}}_j$. We define $P := R_{\mathbb{C}_z \times \mathcal{KM}}/J$ resp. $P^{\mathrm{res}} := R_{\mathbb{C}_z \times \mathcal{KM}}/\mathrm{Quot}$ and denote by $\mathcal{P} = \mathcal{R}_{\mathbb{C}_z \times \mathcal{KM}}/\mathcal{J}$ resp. $\mathcal{P}^{\mathrm{res}} = \mathcal{R}_{\mathbb{C}_z \times \mathcal{KM}}/\mathcal{Q}\mathrm{uot}$ the corresponding $\mathcal{R}_{\mathbb{C}_z \times \mathcal{KM}}$-modules. Notice that we have $\mathcal{J} \subset \mathcal{Q}\mathrm{uot}$, hence there is a canonical surjection $\mathcal{P} \twoheadrightarrow \mathcal{P}^{\mathrm{res}}$.*

*Put $B_\varepsilon^* := \{q \in (\mathbb{C}^*)^r \,|\, 0 < |q| < \varepsilon\} \subset \mathcal{KM}^\circ$, then there is some $\varepsilon$ such that the following diagram is commutative and the horizontal morphisms are isomorphisms of $\mathcal{R}_{\mathbb{C}_z \times B_\varepsilon^*}$-modules.*

$$\begin{array}{ccc}
\mathcal{P}_{|\mathbb{C}_z \times B_\varepsilon^*} & \xrightarrow{\cong} & (\mathrm{id}_{\mathbb{C}_z} \times \mathrm{Mir})^* \left(\mathrm{QDM}(X_\Sigma, \mathcal{E})\right) \\
\downarrow & & \downarrow \xi \\
\mathcal{P}^{\mathrm{res}}_{|\mathbb{C}_z \times B_\varepsilon^*} & \xrightarrow{\cong} & (\mathrm{id}_{\mathbb{C}_z} \times \mathrm{Mir})^* \left(\overline{\mathrm{QDM}}(X_\Sigma, \mathcal{E})\right)
\end{array}$$

Here Mir: $B_\varepsilon^* \to H^0(X_\Sigma) \times U$ is the **mirror map**, as described in [Giv98b, Theorem 0.1] (see also [CG07, Corollary 5 and the remark thereafter]). Recall that $U \subset H^2(X_\Sigma, \mathbb{C})/2\pi i H^2(X_\Sigma, \mathbb{Z}) \cong (\mathbb{C}^*)^r$ is the convergency domain of the twisted quantum product, i.e., the quantum $\mathcal{D}$-modules $\mathrm{QDM}(X_\Sigma, \mathcal{E})$ and $\overline{\mathrm{QDM}}(X_\Sigma, \mathcal{E})$ are defined on $\mathbb{C}_z \times H^0(X_\Sigma, \mathbb{C}) \times U$ (see subsection 4.1).

We now define another quotient $\mathcal{Q}^{res}$ of $\mathcal{P}$ which is better suited to our approach and which turns out to be isomorphic to $\mathcal{P}^{res}$ resp. to $(id_{\mathbb{C}_z} \times \mathrm{Mir})^* \left(\overline{\mathrm{QDM}}(X_\Sigma, \mathcal{E})\right)$ in some neighborhood of $q = 0$.

**Definition 6.10.** *Let $K$ be the following ideal in $R_{\mathbb{C}_z \times \mathcal{KM}}$:*

$$K := \{P \in R_{\mathbb{C}_z \times \mathcal{KM}} \mid \exists\, p \in \mathbb{Z},\, k \in \mathbb{N} \text{ such that } \prod_{i=0}^{k} \widehat{c}_{top}^{\,p+i} P \in J\},$$

*where $\widehat{c}_{top}^{\,i} := \prod_{j=1}^{c}(\widehat{\mathcal{L}}_j + i)$. Define*

$$Q^{res} := R_{\mathbb{C}_z \times \mathcal{KM}}/K$$

*and denote by $\mathcal{Q}^{res}$ be the corresponding $\mathcal{R}_{\mathbb{C}_z \times \mathcal{KM}}$-module.*

**Proposition 6.11.** *Using the notations from above, we have the following isomorphisms:*

$$\mathcal{P}^{res}_{|\mathbb{C}_z \times B_\varepsilon^*} \simeq \mathcal{Q}^{res}_{|\mathbb{C}_z \times B_\varepsilon^*} \simeq (id_{\mathbb{C}_z} \times \mathrm{Mir})^* \left(\overline{\mathrm{QDM}}(X_\Sigma, \mathcal{E})\right).$$

*Proof.* First notice that we have a surjective morphism $\mathcal{P}^{res} \twoheadrightarrow \mathcal{Q}^{res}$ because the generating set $G$ of $Quot$ is contained in the ideal $K$. If we can construct a well-defined morphism

$$\mathcal{Q}^{res}_{|\mathbb{C}_z \times B_\varepsilon^*} \to (id_{\mathbb{C}_z} \times \mathrm{Mir})^* \left(\overline{\mathrm{QDM}}(X_\Sigma, \mathcal{E})\right) \tag{63}$$

such that the following diagram

$$\begin{array}{ccc}
\mathcal{P}^{res}_{|\mathbb{C}_z \times B_\varepsilon^*} & \twoheadrightarrow & \mathcal{Q}^{res}_{|\mathbb{C}_z \times B_\varepsilon^*} \\
\simeq \downarrow & \swarrow & \\
(id_{\mathbb{C}_z} \times \mathrm{Mir})^* \left(\overline{\mathrm{QDM}}(X_\Sigma, \mathcal{E})\right) & &
\end{array}$$

commutes, the proposition follows. In order to construct the morphism (63) we recapitulate the construction from [MM11] of the morphisms

$$\mathcal{P}_{|\mathbb{C}_z \times B_\varepsilon^*} \to (id_{\mathbb{C}_z} \times \mathrm{Mir})^* (\mathrm{QDM}(X_\Sigma, \mathcal{E})) \quad \text{resp.} \quad \mathcal{P}^{res}_{|\mathbb{C}_z \times B_\varepsilon^*} \to (id_{\mathbb{C}_z} \times \mathrm{Mir})^* \left(\overline{\mathrm{QDM}}(X_\Sigma, \mathcal{E})\right).$$

It relies on a certain multivalued section $L^{tw}$ in $\mathrm{End}(\mathrm{QDM}(X_\Sigma, \mathcal{E}))$ having the property that $L^{tw} z^{-\mu} z^{c_1(\mathcal{T}_X) - c_1(\mathcal{E})}$ is a fundamental solution of $\mathrm{QDM}(X_\Sigma, \mathcal{E})$ (see again [Giv98b] and [CG07]). We use the formulation from [MM11, Proposition 2.17]. Moreoverm we also need the multi-valued section $J^{tw}$ having the property that

$$J^{tw} := (L^{tw})^{-1} 1 \quad \text{in} \quad \overline{\mathrm{QDM}}(X_\Sigma, \mathcal{E}).$$

Finally, we are going to use the cohomological multi-valued section

$$I := q^{T/z} \sum_{d \in H_2(X, \mathbb{Z})} q^d A_d(z),$$

where

$$A_d(z) := \prod_{i=1}^{c} \frac{\prod_{m=-\infty}^{d_{L_i}}([L_i] + mz)}{\prod_{m=-\infty}^{0}([L_i] + mz)} \prod_{\theta \in \Sigma(1)} \frac{\prod_{m=-\infty}^{0}([D_\theta] + mz)}{\prod_{m=-\infty}^{d_\theta}([D_\theta] + mz)},$$

$$q^{T/z} := e^{\frac{1}{z} \sum_{a=1}^{r} T_a \log(q_a)},$$

$d_\theta := \int_d D_\theta$ and $d_{L_i} := \int_d c_1(\mathcal{L}_i)$ and which has asymptotic development $I = F(q)1 + O(z^{-1})$. The aforementioned mirror theorem of Givental ([Giv98b, Theorem 0.1] and [CG07, Corollary 5]), which we use it in the version stated in [MM11, Theorem 5.6], says that

$$I(q, z) = F(q) \cdot J^{tw}(\mathrm{Mir}(q), z).$$

Now one defines the following morphism

$$\mathcal{R}_{\mathbb{C}_z \times B_\varepsilon^*} \longrightarrow (id \times \text{Mir})^* \left( \text{QDM}(X_\Sigma, \mathcal{E}) \right), \tag{64}$$
$$P(z, q, zq\partial_q, z^2\partial_z) \mapsto L^{tw}(\text{Mir}(q), z) z^{-\mu} z^{c_1(\mathcal{T}_X) - c_1(\mathcal{E})} P(q, z, z\partial_{q_i}, z^2\partial_z) z^{-c_1(\mathcal{T}_X) + c_1(\mathcal{E})} z^\mu F(q) J^{tw}(\text{Mir}(q), z)$$
$$= L^{tw}(\text{Mir}(q), z) z^{-\mu} z^{c_1(\mathcal{T}_X) - c_1(\mathcal{E})} P(q, z, z\partial_{q_i}, z^2\partial_z) z^{-c_1(\mathcal{T}_X) + c_1(\mathcal{E})} z^\mu I(q, z)$$

the proof of its surjectivity can be found in the proof [MM11, Theorem 5.10].
The morphism above descends to $\mathcal{P}_{|\mathbb{C}_z \times B_\varepsilon^*}$ by the fact that

$$P(q, z, zq\partial_q, z^2\partial_z) z^{-c_1(\mathcal{T}_X) + c_1(\mathcal{E})} z^\mu I = 0 \quad \text{for} \quad P \in \mathcal{J}.$$

If one composes the morphism (64) with the quotient morphism $\xi$, then this descends to a morphism

$$\mathcal{P}^{res}_{|\mathbb{C}_z \times B_\varepsilon^*} \to (id_{\mathbb{C}_z} \times \text{Mir})^* \left( \overline{\text{QDM}}(X_\Sigma, \mathcal{E}) \right), \tag{65}$$

which follows from

$$P(q, z, zq\partial_q, z^2\partial_z) z^{-c_1(\mathcal{T}_X) + c_1(\mathcal{E})} z^\mu I \in ker(m_{c_{top}}) \quad \text{for} \quad P \in \mathcal{Q}uot \tag{66}$$

and the fact that $L^{tw}$ preserves $\ker(m_{c_{top}})$ (cf. [MM11, Lemma 2.31]).

As explained above, the proposition will follow if the morphism (65) descends to $\mathcal{Q}^{res}_{|\mathbb{C}_z \times B_\varepsilon^*}$, i.e. we have to show that

$$P(q, z, zq\partial_q, z^2\partial_z) z^{-c_1(\mathcal{T}_X) + c_1(\mathcal{E})} z^\mu I \in ker(m_{c_{top}}) \quad \text{for} \quad P \in \mathcal{K}. \tag{67}$$

We will adapt the proof of (66) from [MM11, Lemma 5.21] to our situation. First notice that

$$z^{-c_1(\mathcal{T}_X) + c_1(\mathcal{E})} z^\mu I = \sum_{d \in H_2(X, \mathbb{Z})} q^{T+d} z^{-c_1(\mathcal{T}_X) + c_1(\mathcal{E}) - \int_d (c_1(\mathcal{T}_X) - c_1(\mathcal{E}))} A_d(1).$$

Now let $P(q, z, q\partial_q, z^2\partial_z) \in K$ and decompose it:

$$P(q, z, q\partial_q, z^2\partial_z) = \sum_{\substack{d' \in H_2(X, \mathbb{Z}) \\ \text{finite}}} q^{d'} P_{d'}(z, z\partial_q, z\partial_z).$$

This gives

$$P(q, z, q\partial_q, z^2\partial_z) z^{-c_1(\mathcal{T}_X) + c_1(\mathcal{E})} z^\mu I = \sum_{d \in H_2(X, \mathbb{Z})} q^{T+d} z^{-c_1(\mathcal{T}_X) + c_1(\mathcal{E}) - \int_d (c_1(\mathcal{T}_X) - c_1(\mathcal{E}))} B_d(z),$$

where

$$B_d(z) := \sum_{\substack{d' \in H_2(X, \mathbb{Z}) \\ \text{finite}}} P_{d'} \left( z, z(T+d), z(-c_1(\mathcal{T}_X) + c_1(\mathcal{E}) - \int_d (c_1(\mathcal{T}_X) - c_1(\mathcal{E}))) \right) A_{d-d'}(1).$$

Similarly to loc. cit., the statement (67) will follow from the fact that $c_{top} B_d(z) = 0$ for all $d \in H_2(X, \mathbb{Z})$.
Because $P \in K$, there exists $p \in \mathbb{Z}$ and $k \in \mathbb{N}$ such that

$$\left( \prod_{i=0}^{k} \widehat{c}^{\,p+i}_{top} \right) P(q, z, zq\partial_q, z^2\partial_z) z^{-c_1(\mathcal{T}_X \otimes \mathcal{E}^\vee)} z^\mu I = 0,$$

which gives

$$\sum_{d \in H_2(X, \mathbb{Z})} q^{T+d} z^{-c_1(\mathcal{T}_X) + c_1(\mathcal{E}) - d_{\mathcal{T}_X \otimes \mathcal{E}^\vee}} \left( \prod_{i=0}^{k} \prod_{j=1}^{c} z([L_j] + d_{L_j} + p + i) \right) B_d(z) = 0.$$

Notice that the sum above is zero if and only if each summand is zero. For $(z, q) \in \mathbb{C}_z^* \times W$ the term $q^{T+d} z^{-c_1(\mathcal{T}_X) + c_1(\mathcal{E}) - d_{\mathcal{T}_X \otimes \mathcal{E}^\vee}}$ is invertible, so we deduce that

$$\left( \prod_{i=0}^{k} \prod_{j=1}^{c} ([L_j] + d_{L_j} + p + i) \right) B_d(z) = 0 \quad \forall d \in H_2(X, \mathbb{Z}).$$

Let $J_d := \{ j \in \{1, \ldots, c\} \mid \exists i \in \{0, \ldots, k\} \text{ with } d_{L_j} + p + j = 0 \}$ and notice that for every $j$ there is at most one $i \in \{0, \ldots, k\}$ such that $d_{L_j} + p + i = 0$. Because cup-product with $[L_j] + l$ is an automorphism of $H^{2*}(X, \mathbb{C})$ for every $l \neq 0$, we conclude that

$$\left( \prod_{j \in J_d} [L_j] \right) B_d(z) = 0 \quad \forall d \in H_2(X, \mathbb{Z}),$$

which in turn shows that $c_{top} B_d(z) = (\prod_{j=1}^{c} [L_j]) B_d(z) = 0$ for all $d \in H_2(X, \mathbb{Z})$. □

The next proposition compares the $\mathcal{R}$-modules from Theorem 6.9 and Definition 6.10 with ${}_0\mathcal{QM}_{A'}$ and ${}_0\mathcal{QM}_{A'}^{IC}$.

**Proposition 6.12.** *We have isomorphisms of $\mathcal{R}_{\mathbb{C}_z \times \mathcal{KM}^\circ}$-modules*

$$\mathcal{P}_{|\mathbb{C}_z \times \mathcal{KM}^\circ} \cong {}_0\mathcal{QM}_{A'} \qquad \text{and} \qquad \mathcal{Q}^{\mathrm{res}}_{|\mathbb{C}_z \times \mathcal{KM}^\circ} \cong {}_0\mathcal{QM}_{A'}^{IC}.$$

*Proof.* The first isomorphism follows from a similar argument as [RS15, Proposition 3.2], namely, the section

$$\overline{\varrho} = i \circ \varrho' : \mathcal{KM} \hookrightarrow W^*,$$
$$(q_1, \ldots, q_r) \mapsto (\lambda_1 = \underline{q}^{\underline{m}_1}, \ldots, \lambda_m = \underline{q}^{\underline{m}_m}, \lambda_{m+1} = -\underline{q}^{\underline{m}_{m+1}}, \ldots, \lambda_{m+c} = -\underline{q}^{\underline{m}_{m+c}})$$

can be used to construct an isomorphism

$$\theta : F \times \mathcal{KM} \longrightarrow W^*,$$
$$(f_1, \ldots, f_{n+c}, q_1, \ldots, q_r) \mapsto (\underline{q}^{\underline{m}_1} \underline{y}^{\underline{a}_1}, \ldots, \underline{q}^{\underline{m}_m} \underline{y}^{\underline{a}_m}, -\underline{q}^{\underline{m}_{m+1}} \underline{y}^{\underline{a}_{m+1}}, \ldots, -\underline{q}^{\underline{m}_{m+c}} \underline{y}^{\underline{a}_{m+c}})$$

with inverse

$$\theta^{-1} : W^* \longrightarrow F \times \mathcal{KM},$$
$$(\lambda_1, \ldots, \lambda_{m+c}) \mapsto \left( f_j = (-1)^{\sum_{i=m+1}^{m+c} c_{ij}} \underline{\lambda}^{\underline{c}_j}, q_a = (-1)^{\sum_{i=m+1}^{m+c} l_{ia}} \underline{\lambda}^{\underline{l}_a} \right),$$

where $L = (\underline{l}_a)$ resp. $M = (\underline{m}_i)$ are the matrices which were introduced above Definition 6.3 and $C = (\underline{c}_j)$ is a $(m+c) \times (n+c)$-matrix such that the following equations are fulfilled (cf. Section 2.4):

$$M \cdot L = I_r, \quad B \cdot C = I_{n+c}, \quad B \cdot L = 0, \quad M \cdot C = 0, \quad C \cdot B + L \cdot M = I_{m+c}.$$

Under this coordinate change the module ${}_0^* \widehat{\mathcal{N}}_{A'}^{(0,0,0)}$ has the following presentation:

$$\mathcal{R}_{\mathbb{C}_z \times F \times \mathcal{KM}} / ((Q_{\underline{l}})_{\underline{l} \in \mathbb{L}} + (\widehat{E}) + (\widehat{E}'_k)_{k=1,\ldots,n+c})$$

with $Q_{\underline{l}}$ and $\widehat{E}$ as in Definition 6.9 and $\widehat{E}'_k := f_k \partial_k$ for $k \in \{1, \ldots, n+c\}$.

Its module of global sections can be described simply by forgetting $\partial_{f_k}$, i.e. we have the following description

$$\frac{\mathbb{C}[z, f_1^\pm, \ldots, f_{n+c}^\pm, q_1^\pm, \ldots, q_r^\pm] \langle z^2 \partial_z, z \partial_{q_1}, \ldots, z \partial_{q_r} \rangle}{((Q_{\underline{l}})_{\underline{l} \in \mathbb{L}} + (\widehat{E}))}. \tag{68}$$

Notice that the map $\overline{\varrho}$ can be factorized as $\theta \circ i_\theta$ with

$$i_\theta : \mathcal{KM} \longrightarrow F \times \mathcal{KM},$$
$$(q_1, \ldots, q_r) \mapsto (1, \ldots, 1, q_1, \ldots, q_r)$$

Thus the inverse image of (68) with respect to $i_\theta$ is given by

$$\frac{\mathbb{C}[z, q_1^\pm, \ldots, q_r^\pm]\langle z^2\partial_z, z\partial_{q_1}, \ldots, z\partial_{q_r}\rangle}{((Q_{\underline{l}})_{\underline{l} \in \mathbb{L}} + (\widehat{E}))},$$

which is exactly the definition of the module $\mathcal{P}$ from Theorem 6.9.

Concerning the second isomorphism, the associated sub-$R_{\mathbb{C}_z \times F \times \mathcal{KM}}$-module corresponding to $\widehat{\mathcal{K}}_\mathcal{N}$ from Lemma 6.2 can be described by

$$\{P \in \mathbb{C}[z, f_1^\pm, \ldots, f_{n+c}^\pm, q_1^\pm, \ldots, q_r^\pm]\langle z^2\partial_z, z\partial_{q_1}, \ldots, z\partial_{q_r}\rangle \mid \exists p \in \mathbb{Z}, k \in \mathbb{N} \text{ s.t. } \prod_{i=0}^k \widehat{C}_{top}^{p+i} P \in ((Q_{\underline{l}})_{\underline{l} \in \mathbb{L}} + (\widehat{E}))\},$$

where

$$\widehat{C}_{top}^k := \prod_{i=m+1}^{m+c} ((\sum_{j=1}^{n+c} c_{ij} f_j \partial_j + \sum_{a=1}^r l_{ia} q_a \partial_a) + l),$$
$$= \prod_{i=m+1}^{m+c} ((\sum_{j=1}^{n+c} c_{ij} f_j \partial_j + \widehat{\mathcal{D}}_i) + k)$$

for $k \in \mathbb{Z}$. It is easy to see that its inverse image under $(id_{\mathbb{C}_z} \times i_\theta)$ is given by

$$\{P \in \mathbb{C}[z, q_1^\pm, \ldots, q_r^\pm]\langle z^2\partial_z, z\partial_{q_1}, \ldots, z\partial_{q_r}\rangle \mid \exists p \in \mathbb{Z}, k \in \mathbb{N} \text{ s.t. } \prod_{i=0}^k \widehat{c}_{top}^{p+i} P \in ((Q_{\underline{l}})_{\underline{l} \in \mathbb{L}} + (\widehat{E}))\},$$

which is exactly the definition of the ideal $K$ in Definition 6.10. Thus, the second isomorphism follows. □

Combining Proposition 6.12, Theorem 6.9, Lemma 6.4 and 6.6 as well as Proposition 6.7, we obtain the following mirror statement.

**Theorem 6.13.** *Let $X_\Sigma$ and $\mathcal{L}_1, \ldots, \mathcal{L}_c$ be as in Theorem 6.9. Consider the affine resp. non-affine Landau-Ginzburg models $\overline{\pi} = (\widetilde{F}, \underline{q}) : X^{aff} \times \mathcal{KM}^\circ \to \mathbb{C}_{\lambda_0} \times \mathcal{KM}^\circ$, $\pi : S \times \mathcal{KM}^\circ \to \mathbb{C}_{\lambda_0} \times \mathcal{KM}^\circ$ and $\Pi : \mathcal{Z}_X^\circ \hookrightarrow \mathcal{Z}^\circ \xrightarrow{\alpha} \mathbb{C}_{\lambda_0} \times \mathcal{KM}^\circ$ associated to $(X_\Sigma, \mathcal{L}_1, \ldots, \mathcal{L}_c)$. Let $B_\varepsilon^* \subset \mathcal{KM}^\circ$ be the punctured ball from Theorem 6.9. Then there are isomorphisms of $\mathcal{D}_{\mathbb{C}_z \times B_\varepsilon^*}$-modules*

$$\mathrm{FL}^{loc}_{\mathcal{KM}^\circ}\left(\mathcal{H}^0 \pi_\dagger \mathcal{O}_{S \times \mathcal{KM}^\circ}\right)_{|\mathbb{C}_z \times B_\varepsilon^*} \cong (\mathrm{id}_{\mathbb{C}_z} \times \mathrm{Mir})^* \left(\mathrm{QDM}(X_\Sigma, \mathcal{E})\right)(*(\{0\} \times B_\varepsilon^*)),$$

$$\mathrm{FL}^{loc}_{\mathcal{KM}^\circ}\left(\mathcal{H}^0 \alpha_+ \mathcal{M}^{IC}(\mathcal{Z}_X^\circ)\right)_{|\mathbb{C}_z \times B_\varepsilon^*} \cong (\mathrm{id}_{\mathbb{C}_z} \times \mathrm{Mir})^* \left(\overline{\mathrm{QDM}}(X_\Sigma, \mathcal{E})\right)(*(\{0\} \times B_\varepsilon^*))$$

*and an isomorphism of $\mathcal{R}_{\mathbb{C}_z \times B_\varepsilon^*}$-modules*

$$\sigma^* z^n \cdot \left(H^{n+c}(\Omega^\bullet_{X^{aff} \times \mathcal{KM}^\circ / \mathcal{KM}^\circ}(\log D)[z], zd - d\widetilde{F})\right)'_{|\mathbb{C}_z \times B_\varepsilon^*} \cong (\mathrm{id}_{\mathbb{C}_z} \times \mathrm{Mir})^* \left(\mathrm{QDM}(X_\Sigma, \mathcal{E})\right).$$

The following corollary is the promised Hodge theoretic application of the above main theorem.

**Corollary 6.14.** *There exists a variation of non-commutative pure polarized Hodge structures $(\mathcal{F}, \mathcal{L}_\mathbb{Q}, \mathrm{iso}, P)$ on $\mathcal{KM}^\circ$ (see [KKP08], [HS10] or [Sab11] for the definition) such that*

$$\mathcal{F}(*(\{0\} \times B_\varepsilon^*)) \cong (\mathrm{id}_{\mathbb{C}_z} \times \mathrm{Mir})^* \left(\overline{\mathrm{QDM}}(X_\Sigma, \mathcal{E})\right)(*(\{0\} \times B_\varepsilon^*)). \tag{69}$$

*Proof.* Using Theorem 6.13, this is a direct consequence of [Sai88, Théorème 1] and [Sab08, Corollary 3.15]. □

It would of course be desirable to remove the localization with respect to $\{0\} \times B_\varepsilon^*$ from the above theorem. We conjecture that the corresponding statement still holds, however, we cannot give a complete proof of this for the moment as we are not able to control the Hodge filtration on $\mathcal{M}^{IC}(\mathcal{Z}_X^\circ)$. More precisely, we expect the following to be true.

**Conjecture 6.15.** *1. Write $F_\bullet^H \mathcal{H}^0 \alpha_+ \mathcal{M}^{IC}(\mathcal{Z}_X^\circ)$ for the Hodge filtration on $\mathcal{H}^0 \alpha_+ \mathcal{M}^{IC}(\mathcal{Z}_X^\circ)$, which underlies a pure Hodge module due to[Sai88, Théorème 1], and which has weight $n + c + (m - n) = m + c$. Let $F_\bullet^H[\partial_{\lambda_0}^{-1}]$ be the saturation of $F_\bullet^H$ as in the proof of Lemma 6.6 and write $G_\bullet^H$ for the induced filtration on $\mathrm{FL}_{\mathcal{KM}^\circ}(\mathcal{H}^0\alpha_+\mathcal{M}^{IC}(\mathcal{Z}_X^\circ))$. Then under the isomorphism of Proposition 6.7, 2., we have that*
$$G_{\bullet-(m+c)}^H \mathrm{FL}_{\mathcal{KM}^\circ}(\mathcal{H}^0\alpha_+\mathcal{M}^{IC}(\mathcal{Z}_X^\circ)) \cong z^\bullet \cdot {}_0\mathcal{QM}_{A'}^{IC}.$$

*Notice that the bundle $\mathcal{F}$ which was used in the isomorphism from corollary 6.14 is nothing but the object $G_{-(m+c)}^H \mathrm{FL}_{\mathcal{KM}^\circ}(\mathcal{H}^0\alpha_+\mathcal{M}^{IC}(\mathcal{Z}_X^\circ))$.*

*2. The isomorphism (69) holds without localization, i.e., there is an isomorphism of $\mathcal{R}_{\mathbb{C}_z \times B_\varepsilon^*}$-modules*
$$\left(G_{-(m+c)}^H \mathrm{FL}_{\mathcal{KM}^\circ}(\mathcal{H}^0\alpha_+\mathcal{M}^{IC}(\mathcal{Z}_X^\circ))\right)_{|\mathbb{C}_z \times B_\varepsilon^*} \cong (\mathrm{id}_{\mathbb{C}_z} \times \mathrm{Mir})^* \overline{\mathrm{QDM}}(X_\Sigma, \mathcal{E}).$$

*As a consequence, the reduced quantum $\mathcal{D}$-module underlies a variation of non-commutative Hodge structures.*

This conjecture, if proved, should be seen as a first step towards establishing the existence of a very special geometric structure on the cohomology space of the complete intersection subvariety $Y \subset X_\Sigma$, known as $tt^*$-*geometry* (see [CV91, CV93] or [Her03] for a modern account). Its existence is known for the quantum cohomology of nef toric manifolds themselve (this follows from [RS15, Theorem 5.3], see also [Iri09b]). For (non-toric) complete intersections one needs of course to consider its total quantum cohomology, not just the ambient part, but at least on this part the above conjecture would give the desired result.

Comparing Theorem 6.13 with Lemma 6.4 one may wonder whether the module $\mathrm{FL}_{\mathcal{KM}^\circ}^{loc}\left(\mathcal{H}^0\pi_+\mathcal{O}_{S \times \mathcal{KM}^\circ}\right)$ also has an interpretation as a mirror object. This is actually the case, namely, it corresponds to the so-called $Euler^{-1}$-twisted quantum $\mathcal{D}$-module (whereas the object $\mathrm{QDM}(X_\Sigma, \mathcal{E})$ from Definition 4.3 would be the $Euler$-twisted quantum $\mathcal{D}$-module in this terminology). The $Euler^{-1}$-twisted quantum $\mathcal{D}$-module encodes the so-called **local Gromov-Witten invariants** of the dual bundle $\mathcal{E}^\vee$ and is denoted by $\mathrm{QDM}(\mathcal{E}^\vee)$ (see [Giv98a, Theorem 4.2]). There is a non-degenerate pairing between $\mathrm{QDM}(X_\Sigma, \mathcal{E})$ and $\left(\mathrm{id}_{\mathbb{C}_z} \times (h \circ \overline{f})\right)^* \mathrm{QDM}(\mathcal{E}^\vee)$ (this is the non-equivariant limit of the quantum Serre theorem from [CG07, Corollary 2]) where $\overline{f}, h \in \mathbb{C}[[H^*(X_\Sigma, \mathbb{C})^\vee]]^n$ are maps. The existence of this pairing has been proved in the recent paper [IMM16]. However, in the formulation of this result, all objects are defined on the total cohomology space, i.e., correspond to the big (twisted) quantum product. Nevertheless, we are able to obtain a mirror theorem for local Gromov-Witten invariants.

Consider the situation of Theorem 6.13, in particular, let $\mathcal{E} := \oplus_{j=1}^c \mathcal{L}_j$. As $\mathcal{L}_j$ are nef bundles and hence globally generated, also $\mathcal{E}$ is globally generated and therefore *convex*. Let $\mathrm{QDM}(\mathcal{E}^\vee)$ be the $(Euler^{-1})$-twisted quantum $\mathcal{D}$-module governing local Gromov-Witten invariants, that is, integrals over the moduli space $\overline{\mathcal{M}}_{0,l,d}(\mathcal{V}(\mathcal{E}^\vee))$ of stable maps to the total space $\mathcal{V}(\mathcal{E}^\vee)$ (notice that $\overline{\mathcal{M}}_{0,l,d}(\mathcal{V}(\mathcal{E}^\vee))$ is compact unless $d = 0$).

**Theorem 6.16.** *Let again $X_\Sigma$ and $\mathcal{L}_1, \ldots, \mathcal{L}_c$ be as in Theorem 6.9. There is some convergency neighborhood $B_{\varepsilon'}^*$, an isomorphism of $\mathcal{D}_{\mathbb{C}_z \times B_{\varepsilon'}^*}$-modules*
$$\mathrm{FL}_{\mathcal{KM}^\circ}^{loc}\left(\mathcal{H}^0\pi_+\mathcal{O}_{S \times \mathcal{KM}^\circ}\right)_{|\mathbb{C}_z \times B_{\varepsilon'}^*} \cong \left(\mathrm{id}_{\mathbb{C}_z} \times \mathrm{Mir}'\right)^* \left(\mathrm{QDM}(\mathcal{E}^\vee)\right)(*(\{0\} \times B_{\varepsilon'}^*))$$

*and an isomorphism of $\mathcal{R}_{\mathbb{C}_z \times B_{\varepsilon'}^*}$-modules*
$$H^{n+c}(\Omega^\bullet_{X^{\mathrm{aff}} \times \mathcal{KM}^\circ / \mathcal{KM}^\circ}(\log D)[z], zd - d\widetilde{F})_{|\mathbb{C}_z \times B_{\varepsilon'}^*} \cong \sigma^* z^n \cdot \left(\mathrm{id}_{\mathbb{C}_z} \times \mathrm{Mir}'\right)^* \mathrm{QDM}(\mathcal{E}^\vee).$$

*Here* $\mathrm{Mir}'$ *is some base change involving the above mentioned maps $\overline{f}, h$ as well as the base change* $\mathrm{Mir}$.

*Proof.* It is actually sufficient to show the second statement as the first follows by applying the localization functor $(-) \otimes \mathcal{O}_{\mathbb{C}_z \times B_{\varepsilon'}^*}(*(\{0\} \times B_{\varepsilon'}^*))$ (This follows by using Propostion 3.21, Remark 3.16, 2. as well as Proposition 3.3 together with Lemma 3.4).

It follows from [IMM16, Theorem 3.14] that there exists a non-degenerate pairing
$$(\mathrm{QDM}(X_\Sigma, \mathcal{E}))_{|\mathbb{C}_z \times B_{\varepsilon'}^*} \otimes \left(\mathrm{id}_{\mathbb{C}_z} \times (h \circ \overline{f})\right)^*\left(\mathrm{QDM}(\mathcal{E}^\vee)_{|\mathbb{C}_z \times B_{\varepsilon'}^*}\right) \longrightarrow \mathcal{O}_{\mathbb{C}_z \times B_{\varepsilon'}^*},$$

which is compatible with the connection operators induced by the $\mathcal{R}_{\mathbb{C}_z \times B^*_{\varepsilon'}}$-module structures of the objects on the left hand side. As has been pointed out above, this statement is given in loc.cit. for the big quantum $\mathcal{D}$-modules, hence, one has to check that $\left(\mathrm{id}_{\mathbb{C}_z} \times (h \circ \overline{f})\right)^* \left(\mathrm{QDM}(\mathcal{E}^\vee)_{|\mathbb{C}_z \times B^*_{\varepsilon'}}\right)$ is still a vector bundle on $\mathbb{C}_z \times B^*_{\varepsilon'}$. From the definition of the map $h$ (see [IMM16, Proposition 3.11]) it is clear that it restricts to an invertible map $h : H^2(X_\Sigma, \mathbb{C}) \to H^2(X_\Sigma, \mathbb{C})$. We claim that $\overline{f}$ restricts to a map

$$\overline{f} : H^2(X_\Sigma, \mathbb{C}) \longrightarrow H^0(X_\Sigma, \mathbb{C}) \oplus H^2(X_\Sigma, \mathbb{C})$$

so that the pullback $\overline{f}^*\gamma$ of any class $\gamma \in H^2(X_\Sigma, \mathbb{C})$ is still an element of $H^2(X_\Sigma, \mathbb{C})$. This can be seen as follows: From [IMM16, Proof of Lemma 3.2], we know that

$$\overline{f}(\tau) = \sum_{\alpha=0}^{h} \left( \sum_{d \in H_2(X_\Sigma, \mathbb{C}), n \geq 0} \langle T_\alpha, \widetilde{1}, \tau, \underbrace{\tau, \ldots, \tau}_{n-\text{times,}} \rangle_{0, n+3, d} \right) T^\alpha \tag{70}$$

According to the definition 4.1, the correlator $\langle T_\alpha, \widetilde{1}, \tau \underbrace{\tau, \ldots, \tau}_{n-\text{times,}}\rangle_{0,n+3,d}$ is non-zero only if the degree $\deg(T_\alpha) + (n+1) + \deg(e(\mathcal{E}_{0,n+3,d}(2)))$ equals the dimension of the moduli space $[\overline{\mathcal{M}}_{0,n+3,d}(\mathcal{X})]$, i.e., the number $\dim(X_\Sigma) + \int_d c_1(X) + n$. Under the assumption of the theorem, $\mathcal{E}_{0,n+3,d}$ is represented by a vector bundle, which is of rank $\int_d c_1(\mathcal{E}) + \mathrm{rank}(\mathcal{E})$. Hence $\mathcal{E}_{0,n+3,d}(2)$, being the kernel of the map $\mathcal{E}_{0,n+3,d} \to \mathrm{ev}_2^*(\mathcal{E})$ is a bundle of rank $\int_d c_1(\mathcal{E})$, so that we see that $\langle T_\alpha, \widetilde{1}, \tau \underbrace{\tau, \ldots, \tau}_{n-\text{times,}}\rangle_{0,n+3,d} \neq 0$ iff

$$\deg(T_\alpha) + 1 = \dim(X) + \int_d c_1(X_\Sigma) - \int_d c_1(\mathcal{E}) = \dim(X) + \int_d c_1\left(-K_{X_\Sigma} - \sum_{j=1}^{c} \mathcal{L}_j\right) \geq \dim(X)$$

where the last inequality holds due to the assumptions on $X_\Sigma$ and $\mathcal{E}$. We conclude for any class $T^\alpha$ occurring in formula (70) the following holds: either its degree is at most 1 or its coefficient is zero. This means nothing else than $\mathrm{im}(\overline{f}) \subset H^2(X_\Sigma, \mathbb{C}) \oplus H^0(X_\Sigma, \mathbb{C})$.

Hence we can deduce from [IMM16, Theorem 3.14] that there is an isomorphism

$$\left(\mathrm{id}_{\mathbb{C}_z} \times (h \circ \overline{f})\right)^* \left(\mathrm{QDM}(\mathcal{E}^\vee)_{|\mathbb{C}_z \times B^*_{\varepsilon'}}\right) \cong \left((\mathrm{QDM}(X_\Sigma, \mathcal{E}))'\right)_{|\mathbb{C}_z \times B^*_{\varepsilon'}}$$

of $\mathcal{R}_{\mathbb{C}_z \times B^*_{\varepsilon'}}$-modules, and then the desired statement follows from the third line in the displayed formula in Theorem 6.13. $\square$

**Remark:** In view of [Giv98a, corollary 4.3], one may conjecture that $\mathrm{Mir}'$ is the identity if the number $c$ of line bundles defining the bundle $\mathcal{E}$ is strictly bigger than 1. However, at this moment, we do not have any further evidence for this conjecture.

The following consideration shows that the main Theorem 6.13 can also be considered as a generalization of mirror symmetry for Fano manifolds themselves, as presented in our previous paper (see [RS15, Proposition 4.10]). Namely, let us consider the case where the number $c$ of line bundles on the toric variety $X_\Sigma$ is zero. Then we have $A' = A$, and the duality morphism $\phi$ from Definition 5.6 is

$$\phi : \mathcal{M}_{A''}^{-(c+1,\underline{0},\underline{1})} = \mathcal{M}_{A''}^{(-1,\underline{0})} \longrightarrow \mathcal{M}_{A''}^{(0,\underline{0},\underline{0})} = \mathcal{M}_{A''}^{(0,\underline{0})}$$

and is induced by right multiplication by $\partial_{\lambda_0}$. In particular, the induced morphism $\widehat{\phi}$ is simply the identity on $\widehat{\mathcal{M}}_{A'}^{(0,0)}$. In particular, we have that $\mathrm{im}(\widetilde{\phi}) \cong \widehat{\mathcal{M}}_{A'}^{(0,0)}$ so that $\mathcal{QM}_{A'}^{IC} \cong \mathcal{QM}_{A'}$ and $_0\mathcal{QM}_{A'}^{IC} \cong {_0\mathcal{QM}}_{A'}$. On the other hand, the reduced quantum $\mathcal{D}$-module $\overline{\mathrm{QDM}}(X_\Sigma, \mathcal{E})$ is nothing but the quantum $\mathcal{D}$-module of the variety $X_\Sigma$, so that we deduce from Theorem 6.13 that we have an isomorphism of $\mathcal{D}_{\mathbb{C}_z \times B^*_\varepsilon}$-modules

$$\mathrm{FL}_{\mathcal{KM}^\circ}^{loc} \left(\mathcal{H}^0 \pi_+ \mathcal{O}_{S \times \mathcal{KM}^\circ}\right)_{|\mathbb{C}_z \times B^*_\varepsilon} \cong (\mathrm{id}_{\mathbb{C}_z} \times \mathrm{Mir})^* \left(\mathrm{QDM}(X_\Sigma)\right)(*(\{0\} \times B^*_\varepsilon)).$$

One easily sees that we have an even more precise statement, namely, the third assertion of Theorem 6.13 simplifies in this case to an isomorphism of $\mathcal{R}_{\mathbb{C}_z \times B_\varepsilon^*}$-modules

$$H^n(\Omega^\bullet_{S \times \mathcal{KM}^\circ/\mathcal{KM}^\circ}[z], zd - d\widetilde{F})_{|\mathbb{C}_z \times B_\varepsilon^*} \cong (\mathrm{id}_{\mathbb{C}_z} \times \mathrm{Mir})^* \mathrm{QDM}(X_\Sigma, \mathcal{E}).$$

This isomorphism is the restriction of the isomorphism in [RS15, Proposition 4.10] to $\mathbb{C}_z \times B_\varepsilon$ (see also [Iri09a, Proposition 4.8]), notice that the neighborhood $B_\varepsilon$ is called $W_0$ in [RS15]. Hence we see that our main Theorem 6.13 contains in particular the mirror correspondence for smooth toric nef manifolds, at least on the level of $\mathcal{R}_{\mathbb{C}_z \times B_\varepsilon}$-modules.

One may conclude from the above observation that Landau-Ginzburg models, either affine or compactified, appear to be the right point of view to study various type of mirror models of (the quantum cohomology of) smooth projective manifolds, including Calabi-Yau, Fano and more generally nef ones. The preprint [GKR12] where varieties of general types and their mirrors are investigated, also seem to confirm this observation. It would certainly be fruitful to apply our methods to varieties with positive Kodaira dimension to refine the results from loc.cit.

# Index of Notation

## Objects

# Maps and Spaces

# References

[Ado94]  Alan Adolphson, *Hypergeometric functions and rings generated by monomials*, Duke Math. J. **73** (1994), no. 2, 269–290.

[Bat94]  Victor V. Batyrev, *Dual polyhedra and mirror symmetry for Calabi-Yau hypersurfaces in toric varieties*, J. Algebraic Geom. **3** (1994), no. 3, 493–535.

[BB96a]  Victor V. Batyrev and Lev A. Borisov, *Mirror duality and string-theoretic Hodge numbers*, Invent. Math. **126** (1996), no. 1, 183–203.

[BB96b]  \_\_\_\_\_\_, *On Calabi-Yau complete intersections in toric varieties*, Higher-dimensional complex varieties (Trento, 1994), de Gruyter, Berlin, 1996, pp. 39–65.

[BH93]  Winfried Bruns and Jürgen Herzog, *Cohen-Macaulay rings*, Cambridge Studies in Advanced Mathematics, vol. 39, Cambridge University Press, Cambridge, 1993.

[Bry86]  Jean-Luc Brylinski, *Transformations canoniques, dualité projective, théorie de Lefschetz, transformations de Fourier et sommes trigonométriques*, Astérisque (1986), no. 140-141, 3–134, 251, Géométrie et analyse microlocales.

[CG07]  Tom Coates and Alexander Givental, *Quantum Riemann-Roch, Lefschetz and Serre*, Ann. of Math. (2) **165** (2007), no. 1, 15–53.

[CK99]  David A. Cox and Sheldon Katz, *Mirror symmetry and algebraic geometry*, Mathematical Surveys and Monographs, vol. 68, American Mathematical Society, Providence, RI, 1999.

[CLS11]  David A. Cox, John B. Little, and Henry K. Schenck, *Toric varieties*, Graduate Studies in Mathematics, vol. 124, American Mathematical Society, Providence, RI, 2011.

[CV91]  Sergio Cecotti and Cumrun Vafa, *Topological–anti-topological fusion*, Nuclear Phys. B **367** (1991), no. 2, 359–461.

[CV93]  \_\_\_\_\_\_, *On classification of N = 2 supersymmetric theories*, Comm. Math. Phys. **158** (1993), no. 3, 569–644.

[DE03]  Andrea D'Agnolo and Michael Eastwood, *Radon and Fourier transforms for $\mathcal{D}$-modules*, Adv. Math. **180** (2003), no. 2, 452–485.

[Dim92]  Alexandru Dimca, *Singularities and topology of hypersurfaces*, Universitext, Springer-Verlag, New York, 1992.

[Dim04]  \_\_\_\_\_\_, *Sheaves in topology*, Universitext, Springer-Verlag, Berlin, 2004.

[DS03]  Antoine Douai and Claude Sabbah, *Gauss-Manin systems, Brieskorn lattices and Frobenius structures. I*, Ann. Inst. Fourier (Grenoble) **53** (2003), no. 4, 1055–1116.

[Ful93]  William Fulton, *Introduction to toric varieties*, Annals of Mathematics Studies, vol. 131, Princeton University Press, Princeton, NJ, 1993, The William H. Roever Lectures in Geometry.

[Giv98a]  Alexander Givental, *Elliptic Gromov-Witten invariants and the generalized mirror conjecture*, Integrable systems and algebraic geometry (Kobe/Kyoto, 1997) (M.-H. Saito, Y. Shimizu, and K. Ueno, eds.), World Sci. Publ., River Edge, NJ, River Edge, NJ, 1998, pp. 107–155.

[Giv98b]  \_\_\_\_\_\_, *A mirror theorem for toric complete intersections*, Topological field theory, primitive forms and related topics (Kyoto, 1996), Progr. Math., vol. 160, Birkhäuser Boston, Boston, MA, 1998, pp. 141–175.

[GKR12]  Mark Gross, Ludmil Katzarkov, and Helge Ruddat, *Towards mirror symmetry for varieties of general type*, Preprint arxiv:1202.4042, 2012.

[GKZ89]  I. M. Gel′fand, M. M. Kapranov, and A. V. Zelevinsky, *Hypergeometric functions and toric varieties*, Funktsional. Anal. i Prilozhen. **23** (1989), no. 2, 12–26. MR 1011353 (90m:22025)

[GKZ90]  \_\_\_\_\_\_, *Generalized Euler integrals and A-hypergeometric functions*, Adv. Math. **84** (1990), no. 2, 255–271.

[GKZ08]  \_\_\_\_\_\_, *Discriminants, resultants and multidimensional determinants*, Modern Birkhäuser Classics, Birkhäuser Boston Inc., Boston, MA, 2008, Reprint of the 1994 edition.

[GM83]  Mark Goresky and Robert MacPherson, *Intersection homology. II*, Invent. Math. **72** (1983), no. 1, 77–129.

[Gro11]  Mark Gross, *Tropical geometry and mirror symmetry*, CBMS Regional Conference Series in Mathematics, vol. 114, Published for the Conference Board of the Mathematical Sciences, Washington, DC, 2011.

[Her03]  Claus Hertling, *$tt^*$ geometry, Frobenius manifolds, their connections, and the construction for singularities*, J. Reine Angew. Math. **555** (2003), 77–161.

[Hoc72]  M. Hochster, *Rings of invariants of tori, Cohen-Macaulay rings generated by monomials, and polytopes*, Ann. of Math. (2) **96** (1972), 318–337.

[Hot98]  Ryoshi Hotta, *Equivariant $\mathcal{D}$-modules*, Preprint math.RT/9805021, 1998.

[HS10]  Claus Hertling and Christian Sevenheck, *Limits of families of Brieskorn lattices and compactified classifying spaces*, Adv. Math. **223** (2010), no. 4, 1155–1224.

[HTT08]  Ryoshi Hotta, Kiyoshi Takeuchi, and Toshiyuki Tanisaki, *$\mathcal{D}$-modules, perverse sheaves, and representation theory*, Progress in Mathematics, vol. 236, Birkhäuser Boston Inc., Boston, MA, 2008, Translated from the 1995 Japanese edition by Takeuchi.

[IMM16]  Hiroshi Iritani, Etienne Mann, and Thierry Mignon, *Quantum Serre Theorem as a Duality between Quantum D-modules*, International Mathematics Research Notices **2016** (2016), no. 9, 2828–2888.

[Iri09a]  Hiroshi Iritani, *An integral structure in quantum cohomology and mirror symmetry for toric orbifolds*, Adv. Math. **222** (2009), no. 3, 1016–1079.

[Iri09b]  \_\_\_\_\_\_, *$tt^*$-geometry in quantum cohomology*, Preprint arxiv:0906.1307, 2009.

[Iri11]  \_\_\_\_\_\_, *Quantum cohomology and periods*, Preprint arXiv:1112.1552, 2011.

[Kas08]  Masaki Kashiwara, *Equivariant derived category and representation of real semisimple Lie groups*, Representation theory and complex analysis, Lecture Notes in Math., vol. 1931, Springer, Berlin, 2008, Lectures from the C.I.M.E. Summer School held in Venice, June 10–17, 2004, Edited by Enrico Casadio Tarabusi, Andrea D'Agnolo and Massimo Picardello, pp. 137–234.

[KKP08]  L. Katzarkov, M. Kontsevich, and T. Pantev, *Hodge theoretic aspects of mirror symmetry*, From Hodge theory to integrability and TQFT $tt^*$-geometry (Providence, RI) (Ron Y. Donagi and Katrin Wendland, eds.), Proc. Sympos. Pure Math., vol. 78, Amer. Math. Soc., 2008, pp. 87–174.

[Kon95]  Maxim Kontsevich, *Enumeration of rational curves via torus actions*, The moduli space of curves (Texel Island, 1994) (R. Dijkgraaf, C. Faber, and G. van der Geer, eds.), Progr. Math., vol. 129, Birkhäuser Boston, Boston, MA, 1995, Papers from the conference held on Texel Island, April 1994, pp. 335–368.

[KS94]  Masaki Kashiwara and Pierre Schapira, *Sheaves on manifolds*, Grundlehren der Mathematischen Wissenschaften [Fundamental Principles of Mathematical Sciences], vol. 292, Springer-Verlag, Berlin, 1994, With a chapter in French by Christian Houzel, Corrected reprint of the 1990 original.

[KW06]  Frances Kirwan and Jonathan Woolf, *An introduction to intersection homology theory*, second ed., Chapman & Hall/CRC, Boca Raton, FL, 2006.

[MM11]  Etienne Mann and Thierry Mignon, *Quantum $\mathcal{D}$-modules for toric nef complete intersections*, Preprint arXiv:1112.1552, 2011.

[MMW05]  Laura Felicia Matusevich, Ezra Miller, and Uli Walther, *Homological methods for hypergeometric families*, J. Amer. Math. Soc. **18** (2005), no. 4, 919–941 (electronic).

[MS05]  Ezra Miller and Bernd Sturmfels, *Combinatorial commutative algebra*, Graduate Texts in Mathematics, vol. 227, Springer-Verlag, New York, 2005.

[Pan98]  Rahul Pandharipande, *Rational curves on hypersurfaces (after A. Givental)*, Astérisque (1998), no. 252, Exp. No. 848, 5, 307–340, Séminaire Bourbaki. Vol. 1997/98.

[Pha79]  Frédéric Pham, *Singularités des systèmes différentiels de Gauss-Manin*, Progress in Mathematics, vol. 2, Birkhäuser Boston, Mass., 1979, With contributions by Lo Kam Chan, Philippe Maisonobe and Jean-Étienne Rombaldi.

[Pha83]  \_\_\_\_\_\_, *Vanishing homologies and the n variable saddlepoint method*, Singularities, Part 2 (Arcata, Calif., 1981), Proc. Sympos. Pure Math., vol. 40, Amer. Math. Soc., Providence, RI, 1983, pp. 319–333.

[Pha85]  \_\_\_\_\_\_, *La descente des cols par les onglets de Lefschetz, avec vues sur Gauss-Manin*, Astérisque (1985), no. 130, 11–47, Differential systems and singularities (Luminy, 1983).

[Rei14]  Thomas Reichelt, *Laurent polynomials, GKZ-hypergeometric systems and mixed Hodge modules*, Compos. Math. **150** (2014), no. 6, 911–941.

[RS15]  Thomas Reichelt and Christian Sevenheck, *Logarithmic Frobenius manifolds, hypergeometric systems and quantum $\mathcal{D}$-modules*, J. Algebraic Geom. **24** (2015), no. 2, 201–281.

[Sab08]  Claude Sabbah, *Fourier-Laplace transform of a variation of polarized complex Hodge structure.*, J. Reine Angew. Math. **621** (2008), 123–158.

[Sab11] \_\_\_\_\_\_, *Non-commutative Hodge structures*, Ann. Inst. Fourier (Grenoble) **61** (2011), no. 7, 2681–2717.

[Sai88] Morihiko Saito, *Modules de Hodge polarisables*, Publ. Res. Inst. Math. Sci. **24** (1988), no. 6, 849–995 (1989).

[Sai89] \_\_\_\_\_\_, *On the structure of Brieskorn lattice*, Ann. Inst. Fourier (Grenoble) **39** (1989), no. 1, 27–72.

[Sai90] \_\_\_\_\_\_, *Mixed Hodge modules*, Publ. Res. Inst. Math. Sci. **26** (1990), no. 2, 221–333.

[SW09] Mathias Schulze and Uli Walther, *Hypergeometric $\mathcal{D}$-modules and twisted Gauß-Manin systems*, J. Algebra **322** (2009), no. 9, 3392–3409.

[Wal07] Uli Walther, *Duality and monodromy reducibility of A-hypergeometric systems*, Math. Ann. **338** (2007), no. 1, 55–74.

Thomas Reichelt
Mathematisches Institut
Universität Heidelberg
69120 Heidelberg
Germany
treichelt@mathi.uni-heidelberg.de

Christian Sevenheck
Fakultät für Mathematik
Technische Universität Chemnitz
09107 Chemnitz
Germany
christian.sevenheck@mathematik.tu-chemnitz.de